\documentclass[11pt]{article}
\usepackage[T1]{fontenc}
\usepackage[utf8]{inputenc}
\usepackage{lmodern}
\usepackage{amsmath,amssymb,amsthm}
\usepackage{mathtools}
\usepackage{array}
\usepackage{xcolor}
\usepackage[margin=1in]{geometry}
\usepackage{microtype}
\usepackage{parskip}
\usepackage{enumitem}
\usepackage{tikz}
\usepackage{pgfplots}
\pgfplotsset{compat=1.18}
\usepackage[hidelinks,bookmarksnumbered=true]{hyperref}
\hypersetup{pdftitle={Frostman random variables, entropy inequalities, and applications}, pdfauthor={Alex Iosevich, Thang Pham, Nguyen Dac Quan, Steven Senger, and Boqing Xue}}
\newtheorem{theorem}{Theorem}[section]
\newtheorem{proposition}[theorem]{Proposition}
\newtheorem{lemma}[theorem]{Lemma}
\newtheorem{corollary}[theorem]{Corollary}
\newtheorem{conjecture}[theorem]{Conjecture}
\theoremstyle{definition}

\newtheorem{question}[theorem]{Question}
\theoremstyle{remark}
\newtheorem{remark}[theorem]{Remark}

\title{Frostman random variables, entropy inequalities, and applications}
\author{Alex Iosevich \and Thang Pham \and Nguyen D. Quan \and Steven Senger \and Boqing Xue}
\date{}

\begin{document}
\maketitle
\begin{abstract}
We introduce jointly and conditionally Frostman random variables and study the growth of the dyadic entropies of $X+Y$ and $f(X,Y)$, where $f(x,y)$ is any real polynomial that is not a function of a single affine coordinate. This restriction cannot be removed, since entropy expansion fails for such functions even when the inputs are independent and Frostman. Under conditional Frostman bounds, we obtain, for every $0<s<1$, an explicit gain over $sn$ in the larger of the two output entropies. Under joint Frostman bounds, both output entropies may equal $sn+O(1)$, and we prove expansion once the average marginal entropy is included in the maximum. Uniformly over this class of polynomials, the sharp leading coefficients are $\frac{1+3s}{4}$ for $\frac{3}{4}\le s<1$ in the joint setting and $\frac{1+2s}{3}$ for $\frac{2}{3}\le s<1$ in the conditional setting, and we classify the polynomials for which the sharp conditional coefficient improves to $\frac{1+s}{2}$ when $\frac{3}{4}\le s<1$. As applications, we recover the known discretized sum--product exponents under uniform Frostman nonconcentration and extend the corresponding entropy estimates to pairs of polynomial images.
\end{abstract}
\tableofcontents

\section{Introduction}\label{sec:main}

A recurring principle in additive combinatorics and geometric measure theory is that a nonconcentrated set or measure cannot stay small under two sufficiently transverse algebraic maps. At a fixed scale, the Frostman condition is the measure-theoretic counterpart of the discretized nonconcentration conditions of Katz and Tao \cite{KatzTao}, and Bourgain's discretized sum--product and projection theorem \cite{BourgainProjection} established the basic expansion principle for such sets. Nonlinear versions of the projection theorem were later developed by Shmerkin \cite{ShmerkinNonlinear} and by Raz and Zahl \cite{RZ}, while the Elekes--R\'onyai dichotomy \cite{ER} identifies the rigid additive and multiplicative polynomial forms responsible for the exceptional cases.

Dyadic entropy is a convenient, scale-sensitive language for the same phenomenon. Local entropy averages have been very effective in the study of projections and differentiable images of fractal measures \cite{HochmanShmerkin}, and inverse theorems for the growth of entropy under convolution \cite{HochmanEntropy} are another influential part of this circle of ideas. More recently, M\'ath\'e and O'Regan \cite{MatheORegan} used Shannon-type information inequalities to obtain quantitative discretized sum--product estimates. None of these results is used as a black box below, but they shaped the methods of this paper.

The present paper develops the entropy framework of our earlier preprint, whose fifth and sixth versions are \cite{IPQSX} and \cite{IPQSXv6}, and subsumes it. Theorems~1.2 and~1.3 of \cite{IPQSX} treat jointly Frostman laws for diagonal polynomials of a common degree and conditionally Frostman laws for polynomials obtained from such diagonal forms by an invertible rational linear change of variables, respectively. In \cite{IPQSXv6}, the polynomial may be arbitrary outside the family $H(ax+by+c)$, with no rationality assumption on its coefficients, and explicit estimates are obtained in some ranges. Here we establish explicit entropy gains throughout the range $0<s<1$ for every real polynomial that is not a function of one affine coordinate, and we determine the sharp leading coefficients for $\frac{3}{4}\le s<1$ under joint Frostman bounds and for $\frac{2}{3}\le s<1$ under conditional Frostman bounds.

\subsection{Frostman hypotheses and entropy expansion}

All logarithms are to base two, and $\ln$ denotes the natural logarithm. The Shannon entropy of a discrete law with atom probabilities $p_i$ is $H=-\sum_i p_i\log p_i$, with the convention $0\log 0=0$. We write $\mathcal L(Z)$ for the law of a random variable $Z$ and use the dyadic quantization at level $n$, defined coordinatewise by $D_nz=2^{-n}\lfloor 2^nz\rfloor$, to set $H_n(Z)=H(D_nZ)$. Similarly, $H_n(\mu)$ denotes the dyadic entropy of a law $\mu$, the conditional dyadic entropy is $H_n(Z\mid W)=H(D_nZ\mid W)$, and the dyadic mutual information is
$${\mathcal{I}_n}(U;V)=H_n(U)+H_n(V)-H_n(U,V).$$
Without a subscript, ${\mathcal{I}}$ denotes the usual mutual information of discrete random variables.

Let $0<s\le1$ and $C>0$. A real random variable $X$ is $(s,C)$-Frostman if $\mathbb P(X\in I)\le C|I|^s$ for every interval $I$ of length at most one. A pair $(X,Y)$ is jointly $(s,s,C)$-Frostman if
$$\mathbb P(X\in I,Y\in J)\le C|I|^s|J|^s$$
for all intervals $I,J$ of length at most one, and it is conditionally $(s,C_1;s,C_2)$-Frostman if
$$\begin{aligned}
\mathbb P(X\in I,Y\in E)&\le C_1|I|^{s}\mathbb P(Y\in E),\\
\mathbb P(Y\in J,X\in E)&\le C_2|J|^{s}\mathbb P(X\in E)
\end{aligned}$$
for all such intervals and all Borel sets $E$. Equivalently, $\mathcal L(X\mid Y=y)$ is $(s,C_1)$-Frostman for $\mathcal L(Y)$-almost every $y$, and $\mathcal L(Y\mid X=x)$ is $(s,C_2)$-Frostman for $\mathcal L(X)$-almost every $x$. Independent $(s,C_1)$- and $(s,C_2)$-Frostman variables form a conditionally $(s,C_1;s,C_2)$-Frostman pair, and every conditionally $(s,C_1;s,C_2)$-Frostman pair is jointly $(s,s,C_1C_2)$-Frostman.

We call a polynomial $f\in\mathbb R[x,y]$ admissible if it cannot be written as
\begin{equation}\label{eq:excluded-affine-form}
f(x,y)=h(ax+by+c)
\end{equation}
with $h\in\mathbb R[t]$ and $a,b,c\in\mathbb R$, that is, if $f$ is not a polynomial in a single affine coordinate. No rationality assumption is imposed on the coefficients of $f$. We set $S=X+Y$ and $F=f(X,Y)$, and we study the quantity $M_n=\max\{H_n(S),H_n(F)\}$ together with
$$M_{\mathrm j,n}=\max\left\{\frac{H_n(X)+H_n(Y)}{2},M_n\right\}.$$
Throughout the paper, the support and the Frostman constants are fixed, and all error terms hold for every $n\ge1$, with implied constants independent of $n$ and of the particular law of $(X,Y)$. These constants may depend on the polynomial, the exponents, the support and the Frostman constants, and also on $\varepsilon$ in estimates involving this parameter.

For $d\ge2$, we define
\begin{equation}\label{eq:joint-main-gain}
\Delta_{\mathrm j}(s,d)=
\begin{cases}
\frac{s^2}{64d}, & 0<s\le\frac{1}{2},\\
\frac{s}{16}, & \frac{1}{2}<s\le\frac{2}{3},\\
\frac{2s-1}{8}, & \frac{2}{3}\le s\le\frac{3}{4},\\
\frac{1-s}{4}, & \frac{3}{4}\le s<1.
\end{cases}
\end{equation}

\begin{theorem}\label{thm:joint-main}
Let $0<s<1$, let $f$ be an admissible polynomial of degree $d$, and fix $R,C>0$. Every jointly $(s,s,C)$-Frostman pair supported in $[-R,R]^2$ satisfies
\begin{equation}\label{eq:joint-main}
M_{\mathrm j,n}\ge(s+\Delta_{\mathrm j}(s,d)-\varepsilon)n-O_\varepsilon(1)
\end{equation}
for every $\varepsilon>0$ and every integer $n\ge1$. In particular,
$$M_{\mathrm j,n}\ge\left(\frac{1+3s}{4}-\varepsilon\right)n-O_\varepsilon(1)$$
for $\frac{3}{4}\le s<1$, and for $\frac{5}{6}\le s<1$, this improves to
\begin{equation}\label{eq:joint-main-endpoint}
M_{\mathrm j,n}\ge\frac{1+3s}{4}n-O(\log(n+2)),
\end{equation}
where the error term can be replaced by $O(1)$ when $\frac{9}{10}<s<1$. Moreover, throughout the range $0<s<1$,
\begin{equation}\label{eq:joint-main-bounded}
M_{\mathrm j,n}\ge\left(s+\frac{s\min\{s,2-2s\}}{64d}\right)n-O(1).
\end{equation}
\end{theorem}

The leading coefficient $\frac{1+3s}{4}$ cannot be improved for $\frac{3}{4}\le s<1$ in a bound valid for all admissible polynomials, as the examples of Proposition~\ref{prop:joint-high-mixture-upper} show. Here and below, sharpness always refers to the leading coefficient, and we make no claim about the optimality of the error terms. The gain $\Delta_{\mathrm j}(s,d)$ is positive for every $0<s<1$, and it does not depend on $d$ when $s>\frac{1}{2}$. The marginal term $\frac{H_n(X)+H_n(Y)}{2}$ cannot be dropped from $M_{\mathrm j,n}$. Indeed, if $T$ is uniform on $[0,1]$, then $(T,-T)$ is jointly $(\frac{1}{2},\frac{1}{2},1)$-Frostman, while both outputs vanish when $f(x,y)=x^2-y^2$.

For the conditional estimate, we set
$$\begin{aligned}
s_{\mathrm a}&=\frac{203-\sqrt{2149}}{270}=0.580158\ldots,\\
s_{\mathrm b}&=\frac{45-\sqrt{65}}{60}=0.615629\ldots,
\end{aligned}$$
and we define the continuous gain
\begin{equation}\label{eq:conditional-main-gain}
\Delta_{\mathrm c}(s)=
\begin{cases}
\frac{4s}{31}, & 0<s\le s_{\mathrm a},\\
\frac{23s-15s^2-7}{26-15s}, & s_{\mathrm a}\le s\le s_{\mathrm b},\\
\frac{(1-s)(15s-7)}{19-15s}, & s_{\mathrm b}\le s\le\frac{2}{3},\\
\frac{1-s}{3}, & \frac{2}{3}\le s<1.
\end{cases}
\end{equation}

\begin{theorem}\label{thm:conditional-main}
Let $0<s<1$, let $f$ be admissible, and fix $R,C_1,C_2>0$. Every conditionally $(s,C_1;s,C_2)$-Frostman pair supported in $[-R,R]^2$ satisfies
\begin{equation}\label{eq:conditional-main}
M_n\ge(s+\Delta_{\mathrm c}(s)-\varepsilon)n-O_\varepsilon(1)
\end{equation}
for every $\varepsilon>0$. In particular,
$$M_n\ge\left(\frac{1+2s}{3}-\varepsilon\right)n-O_\varepsilon(1)$$
for $\frac{2}{3}\le s<1$, and for $\frac{5}{6}\le s<1$, this improves to
\begin{equation}\label{eq:conditional-main-endpoint}
M_n\ge\frac{1+2s}{3}n-O(\log(n+2)),
\end{equation}
where the error term can be replaced by $O(1)$ when $\frac{5}{6}<s<1$.
\end{theorem}

\begin{figure}[t]
\centering
\begin{tikzpicture}
\begin{axis}[
  width=0.5\textwidth, height=0.42\textwidth,
  xmin=0, xmax=1, ymin=0, ymax=0.24,
  xlabel={$s$},
  xtick={0,0.25,0.5,0.75,1},
  xticklabels={$0$,$\frac{1}{4}$,$\frac{1}{2}$,$\frac{3}{4}$,$1$},
  ytick={0,0.05,0.1,0.15,0.2},
  yticklabel style={/pgf/number format/fixed},
  title={conditional Frostman bounds},
  tick label style={font=\small}, label style={font=\small}, title style={font=\small}]
\addplot[thick, domain=0:0.580158, samples=2] {4*x/31};
\addplot[thick, domain=0.580158:0.615629, samples=30] {(23*x-15*x^2-7)/(26-15*x)};
\addplot[thick, domain=0.615629:0.666667, samples=30] {(1-x)*(15*x-7)/(19-15*x)};
\addplot[thick, domain=0.666667:1, samples=2] {(1-x)/3};
\addplot[dashed, domain=0:0.333333, samples=40] {x-x^2};
\addplot[dashed, domain=0.333333:0.666667, samples=2] {(1-x)/3};
\addplot[densely dotted, very thick, domain=0.75:1, samples=2] {(1-x)/2};
\end{axis}
\end{tikzpicture}\hfill
\begin{tikzpicture}
\begin{axis}[
  width=0.5\textwidth, height=0.42\textwidth,
  xmin=0, xmax=1, ymin=0, ymax=0.14,
  xlabel={$s$},
  xtick={0,0.25,0.5,0.75,1},
  xticklabels={$0$,$\frac{1}{4}$,$\frac{1}{2}$,$\frac{3}{4}$,$1$},
  ytick={0,0.04,0.08,0.12},
  yticklabel style={/pgf/number format/fixed},
  title={joint Frostman bounds, $d=2$},
  tick label style={font=\small}, label style={font=\small}, title style={font=\small}]
\addplot[thick, domain=0:0.5, samples=30] {x^2/128};
\addplot[thick, domain=0.5:0.666667, samples=2] {x/16};
\addplot[thick, domain=0.666667:0.75, samples=2] {(2*x-1)/8};
\addplot[thick, domain=0.75:1, samples=2] {(1-x)/4};
\addplot[dashed, domain=0:0.5, samples=2] {x/5};
\addplot[dashed, domain=0.5:0.75, samples=2] {(1-x)/4};
\end{axis}
\end{tikzpicture}
\caption{Gains over $s$ in the leading coefficient. On the left, the solid curve is the gain $\Delta_{\mathrm c}(s)$ of Theorem~\ref{thm:conditional-main}, and the dashed curve is the upper bound $\min\{s-s^2,\frac{1-s}{3}\}$ for any gain valid for all admissible polynomials, coming from Propositions~\ref{prop:conditional-low-upper} and~\ref{upper-modular:prop}. The two agree on $[\frac{2}{3},1)$. The dotted line is the sharp gain $\frac{1-s}{2}$ for polynomials outside the class $\mathfrak E_S$ of Theorem~\ref{thm:conditional-classification}. On the right, the solid curve is the gain $\Delta_{\mathrm j}(s,2)$ of Theorem~\ref{thm:joint-main} for quadratic polynomials, and the dashed curve is the upper bound for quadratic polynomials coming from Proposition~\ref{prop:joint-high-mixture-upper}. The two agree on $[\frac{3}{4},1)$.}
\label{fig:gains}
\end{figure}
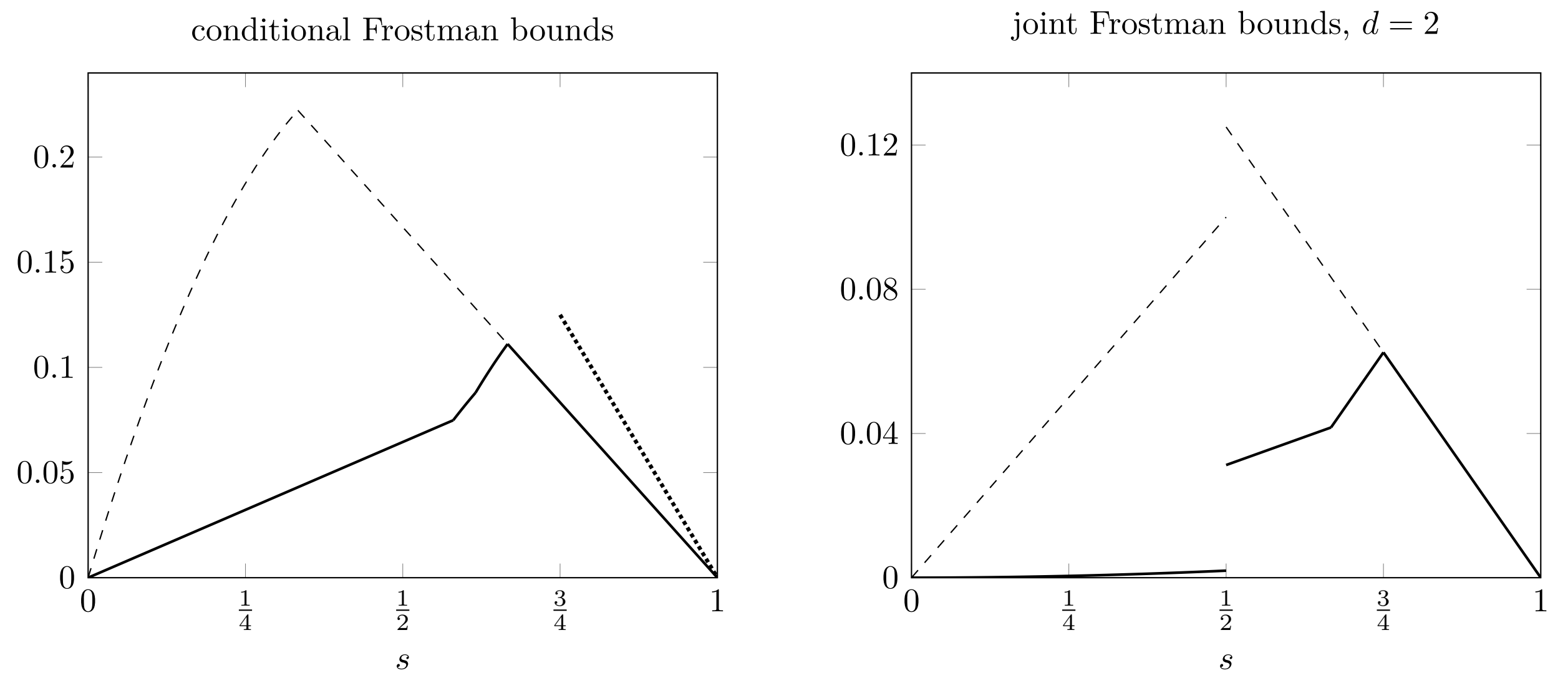

Again, the coefficient $\frac{1+2s}{3}$ is sharp for $\frac{2}{3}\le s<1$ uniformly over admissible polynomials, by the construction in Proposition~\ref{upper-modular:prop}, and Theorem~\ref{thm:conditional-classification} below identifies the sharp coefficient for each fixed polynomial in the relevant ranges. The gain $\Delta_{\mathrm c}(s)$ depends only on $s$ and is positive throughout $0<s<1$. Figure~\ref{fig:gains} compares the gains in Theorems~\ref{thm:joint-main} and~\ref{thm:conditional-main} with the upper bounds furnished by the constructions of Section~\ref{sec:sharpness-counterexamples}.

Neither theorem holds for arbitrary polynomials. Let $X$ and $Y$ be independent and uniformly distributed on two arithmetic progressions with the same step and about $2^{sn}$ terms each, smoothed by uniform noise at scale $2^{-n}$. Lemma~\ref{lem:progressions-mixtures} shows that such variables can be chosen in a fixed bounded rectangle so that they satisfy both the conditional and the joint Frostman hypotheses with constants independent of $n$, and for $f(x,y)=(x+y)^2$, they satisfy
$$M_n\le M_{\mathrm j,n}\le sn+O(1)$$
with an implied constant that is also independent of $n$.

\subsection{Sharp conditional coefficients and joint output bounds}

Under conditional Frostman bounds, the sharp leading coefficient of $M_n$ depends, in the ranges below, on whether $f$ belongs to the class
\begin{equation}\label{eq:ES-main}
\mathfrak E_S=\{h(ax+P(x+y)):h,P\in\mathbb R[t]\text{ and }a\in\mathbb R\}.
\end{equation}
For $f\in\mathfrak E_S$, the lower bound in the next theorem is already contained in Theorem~\ref{thm:conditional-main}, while for $f\notin\mathfrak E_S$, the coefficient improves from $\frac{1+2s}{3}$ to $\frac{1+s}{2}$ when $\frac{3}{4}\le s<1$. In both cases, matching constructions show that the coefficients are sharp for each admissible polynomial in the stated ranges.

\begin{theorem}\label{thm:conditional-classification}
Let $f$ be admissible and fix $R,C_1,C_2>0$. Every conditionally $(s,C_1;s,C_2)$-Frostman pair supported in $[-R,R]^2$ satisfies, for every $\varepsilon>0$ and every integer $n\ge1$,
$$M_n\ge
\begin{cases}
\left(\frac{1+2s}{3}-\varepsilon\right)n-O_\varepsilon(1) & \text{if } f\in\mathfrak E_S\text{ and }\frac{2}{3}\le s<1,\\
\left(\frac{1+s}{2}-\varepsilon\right)n-O_\varepsilon(1) & \text{if } f\notin\mathfrak E_S\text{ and }\frac{3}{4}\le s<1.
\end{cases}$$
For $\frac{5}{6}\le s<1$, these estimates improve to
$$M_n\ge
\begin{cases}
\frac{1+2s}{3}n-O(\log(n+2)) & \text{if } f\in\mathfrak E_S,\\
\frac{1+s}{2}n-O(\log(n+2)) & \text{if } f\notin\mathfrak E_S,
\end{cases}$$
and both error terms can be replaced by $O(1)$ when $\frac{5}{6}<s<1$.

Conversely, for each $s$ in the corresponding range, there is a constant $C_0=C_0(s,f,R)<\infty$ such that for every integer $n\ge1$, some conditionally $(s,C_0;s,C_0)$-Frostman pair $(X,Y)$ supported in $[-R,R]^2$ has each of the two entropies $H_n(X+Y)$ and $H_n(f(X,Y))$ equal to
$$\begin{cases}
\frac{1+2s}{3}n+O(1) & \text{if } f\in\mathfrak E_S\text{ and }\frac{2}{3}\le s<1,\\
\frac{1+s}{2}n+O(1) & \text{if } f\notin\mathfrak E_S\text{ and }\frac{3}{4}\le s<1.
\end{cases}$$
The constant $C_0$ and the constants in these error terms are independent of $n$. If $f\notin\mathfrak E_S$, a single law can be chosen for all $n$, while if $f\in\mathfrak E_S$, the law may depend on $n$.
\end{theorem}

The matching laws come from Proposition~\ref{upper-modular:prop} when $f\in\mathfrak E_S$ and from Proposition~\ref{prop:balanced-conditional-source} when $f\notin\mathfrak E_S$. The coefficients are the same for all polynomials in each class, while the error constants and the matching laws may depend on $f$, and the constant $C_0$ of the construction need not equal either of the constants $C_1,C_2$ in the lower bounds.

The class $\mathfrak E_S$ has a transparent description in the coordinates $(X,S)$ and $(Y,S)$. Following Raz and Zahl \cite{RZ}, we say that a real-analytic function $\Phi(x,y)$ is analytic-special if, on each connected component of the set where $\Phi_x\Phi_y\ne0$, it can be written as $\Phi(x,y)=g(h(x)+k(y))$ with analytic functions $g,h,k$ of one variable, and we call $\Phi$ nonspecial otherwise. Where $\Phi_x\Phi_y\ne0$, the function
$$\partial_x\partial_y\ln\left|\frac{\Phi_x}{\Phi_y}\right|$$
is called the curvature expression of $\Phi$, and it vanishes identically exactly when $\Phi$ is analytic-special \cite[Lemma~3.2]{RZ}. Writing $f$ in the coordinates $(x,z)$ and $(y,z)$ with $z=x+y$, we obtain the two sum-coordinate representations $G_X(x,z)=f(x,z-x)$ and $G_Y(y,z)=f(z-y,y)$, so that $F=G_X(X,S)=G_Y(Y,S)$. Theorem~\ref{thm:two-transition} below shows that a a nonconstant polynomial $f$ belongs to $\mathfrak E_S$ if and only if both $G_X$ and $G_Y$ are analytic-special.

For example, $x+(x+y)^2$ has coefficient $\frac{1+2s}{3}$ for $s\ge\frac{2}{3}$, while the quadratic forms $x^2+y^2$ and $xy$, which lie outside $\mathfrak E_S$, have coefficient $\frac{1+s}{2}$ for $s\ge\frac{3}{4}$. Joint Frostman bounds alone lead to a different picture, since for these three polynomials the sharp coefficient of $M_n$ is then $s$ throughout $0<s\le1$. When $s>\frac{1}{2}$, the same coefficient $s$ is sharp for every admissible polynomial, as the next proposition shows.

\begin{proposition}\label{thm:joint-models}
Let $0<s\le1$ and fix $R,C>0$. Every jointly $(s,s,C)$-Frostman pair supported in $[-R,R]^2$ satisfies the following estimates for every integer $n\ge1$, where the positive constants $K_{a,P,R}$ and $K_R$ depend only on the indicated parameters.
\begin{enumerate}[label=\textup{(\roman*)},nosep]
\item If $f(x,y)=ax+P(x+y)$ with $a\ne0$ and $\deg P\ge2$, then
\begin{equation}\label{joint-models:triangular-lower}
M_n\ge sn-\frac{1}{2}\log C-K_{a,P,R}.
\end{equation}
\item If $f(x,y)=x^2+y^2$ or $f(x,y)=xy$, then
\begin{equation}\label{joint-models:radial-lower}
H_n(S)+H_n(F)\ge 2sn-K_R(\log(C+2)+\log(n+2)).
\end{equation}
\end{enumerate}
If $s>\frac{1}{2}$, then
\begin{equation}\label{joint-models:general-recovery}
H_n(S,F)=H_n(X,Y)+O_{s,f,R,C}(1)\ge 2sn-O_{s,f,R,C}(1)
\end{equation}
for every polynomial $f$ with $f_y-f_x\not\equiv0$, and in particular for every admissible polynomial. Conversely, for every polynomial $f$ with $f_y-f_x\not\equiv0$ and every $0<s\le1$, there is a single jointly $(s,s,C_0)$-Frostman law supported in $[-R,R]^2$, with $C_0=C_0(s,f,R)<\infty$, such that
\begin{equation}\label{joint-models:fixed-upper}
H_n(S)=sn+O(1)\text{ and }H_n(F)=sn+O(1)\text{ for all }n\ge1.
\end{equation}
The law, the constant $C_0$ and the bounded errors in \eqref{joint-models:fixed-upper} do not depend on $n$.
\end{proposition}

\subsection{Main ideas and sketch of proofs}

Two estimates drive the proofs. The first, Lemma~\ref{lem:copied-output-incidence}, turns a rectangular incidence estimate for curves into an entropy inequality for conditional copies, with explicit terms that measure the dependence between the inputs. The second, Theorem~\ref{thm:local-multiscale-profile}, is an entropy inequality on regular charts obtained from a multiscale Fourier argument. Let us explain how they lead to the three main theorems.

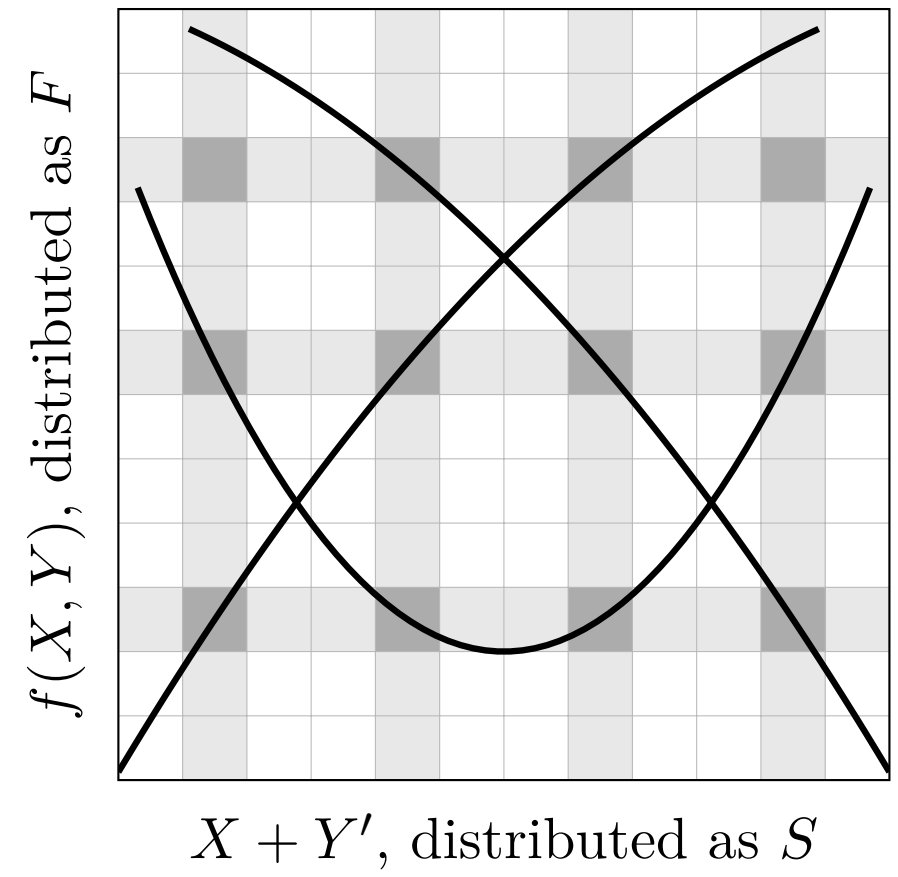
\begin{figure}[t]
\centering
\begin{tikzpicture}[scale=0.85]
\foreach \c in {1,4,7,10} {\fill[gray!18] (\c*0.5,0) rectangle (\c*0.5+0.5,6);}
\foreach \r in {2,6,9} {\fill[gray!18] (0,\r*0.5) rectangle (6,\r*0.5+0.5);}
\foreach \c in {1,4,7,10} {\foreach \r in {2,6,9} {\fill[gray!65] (\c*0.5,\r*0.5) rectangle (\c*0.5+0.5,\r*0.5+0.5);}}
\draw[step=0.5, gray!55, very thin] (0,0) grid (6,6);
\draw (0,0) rectangle (6,6);
\begin{scope}
\clip (0,0) rectangle (6,6);
\draw[very thick, domain=0:5.45, samples=60] plot (\x, {-\x*\x/9+5*\x/3+1/16});
\draw[very thick, domain=0.55:6, samples=60] plot (\x, {-(6-\x)*(6-\x)/9+5*(6-\x)/3+1/16});
\draw[very thick, domain=0.15:5.85, samples=60] plot (\x, {4*(\x-3)*(\x-3)/9+1});
\end{scope}
\node[below] at (3,-0.1) {$X+Y'$, distributed as $S$};
\node[rotate=90, above] at (-0.1,3) {$f(X,Y)$, distributed as $F$};
\end{tikzpicture}
\caption{The incidence step in the proof of Theorem~\ref{thm:joint-main}. The light columns and rows are the dyadic cells occupied by $S$ and by $F$. Since $X+Y'$ has the law of $S$ and $f(X,Y)$ has the law of $F$, every point $(X+Y',f(X,Y))$ lies in one of the dark squares, and it also lies on the curve $t\mapsto f(t-Y',Y)$. If both outputs occupy few cells, many curves of a transversal family must pass through a small set of squares, and the incidence theorem limits how often this can happen.}
\label{fig:incidences}
\end{figure}

In Theorem~\ref{thm:joint-main}, we must show that one of the outputs or the average marginal entropy grows. After interchanging $X$ and $Y$ if necessary, we may choose an orientation in which the family of curves described below is nondegenerate away from a proper algebraic set. Let $Y$ and $Y'$ be conditionally independent given $X$, each having the prescribed joint law with $X$. The point $(X+Y',f(X,Y))$ lies on the curve $t\mapsto f(t-Y',Y)$, and its two coordinates have the laws of $S$ and $F$, respectively. If both $S$ and $F$ are concentrated on few dyadic intervals, these points lie in few dyadic squares, and we obtain many incidences between dyadic squares and curves, as in Figure~\ref{fig:incidences}. After grouping cells of comparable probability, we can apply the quantitative incidence theorem of O'Regan, Wu and Yi~\cite{OWY}. The main difficulty is the possible concentration of the law near the algebraic set where the family of curves degenerates, and we control the resulting logarithmic loss by the excess of the average marginal entropy over $sn$. Either this excess already provides the required growth, or the incidence estimate forces one of the outputs to grow. The multiplicities of the polynomials involved account for the dependence on the degree, and when $s>\frac{1}{2}$, thin neighborhoods of the critical algebraic curves carry uniformly small mass, so that the leading coefficient becomes independent of the degree.

For the stronger bounds, we use the algebraic classification of Corollary~\ref{cor:three-coordinate}, according to which every admissible polynomial has a nonspecial representation $F=G(U,V)$ with $(U,V)$ equal to $(X,Y)$, $(X,S)$ or $(Y,S)$. On regular charts, that is, on small product rectangles where $G_u$, $G_v$ and the curvature expression of $G$ do not vanish, a Fourier estimate for the associated Radon-type operator is available. On each class of cells of comparable probability, we smooth the restricted input law at the input scale, form the product of its two coordinate marginals, apply $G$ to both laws and add the same independent noise to each output. The Fourier estimate bounds the probability that independent samples from the two output laws fall in the same dyadic cell, data processing bounds the relative entropy of the outputs by the mutual information of the smoothed input law, and Jensen's inequality converts the collision bound into a lower bound for the output entropy plus this information. Smoothing and noise change the dyadic output entropy by only $O(1)$, and iterating the resulting inequality and averaging over the regular charts leads to
$$H_n(F)+{\mathcal{I}_n}(U;V)\ge\left(\min\left\{1,2s-\frac{1}{2}\right\}-\varepsilon\right)n-O_\varepsilon(1)$$
for $\frac{1}{2}<s<1$. Together with the Frostman lower bounds for the input entropies, this yields the coefficient $\frac{1+3s}{4}$ for $\frac{3}{4}\le s<1$, up to an arbitrarily small loss.

In Theorem~\ref{thm:conditional-main}, the exact conditional Frostman bounds provide additional entropy estimates, which allow us to eliminate the input entropies. For $0<s\le\frac{2}{3}$, we apply the copied incidence estimate both in the original coordinates and on regular local inverse branches of $(x,y)\mapsto(x+y,f(x,y))$. On a common restricted law, the output mutual information terms of the forward and inverse inequalities cancel, and we obtain the coefficient $\frac{35s}{31}$ up to an arbitrarily small loss. Localization requires some care here, because restricting a conditionally Frostman law to an event need not preserve uniform conditional bounds. We therefore split the restricted law into pieces on which both conditional bounds are under control, eliminate the copied variables on each piece, and only then average, keeping track of the normalization and labeling costs. For larger exponents, the copy inequalities strengthen the Fourier iteration by retaining the output mutual information ${\mathcal{I}_n}(S;F)$, and this already gives the coefficient $\frac{1+2s}{3}$ for $\frac{2}{3}\le s<1$. Keeping track, in addition, of the probabilities of cells at the final observation scale produces one more entropy inequality, and two different ways of eliminating its mutual information term lead to the two middle branches of \eqref{eq:conditional-main-gain}. In every case, the algebraic classification guarantees that one of these inequalities applies, either in the original or in the sum coordinates. The bounds with error $O(\log(n+2))$ or $O(1)$ in both theorems come from separate, direct Fourier estimates.

Theorem~\ref{thm:conditional-classification} combines these lower bounds with matching examples. Both sum-coordinate representations are special exactly when $f\in\mathfrak E_S$, and in this class, the coefficient $\frac{1+2s}{3}$ of Theorem~\ref{thm:conditional-main} is sharp for $\frac{2}{3}\le s<1$. Outside $\mathfrak E_S$, one of the sum-coordinate representations is nonspecial, the stronger entropy inequality applies to it, and the coefficient improves to $\frac{1+s}{2}$ for $\frac{3}{4}\le s<1$. The examples showing sharpness outside $\mathfrak E_S$ are mixtures of two fixed laws spread along families of level curves, as in Figure~\ref{fig:level-sets}, whose output entropy pairs $(H_n(S),H_n(F))$ are $(sn,n)+O(1)$ and $(n,sn)+O(1)$. Their equal mixture keeps both exact conditional Frostman bounds, and each of its output entropies equals $\frac{1+s}{2}n+O(1)$ at every scale. Inside $\mathfrak E_S$, we use a modular grid construction that places each output in $O(2^{\frac{(1+2s)n}{3}})$ dyadic cells at level $n$. Polynomial discrepancy estimates and a random selection produce comparable row and column degrees, together with uniform bounds on the number of neighbors in every interval of consecutive indices, and these counts guarantee that both exact conditional Frostman bounds hold at every scale after smoothing. The smoothed law has joint output entropy $\frac{2(1+2s)}{3}n+O(1)$, so subadditivity and the scalar upper bounds force each output entropy to equal $\frac{1+2s}{3}n+O(1)$. This law depends on $n$, but its Frostman constants do not.

The same two estimates lead to the results of the next subsection. For independent inputs, the dependence terms in Lemma~\ref{lem:copied-output-incidence} vanish, which yields Theorem~\ref{thm:independent-polynomial-pair} and the dimension bound for Cartesian products in Corollary~\ref{cor:polynomial-image-product}, while adding the joint input entropy to the inequality of Theorem~\ref{thm:local-multiscale-profile} leads to the three-image bound of Corollary~\ref{cor:sharp-three-images}.

\subsection{Some consequences of the method}

For a pair of polynomial maps, the nonconstancy condition in the next theorem guarantees that, after interchanging the inputs if necessary, the family of curves in Lemma~\ref{lem:copied-output-incidence} is nondegenerate on suitable charts. For independent inputs, that lemma then yields the following bound.

\begin{theorem}\label{thm:independent-polynomial-pair}
Let $0<s<1$, and let $P,Q\in\mathbb R[x,y]$ depend on both variables. Suppose that the rational function
$$\frac{P_xQ_y}{P_yQ_x}$$
is nonconstant, and fix $R,C_1,C_2>0$. If $X$ and $Y$ are independent, supported in $[-R,R]$, and have $(s,C_1)$- and $(s,C_2)$-Frostman laws, respectively, then, for every $\varepsilon>0$ and every integer $n\ge1$,
\begin{equation}\label{eq:independent-polynomial-pair}
H_n(P(X,Y))+H_n(Q(X,Y))\ge\left(\min\left\{\frac{5s}{2},1+s\right\}-\varepsilon\right)n-O_\varepsilon(1).
\end{equation}
The implied constant is uniform in the laws and in $n\ge1$, and it may depend on $s,P,Q,R,C_1,C_2$ and $\varepsilon$.
\end{theorem}

Entropy is at most the logarithm of the number of occupied cells, so the entropy bound implies a lower bound for the number of dyadic cells meeting the Cartesian product of the two image sets. Choosing independent inputs with a common Frostman law, we obtain the following dimension estimate.

\begin{corollary}\label{cor:polynomial-image-product}
Let $P,Q\in\mathbb R[x,y]$ satisfy the hypotheses of Theorem~\ref{thm:independent-polynomial-pair}, and let $A\subset\mathbb R$ be a bounded Borel set with $\dim_{\mathrm H}A=s\in(0,1]$. Then
\begin{equation}\label{eq:polynomial-image-product}
\underline{\dim}_{\mathrm B}(P(A,A)\times Q(A,A))\ge\min\left\{\frac{5s}{2},1+s\right\},
\end{equation}
where $P(A,A)=\{P(a,b):a,b\in A\}$ and $\underline{\dim}_{\mathrm B}$ denotes the lower box dimension.
\end{corollary}

For $P(x,y)=x+y$ and $Q(x,y)=xy$, the bound in Corollary~\ref{cor:polynomial-image-product} already follows from the work of Ren and Wang~\cite[Theorem~1.1]{RenWang}, even with Hausdorff dimension in place of lower box dimension, and an explicit derivation can be found in Orponen's survey~\cite[Section~4.1]{OrponenIncidenceSurvey}. The case $P(x,y)=x+y$ and $Q(x,y)=f(x,y)$ has a finite-set analogue in Shen's estimate $|A+A|\cdot|f(A,A)|\gtrsim_f|A|^{\frac{5}{2}}$, valid for every admissible $f$ and every finite set $A\subset\mathbb R$~\cite[Theorem~3.4]{ShenAlgebraicSumProduct}.

Corollary~\ref{cor:polynomial-image-product} concerns the full image sets, whereas Theorem~\ref{thm:independent-polynomial-pair} concerns the distributions of the images under a prescribed product law. The entropy estimate implies the covering bound, but a covering bound alone does not control the entropy of a given law, even for independent inputs, and a dependent law may have smaller output entropies. For instance, the jointly $(s,s,C)$-Frostman laws in Proposition~\ref{thm:joint-models} have output entropies adding up to only $2sn+O(1)$. This does not contradict \eqref{eq:polynomial-image-product}, since the support of such a law is not a Cartesian product of the form $A\times A$.

The entropy inequality of Theorem~\ref{thm:local-multiscale-profile}, which is used in the proofs of Theorems~\ref{thm:joint-main}--\ref{thm:conditional-classification}, also gives a sharp result for three images of a planar law, namely the two coordinate projections $u$, $v$ and a real-analytic function $g(u,v)$. The lower bound below follows by combining Theorem~\ref{thm:local-multiscale-profile} with the input entropy provided by the ball bound, and the fixed law constructed in Proposition~\ref{prop:three-image-fixed-sharpness} shows that it is sharp.

\begin{corollary}\label{cor:sharp-three-images}
Let $B_0,B_1\subset\mathbb R^2$ be compact product rectangles with nonempty interiors and $B_0\subset\operatorname{int}B_1$. Let $g$ be real analytic near $B_1$, with
$$g_ug_v\partial_u\partial_v\ln\left|\frac{g_u}{g_v}\right|\ne0\text{ on }B_1.$$
Fix $\frac{3}{2}\le\alpha\le2$ and $C\ge1$. Every law of $(U,V)$ supported in $B_0$ that satisfies the Frostman bound
$$\mathbb P((U,V)\in B(z,r))\le Cr^\alpha\text{ for all }z\in\mathbb R^2\text{ and }0<r\le1$$
obeys, for every $\varepsilon>0$ and $n\ge1$,
\begin{equation}\label{eq:sharp-three-image-lower}
\max\{H_n(U),H_n(V),H_n(g(U,V))\}\ge\left(\frac{\alpha+1}{3}-\varepsilon\right)n-O_{\alpha,\varepsilon,g,B_0,B_1}(1+\log C).
\end{equation}
The coefficient is optimal. For some finite $C_0$, there is a single law supported in $B_0$ and satisfying the same ball bound with constant $C_0$ such that
\begin{equation}\label{eq:sharp-three-image-upper}
H_n(Z)=\frac{\alpha+1}{3}n+O(1)\text{ for }Z\in\{U,V,g(U,V)\}\text{ and all }n\ge1.
\end{equation}
Here $C_0$ and the constants in the $O(1)$ terms of \eqref{eq:sharp-three-image-upper} depend only on $\alpha,g,B_0,B_1$.
\end{corollary}

A related estimate for sets was proved by Raz and Zahl~\cite[Theorem~1.11]{RZ}. With $B_0$, $B_1$ and $g$ as in Corollary~\ref{cor:sharp-three-images}, and with $\pi_1(u,v)=u$ and $\pi_2(u,v)=v$, their result says that for each $0<\alpha<2$, there is $c=c(\alpha)>0$ such that, for all sufficiently small $\delta>0$, every union $E\subset B_0$ of squares of side $\delta$ satisfying
$$|E|\ge\delta^{2-\alpha+c}\text{ and }|E\cap B|\le r^\alpha\delta^{2-\alpha-c}\text{ for every ball }B\text{ of diameter }r\ge\delta$$
obeys
$$\max\{|\pi_1(E)|,|\pi_2(E)|,|g(E)|\}\ge\delta^{1-\frac{\alpha}{2}-c}.$$
Here $|\cdot|$ is Lebesgue measure in the relevant dimension, and the threshold for $\delta$ may depend on $\alpha$, $g$, $B_0$ and $B_1$. In other words, one of the three images cannot be covered by fewer than $\delta^{-\frac{\alpha}{2}-c}$ intervals of length $\delta$. This is a statement about the size of image sets, whereas Corollary~\ref{cor:sharp-three-images} concerns Shannon entropies under a prescribed planar law and gives the sharp coefficient $\frac{\alpha+1}{3}$ for $\frac{3}{2}\le\alpha\le2$.

For independent inputs, the mutual information term in Theorem~\ref{thm:local-multiscale-profile} vanishes, and after localization, we obtain the following bound for a single polynomial image.

\begin{corollary}\label{cor:independent-single-image}
Let $f\in\mathbb R[x,y]$ be a polynomial that is not analytic-special. Fix $R>0$, exponents $0<s_1,s_2\le1$ with $s_1+s_2>1$, and constants $C_1,C_2\ge1$. If $X$ and $Y$ are independent, supported in $[-R,R]$, and have $(s_1,C_1)$- and $(s_2,C_2)$-Frostman laws, respectively, then, for every $\varepsilon>0$ and every integer $n\ge1$,
\begin{equation}\label{eq:independent-single-image}
H_n(f(X,Y))\ge\left(\min\left\{1,s_1+s_2-\frac{1}{2}\right\}-\varepsilon\right)n-O_\varepsilon(1).
\end{equation}
The implied constant is uniform in the input laws and in $n$, and it may depend on $f,R,s_1,s_2,C_1,C_2$ and $\varepsilon$. In particular, if $s_1+s_2\ge\frac{3}{2}$, then
$$\lim_{n\to\infty}\frac{H_n(f(X,Y))}{n}=1.$$
\end{corollary}

\begin{figure}[ht]
\centering
\begin{tikzpicture}[
box/.style={draw,rounded corners=2pt,align=center,font=\small,inner sep=4pt},
res/.style={draw,align=center,font=\small,inner sep=4pt,text width=2.15cm,minimum height=1.1cm},
arr/.style={->,>=stealth},
darr/.style={->,>=stealth,dashed}
]
\node[box,text width=8.2cm] (N2) at (0,0) {Section~\ref{sec:algebraic-entropy}: classification in sum coordinates (Theorem~\ref{thm:two-transition}) and entropy under changes of coordinates (Theorem~\ref{thm:coordinate-entropy})};
\node[box,text width=4.9cm] (N3) at (-3.35,-2.3) {Section~\ref{sec:incidence-copies}: incidences for conditional copies (Lemma~\ref{lem:copied-output-incidence}, Theorem~\ref{thm:joint-incidence-bound}, Proposition~\ref{prop:conditional-inverse-output})};
\node[box,text width=4.9cm] (N4) at (3.35,-2.3) {Section~\ref{sec:fourier-comparisons}: Radon estimates (Proposition~\ref{prop:entropy-radon-shell}, Lemma~\ref{lem:output-profile-transfer}, Proposition~\ref{prop:joint-mass-classes})};
\node[box,text width=5.2cm] (N5) at (0,-4.5) {Section~\ref{sec:multiscale-estimates}: multiscale entropy inequality (Theorem~\ref{thm:local-multiscale-profile}) and conditional bounds};
\node[res] (R14) at (-6.75,-6.7) {Proposition\\\ref{thm:joint-models}};
\node[res] (R15) at (-4.05,-6.7) {Theorem~\ref{thm:independent-polynomial-pair}, Corollary~\ref{cor:polynomial-image-product}};
\node[res] (R11) at (-1.35,-6.7) {Theorem~\ref{thm:joint-main}};
\node[res] (R17) at (1.35,-6.7) {Corollaries~\ref{cor:sharp-three-images} and~\ref{cor:independent-single-image}};
\node[res] (R12) at (4.05,-6.7) {Theorem~\ref{thm:conditional-main}};
\node[res] (R13) at (6.75,-6.7) {Theorem~\ref{thm:conditional-classification}};
\node[box,text width=8.2cm] (N8) at (0,-8.9) {Section~\ref{sec:sharpness-counterexamples}: constructions and counterexamples};
\draw[arr] (N3 |- N2.south) -- (N3.north);
\draw[arr] (N4 |- N2.south) -- (N4.north);
\draw[arr] (N3.south) -- (N5.north west);
\draw[arr] (N4.south) -- (N5.north east);
\draw[arr] (N2.west) -| (R14.north);
\draw[arr] (R15 |- N3.south) -- (R15.north);
\draw[arr] (N5.south) -- (R11.north);
\draw[arr] (N5.south) -- (R17.north);
\draw[arr] (N5.south) -- (R12.north);
\draw[arr] (N5.south) -- (R13.north);
\draw[arr] (R12.east) -- (R13.west);
\foreach \r in {R11,R17,R12}{\draw[darr] (\r |- N8.north) -- (\r.south);}
\draw[darr] (N8.north west) -- (R14.south);
\draw[darr] (N8.north east) -- (R13.south);
\end{tikzpicture}
\caption{The structure of the paper. A solid arrow indicates that results of one part are used in the proofs of the next, and dashed arrows indicate the constructions that show the sharpness of the corresponding results. The lower bounds in Proposition~\ref{thm:joint-models} and Corollary~\ref{cor:polynomial-image-product} are proved in Section~\ref{sec:joint-models}, and Theorem~\ref{thm:conditional-classification} is proved in Section~\ref{sec:sharpness}.}
\label{fig:roadmap}
\end{figure}
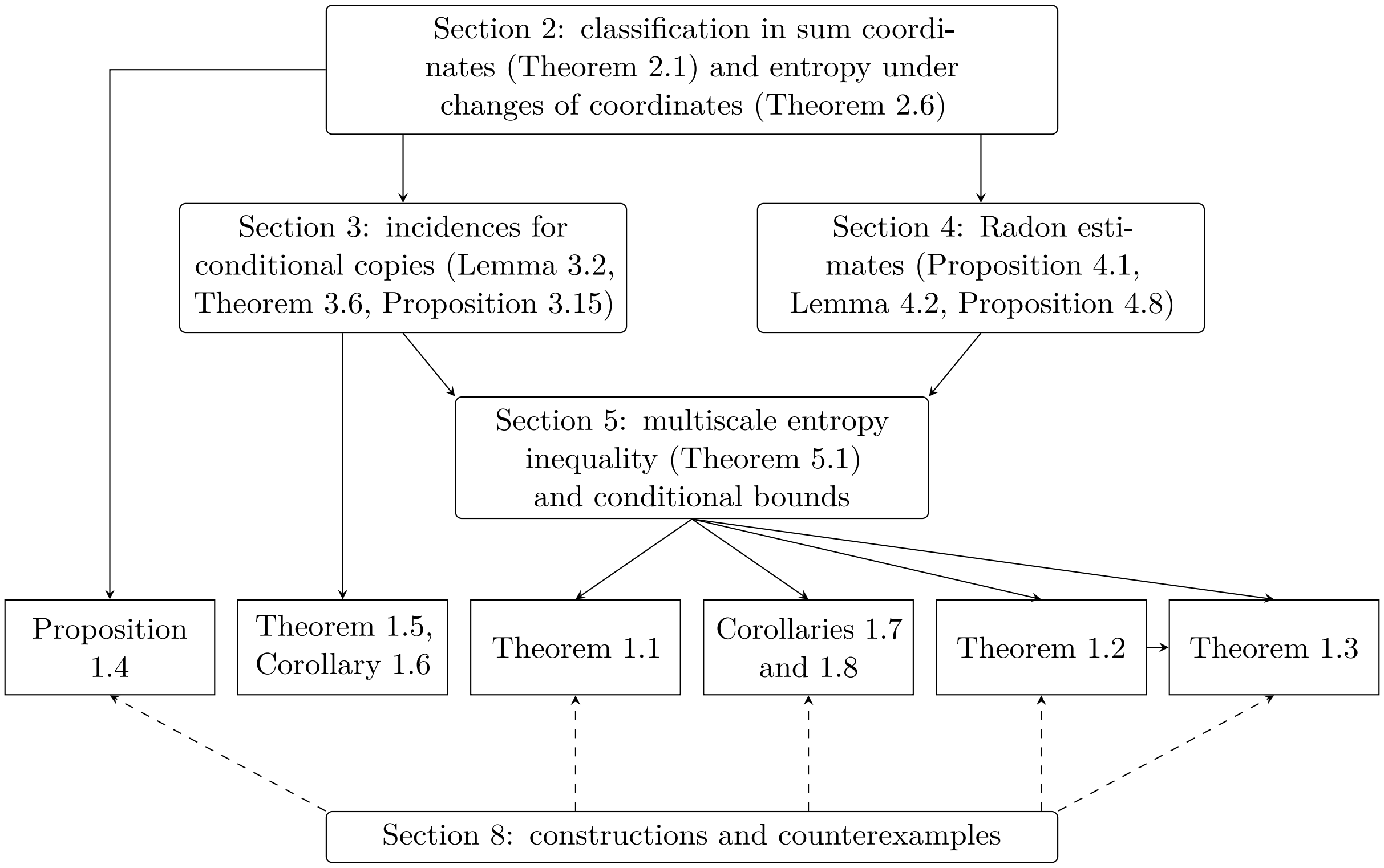

The paper is organized as follows. Section~\ref{sec:algebraic-entropy} contains the polynomial classification and the basic entropy and localization estimates, and Section~\ref{sec:incidence-copies} develops the incidence estimates for conditional copies and proves Theorem~\ref{thm:independent-polynomial-pair}. The Fourier estimates are developed in Section~\ref{sec:fourier-comparisons}, and in Section~\ref{sec:multiscale-estimates}, we combine these ingredients in a multiscale entropy argument, which proves Theorems~\ref{thm:joint-main} and~\ref{thm:conditional-main}, the lower bound in Corollary~\ref{cor:sharp-three-images} and Corollary~\ref{cor:independent-single-image}. Section~\ref{sec:sharpness} proves Theorem~\ref{thm:conditional-classification}, and Section~\ref{sec:joint-models} proves the lower bounds in Proposition~\ref{thm:joint-models} together with Corollary~\ref{cor:polynomial-image-product}. All constructions are collected in Section~\ref{sec:sharpness-counterexamples}, in the order of the results of this introduction, together with the matching parts of Proposition~\ref{thm:joint-models} and Corollary~\ref{cor:sharp-three-images} and the arithmetic counterexamples. We conclude in Section~\ref{sec:conjectures-support} with some questions about the optimal coefficients. Figure~\ref{fig:roadmap} shows how these parts depend on one another.

\section{Polynomial classification and entropy estimates}
\label{sec:algebraic-entropy}

We first determine which polynomial coordinate systems are available for the analytic arguments, and then we collect the entropy and localization estimates used throughout the paper.

\subsection{Polynomial classification in sum coordinates}
\label{sec:polynomial-structure}

A polynomial $f$ is of special form if
\begin{equation}\label{eq:polynomial-special}
f(x,y)=h_0(A(x)+B(y))\text{ or }f(x,y)=h_0(A(x)B(y))
\end{equation}
for some univariate real polynomials $h_0,A,B$. The regular set of an analytic function $\Phi(x,y)$ is the set where both first partial derivatives are nonzero, and, as in the introduction, $\Phi$ is analytic-special if, on each connected component of the regular set, it has the form
\begin{equation}\label{eq:analytic-special}
\Phi(x,y)=g(h(x)+k(y))
\end{equation}
with analytic functions of one variable. Multiplicative forms are included by taking logarithms on the components where the factors have constant signs. For polynomials, the two notions agree. If a first partial derivative vanishes identically, the polynomial depends on one variable only, and if the mixed partial derivative vanishes identically, it is a sum of a function of $x$ and a function of $y$. Both cases have the special forms above, and in the remaining case, the analytic and the polynomial classes are identified in \cite[Lemmas~3.2 and~3.3]{RZ}. On a regular product rectangle, \eqref{eq:analytic-special} implies
\begin{equation}\label{eq:separated-derivative-ratio}
\frac{\Phi_x}{\Phi_y}=\frac{h'(x)}{k'(y)}.
\end{equation}
For a polynomial, the separated factors can be chosen rational, as explained in the proof below, and clearing denominators then extends this identity from a nonempty regular rectangle to an identity between rational functions. Functions that are not analytic-special are called nonspecial. By \cite[Lemma~3.2]{RZ}, $\Phi$ is analytic-special exactly when its curvature expression $\partial_x\partial_y\ln\left|\frac{\Phi_x}{\Phi_y}\right|$ vanishes throughout the regular set, so a nonspecial polynomial has a regular product rectangle on which the curvature expression is bounded away from zero.

Recall the class $\mathfrak E_S$ from \eqref{eq:ES-main} and the two sum-coordinate representations
\begin{equation}\label{eq:sum-transitions}
G_X(x,z)=f(x,z-x)\text{ and }G_Y(y,z)=f(z-y,y).
\end{equation}

\begin{theorem}\label{thm:two-transition}
For a nonconstant real polynomial $f$, both $G_X$ and $G_Y$ are analytic-special if and only if $f\in\mathfrak E_S$.
\end{theorem}

\begin{proof}
If $f=h_0(ax+P(x+y))$, then $G_X(x,z)=h_0(ax+P(z))$ and $G_Y(y,z)=h_0(-ay+P(z)+az)$, so both representations are special.

Conversely, suppose that both representations are special, and write $p=f_x$ and $q=f_y$. If $p$, $q$ or $p-q$ vanishes identically, then $f$ is a polynomial in $y$, in $x$ or in $x+y$, respectively, and all of these belong to $\mathfrak E_S$, where for a polynomial in $y$, we use $y=-x+(x+y)$. We may therefore assume that $p$, $q$ and $p-q$ are nonzero, and we set $r=\frac{p}{q}$. Our first aim is to show that $r-1$ depends only on $z=x+y$, which will give a first-order equation for $G_X$ that can be integrated.

The derivative ratios of the two representations are
$$\begin{aligned}
\frac{(G_X)_x}{(G_X)_z}&=r-1,\\
\frac{(G_Y)_y}{(G_Y)_z}&=\frac{1-r}{r},
\end{aligned}$$
and applying \eqref{eq:separated-derivative-ratio} to each of them, with $z=x+y$, we obtain factorizations
$$\begin{aligned}
r-1&=A(x)B(z),\\
\frac{1-r}{r}&=C(y)D(z)
\end{aligned}$$
with nonzero one-variable factors, which may be chosen rational. To see this for the first ratio, choose a regular base point $(x_*,z_*)$ at which it does not vanish. Fixing $z=z_*$ and dividing by the value at the base point produces one normalized factor, fixing $x=x_*$ produces the other, and both are rational functions. The same argument in the coordinates $(y,z)$ applies to the second ratio. Each factorization holds on a product rectangle, so clearing denominators turns it into an identity between rational functions.

The two factorizations constrain the translates of a single rational function of one variable. Set $P_0=\frac{1}{A}$ and $Q_0=\frac{1}{C}$. Since
$$\begin{aligned}
r&=1+\frac{B(z)}{P_0(x)},\\
\frac{1-r}{r}&=-\frac{B(z)}{P_0(x)+B(z)}=\frac{D(z)}{Q_0(y)},
\end{aligned}$$
eliminating $r$ yields
\begin{equation}\label{eq:translation-rank-prelim}
P_0(z-y)=R_0(z)+R_1(z)Q_0(y),\text{ where }R_0=-B\text{ and }R_1=-\frac{B}{D}\ne0.
\end{equation}
As functions of $y$, the translates $P_0(z-y)$ therefore lie in the span of $1$ and $Q_0$, a space of dimension at most two. A complex pole of $P_0$ would produce a pole of $P_0(z-y)$ whose location moves with $z$, which is impossible because every function in this fixed space has its poles among the finitely many poles of $Q_0$. Hence, $P_0$ has no finite complex poles and is a polynomial. Expanding its translates in powers of the translation parameter shows that they span a space whose dimension is one more than the degree of $P_0$, so $P_0$ is constant or affine.

The affine case cannot occur, since it would force $f$ to be constant. Indeed, suppose that $P_0$ is nonconstant. Fixing in \eqref{eq:translation-rank-prelim} a regular value of $z$ with $R_1(z)\ne0$ shows that $Q_0$ is a nonconstant affine function as well, so we may write $P_0(x)=a_1x+a_0$ and $Q_0(y)=b_1y+b_0$ with $a_1b_1\ne0$. Comparing the coefficients of $y$ and the constant terms in \eqref{eq:translation-rank-prelim}, we find
$$\begin{aligned}
\frac{B}{D}&=\frac{a_1}{b_1},\\
B(z)&=-a_1z-a_0-\frac{a_1b_0}{b_1},
\end{aligned}$$
and substituting into $r=1+\frac{B}{P_0}$ yields
$$r(x,y)=-\frac{y-y_0}{x-x_0},\text{ where }x_0=-\frac{a_0}{a_1}\text{ and }y_0=-\frac{b_0}{b_1}.$$
Hence, $(x-x_0)f_x+(y-y_0)f_y=0$. By Euler's identity, this operator multiplies each homogeneous part of $f(x+x_0,y+y_0)$ by its degree, so every homogeneous part of positive degree vanishes, contradicting the assumption that $f$ is nonconstant.

Therefore, $P_0$ is constant, and $r-1=b(z)$ for a rational function $b$, which is not identically zero because $p-q\not\equiv0$. In other words, $(G_X)_x=b(z)(G_X)_z$. Choose an interval on which $b$ is finite and nonzero, and let $P_1$ be a primitive of $\frac{1}{b}$ there. In the coordinates $(x+P_1(z),z)$, the derivative in the second coordinate with the first one held fixed is
$$(G_X)_z-P_1'(z)(G_X)_x=(G_X)_z-\frac{(G_X)_x}{b(z)}=0,$$
and consequently $G_X(x,z)=h_0(x+P_1(z))$ locally. It remains to see that this local representation involves polynomials. Fixing $z$ shows that $h_0$ agrees on its local interval with a translate of a polynomial in $x$, and we may use that polynomial for $h_0$. Its degree $m$ is at least one, because $G_X$ is nonconstant. If $c_m\ne0$ and $c_{m-1}$ are its two leading coefficients, then the coefficient of $x^{m-1}$ in $G_X$ equals $mc_mP_1(z)+c_{m-1}$. This coefficient is a polynomial in $z$ and $mc_m\ne0$, so $P_1$ also agrees with a polynomial on the interval. The resulting identity between polynomials holds on an open set and therefore everywhere, and substituting $z=x+y$, we conclude that $f=h_0(x+P_1(x+y))\in\mathfrak E_S$.
\end{proof}

\begin{corollary}\label{cor:three-coordinate}
The three representations $f$, $G_X$ and $G_Y$ are all analytic-special if and only if $f$ is a polynomial in one affine coordinate. Consequently, every admissible polynomial has a nonspecial representation among these three.
\end{corollary}

\begin{proof}
Suppose that all three representations are special. A constant polynomial has the required form, so let $f$ be nonconstant. Theorem~\ref{thm:two-transition}, applied to the two sum-coordinate representations, shows that $f=h_0(ax+P(x+y))$. If $a=0$, if $P'\equiv0$, or if $a+P'\equiv0$, then $f$ is a polynomial in $x+y$, in $x$ or in $y$, respectively, and again has the required form. Otherwise, $a\ne0$, and the numerator and the denominator of
$$\frac{f_x}{f_y}=\frac{a+P'(x+y)}{P'(x+y)}$$
are nonzero polynomials. Since $f$ is special, this rational function of $x+y$ factors as $A(x)B(y)$, so that, as $y$ ranges over regular values, all its translates in $x$ are scalar multiples of $A(x)$. The argument with translated poles from the preceding proof shows that the ratio has no finite poles and is therefore a polynomial. The translates of a nonconstant polynomial span a space of dimension at least two, so the ratio is constant, and since $a\ne0$, the identity $\frac{a+P'}{P'}=1+\frac{a}{P'}$ forces $P'$ to be constant. Hence, $ax+P(x+y)$ is affine, and $f$ has the required form.

Conversely, an invertible shear maps polynomials in one affine coordinate to polynomials in one affine coordinate, so for such $f$, all three representations are special, including when $f$ is constant. The statement about admissible polynomials is the contrapositive.
\end{proof}

\subsection{Frostman bounds and entropy under coordinate changes}
\label{sec:preliminaries}

We now collect the entropy inequalities and the estimates near critical sets that are used throughout the paper. For unequal exponents, conditional $(a,C_1;b,C_2)$ regularity means that $X\mid Y$ is $(a,C_1)$-Frostman and $Y\mid X$ is $(b,C_2)$-Frostman, and joint $(a,b,C)$ regularity means that $\mathbb P(X\in I,Y\in J)\le C|I|^a|J|^b$ whenever $|I|,|J|\le1$. We sometimes suppress fixed constants from the notation.

\subsubsection{Dyadic entropy and sum shears}

The relative entropy of two probability measures is $D(\mu\Vert\nu)=\int\log\frac{d\mu}{d\nu}d\mu$ if $\mu\ll\nu$, and $D(\mu\Vert\nu)=\infty$ otherwise. For probability vectors, it equals $\sum_ip_i\log\frac{p_i}{q_i}$, where the terms with $p_i=0$ are omitted. Given a discrete label $L$, the conditional dyadic mutual information is
$${\mathcal{I}_n}(U;V\mid L)=H(D_nU\mid L)+H(D_nV\mid L)-H(D_nU,D_nV\mid L).$$

Let $Z$ be a bounded random vector in $\mathbb R^k$, let $m\ge0$ be an integer, let $T$ be a fixed invertible linear map and let $\varphi$ be a fixed Lipschitz map. Then
\begin{equation}\label{eq:entropy-refinement}
H_n(Z)\le H_{n+m}(Z)\le H_n(Z)+km,
\end{equation}
\begin{equation}\label{eq:linear-entropy}
H_n(TZ)=H_n(Z)+O_T(1),
\end{equation}
and
\begin{equation}\label{eq:Lipschitz-entropy}
H_n(\varphi(Z))\le H_n(Z)+O_\varphi(1).
\end{equation}
Indeed, the cell at level $n+m$ determines the cell at level $n$, and each cell at level $n$ has exactly $2^{km}$ descendants at level $n+m$, so \eqref{eq:entropy-refinement} follows from the chain rule. The image under $\varphi$ of a cell at level $n$ meets a bounded number of output cells, so the conditional entropy of the output given the input cell is $O_\varphi(1)$, which proves \eqref{eq:Lipschitz-entropy}, and applying this to $T$ and $T^{-1}$, we obtain \eqref{eq:linear-entropy}. In the same way, if two random vectors can be coupled so that they differ by $O(2^{-n})$, then each cell at level $n$ of one of them leaves only boundedly many possible cells for the other, and their entropies at level $n$ differ by $O(1)$.

\begin{lemma}\label{lem:Frostman-entropy}
An $(s,C)$-Frostman random variable satisfies $H_n(X)\ge sn-\log C$, and joint $(s_1,s_2,C)$ regularity implies $H_n(X,Y)\ge(s_1+s_2)n-\log C$. Under conditional $(s_1,C_1;s_2,C_2)$ regularity,
$$\begin{aligned}
H(D_nX\mid D_nY)&\ge s_1n-\log C_1,\\
H(D_nY\mid D_nX)&\ge s_2n-\log C_2.
\end{aligned}$$
\end{lemma}

\begin{proof}
If all atoms of a discrete law have probability at most $p_*$, then
$$H=-\sum_ip_i\log p_i\ge-\sum_ip_i\log p_*=-\log p_*.$$
A dyadic interval at level $n$ has length $2^{-n}$, and a dyadic square at level $n$ has both sides of length $2^{-n}$, so the Frostman hypotheses allow us to take $p_*=C2^{-sn}$ and $p_*=C2^{-(s_1+s_2)n}$, respectively, which proves the first two assertions. For the first conditional assertion, fix a bin of $D_nY$ of positive probability. The law of $D_nX$ conditioned on this bin is a mixture of the laws conditioned on the exact values $Y=y$, and for almost every $y$, the atoms of the latter have masses at most $C_12^{-s_1n}$. Integrating over the conditioning bin preserves this bound, so the conditional entropy given the bin is at least $s_1n-\log C_1$, and averaging over the bins of $D_nY$ proves the claim. The second conditional assertion follows by interchanging $X$ and $Y$.
\end{proof}

For finite-valued random variables $X,Y,Z$ and a finite or countable label $L$ with finite entropy, we have
\begin{equation}\label{eq:label-entropy}
H(Z\mid L)\le H(Z)\le H(Z\mid L)+H(L),
\end{equation}
\begin{equation}\label{eq:label-MI}
{\mathcal{I}}(X;Y\mid L)\le {\mathcal{I}}(X;Y)+H(L),
\end{equation}
and
\begin{equation}\label{eq:mixture-concavity}
H(Z)\ge\sum_j\mathbb P(L=j)H(Z\mid L=j).
\end{equation}
The entropy bounds follow from the chain rule and the nonnegativity of conditional entropy, and the chain rule also shows that
$${\mathcal{I}}(X;Y\mid L)-{\mathcal{I}}(X;Y)={\mathcal{I}}(X;L\mid Y)-{\mathcal{I}}(X;L)\le H(L),$$
which is \eqref{eq:label-MI}. For a countable label, we first merge its tail into a single value and then pass to the limit, which is justified because $X$, $Y$, $Z$ and $L$ have finite entropy. We shall often use the fact that the average of the restriction cost $\log\frac{1}{\mathbb P(L=j)}$ is exactly $H(L)$.

\begin{lemma}\label{lem:conditional-shear}
If $(X,Y)$ is conditionally $(s_1,C_1;s_2,C_2)$-Frostman and $S=X+Y$, then $(X,S)$ is jointly $(s_1,s_2,C_1C_2)$-Frostman and $(Y,S)$ is jointly $(s_2,s_1,C_1C_2)$-Frostman. Moreover,
\begin{equation}\label{eq:shear-cond-X}
H(D_nS\mid D_nX)\ge s_2n-O(1)
\end{equation}
and
\begin{equation}\label{eq:shear-cond-Y}
H(D_nS\mid D_nY)\ge s_1n-O(1).
\end{equation}
\end{lemma}

\begin{proof}
Integrating the exact conditional bounds shows that $X$ is $(s_1,C_1)$-Frostman and $Y$ is $(s_2,C_2)$-Frostman. To prove the rectangle bound for $(X,S)$, we condition on $X$ first. For intervals $I,J$ of length at most one,
$$\mathbb P(X\in I,S\in J)=\int_I\mathbb P(Y\in J-x\mid X=x)d\mathcal L(X)(x)\le C_1C_2|I|^{s_1}|J|^{s_2},$$
because the translate $J-x$ has the same length as $J$, so that the bound for $Y\mid X$ controls the integrand for almost every $x$, and the marginal bound for $X$ then controls the integral over $I$. The same translation applied to a dyadic interval $J$ at level $n$ shows that $\mathbb P(S\in J\mid X=x)\le C_22^{-s_2n}$. This atom bound survives averaging over a bin of $D_nX$ of positive probability, so the conditional entropy of $D_nS$ given such a bin is at least $s_2n-\log C_2$, and averaging over the bins proves \eqref{eq:shear-cond-X}. Interchanging the roles of $X$ and $Y$, together with those of $(s_1,C_1)$ and $(s_2,C_2)$, proves the assertion for $(Y,S)$ and \eqref{eq:shear-cond-Y}.
\end{proof}

\subsubsection{Polynomial sublevel bounds and entropy recovery}

\begin{lemma}\label{lem:polynomial-sublevel}
Let $(X,Y)$ be a bounded conditionally $(s_1,C_1;s_2,C_2)$-Frostman pair with $s_1,s_2>0$, and let $P$ be a fixed nonzero polynomial. There are constants $A<\infty$ and $\theta>0$, uniform over such laws, for which
\begin{equation}\label{eq:polynomial-sublevel}
\mathbb P(|P(X,Y)|\le r)\le Ar^\theta\text{ for }0<r\le1.
\end{equation}
\end{lemma}

\begin{proof}
The idea is to factor $P$ in $x$ on the fibers where its leading coefficient is not too small and to control the remaining fibers by the marginal Frostman bound in $y$. A nonzero constant polynomial has an empty sublevel set for small $r$, and the assertion follows after increasing $A$. If $P$ depends only on $y$ and has degree $e\ge1$, then factoring over $\mathbb C$ shows that $|P(y)|\le r$ forces one of the root factors to have modulus $O_P(r^{\frac{1}{e}})$. Projecting the roots to the real axis, we cover the sublevel set by at most $e$ real intervals of this length, and the marginal $(s_2,C_2)$-Frostman bound for $Y$ proves the claim for small $r$, while increasing $A$ takes care of the remaining values.

Suppose now that $P$ has degree $d\ge1$ in $x$, and let $a(y)$ be its leading coefficient, a nonzero polynomial. We split the fibers according to whether $|a(y)|\le r^{\frac{1}{2}}$, a threshold that leaves a positive power of $r$ in both of the estimates below. If $|a(y)|>r^{\frac{1}{2}}$ and $|P(x,y)|\le r$, the factorization of $P(\cdot,y)$ shows that one of the factors has modulus at most $\left(\frac{r}{|a(y)|}\right)^{\frac{1}{d}}\le r^{\frac{1}{2d}}$. Projecting the corresponding complex roots to the real axis, we cover the sublevel set of the fiber by at most $d$ real intervals of length $O(r^{\frac{1}{2d}})$, which have length at most one when $r$ is small. Applying the exact $(s_1,C_1)$-Frostman bound to $X\mid Y=y$ and integrating in $y$, we obtain
$$\mathbb P(|P(X,Y)|\le r)\le O_P(C_1r^{\frac{s_1}{2d}})+\mathbb P(|a(Y)|\le r^{\frac{1}{2}}).$$
If $a$ is nonconstant of degree $e$, the univariate argument above covers the exceptional set $\{|a(y)|\le r^{\frac{1}{2}}\}$ by at most $e$ intervals of length $O_P(r^{\frac{1}{2e}})$, whose probability is $O_P(C_2r^{\frac{s_2}{2e}})$, so we may take $\theta=\min\{\frac{s_1}{2d},\frac{s_2}{2e}\}$. If $a$ is a nonzero constant, the exceptional set is empty for small $r$ and $\theta=\frac{s_1}{2d}$ suffices. In both cases, increasing $A$ covers the remaining values $r\le1$, and all choices depend only on the polynomial and the Frostman data.
\end{proof}

Let $P$ be a fixed nonzero polynomial with zero set $Z(P)=\{(x,y):P(x,y)=0\}$. Inside a fixed bounded box, after a subdivision according to the direction of the tangent, $Z(P)$ splits into finitely many points and finitely many Lipschitz graphs with bounded coordinate ranges and Lipschitz constants. The $r$-neighborhood of each graph is covered by $O_P(r^{-1})$ squares of side $O_P(r)$, and under joint $(s,s,C)$ regularity, each of these squares has mass $O_{P,s,C}(r^{2s})$. Summing over the squares, we find for $s>\frac{1}{2}$ that
\begin{equation}\label{eq:joint-tube}
\mathbb P(\operatorname{dist}((X,Y),Z(P))\le r)\lesssim_P r^{2s-1}.
\end{equation}

\begin{theorem}\label{thm:coordinate-entropy}
Let $\mathsf L(x,y)=ax+by\ne0$ and let $P\in\mathbb R[x,y]$ satisfy
\begin{equation}\label{eq:coordinate-J}
J_{\mathsf L,P}=aP_y-bP_x\not\equiv0.
\end{equation}
Every bounded conditionally $(s_1,C_1;s_2,C_2)$-Frostman pair with positive exponents satisfies
\begin{equation}\label{eq:coordinate-recovery}
|H_n(X,Y)-H_n(\mathsf L(X,Y),P(X,Y))|\le O(1),
\end{equation}
where the implied constant depends only on the fixed data and not on the law or on $n$. Consequently,
\begin{equation}\label{eq:coordinate-baseline}
\max\{H_n(\mathsf L(X,Y)),H_n(P(X,Y))\}\ge\frac{(s_1+s_2)n}{2}-O(1).
\end{equation}
The comparison \eqref{eq:coordinate-recovery} also holds for every bounded pair satisfying
\begin{equation}\label{eq:coordinate-isotropic-hypothesis}
\mathbb P((X,Y)\in B)\le Ch^t\text{ for every axis-parallel square }B\text{ of side }0<h\le1,
\end{equation}
where $1<t\le2$. In this case, the maximum of the two output entropies is at least $\frac{tn}{2}-O(1)$, with constants depending only on $t$, $C$, $\mathsf L$, $P$ and a fixed support box. For $\mathsf L=x+y$, the excluded case $J_{\mathsf L,P}\equiv0$ occurs exactly when $P(x,y)=q(x+y)$ for a real univariate polynomial $q$.
\end{theorem}

\begin{proof}
By the Lipschitz bound \eqref{eq:Lipschitz-entropy}, $H_n(\mathsf L(X,Y),P(X,Y))\le H_n(X,Y)+O(1)$. By the chain rule, the reverse inequality amounts to a bound for the conditional entropy of the input cells given the output cells. We shall count the input cells on which the Jacobian has a prescribed size and then average the logarithm of this count with the help of a polynomial sublevel bound.

We use the fixed linear change of coordinates $(x,y)\mapsto(u,v)=(ax+by,-bx+ay)$, whose determinant is $a^2+b^2>0$, write $\widetilde P$ for $P$ in these coordinates, and denote the transformed random pair by $(U,V)$. The relevant map is then
$$\mathcal T(u,v)=(u,\widetilde P(u,v)),\text{ with }\det D\mathcal T=\partial_v\widetilde P\not\equiv0,$$
because $\partial_v\widetilde P(u,v)=\frac{J_{\mathsf L,P}(x,y)}{a^2+b^2}$ at corresponding points. Since $U=\mathsf L(X,Y)$, we have $\mathcal T(U,V)=(\mathsf L(X,Y),P(X,Y))$, and since the linear change of input coordinates preserves entropy up to $O(1)$ by \eqref{eq:linear-entropy}, it suffices to prove that
$$H(D_n(U,V)\mid D_n\mathcal T(U,V))=O(1).$$
All geometric bounds below are taken on a fixed neighborhood of the transformed support that contains every input cell we use.

The count of cells relies on the area formula, which requires a bound on the multiplicity of $\mathcal T$ almost everywhere. Let $d=\deg_v\widetilde P\ge1$. For every $u$ at which the coefficients of the positive powers of $v$ in $\widetilde P(u,v)$ do not all vanish, each target value $w$ leads to an equation $\widetilde P(u,v)=w$ with at most $d$ solutions. At least one of these coefficients is a nonzero polynomial in $u$, so the exceptional values of $u$ form a finite set, and $\mathcal T$ has multiplicity at most $d$ outside finitely many vertical lines in the target, which is all the area formula needs.

To average the possible loss in inverting $\mathcal T$, we record the size of the Jacobian until it reaches the observation scale $2^{-n}$. Set $K=0$ when $|\partial_v\widetilde P(U,V)|>\frac{1}{2}$, set $K=k$ when
$$2^{-k-1}<|\partial_v\widetilde P(U,V)|\le2^{-k}\text{ and }1\le k<n,$$
and set $K=n$ when $|\partial_v\widetilde P(U,V)|\le2^{-n}$. These regions cover the source, including the zero set of the Jacobian, which belongs to the last region. We next prove a tail bound for $K$ under either hypothesis on the source, using that the Jacobian at $(U,V)$ is a fixed nonzero multiple of $J_{\mathsf L,P}(X,Y)$. Under the conditional hypothesis, we simply apply Lemma~\ref{lem:polynomial-sublevel} to this nonzero polynomial in the original variables. Under \eqref{eq:coordinate-isotropic-hypothesis}, suppose first that the zero set of $J_{\mathsf L,P}$ meets the fixed neighborhood of the original support. The covering argument behind \eqref{eq:joint-tube}, which now assigns mass $O(r^t)$ to each of $O(r^{-1})$ squares, yields
$$\mathbb P(\operatorname{dist}((X,Y),Z(J_{\mathsf L,P}))\le r)\lesssim_{\mathsf L,P,t,C}r^{t-1}$$
on this neighborhood. The functions $|J_{\mathsf L,P}|$ and $\min\{1,\operatorname{dist}(\cdot,Z(J_{\mathsf L,P}))\}$ are continuous and semialgebraic on a compact box and have the same zero set there, so the the polynomial {\L}ojasiewicz inequality \cite[Corollary~2.6.7]{BCR} provides fixed $c,m>0$ with
$$|J_{\mathsf L,P}(x,y)|\ge c\min\{1,\operatorname{dist}((x,y),Z(J_{\mathsf L,P}))\}^m$$
on the box. Hence, $\mathbb P(|J_{\mathsf L,P}(X,Y)|\le\tau)\lesssim\tau^{\frac{t-1}{m}}$ for all sufficiently small $\tau>0$. If the polynomial has no zeros near the support, then its absolute value has a positive minimum there, and the same conclusion holds after adjusting the constants. Under either hypothesis, we therefore find fixed $A<\infty$ and $\theta>0$ with
\begin{equation}\label{eq:Jacobian-shell-tail}
\mathbb P(K\ge k)\le A2^{-\theta k}\text{ for }1\le k\le n.
\end{equation}
Summing this tail shows that $\mathbb EK=O(1)$ uniformly in $n$, and comparing the law of $K$ with the geometric probabilities $2^{-k-1}$, using the nonnegativity of relative entropy, we get $H(K)\le1+\mathbb EK=O(1)$. These are the two quantities that enter when we average the counts of cells. Note that in both cases, the sublevel estimates were applied in the original variables, so we never needed a Frostman hypothesis after the linear change of coordinates.

It remains to count the input cells lying over a fixed dyadic output square. Fix $k<n$, and let $\mathcal Q$ be the set of input squares that contain a point of the shell $K=k$ whose image lies in that output square. We claim that $|\mathcal Q|=O(2^k)$. Since the derivatives are bounded, the Jacobian varies by $O(2^{-n})$ on each input square, and we choose a fixed integer $C_0$ so that this variation is less than half of the lower bound $2^{-k-1}$ at the selected point whenever $k\le n-C_0$. In this range, $|\det D\mathcal T|\gtrsim2^{-k}$ throughout every square of $\mathcal Q$, and the Lipschitz bound places the images of these squares in a fixed enlargement of the output square, so the area formula and the multiplicity bound yield
$$|\mathcal Q|2^{-2n-k}\lesssim\int_{\bigcup_{B\in\mathcal Q}B}|\det D\mathcal T|\lesssim d2^{-2n},$$
which proves the claim. In the remaining range $n-C_0<k<n$, the first output coordinate confines $u$ to $O(1)$ columns of dyadic squares, and since the support is bounded and $n-k<C_0$, these columns contain $O(2^n)=O(2^k)$ input squares. On the last shell $K=n$, we use the trivial count $O(2^{2n})$, since the probability of this shell is already controlled by \eqref{eq:Jacobian-shell-tail}.

Entropy is at most the logarithm of the size of the support, so conditioned on the output cell and on $K=k$, the entropy of the input cell is at most $k+O(1)$ for $k<n$ and at most $2n+O(1)$ for $k=n$. The chain rule and averaging over $K$ lead to
$$H(D_n(U,V)\mid D_n\mathcal T(U,V))\le H(K)+\mathbb E[K1_{\{K<n\}}]+2n\mathbb P(K=n)+O(1)=O(1),$$
where the last term is bounded because $n\mathbb P(K=n)\le An2^{-\theta n}=O(1)$. This proves \eqref{eq:coordinate-recovery} under both hypotheses. For \eqref{eq:coordinate-baseline}, we integrate the two exact conditional Frostman bounds to obtain the joint atom bound with constant $C_1C_2$, so that by subadditivity
$$(s_1+s_2)n-O(1)\le H_n(X,Y)\le H_n(\mathsf L(X,Y))+H_n(P(X,Y))+O(1),$$
and \eqref{eq:coordinate-baseline} follows. Under \eqref{eq:coordinate-isotropic-hypothesis}, the input atom bound implies $H_n(X,Y)\ge tn-\log C$ instead, and the output maximum is at least $\frac{tn}{2}-O(1)$. Finally, for $\mathsf L=x+y$, the identity $P_y-P_x=0$ says that $P$ does not depend on $x-y$ in the coordinates $(x+y,x-y)$, so that $P$ is a polynomial $q(x+y)$, and the converse follows by differentiation.
\end{proof}

\subsubsection{Frostman bounds for regular outputs}

By a buffered chart we mean a compact core contained in a larger chart on which the same regularity bounds hold, so that sufficiently fine dyadic cells meeting the core remain inside the larger chart.

\begin{lemma}\label{lem:regular-output}
Let $g$ be smooth near a fixed product rectangle $B$, with $|g_y|$ bounded away from zero, and let $E\subseteq\{(X,Y)\in B\}$ be an event of probability $w>0$. For intervals $J$ with $|J|\le1$, conditional $(s,C_1;s,C_2)$ regularity implies
\begin{equation}\label{eq:localized-output-cap}
\mathbb P(g(X,Y)\in J\mid E)\lesssim_{g,B}C_2w^{-1}|J|^s,
\end{equation}
and joint $(s,s,C)$ regularity with $s>\frac{1}{2}$ gives
\begin{equation}\label{eq:localized-output-cap-joint}
\mathbb P(g(X,Y)\in J\mid E)\lesssim_{g,B}Cw^{-1}|J|^{2s-1}.
\end{equation}
If $|g_x|$, rather than $|g_y|$, is bounded away from zero, then \eqref{eq:localized-output-cap} holds with $C_1$ in place of $C_2$.
\end{lemma}

\begin{proof}
We first bound, under the original law, the mass of the event that $(X,Y)\in B$ and $g(X,Y)\in J$. Intersecting with $E$ can only decrease this mass, and dividing by $w$, we then obtain the asserted bounds. Proceeding in this order allows us to use the exact conditional hypotheses for the original law, without assuming that they survive the restriction to $E$.

Let $c=\inf_B|g_y|>0$. Since $B$ is connected, $g_y$ has a fixed sign, and for each fixed $x$, the mean value theorem places the set of $y$ with $(x,y)\in B$ and $g(x,y)\in J$ in an interval of length at most $c^{-1}|J|$. Integrating the exact conditional Frostman bound for $Y\mid X$, we obtain the unnormalized estimate $O_{g,B}(C_2|J|^s)$. If the interval is longer than one, we cover it by boundedly many intervals of length at most one, and the resulting constant is absorbed because $|J|>c$ in this case.

For the joint assertion, the available information concerns rectangles rather than conditional fibers. Let $r=|J|>0$ and divide the $x$-range of $B$ into $O_B(r^{-1})$ strips of width at most $r$. If two points $(x,y)$ and $(x',y')$ of the same strip both have output values in $J$, then
$$c|y-y'|\le|g(x,y)-g(x,y')|\le|J|+\|g_x\|_{L^\infty(B)}|x-x'|\lesssim_{g,B}r,$$
so the relevant part of each strip is covered by $O_{g,B}(1)$ squares of side $r$, each of mass at most $Cr^{2s}$. Summing over the strips leads to the unnormalized bound $O_{g,B}(Cr^{2s-1})$. Restricting both unnormalized estimates to $E$ and dividing by $w$ proves the displayed bounds, and intervals of length zero are handled by shrinking intervals of positive length to a point, since both exponents are positive. Finally, if $|g_x|$ is bounded away from zero instead, we interchange the coordinates and use the bound for $X\mid Y$, whose constant is $C_1$.
\end{proof}

\section{Incidence estimates and conditional copies}
\label{sec:incidence-copies}

A single weighted incidence calculation underlies both the joint and the conditional arguments. We state it once and then specialize its geometry and its information terms to the copied variables needed later. Throughout the paper, for $0<s<1$, we write
$$t_s=\min\{s,2-2s\}\text{ and }\beta_s=s+\frac{t_s}{2}=\min\left\{\frac{3s}{2},1\right\}.$$

\subsection{Rectangular incidences and algebraic concentration}
\label{sec:joint-incidence-inputs}

We first convert a rectangular incidence estimate into an entropy inequality for conditional copies, and then we control the concentration of the law near an algebraic curve in terms of the marginal entropies.

\subsubsection{Entropy bounds from rectangular incidences}

A family $\mathcal G\subset C^2([-2,2])$ is $T$-transversal if
$$\|g-h\|_{C^2}\le T\inf_t(|g(t)-h(t)|+|g'(t)-h'(t)|)\text{ for all }g,h\in\mathcal G,$$
where $\|g\|_{C^2}=\sup_{t\in[-2,2]}\sum_{k=0}^2|g^{(k)}(t)|$, and the separation of curves is measured in the same norm. For fixed $t$, the value--slope map sends $g$ to the pair $(g(t),g'(t))$. We use the following consequence of the incidence theorem of O'Regan, Wu and Yi \cite[Definition~4.3 and Theorem~4.4, pp.~29--30]{OWY}, in which the cardinality conditions below imply the covering number conditions required there.

\begin{theorem}\label{thm:rectangular-incidence}
Let $0<s<1$, $T\ge1$ and $\eta>0$, and let $0<\delta<\delta_0(s,T,\eta)$. Let $\mathcal G$ be a finite, $\delta$-separated, $T$-transversal family in the unit ball of $C^2([-2,2])$. For each $g\in\mathcal G$, let $\mathcal P(g)$ be a set of dyadic $\delta$-squares in $[-1,1]^2$ that meet the graph of $g$, and set $\mathcal P=\bigcup_g\mathcal P(g)$. Suppose that $C_1,C_2\ge1$ and that
$$\begin{aligned}
|\{g:(g(t),g'(t))\in B\}|&\le C_1\left(\frac{\sqrt{ab}}{\delta}\right)^{2s},\\
|\{p\in\mathcal P(g):p\cap B_a\ne\varnothing\}|&\le C_2\left(\frac{a}{\delta}\right)^\sigma,\text{ where }\sigma=\min\{2s,2-2s\},
\end{aligned}$$
for every $t$, every $g\in\mathcal G$, every axis-parallel rectangle $B$ with sides $a,b\ge\delta$, and every ball $B_a$ of radius $a\ge\delta$. Then
\begin{equation}\label{eq:rectangular-incidence}
\sum_g|\mathcal P(g)|\le C_1^{\frac{2}{3}}C_2^{\frac{1}{3}}\delta^{-\frac{2s}{3}-\eta}|\mathcal P|^{\frac{2}{3}}|\mathcal G|^{\frac{1}{3}}.
\end{equation}
The threshold $\delta_0$ does not depend on $C_1$ and $C_2$.
\end{theorem}

Since the threshold does not depend on $C_1$ and $C_2$, these constants may depend on the working scale and on the probability class under consideration. Fixed affine normalizations and fixed enlargements of the squares only change the fixed factors in the hypotheses. The counting hypotheses imply the covering bounds of \cite{OWY} directly, because a finite set of value--slope pairs is covered by one $\delta$-ball per point and each dyadic $\delta$-square fits inside a ball of radius $\delta$, so that counting the squares that meet a ball also bounds the covering number of their union inside that ball.

\begin{lemma}\label{lem:copied-output-incidence}
Fix $0<s<1$ and let $\delta=2^{-n}$. Let $U$ and $V$ take values in a $\delta$-grid in fixed bounded intervals, and let $V$ and $V'$ be independent conditionally on $U$, with the same conditional law. Suppose that $\mathbb P(V\in I)\le C_V|I|^s$ and $\mathbb P(U\in I\mid V=v)\le K_v|I|^s$ for every interval $I$ with $\delta\le|I|\le1$, where $C_V\ge1$ and $K_v\ge1$. Let $\Phi$ and $\Psi$ be fixed smooth functions, and set
$$W=(\Phi(U,V'),\Psi(U,V)).$$
Assume that $|\Phi_u|$ is bounded below on fixed product boxes, and that the curves
$$g_{v,z}(t)=\Psi(u(t,z),v),\text{ where }\Phi(u(t,z),z)=t,$$
have bounded $C^2$ norms and uniformly nonsingular value--slope maps $(v,z)\mapsto(g_{v,z}(t),g'_{v,z}(t))$. Assume also that each family of curves has a common buffered interval of values of $t$ containing its output values, on which $u(t,z)$ is defined for every parameter $z$ in the buffer, and that after a subdivision of the boxes indexed by a fixed finite label these maps are uniformly bi-Lipschitz and the families are uniformly transversal. Then, for every $\varepsilon>0$,
\begin{equation}\label{eq:copied-output-incidence}
\begin{aligned}
H_n(W)\ge{}&H(U)-\frac{5}{2}{\mathcal{I}}(U;V)+\frac{1}{2}{\mathcal{I}}(V;V')+\beta_sn\\
&-\frac{1}{2}\mathbb E\log K_V-\varepsilon n-O_\varepsilon(1+\log C_V),
\end{aligned}
\end{equation}
where the constants depend on the fixed geometry but not on $C_V$ or on the $K_v$.
\end{lemma}

\begin{proof}
We apply Theorem~\ref{thm:rectangular-incidence} on classes of copied triples whose relevant probabilities are comparable. The regularity required by the incidence theorem comes from reference laws, and the ratios between actual and reference probabilities measure the entropy lost in the process.

Write $p_{ij}=\mathbb P(U=u_i,V=v_j)$, and let $p_i$ and $q_j$ be the marginal probabilities. For indices of positive probability, we define the law of the copied triple and the reference probabilities by
$$\begin{aligned}
p_{ijj'}&=\frac{p_{ij}p_{ij'}}{p_i}, & \nu_{jj'}&=\sum_ip_{ijj'},\\
\pi_{i\mid jj'}&=\frac{p_{ijj'}}{\nu_{jj'}}, & a_{i\mid j}&=\frac{p_{ij}}{q_j},
\end{aligned}$$
so that $\nu$ is the law of $(V,V')$, $\pi$ is the conditional law of $U$ given this pair, and $a$ is the conditional law of $U$ given $V$. These laws are formed before any class labels are attached, and all information quantities in \eqref{eq:joint-copy-identities} below refer to the copied law and its marginals. The restrictions to classes introduced later need not preserve conditional independence, but averaging the class inequalities brings us back to the information quantities of the copied law, to which the identities apply. Conditional independence and the chain rule imply
\begin{equation}\label{eq:joint-copy-identities}
\begin{aligned}
H(U\mid V,V')&=H(U)-2{\mathcal{I}}(U;V)+{\mathcal{I}}(V;V'),\\
D(\nu\Vert q\otimes q)&={\mathcal{I}}(V;V'),\\
\mathbb ED(\pi_{\cdot\mid VV'}\Vert a_{\cdot\mid V})&={\mathcal{I}}(U;V)-{\mathcal{I}}(V;V').
\end{aligned}
\end{equation}
For the last identity, observe that its left-hand side equals ${\mathcal{I}}(U;V'\mid V)$, which the chain rule for ${\mathcal{I}}((U,V);V')$ identifies with ${\mathcal{I}}(U;V)-{\mathcal{I}}(V;V')$. In particular, $0\le {\mathcal{I}}(V;V')\le {\mathcal{I}}(U;V)$. If $\mathbb E\log K_V\ge2H(U)+2\beta_sn$, then the right-hand side of \eqref{eq:copied-output-incidence} is nonpositive and there is nothing to prove, so we may assume that $\mathbb E\log K_V=O(n+1)$, since $H(U)\le n+O(1)$ on the fixed grid.

Before attaching probability labels, we check the geometric regularity of the product reference $q\otimes q$. On the fixed small parameter boxes, the inverses of the value--slope maps are uniformly Lipschitz. Consider an axis-parallel rectangle in the value--slope plane with sides $r_1$ and $r_2$. The component of the map in the direction of the shorter side has a partial derivative of fixed sign bounded away from zero, so the fibers of the preimage have length $O(\min\{r_1,r_2\})$, while the inverse Lipschitz bound confines the projection of the preimage to the other parameter coordinate to an interval of length $O(\max\{r_1,r_2\})$. Integrating the two marginal bounds, we see that $q\otimes q$ assigns mass $O(C_V^2(r_1r_2)^s)$ to the preimage. The same inverse Lipschitz bound shows that the curves are separated in the $C^2$ norm once both indices of the parameter grid are colored by finitely many residue classes. Fixed normalizations place the curves in the unit ball of $C^2$ on the required common interval, transversality is part of the geometric hypotheses, and the labels recording the charts and the colors have bounded entropy.

The remaining labels make the masses of curves, fibers and outputs, and their ratios to the references, comparable within each class. We partition the copied triples according to the integer parts of the six quantities
$$\begin{aligned}
&-\log\nu_{jj'}, & &-\log\pi_{i\mid jj'}, & &-\log\mathbb P(D_nW=D_nW(i,j,j')),\\
&\log\frac{\nu_{jj'}}{q_jq_{j'}}, & &\log\frac{\pi_{i\mid jj'}}{a_{i\mid j}}, & &\log K_{v_j},
\end{aligned}$$
where $W(i,j,j')=(\Phi(u_i,v_{j'}),\Psi(u_i,v_j))$, and we denote the first five levels by $b_\nu$, $c$, $\zeta$, $h$ and $\lambda$. Let $L$ be the complete label, which also records the charts and colors, and write $w_\ell=\mathbb P(L=\ell)$ under the copied law. Since the integer parts change the quantities by bounded amounts, the expectations of the levels under the copied law are
$$\begin{aligned}
\mathbb Ec&=H(U\mid V,V')+O(1), & \mathbb E\zeta&=H_n(W)+O(1),\\
\mathbb Eh&={\mathcal{I}}(V;V')+O(1), & \mathbb E\lambda&={\mathcal{I}}(U;V)-{\mathcal{I}}(V;V')+O(1),\\
\mathbb E\log K_V&=\sum_jq_j\log K_{v_j}, & \mathbb E\log w_L&=-H(L).
\end{aligned}$$
The level $b_\nu$ of the curve masses will cancel when the incidence inequality is solved for $\zeta$.

Consider a class $L=\ell$ and write $w=w_\ell$. Each triple in the class has mass $O(2^{-b_\nu-c})$, and since $|\Phi_u|$ is bounded below, only boundedly many input grid points produce the same pair of a curve and a square. Hence, the class contains at least $c_0w2^{b_\nu+c}$ incidences between curves and squares, at most $C_02^{b_\nu}$ curves, and at most $C_02^\zeta$ squares. The reference masses of the parameters and of the fibers are comparable to $2^{-b_\nu-h}$ and $2^{-c-\lambda}$, respectively, and the last label makes $K_{v_j}$ comparable within the class, so a single representative value can be used in the incidence bound. The reference estimates above therefore provide the incidence constants
$$C_1=C_0C_V^22^{b_\nu+h}\delta^{2s}\text{ and }C_2=C_0K_{v_j}2^{c+\lambda}\delta^{t_s}.$$
For the second of these, the lower bound for $|\Phi_u|$ confines the input points whose outputs lie in a ball of radius $r$ to an interval of length $O(r)$, and the exponent required by the incidence theorem follows from the inequality $r^s\lesssim\delta^{t_s}\left(\frac{r}{\delta}\right)^\sigma$, valid for $\delta\le r\le O(1)$ with $\sigma=\min\{2s,2-2s\}$. The bounds $q_jq_{j'}\lesssim C_V^2\delta^{2s}$ and $a_{i\mid j}\lesssim K_{v_j}\delta^s$ for the selected atoms show that $C_1,C_2\ge1$ after increasing $C_0$.

We now apply Theorem~\ref{thm:rectangular-incidence} to this family of curves and its selected output squares. Taking logarithms, we get
$$\log w+b_\nu+c\le\frac{2}{3}\log C_1+\frac{1}{3}\log C_2+\left(\frac{2s}{3}+\eta\right)n+\frac{2}{3}\zeta+\frac{1}{3}b_\nu+O(1),$$
and substituting the values of $C_1$ and $C_2$, multiplying by $\frac{3}{2}$, using $s+\frac{t_s}{2}=\beta_s$ and choosing $\frac{3\eta}{2}\le\varepsilon$, we arrive at
\begin{equation}\label{eq:joint-class-incidence}
\zeta\ge c+\beta_sn-h-\frac{\lambda}{2}-\frac{1}{2}\log K_{v_j}+\frac{3}{2}\log w-2\log C_V-\varepsilon n-O_\varepsilon(1).
\end{equation}
The incidence threshold does not depend on the reference constants, and the values of $n$ below this fixed threshold are absorbed into the error term.

It remains to average \eqref{eq:joint-class-incidence} without losing a fixed proportion of $n$, and for this, it suffices to bound the entropy of the complete label $L$. The nonnegative levels $b_\nu$, $c$ and $\zeta$ have means $O(n+1)$, and so does $\log K_V$ by our initial reduction. For the signed levels, we use the inequalities
$$\begin{aligned}
\mathbb E_{\mu_1}\left[\left(-\log\frac{d\mu_1}{d\mu_2}\right)_+\right]&\le\frac{1}{\ln2},\\
\mathbb E_{\mu_1}\left|\log\frac{d\mu_1}{d\mu_2}\right|&\le D(\mu_1\Vert\mu_2)+\frac{2}{\ln2},
\end{aligned}$$
valid for probability laws $\mu_1\ll\mu_2$. The relevant divergences are ${\mathcal{I}}(V;V')$ and ${\mathcal{I}}(U;V)-{\mathcal{I}}(V;V')$, both of size $O(n+1)$. An integer-valued random variable whose absolute value has mean $m$ has entropy $O(\log(m+2))$, so the complete label has entropy $O(\log(n+2))$, and the term $\frac{3}{2}\mathbb E\log w_L=-\frac{3}{2}H(L)$ costs only this logarithm. By \eqref{eq:joint-copy-identities}, the remaining information terms average to
$$\begin{aligned}
\mathbb E\left(c-h-\frac{\lambda}{2}\right)&=H(U\mid V,V')-{\mathcal{I}}(V;V')-\frac{1}{2}({\mathcal{I}}(U;V)-{\mathcal{I}}(V;V'))+O(1)\\
&=H(U)-\frac{5}{2}{\mathcal{I}}(U;V)+\frac{1}{2}{\mathcal{I}}(V;V')+O(1).
\end{aligned}$$
Therefore,
$$\begin{aligned}
H_n(W)\ge{}&H(U)-\frac{5}{2}{\mathcal{I}}(U;V)+\frac{1}{2}{\mathcal{I}}(V;V')+\beta_sn\\
&-\frac{1}{2}\mathbb E\log K_V-\varepsilon n-O_\varepsilon(\log(n+2)+\log C_V).
\end{aligned}$$
Running the incidence argument with a smaller rate tolerance leaves room to absorb $O_\varepsilon(\log(n+2))$ into the loss $\varepsilon n$ and a constant, and this proves \eqref{eq:copied-output-incidence} with constants independent of $C_V$ and of the $K_v$.
\end{proof}

\begin{proof}[Proof of Theorem~\ref{thm:independent-polynomial-pair}]
We apply Lemma~\ref{lem:copied-output-incidence} on fixed product rectangles on which its family of curves is regular. Independence removes both mutual information terms, and it is preserved when we restrict the law to such a rectangle. We first identify a nonzero polynomial that controls the regularity, and then we choose the rectangles so that the discarded mass and the cost of averaging are uniform in the input laws.

Set
$$\rho=\frac{P_xQ_y}{P_yQ_x}\text{ and }A=P_xP_yQ_xQ_y.$$
By hypothesis, all four first partial derivatives of $P$ and $Q$ are nonzero polynomials, and the derivatives of $\rho$ will decide which input we copy. To express the corresponding nondegeneracy conditions in polynomial form, we define
$$\begin{aligned}
N_x&=P_xP_y(Q_xQ_{xy}-Q_yQ_{xx})+Q_xQ_y(P_yP_{xx}-P_xP_{xy}),\\
N_y&=P_xP_y(Q_yQ_{xy}-Q_xQ_{yy})+Q_xQ_y(P_xP_{yy}-P_yP_{xy}),
\end{aligned}$$
so that, as identities between rational functions,
$$N_x=A\frac{\partial_x\rho}{\rho}\text{ and }N_y=-A\frac{\partial_y\rho}{\rho}.$$
A nonconstant rational function in characteristic zero has a nonzero partial derivative, so at least one of $N_x$ and $N_y$ is nonzero. If $N_x\equiv0$, we interchange $X$ and $Y$ and replace $P(x,y)$ and $Q(x,y)$ by $P(y,x)$ and $Q(y,x)$. The new $N_x$ is the old $N_y$ with its variables interchanged, hence nonzero, and this change preserves both output laws and all the hypotheses, with $C_1$ and $C_2$ interchanged. We may therefore assume that $N_x\not\equiv0$.

Next, we verify the geometric hypotheses of Lemma~\ref{lem:copied-output-incidence} with $\Phi=P$ and $\Psi=Q$. Fix $(a,b)$ with $P_x(a,b)N_x(a,b)\ne0$ and write $t_0=P(a,b)$. The implicit function theorem provides a smooth function $u(t,z)$ near $(t_0,b)$ with $P(u(t,z),z)=t$. We choose small intervals $I\ni a$ and $J\ni b$ and a common buffered interval $T$ containing $P(I,J)$, so that $u(t,z)$ is defined on a neighborhood of $T\times J$, and for $v,z\in J$, we set
$$g_{v,z}(t)=Q(u(t,z),v)\text{ for }t\in T.$$
The nondegeneracy we need concerns the value--slope map $(v,z)\mapsto(g_{v,z}(t),g'_{v,z}(t))$. Differentiating the implicit identity, we find $u_t=\frac{1}{P_x(u,z)}$ and $u_z=-\frac{P_y(u,z)}{P_x(u,z)}$, and at $(t,v,z)=(P(x,y),y,y)$, it follows that
$$\begin{aligned}
\partial_vg_{v,z}(t)&=Q_y, & \partial_zg_{v,z}(t)&=-\frac{Q_xP_y}{P_x},\\
\partial_vg'_{v,z}(t)&=\frac{Q_{xy}}{P_x}, & \partial_zg'_{v,z}(t)&=\frac{Q_x(P_yP_{xx}-P_xP_{xy})-P_xP_yQ_{xx}}{P_x^3},
\end{aligned}$$
where all derivatives on the right are evaluated at $(x,y)$. Taking the determinant, we find
\begin{equation}\label{eq:independent-pair-jet}
\det D_{(v,z)}(g_{v,z}(t),g'_{v,z}(t))=\frac{N_x(x,y)}{P_x(x,y)^3}\text{ at }(t,v,z)=(P(x,y),y,y),
\end{equation}
and in particular the derivative matrix at $(t_0,b,b)$ is invertible. We shrink the intervals and their buffers until the derivative differs from this fixed matrix by less than half of its smallest singular value. Integrating along the segment joining two parameter pairs, we then obtain some $c>0$ with
$$|g_{v,z}(t)-g_{\widetilde v,\widetilde z}(t)|+|g'_{v,z}(t)-g'_{\widetilde v,\widetilde z}(t)|\ge c|(v,z)-(\widetilde v,\widetilde z)|\text{ for }t\in T,$$
while the bounded derivatives on the buffered domain provide a fixed $C<\infty$ with
$$\|g_{v,z}-g_{\widetilde v,\widetilde z}\|_{C^2(T)}\le C|(v,z)-(\widetilde v,\widetilde z)|.$$
The value--slope maps are therefore uniformly bi-Lipschitz, and the two inequalities together bound the $C^2$ distance by a fixed multiple of the value--slope distance at every $t$, which is uniform transversality. The lower bound for $|P_x|$ and the common buffer take care of the remaining hypotheses of Lemma~\ref{lem:copied-output-incidence}.

To apply this local geometry with constants uniform in the input laws, we choose a fixed finite partition. It suffices to consider $0<\varepsilon<s+\beta_s$, since otherwise the asserted lower bound follows from the nonnegativity of entropy. We set $\eta=\frac{\varepsilon}{4(s+\beta_s)}$, so that discarding mass at most $\eta$ costs at most $\frac{\varepsilon n}{4}$ in the leading term. By independence, $(X,Y)$ is conditionally $(s,C_1;s,C_2)$-Frostman, and Lemma~\ref{lem:polynomial-sublevel}, applied to the nonzero polynomial $P_xN_x$, provides $r>0$, depending only on the fixed data and on $\eta$, with
$$\mathbb P(|P_x(X,Y)N_x(X,Y)|<r)\le\eta.$$
The compact set
$$\Omega=\{(x,y)\in[-R,R]^2:|P_x(x,y)N_x(x,y)|\ge r\}$$
is covered by the regular neighborhoods constructed above. We choose a finite subcover and a fixed product grid in $[-R,R]^2$, fine enough that every grid rectangle meeting $\Omega$ lies, together with the required buffers, in one of these neighborhoods. We call these rectangles good. The grid and all the geometric constants depend only on the fixed data and on $\varepsilon$, and the good rectangles have total probability at least $1-\eta$.

Let $I_i$ and $J_j$ be the intervals of the grid, chosen disjoint and half-open, with the outer endpoints assigned to the adjacent intervals, and set $p_i=\mathbb P(X\in I_i)$ and $q_j=\mathbb P(Y\in J_j)$. Conditioning on a rectangle $I_i\times J_j$ of positive probability preserves independence, and the two conditional laws are the normalized restrictions of $X$ and $Y$, with Frostman constants $\frac{C_1}{p_i}$ and $\frac{C_2}{q_j}$. On a good rectangle, we introduce an independent copy of the restricted $Y$, quantize all three inputs, and take $U$, $V$ and $V'$ in Lemma~\ref{lem:copied-output-incidence} to be the quantized versions of $X$, $Y$ and the copy of $Y$, respectively. These three variables are independent, so ${\mathcal{I}}(U;V)={\mathcal{I}}(V;V')=0$. After fixed enlargements of the test intervals to account for the quantization, all the conditional reference constants $K_v$ can be taken equal to a constant $C_{X,i}$, and the marginal reference constant $C_V$ can be taken equal to $C_{Y,j}$, where
$$C_{X,i}\le\frac{C(1+C_1)}{p_i}\text{ and }C_{Y,j}\le\frac{C(1+C_2)}{q_j}$$
and both constants are at least one. The atom bound for the input implies $H(U)\ge sn-\log C_{X,i}$, so the lemma with rate tolerance $\frac{\varepsilon}{2}$ produces the coefficient $s+\beta_s$ with an error of at most $\frac{3}{2}\log C_{X,i}+O_\varepsilon(1+\log C_{Y,j})$.

By subadditivity, the entropy of the copied output vector is at most the sum of the entropies of its two coordinates, which have the laws of the two polynomial outputs of the restricted and quantized input pair. Since the derivatives are bounded, these differ from the outputs of the original restricted pair by $O(1)$ in entropy, and consequently
\begin{equation}\label{eq:independent-pair-cell}
\begin{aligned}
&H_n(P(X,Y)\mid X\in I_i,Y\in J_j)+H_n(Q(X,Y)\mid X\in I_i,Y\in J_j)\\
&\ge\left(s+\beta_s-\frac{\varepsilon}{2}\right)n-C_\varepsilon\left(1+\log\frac{1}{p_i}+\log\frac{1}{q_j}\right).
\end{aligned}
\end{equation}
At sufficiently fine scales, all representatives lie in the fixed buffers, and since the incidence threshold does not depend on $C_{X,i}$ and $C_{Y,j}$, the finitely many remaining scales are absorbed into the same error by the nonnegativity of entropy. There are finitely many fixed charts, so $C_\varepsilon$ is uniform in $i$, $j$, the input laws and $n$.

Finally, let $L_X=i$ when $X\in I_i$ and $L_Y=j$ when $Y\in J_j$. Entropy is at least the conditional entropy given these labels, the rectangle $I_i\times J_j$ has probability $p_iq_j$, and the normalization costs satisfy
$$\sum_{i,j}p_iq_j\left(\log\frac{1}{p_i}+\log\frac{1}{q_j}\right)=H(L_X)+H(L_Y)\le\log\#\{I_i\}+\log\#\{J_j\},$$
where terms of probability zero are omitted. Applying \eqref{eq:independent-pair-cell} on the good rectangles and the nonnegativity of entropy on the others, we obtain
$$H_n(P(X,Y))+H_n(Q(X,Y))\ge(1-\eta)\left(s+\beta_s-\frac{\varepsilon}{2}\right)n-O_\varepsilon(1)\ge(s+\beta_s-\varepsilon)n-O_\varepsilon(1),$$
because $\eta(s+\beta_s)=\frac{\varepsilon}{4}$. Since $s+\beta_s=\min\{\frac{5s}{2},1+s\}$ and all partitions, buffers and constants were chosen from the fixed data, this proves \eqref{eq:independent-polynomial-pair} with an error that is uniform in the input laws and in $n$.
\end{proof}

\subsubsection{Concentration near polynomial zero sets}

For a nonzero polynomial $P$, we define the truncated logarithmic depth
\begin{equation}\label{eq:joint-depth-definition}
K_n(P)=\min\left\{n,\left\lfloor\log_+\frac{1}{|P(X,Y)|}\right\rfloor\right\},\text{ where }\log_+u=\max\{0,\log u\},
\end{equation}
with $K_n(P)=n$ on the event $\{P(X,Y)=0\}$. A factor of a polynomial is called mixed if it depends on both variables.

\begin{lemma}\label{lem:joint-multiplicity-charge}
Let $(X,Y)$ be jointly $(s,s,C)$-Frostman in a fixed bounded box, where $0<s<1$. For a fixed nonzero real polynomial $P$, let $r$ be the largest multiplicity of a mixed irreducible factor of $P$, with $r=0$ if there is no such factor. Then
\begin{equation}\label{eq:joint-multiplicity-charge}
\mathbb EK_n(P)\le\frac{r}{s}\left(\frac{H_n(X)+H_n(Y)}{2}-sn\right)+O(\log(n+2)),
\end{equation}
and if $r=0$, then $\mathbb EK_n(P)=O(1)$.
\end{lemma}

\begin{proof}
Assume first that $P$ is squarefree and has no univariate factor. We bound the probability that $X$ lies in a given dyadic interval while the depth takes a given value, and then taking negative logarithms and averaging controls the expected depth in terms of $H_n(X)$. The same argument in the other coordinate proves the symmetric estimate.

Since $P$ is squarefree and has no factor depending only on $x$, the polynomials $P$ and $P_y$ are relatively prime over $\mathbb R(x)$. Clearing denominators in B\'ezout's identity, we obtain polynomials $A$, $B$ and $\Delta_P$ with
$$AP+BP_y=\Delta_P(x)\text{ and }\Delta_P\not\equiv0.$$
This identity controls the lengths of the fibers of sublevel sets. For fixed $x$ with $\Delta_P(x)\ne0$, the set of $y$ in the support with $|P(x,y)|\le u$ is covered by $O_P(1)$ intervals of length $O_P(\frac{u}{|\Delta_P(x)|})$. Indeed, if $u\le c_P|\Delta_P(x)|$, the identity and the boundedness of $A$ and $B$ on the support give $|P_y|\ge c'_P|\Delta_P(x)|$ on the sublevel set, so $P_y$ has a constant sign on each component of the sublevel set, whose length is then bounded by the mean value theorem, while the degree bounds the number of components. If $u>c_P|\Delta_P(x)|$, the support interval itself satisfies the required length bound.

We apply this bound to the dyadic intervals of the input. Let $\delta=2^{-n}$, let $I_i=[i\delta,(i+1)\delta)$, and set $b_i=\min\{1,\max\{\delta,|\Delta_P(i\delta)|\}\}$. If $X\in I_i$ and $K_n(P)=k$, the bounded derivatives of $P$ imply $|P(i\delta,Y)|\le C_P2^{-k}$. Covering the corresponding sublevel set at $x=i\delta$ as above and using the joint rectangle bound, we obtain
\begin{equation}\label{eq:joint-depth-atom}
\mathbb P(X\in I_i,K_n(P)=k)\le C_P\delta^s2^{-sk}b_i^{-s}\text{ for }0\le k\le n.
\end{equation}
When $k=0$, when the covering intervals are longer than the support, or when $|\Delta_P(i\delta)|<\delta$, the right-hand side is at least a fixed multiple of $\delta^s$, and the marginal bound yields the same estimate.

The factor $b_i^{-s}$ is the only loss in \eqref{eq:joint-depth-atom}, and its mean logarithmic cost is bounded. Indeed, if $\Delta_P$ is nonconstant, factorization in one variable and the marginal bound imply
$$\mathbb P(|\Delta_P(D_nX)|\le u)\lesssim_{\Delta_P}u^{\frac{s}{\deg\Delta_P}}+\delta^s\text{ for }0<u\le1,$$
where the term $\delta^s$ accounts for the dyadic intervals at the boundary of the finitely many sublevel intervals. Summing over $u=2^{-j}$ for $0\le j\le n$, and noting that $n\delta^s$ is bounded, we get $\mathbb E\log\frac{1}{b_{\lfloor 2^nX\rfloor}}=O(1)$, which is immediate when $\Delta_P$ is constant. Taking negative logarithms in \eqref{eq:joint-depth-atom} and averaging now yields
$$H_n(X)+H(K_n(P))\ge H(D_nX,K_n(P))\ge sn+s\mathbb EK_n(P)-O(1).$$
Since $P$ has no univariate factor in either variable, the same argument applies with $X$ and $Y$ interchanged. Averaging the two inequalities and paying for the depth label with $H(K_n(P))\le\log(n+1)$, we obtain \eqref{eq:joint-multiplicity-charge} in the squarefree case.

For general $P$, we write $P=a(x)b(y)\prod_jP_j^{r_j}$ with distinct mixed irreducible factors $P_j$ and set $Q=\prod_jP_j$. Since all these factors are bounded on the support,
$$K_n(P)\le rK_n(Q)+\log_+\frac{1}{|a(X)|}+\log_+\frac{1}{|b(Y)|}+O(1),$$
and the last two terms have bounded expectation by univariate factorization and the marginal Frostman bounds. If $r>0$, then $Q$ is squarefree and has no univariate factor, so the estimate already proved applies to $Q$ and yields the coefficient $\frac{r}{s}$. If $r=0$, only the two bounded univariate terms remain, and $\mathbb EK_n(P)=O(1)$.
\end{proof}

\subsection{Joint expansion from conditional copies}
\label{sec:joint-forward}

The incidence inequality applies away from a polynomial critical set. We keep track of the logarithmic cost of approaching this set, and we bound this cost in terms of the marginal entropies by means of Lemma~\ref{lem:joint-multiplicity-charge}. Set
$$\mathcal W_x=f_yf_{xx}-f_xf_{xy}\text{ and }\mathcal W_y=f_xf_{yy}-f_yf_{xy}.$$

\begin{lemma}\label{lem:joint-critical-multiplicity}
If $f$ has degree $d$ and is not of the form \eqref{eq:excluded-affine-form}, then at least one of $\mathcal W_x$ and $\mathcal W_y$ is nonzero, and every mixed irreducible factor of a nonzero one among them has multiplicity at most $2d-4$.
\end{lemma}

\begin{proof}
If one of the first partial derivatives of $f$ vanished identically, $f$ would depend on one variable only and would have the excluded form, so $f_x,f_y\not\equiv0$. Since
$$\partial_x\frac{f_x}{f_y}=\frac{\mathcal W_x}{f_y^2}\text{ and }\partial_y\frac{f_x}{f_y}=-\frac{\mathcal W_y}{f_y^2},$$
the vanishing of both critical polynomials would make $\frac{f_x}{f_y}$ equal to a constant $c$. Then $f_x=cf_y$, and an invertible linear change of variables would make $f$ independent of one coordinate, contrary to the hypothesis. Hence, at least one of $\mathcal W_x$ and $\mathcal W_y$ is nonzero.

For the multiplicity bound, it suffices to bound the degree in a variable on which every mixed factor depends. We regard $\mathcal W_x=f_y(f_x)_x-f_x(f_y)_x$ as the Wronskian of two polynomials in $x$ with coefficients in $\mathbb R[y]$, both of degree at most $d-1$ in $x$. If the two degrees agree, the leading terms cancel and the Wronskian has degree at most $2d-4$ in $x$, and if they differ, the sum of the degrees is at most $2d-3$ and the differentiation again leads to degree at most $2d-4$. Every mixed irreducible factor has positive degree in $x$, so its multiplicity in $\mathcal W_x$, when $\mathcal W_x\ne0$, is at most $2d-4$. Interchanging $x$ and $y$, we obtain the same bound for $\mathcal W_y$.
\end{proof}

\begin{proposition}\label{prop:joint-depth-seed}
Let $(X,Y)$ be jointly $(s,s,C)$-Frostman in a fixed bounded box, where $0<s<1$, and suppose that $\mathcal W_x\not\equiv0$. For every $\varepsilon>0$,
\begin{equation}\label{eq:joint-depth-seed}
\begin{aligned}
&H_n(S)+H_n(F)+\frac{3}{2}H_n(X)+3H_n(Y)-\frac{5}{2}H_n(X,Y)\\
&\ge\frac{3s+t_s}{2}n-16\mathbb EK_n(\mathcal W_x)-\varepsilon n-O_\varepsilon(\log(n+2)).
\end{aligned}
\end{equation}
If $\mathcal W_y\not\equiv0$, the same holds with $X$ and $Y$ interchanged and with $K_n(\mathcal W_y)$ in place of $K_n(\mathcal W_x)$.
\end{proposition}

\begin{proof}
We first prove the inequality for quantized inputs, keeping the logarithmic cost of the critical set in every incidence class, and at the end, a comparison of polynomial values returns the bound to the original inputs. Let $\delta=2^{-n}$ and replace each input by the center of its dyadic interval of length $\delta$. This does not change the input entropies, and since the two output maps have bounded derivatives, it changes the output entropies by $O_f(1)$. Let $k$ denote the truncated depth of $\mathcal W_x$ at these centers. We begin with the part of the law where $k\le\frac{n}{2}$, conditioning the original law on this event before forming any copies. This cutoff guarantees that all the normalized meshes below tend to zero uniformly. Conditioning on an event of probability $p$ enlarges the rectangle and marginal constants by a factor of at most $\frac{1}{p}$, and we keep this logarithmic cost until the two parts are averaged at the end.

For now, all entropies and reference probabilities refer to the conditional law on $\{k\le\frac{n}{2}\}$. We sample $Y$ and $Y'$ independently given the quantized $X$ under this law and use the class calculation of Lemma~\ref{lem:copied-output-incidence} with
$$W=(X+Y',f(X,Y))\text{ and }g_{y,z}(t)=f(t-z,y).$$
The first coordinate of $W$ is invertible in $X$ with derivative one, and the laws of the two coordinates of $W$ coincide with those of $S$ and $F$, so
\begin{equation}\label{eq:joint-copy-output}
H_n(W)\le H_n(S)+H_n(F)+O_f(1).
\end{equation}
To bound the entropy of $W$ from below, we have to make the incidence geometry uniform as $k$ varies and to account for the necessary labels. This leads to \eqref{eq:joint-good-depth} below, with a coefficient of the depth at most $13$, which we then combine with a trivial estimate on the remaining part of the law.

Fix a depth $k<n$ in the range under consideration, write $\rho=2^{-k}$, and partition the copied triples into product cells of side $\ell=2^{-a_f}\rho$. By the definition of the depth, $|\mathcal W_x|\ge\frac{\rho}{2}$ at the center of a triple in the class, and choosing the fixed integer $a_f$ large enough in terms of the derivative bounds of $f$, we ensure that $|\mathcal W_x|\ge\frac{\rho}{4}$ on every occupied cell and on its required enlargements. For fixed $t$, the value--slope map
$$\Phi_t(y,z)=(f(t-z,y),f_x(t-z,y))$$
has determinant $-\mathcal W_x$, so its smallest singular value is at least $c_f\rho$, while its derivative varies by $O_f(\ell)$ across the cell. Comparing the integral of the derivative along a segment with the derivative at one endpoint, we obtain
\begin{equation}\label{eq:joint-jet-inverse}
|\Phi_t(y,z)-\Phi_t(y',z')|\ge c_f\rho|(y,z)-(y',z')|
\end{equation}
for all parameter points $(y,z)$ and $(y',z')$ in the cell.

Next, we normalize each cell so that the incidence theorem can be used with a transversality constant independent of $k$. We rescale the horizontal and the vertical coordinates by $\ell$, subtract a common vertical constant and apply fixed further dilations. The curves then lie in the unit ball of $C^2$ on a common buffered interval $[-2,2]$, at mesh $\delta'=\frac{\delta}{\ell}$. If $\mathfrak h$ is the difference of two normalized curves with parameters $(y,z)$ and $(y',z')$, then
$$|\mathfrak h''|\le C_f\ell|(y,z)-(y',z')|\le C_f\ell\rho^{-1}(|\mathfrak h|+|\mathfrak h'|),$$
and after increasing $a_f$, the last coefficient is at most one. Gronwall's inequality for $(\mathfrak h,\mathfrak h')$ now provides a transversality constant $T$ independent of $k$. Moreover, our cutoff implies that $\delta'\le2^{a_f-\frac{n}{2}}<\delta_0(s,T,\eta)$ for all sufficiently large $n$, where the tolerance $\eta$ is chosen below.

We also need a bound for the mass that the parameter reference assigns to rectangles after this normalization. The derivative of either component of the normalized value--slope map has the form $\xi D\Phi_t$, where, up to fixed factors, $\xi=(\ell^{-1},0)$ or $\xi=(0,1)$. Its norm is at least $c_f\rho|\xi|$ and its variation over the cell is $O_f(\ell|\xi|)$, so one of its two partial derivatives keeps its sign and has magnitude at least $c_f\rho$ throughout the cell, since $|\xi|\ge1$. For an axis-parallel rectangle in the value--slope plane with sides $a$ and $b$, the inverse bound \eqref{eq:joint-jet-inverse} confines the projection of the preimage to an interval of length $O_f(\rho^{-1}\max\{a,b\})$, while the derivative bound produces fibers of length $O_f(\rho^{-1}\min\{a,b\})$. A product of marginal $s$-Frostman references with constants $C_0$ therefore assigns the preimage mass at most
\begin{equation}\label{eq:joint-jet-reference}
C_fC_0^2\rho^{-2s}(ab)^s.
\end{equation}
The references on the grid satisfy their bounds down to scale $\delta$, and since $a,b\ge\delta'$, enlarging the corresponding parameter intervals by $\delta$ changes this estimate only by a fixed factor.

The last geometric requirement is the separation of the curves at the normalized mesh. We color the two indices of the parameter grid modulo an integer $N_k\le C_f2^{2k}$ chosen so that $N_k\delta\ge\frac{C_f\delta}{\ell\rho}$, and then, by \eqref{eq:joint-jet-inverse}, distinct parameter pairs of the same color give $\delta'$-separated curves. At depth $k$, there are $O_f(2^{4k})$ colors and $O_f(2^{3k})$ product cells, and the combined label $L$ recording the depth, the cell and the color, which is attached after copying, satisfies
\begin{equation}\label{eq:joint-geometric-label}
H(L)\le7\mathbb Ek+O_f(\log(n+2)).
\end{equation}
Indeed, at a fixed depth, the logarithm of the number of choices of a cell and a color is at most $7k+O_f(1)$, and averaging this bound and adding the entropy of the depth, which takes at most $n+1$ values, yields \eqref{eq:joint-geometric-label}.

With the geometric hypotheses of Theorem~\ref{thm:rectangular-incidence} at mesh $\delta'$ in place, we turn to the probability weights. We use the levels $b_\nu,c,\zeta,h,\lambda$ from the proof of Lemma~\ref{lem:copied-output-incidence}, now computed for the conditional law under consideration, and we write $y_j$ for the centers of the dyadic intervals of $Y$ and $q_j=\mathbb P(Y=y_j)$ for their probabilities. The joint rectangle bound provides the fiber reference constants
$$K_{y_j}=\frac{C_0\delta^s}{q_j}\ge1,\text{ for which }\mathbb E\log K_Y=H_n(Y)-sn+\log C_0,$$
where the factor $C_0$ includes the cost of the initial restriction. On a class of probability $w$ and depth $k$, the parameter and fiber estimates imply
$$C_1=C_fC_0^22^{2sk}2^{b_\nu+h}(\delta')^{2s}\text{ and }C_2=C_fK_{y_j}\ell^s2^{c+\lambda}(\delta')^{t_s},$$
and both constants are at least one by the bounds for the selected atoms in the proof of Lemma~\ref{lem:copied-output-incidence}. The counts of curves and incidences from that proof are unchanged, and the count of squares remains $O_f(2^\zeta)$ because a $\delta$-square of the original grid meets only boundedly many normalized $\delta'$-squares.

We apply Theorem~\ref{thm:rectangular-incidence} at mesh $\delta'$, with its tolerance decreased so that the final linear loss is at most $\varepsilon n$. Substituting the constants into the logarithmic form of the incidence bound, we get the class inequality
$$\begin{aligned}
\zeta\ge{}&c+\left(s+\frac{t_s}{2}\right)n-h-\frac{\lambda}{2}-\frac{1}{2}\log K_{y_j}+\frac{3}{2}\log w\\
&-\frac{5s+t_s}{2}k-2\log C_0-\varepsilon n-O_{f,\varepsilon}(1).
\end{aligned}$$
It contains two explicit costs of the depth, namely $\frac{(s+t_s)k}{2}$ from the rescaling of the mesh and the fibers, and $2sk$ from the rectangle bound \eqref{eq:joint-jet-reference}. We average this inequality using \eqref{eq:joint-copy-identities}, \eqref{eq:joint-copy-output} and \eqref{eq:joint-geometric-label}. The labels that do not come from the geometry cost $O(\log(n+2))$, exactly as in the proof of Lemma~\ref{lem:copied-output-incidence}, because the shift in $\log K_{y_j}=\log C_0-sn-\log q_j$ is deterministic. Before we discard the nonnegative mutual information of the copies, the information terms are
$$H_n(X)-\frac{5}{2}{\mathcal{I}}(X;Y)+\frac{1}{2}{\mathcal{I}}(Y;Y')-\frac{1}{2}(H_n(Y)-sn+\log C_0),$$
and the geometric label contributes the further amount $\frac{3}{2}\cdot7\mathbb Ek$ through the term $\frac{3}{2}\log w$. Expanding ${\mathcal{I}}(X;Y)$ and discarding $\frac{1}{2}{\mathcal{I}}(Y;Y')\ge0$, we arrive at
\begin{equation}\label{eq:joint-good-depth}
\begin{aligned}
&H_n(S)+H_n(F)+\frac{3}{2}H_n(X)+3H_n(Y)-\frac{5}{2}H_n(X,Y)\\
&\ge\frac{3s+t_s}{2}n-\left(\frac{5s+t_s}{2}+\frac{21}{2}\right)\mathbb Ek-\varepsilon n-O_\varepsilon(\log(n+2)).
\end{aligned}
\end{equation}
For the restricted law, the suppressed constant contains the additional cost $O(\log\frac{1}{p})$, whose average with respect to the probabilities of the two parts is bounded, since a binary label has entropy at most one.

This proves the required estimate on the part where $k\le\frac{n}{2}$. On the complementary part, we use the conditional law on $\{k>\frac{n}{2}\}$, for which the inequality $H_n(X,Y)\le H_n(X)+H_n(Y)$ and the nonnegativity of entropy yield
$$H_n(S)+H_n(F)+\frac{3}{2}H_n(X)+3H_n(Y)-\frac{5}{2}H_n(X,Y)\ge-H_n(X)\ge-n-O_f(1).$$
Since $\frac{3s+t_s}{2}\le\frac{3}{2}$, the loss compared with the target is at most $\frac{5}{2}n+O_f(1)\le5\mathbb Ek+O_f(1)$ on this part, while the coefficient of $\mathbb Ek$ in \eqref{eq:joint-good-depth} is at most $13$ throughout $0<s<1$. Each entropy differs from its average over the two parts by at most the entropy of the binary label, so the same is true, up to $O(1)$, for our fixed signed combination of entropies. Averaging the two inequalities therefore proves the asserted bound for the centers, since both coefficients of the depth are at most $16$, and the finitely many small values of $n$ are absorbed into the constant.

It remains to replace the depth at the centers by the depth in the statement. The bounded derivatives of $\mathcal W_x$ imply
$$k\le K_n(\mathcal W_x)+O_f(1).$$
Indeed, if the value of $\mathcal W_x$ at the center exceeds $2C_f\delta$ in absolute value, then its absolute value at the actual point is at most twice that at the center, and otherwise the actual value has absolute value at most $3C_f\delta$, so that its truncated depth is at least $n-O_f(1)$. This returns the estimate to the original law. Interchanging $X$ and $Y$ preserves the rectangle bound and leads to the assertion for $\mathcal W_y$.
\end{proof}

\begin{theorem}\label{thm:joint-incidence-bound}
Let $0<s<1$, and let $f$ have degree $d$ and satisfy the hypothesis of Lemma~\ref{lem:joint-critical-multiplicity}. Fix the support box and the joint $(s,s,C)$-Frostman constant. Choose a nonzero polynomial among $\mathcal W_x$ and $\mathcal W_y$, and let $r$ be the largest multiplicity of a mixed irreducible factor of it, with $r=0$ if it has no such factor. For every $\varepsilon>0$,
\begin{equation}\label{eq:joint-multiplicity-gain}
M_{\mathrm j,n}\ge\left(s+\frac{s\min\{s,2-2s\}}{16s+32r}-\varepsilon\right)n-O_\varepsilon(1).
\end{equation}
In particular, for every $0<s<1$,
\begin{equation}\label{eq:joint-all-exponents}
M_{\mathrm j,n}\ge\left(s+\frac{s\min\{s,2-2s\}}{64d}\right)n-O(1),
\end{equation}
and for $s>\frac{1}{2}$, we also have
\begin{equation}\label{eq:joint-regular-gain}
M_{\mathrm j,n}\ge\left(s+\frac{\min\{s,2-2s\}}{16}-\varepsilon\right)n-O_\varepsilon(1).
\end{equation}
All constants are uniform over laws with the fixed data.
\end{theorem}

\begin{proof}
We replace the loss caused by the depth in Proposition~\ref{prop:joint-depth-seed} by the excess of the marginal entropies, and then we bound the resulting weighted inequality from above by a positive multiple of $M_{\mathrm j,n}$. Lemma~\ref{lem:joint-critical-multiplicity} provides a nonzero orientation. We treat the case where the chosen polynomial is $\mathcal W_x$, since interchanging $X$ and $Y$ preserves both the Frostman hypothesis and $M_{\mathrm j,n}$ and therefore covers the case of $\mathcal W_y$.

We apply Lemma~\ref{lem:joint-multiplicity-charge} to $\mathcal W_x$, whose mixed factors have multiplicity at most $r$, and insert the resulting bound into \eqref{eq:joint-depth-seed}. Since $H_n(X,Y)\ge2sn-O(1)$ by Lemma~\ref{lem:Frostman-entropy}, we obtain
$$\begin{aligned}
&H_n(S)+H_n(F)+\left(\frac{3}{2}+\frac{8r}{s}\right)H_n(X)+\left(3+\frac{8r}{s}\right)H_n(Y)\\
&\ge\left(\frac{13s+t_s}{2}+16r\right)n-\varepsilon n-O_\varepsilon(\log(n+2)).
\end{aligned}$$
The two marginal coefficients are different, while $M_{\mathrm j,n}$ only controls the unweighted sum of the marginal entropies. To use this control without increasing the final denominator, we write the left-hand side as
$$H_n(S)+H_n(F)+\left(3+\frac{8r}{s}\right)(H_n(X)+H_n(Y))-\frac{3}{2}H_n(X).$$
The inequalities $H_n(S)+H_n(F)\le2M_{\mathrm j,n}$, $H_n(X)+H_n(Y)\le2M_{\mathrm j,n}$ and $H_n(X)\ge sn-O(1)$ then imply
$$\left(8+\frac{16r}{s}\right)M_{\mathrm j,n}\ge\left(8s+16r+\frac{t_s}{2}\right)n-\varepsilon n-O_\varepsilon(\log(n+2)).$$
The coefficient of $M_{\mathrm j,n}$ is positive, and dividing by it, after starting with a smaller tolerance to absorb the term $O(\log(n+2))$, proves \eqref{eq:joint-multiplicity-gain}.

The bound \eqref{eq:joint-all-exponents} with a bounded error follows by leaving a fixed gap in the leading coefficient. By the multiplicity bound of Lemma~\ref{lem:joint-critical-multiplicity},
$$16s+32r\le16s+32(2d-4)<64d,$$
so, since $st_s>0$, the gain in \eqref{eq:joint-multiplicity-gain} is strictly larger than $\frac{st_s}{64d}$. Choosing the tolerance smaller than this fixed positive difference, we obtain \eqref{eq:joint-all-exponents} with a constant independent of $n$ and of the law.

Finally, we prove the rate \eqref{eq:joint-regular-gain}, which does not depend on the degree. When $s>\frac{1}{2}$, the critical set can be removed at an arbitrarily small cost in probability. The tube estimate \eqref{eq:joint-tube}, applied to the fixed algebraic set $\{\mathcal W_x=0\}$, bounds the mass of its $h$-neighborhood by $O_f(h^{2s-1})$, and on the complement of such a fixed neighborhood, compactness keeps $|\mathcal W_x|$ away from zero, so that $K_n(\mathcal W_x)=O_{f,h}(1)$. We apply Proposition~\ref{prop:joint-depth-seed} to the normalized restriction of the law to this complement. On the removed part, the left-hand side of \eqref{eq:joint-depth-seed} is at least $-n-O_f(1)$, as shown in the preceding proof. Choosing $h>0$ so small that the removed mass contributes less than the prescribed linear tolerance and averaging, at the cost of the bounded entropy of a binary label, we obtain the inequality of Proposition~\ref{prop:joint-depth-seed} without any linear loss from the depth. The calculation in the first part of this proof then applies with $r=0$, and after adjusting the tolerance, it yields \eqref{eq:joint-regular-gain}.
\end{proof}

\subsubsection{Joint bounds from the multiscale entropy inequality}

For exponents above two thirds, the multiscale estimate leads to a stronger gain. The following deduction shows how the marginal term in $M_{\mathrm j,n}$ enters.

\begin{corollary}\label{cor:joint-multiscale-gain}
Under the hypotheses of Theorem~\ref{thm:joint-incidence-bound}, let $\frac{1}{2}<s<1$ and $\kappa_s=\min\{1,2s-\frac{1}{2}\}$. For every $\varepsilon>0$,
\begin{equation}\label{eq:joint-multiscale-gain}
M_{\mathrm j,n}\ge\left(\frac{3s+\kappa_s}{4}-\varepsilon\right)n-O_\varepsilon(1).
\end{equation}
The gain over $s$ equals $\frac{2s-1}{8}$ for $s\le\frac{3}{4}$ and $\frac{1-s}{4}$ for $s\ge\frac{3}{4}$.
\end{corollary}

\begin{proof}
By Corollary~\ref{cor:three-coordinate}, one of $f$, $G_Y$ and $G_X$ is nonspecial. We first treat a nonspecial sum-coordinate representation, where the input marginal that does not appear accounts for the denominator four in the conclusion. If $G_Y$ is nonspecial, Corollary~\ref{cor:global-multiscale-profile} in the coordinates $(Y,S)$ shows that
$$H_n(F)+{\mathcal{I}_n}(Y;S)\ge(\kappa_s-\varepsilon)n-O_\varepsilon(1).$$
The shear $(X,Y)\mapsto(Y,S)$ changes the joint dyadic entropy by $O(1)$, so $H_n(Y,S)\ge2sn-O(1)$, and adding this to the preceding inequality and using the definition of mutual information, we find
$$H_n(Y)+H_n(S)+H_n(F)\ge(2s+\kappa_s-\varepsilon)n-O_\varepsilon(1).$$
The left-hand side is at most $4M_{\mathrm j,n}-H_n(X)\le4M_{\mathrm j,n}-sn+O(1)$, and dividing by four, we obtain the required coefficient. If $G_X$ is nonspecial, we interchange $X$ and $Y$ in these steps, and the same maximum and the same lower bound for the missing marginal apply.

If neither sum-coordinate representation is nonspecial, the classification shows that $f$ itself is nonspecial, and Corollary~\ref{cor:global-multiscale-profile} in the original coordinates yields $H_n(F)+{\mathcal{I}_n}(X;Y)\ge(\kappa_s-\varepsilon)n-O_\varepsilon(1)$. Adding $H_n(X,Y)\ge2sn-O(1)$, we now obtain
$$3M_{\mathrm j,n}\ge H_n(X)+H_n(Y)+H_n(F)\ge(2s+\kappa_s-\varepsilon)n-O_\varepsilon(1),$$
and since $\kappa_s\ge s$ implies $\frac{2s+\kappa_s}{3}\ge\frac{3s+\kappa_s}{4}$, this is even stronger than the desired bound. Adjusting the initial tolerance proves \eqref{eq:joint-multiscale-gain} in all cases. Finally, $\kappa_s=2s-\frac{1}{2}$ for $s\le\frac{3}{4}$ and $\kappa_s=1$ for $s\ge\frac{3}{4}$, and substituting these values, we obtain the two stated gains, which agree at $s=\frac{3}{4}$.
\end{proof}

\subsection{Conditional copies and inverse coordinates}
\label{sec:conditional-copies-new}

In this subsection, we derive three local entropy inequalities from conditional copies and the incidence estimate. Restrictions to events and the passage from regular charts to the full law are handled once and for all, so that the same inequalities can be used in the multiscale argument. Throughout, both exact conditional $s$-Frostman constants of the original law are at most $C\ge1$.

\subsubsection{Restrictions and changes of coordinates}

\begin{lemma}\label{lem:event-conditional-entropy}\label{lem:conditional-cores}
If finite-valued random variables satisfy $\mathbb P(U=u\mid V=v)\le2^{-h}$, then every event $A$ of probability $a>0$ satisfies $H(U\mid V,A)\ge h+\log a$. In particular, both $H(D_jX\mid D_jY,A)$ and $H(D_jY\mid D_jX,A)$ are at least $sj-\log C+\log a$. Moreover, $A$ has a countable partition into conditional cores $A_r$ such that the label of the core has bounded entropy under the law conditioned on $A$, and the two exact conditional constants $C_r$ of the cores satisfy
$$\sum_r\mathbb P(A_r\mid A)\log C_r\le\log\frac{C}{a}+O(1).$$
Consequently, suppose that a fixed linear combination of entropies of fixed functions of $(X,Y)$ satisfies a lower bound with error $O_\varepsilon(1+\log C)$ for every law on a chart whose two exact conditional constants are at most $C$, with coefficients and rate independent of $C$. Then the same lower bound holds with error $O_\varepsilon(1+\log\frac{C}{a})$ for the normalized restriction of the original law to every event $A$ of probability $a$ in the chart.
\end{lemma}

\begin{proof}
We first measure the loss of entropy under restriction. For a conditioning value $v$ of positive probability, let $r_v=\mathbb P(A\mid V=v)$. The values with $r_v=0$ have probability zero under the law conditioned on $A$ and can be ignored, and under this law the conditional atoms of $U$ are at most $\frac{2^{-h}}{r_v}$. Therefore,
$$H(U\mid V,A)\ge h+\mathbb E[\log r_V\mid A]=h+\log a+D(\mathcal L(V\mid A)\Vert\mathcal L(V))\ge h+\log a.$$
The exact conditional Frostman bounds, averaged over the conditioning cells, show that the atoms of $D_jX\mid D_jY$ and of $D_jY\mid D_jX$ have size at most $C2^{-sj}$, and the finite assertion in these two orientations yields the dyadic entropy bounds.

Next, we partition $A$ into pieces on which both exact conditional bounds hold uniformly, measuring all masses in this extraction under the original law. Starting from an event of original mass $w>0$, we alternately delete all $X$-fibers and all $Y$-fibers whose remaining acceptance, measured under the original conditional law, is less than $\frac{w}{4}$. At each stage, the entire remaining portion of such a fiber is removed, so a deleted coordinate value never returns.The newly deleted coordinate values in each orientation are disjoint, and integrating their acceptance probabilities shows that each orientation removes mass at most $\frac{w}{4}$. By conditional monotone convergence, what remains is a measurable core of mass at least $\frac{w}{2}$, and up to the usual null sets, every surviving fiber has original acceptance at least $\frac{w}{4}$. For instance,
$$\mathbb P(X\in I\mid Y=y,\text{core})\le\frac{\mathbb P(X\in I\mid Y=y)}{\mathbb P(\text{core}\mid Y=y)}\le\frac{4C}{w}|I|^s,$$
and the other orientation satisfies the same bound.

We repeat this extraction on the remainder of $A$. If $v_r$ is the mass of the core $A_r$ and $w_r$ is the mass remaining before its extraction, then $\frac{w_r}{2}\le v_r\le w_r$ and $w_{r+1}\le\frac{w_r}{2}$. The remainders have masses tending to zero, so the cores partition $A$ up to a null set. We index the cores by $r\ge1$ and let $L=r$ on $A_r$. Then $\mathbb P(L\ge r\mid A)=\frac{w_r}{a}\le2^{1-r}$, so that $\mathbb E[L\mid A]\le2$, and comparing the law of $L$ with the probabilities $2^{-r}$ shows that $H(L\mid A)\le\mathbb E[L\mid A]\le2$. Since $C_r\le\frac{4C}{w_r}\le\frac{4C}{v_r}$, we obtain
$$\sum_r\frac{v_r}{a}\log C_r\le\log\frac{4C}{a}+\sum_r\frac{v_r}{a}\log\frac{a}{v_r}=\log\frac{C}{a}+O(1).$$

Finally, for every finite quantization $Z$ of a fixed source observable,
$$0\le H(Z\mid A)-\sum_r\mathbb P(A_r\mid A)H(Z\mid A_r)\le H(L\mid A)=O(1),$$
so passing from the cores to $A$ changes a fixed signed combination of entropies by at most the sum of the absolute values of its coefficients times $H(L\mid A)$. We apply the given inequality on each core and average its logarithmic constant by the preceding estimate. Copies may be formed separately on the cores, as long as they are eliminated before this final averaging of the entropies of fixed source observables.
\end{proof}

Let us stress the difference between the two averaging procedures used below. Restricting the law to an arbitrary event $A$ preserves the conditional dyadic entropy bounds of Lemma~\ref{lem:conditional-cores} with the loss $\log\frac{C}{\mathbb P(A)}$, but the restricted law need not satisfy uniform exact conditional Frostman bounds with constant $\frac{C}{\mathbb P(A)}$. When exact conditional bounds are needed, we therefore form copies separately under the normalized laws on the cores $A_r$, and we eliminate these copies before averaging the entropies of fixed source observables. Since the core label has bounded entropy, signed coefficients are allowed in this step. In Lemma~\ref{lem:conditional-chart-averaging}, by contrast, the exceptional part of the chart is discarded, and this requires nonnegative coefficients.

\begin{lemma}\label{lem:conditional-chart-averaging}
Suppose that, outside the zero sets of finitely many nonzero polynomials, each point has a fixed buffered chart with the following property. For every $\eta>0$, every event $A$ of positive original probability contained in that chart obeys
$$\sum_ib_iH_n(Z_i\mid A)\ge(\gamma-\eta)n-O_{\eta,\mathrm{chart}}\left(1+\log\frac{C}{\mathbb P(A)}\right),$$
where the coefficients $b_i$ are nonnegative. Then $\sum_ib_iH_n(Z_i)\ge(\gamma-\varepsilon)n-O_\varepsilon(1)$ for every $\varepsilon>0$, uniformly over the prescribed source laws.
\end{lemma}

\begin{proof}
Fix $\varepsilon>0$ and choose $\eta,\rho>0$ with $\eta+\rho(|\gamma|+\eta)<\varepsilon$. These parameters stay fixed while we choose a finite family of charts, so that both the discarded mass and all the constants of the charts are uniform in the law. Lemma~\ref{lem:polynomial-sublevel}, applied to the product of the nonzero critical polynomials, provides a fixed neighborhood of their zero sets of probability at most $\rho$, uniformly over the prescribed laws. The rest of the fixed support is compact, so we can cover it by finitely many buffered charts and partition it into disjoint events $A_1,\ldots,A_m$ subordinate to this cover. Let $A_0$ be the discarded event and $p_k=\mathbb P(A_k)$, and use $\eta$ as the local tolerance. Since there are finitely many charts, their error coefficients have a common bound $K_\eta$, and the concavity of entropy and the nonnegativity of the $b_i$ imply
$$\begin{aligned}
\sum_ib_iH_n(Z_i)&\ge\sum_{k=1}^mp_k\sum_ib_iH_n(Z_i\mid A_k)\\
&\ge(1-p_0)(\gamma-\eta)n-K_\eta\sum_{k=1}^mp_k\left(1+\log C+\log\frac{1}{p_k}\right),
\end{aligned}$$
where events of probability zero do not contribute. The last sum is at most $1+\log C+\sum_{k=0}^mp_k\log\frac{1}{p_k}\le1+\log C+\log(m+1)$, and moreover
$$(1-p_0)(\gamma-\eta)\ge\gamma-\eta-\rho(|\gamma|+\eta),$$
so the loss in the rate is smaller than $\varepsilon$ by the choice of $\rho$ and $\eta$. The cover, the number $m$ and the constant $K_\eta$ depend only on these choices and on the fixed data, so the remaining error is $O_\varepsilon(1)$ uniformly in the law and in $n$.
\end{proof}

\begin{lemma}\label{lem:conditional-jet-transport}
On a sufficiently small rectangle, let $\Phi$ be a fixed $C^1$ map with invertible derivative and sufficiently small variation of the derivative. If a law satisfies both exact conditional $s$-Frostman bounds with constant $C$, then its unnormalized image under $\Phi$ assigns to every rectangle with side lengths $a,b\le1$, rotated or not, mass at most $C_\Phi C^2a^sb^s$.
\end{lemma}

\begin{proof}
The statement is symmetric in the two side lengths, so we may assume that $a\le b$. The idea is to use the exact conditional bound in a coordinate with short fibers and the marginal bound in the other coordinate, which is why both directions of conditional regularity are needed.

Let $w$ be a unit normal to the long sides of the output rectangle. Since $D\Phi$ is invertible at the center of the chart, at least one coordinate derivative of $w\cdot\Phi$ has absolute value bounded below by a fixed positive constant, which can be chosen uniformly in $w$ by means of the smallest singular value of the derivative at the center. The small variation of the derivative preserves the sign and the lower bound on the source rectangle, so the preimage of the thin output strip has fibers of length $O_\Phi(a)$ in that coordinate. The same smallness assumption also provides a uniform inverse Lipschitz bound, because along a segment in the source rectangle the increment of $\Phi$ differs from that of its invertible linearization at the center by at most half the smallest singular value times the length of the segment. Since the output rectangle has diameter $O(b)$, the projection of its preimage to the other source coordinate has diameter $O_\Phi(b)$.

Suppose, for instance, that the selected fiber coordinate is $X$. Then the preimage lies over an interval $J$ of values of $Y$ of length $O_\Phi(b)$, and its fiber over $y$ is contained in an interval $I_y$ of length $O_\Phi(a)$. The exact conditional bound allows $I_y$ to depend on $y$, and we obtain
$$\begin{aligned}
&\mathbb P(\Phi(X,Y)\text{ lies in the rectangle and }(X,Y)\text{ lies in the chart})\\
&\le\int_J\mathbb P(X\in I_y\mid Y=y)d\mathcal L(Y)(y)\le C_\Phi C^2a^sb^s,
\end{aligned}$$
where the marginal $s$-Frostman bound for $Y$ comes from integrating the exact bound for $Y\mid X$. If a covering interval has length greater than one, the trivial bound by one and the fixed lower bound for its length give the same estimate after increasing $C_\Phi$. If the selected coordinate is $Y$, we interchange $X$ and $Y$ in this calculation. All derivative bounds are uniform in $w$, so the estimate also holds for rotated rectangles.
\end{proof}

\subsubsection{Entropy inequalities for conditional copies}

\begin{lemma}\label{lem:replica-ruzsa}
Let $U,U',V$ be finite-valued random variables in an abelian group, where $U$ and $U'$ are conditionally independent copies given $V$. Then
$$H(U-U')\le2H(U+V)-H(V)+2{\mathcal{I}}(U;V).$$
\end{lemma}

\begin{proof}
The difference can be recovered from the two sums, so we compare their joint entropy before using conditional independence. Let $Z=U+V$. The pairs $(U+V,U'+V)$ and $(Z,U-U')$ determine each other, hence
$$\begin{aligned}
H(U-U')&=H(Z,U-U')-H(Z\mid U-U')\\
&\le2H(U+V)-H(Z\mid U-U',U)=2H(U+V)-H(V\mid U,U').
\end{aligned}$$
The inequality uses subadditivity, the fact that $U+V$ and $U'+V$ have the same law, and the fact that conditioning on $U$ can only decrease $H(Z\mid U-U')$. The last equality holds because $(U-U',U)$ determines $(U,U')$ and because, once $U$ is fixed, $Z$ determines $V$. By conditional independence and the equality of the conditional laws,
$$H(U,U'\mid V)=2H(U\mid V)\text{ and }H(U,U')=2H(U)-{\mathcal{I}}(U;U'),$$
and subtracting, we find that ${\mathcal{I}}(V;(U,U'))=2{\mathcal{I}}(U;V)-{\mathcal{I}}(U;U')\le2{\mathcal{I}}(U;V)$. Therefore, $H(V\mid U,U')\ge H(V)-2{\mathcal{I}}(U;V)$, and substituting this bound proves the lemma.
\end{proof}

\begin{lemma}\label{lem:retained-copy-hereditary}
Fix a small buffered product chart on which $f_x\mathcal W_x\ne0$. Every event $A$ of original probability $a>0$ in this chart satisfies
\begin{equation}\label{eq:retained-copy-hereditary-new}
\begin{aligned}
&H_j(F\mid A)+2H_j(S\mid A)+2H_j(X\mid A)+H_j(Y\mid A)-2H_j(X,Y\mid A)\\
&\ge(s+\beta_s-\varepsilon)j-O_{\varepsilon,\mathrm{chart}}\left(1+\log\frac{C}{a}\right).
\end{aligned}
\end{equation}
The same holds with $X$ and $Y$ interchanged on a chart where $f_y\mathcal W_y\ne0$.
\end{lemma}

\begin{proof}
We first consider a law supported on the chart that satisfies both exact conditional bounds with constant $C$. The main step is the lower bound \eqref{eq:inverse-copy-retained-fibre} for the entropy of $(F,X-X')$, where $X'$ is a conditional copy of $X$, after which the Ruzsa inequality of Lemma~\ref{lem:replica-ruzsa} eliminates $X'$. Only after this elimination do we pass to arbitrary events by means of conditional cores. All intermediate rate losses are chosen sufficiently small in terms of $\varepsilon$, while the chart and the parameters chosen from it remain fixed.

Let $\ell=j+r$ and $\delta=2^{-\ell}$, where the nonnegative integer $r$ is fixed by the chart and chosen so large that the filled cells described next stay inside the buffer of the chart. We fill each cell of the input mesh uniformly, and independently in the two coordinates. By the fixed Lipschitz bounds, this changes the entropies of the observed source and output variables by $O(1)$, and it enlarges the exact conditional Frostman constants by a fixed factor only. Above scale $\delta$, this follows by averaging the original bounds over the conditioning cells and enlarging the test interval by $O(\delta)$. Below that scale, an interval of length $h\le\delta$ meets at most two cells, each of conditional mass $O(C\delta^s)$, and the uniform filling bounds its mass by $O(Ch\delta^{s-1})\le O(Ch^s)$.

We keep writing $X$ and $Y$ for the filled variables, whose conditional laws depend only on the conditioning cell at level $\ell$. Sampling $X$ and $X'$ independently given $Y$, we obtain a pair $(X',Y)$ with exactly the filled source law, and $D_\ell X$ and $D_\ell X'$ are independent and identically distributed given $D_\ell Y$, as required for the finite Ruzsa inequality below.

Since $f_x$ is bounded away from zero, the equation $z=f(x,y)$ can be solved as $x=\phi(y,z)$, and when choosing the buffered chart, we take the interior source rectangle so small that this inverse has a common domain in $(y,z)$ containing all sampled values of $(Y,F)$. The point $(F,X-X')$ lies on the curve
$$h_{x',y}(z)=\phi(y,z)-x'$$
with parameters $(x',y)=(X',Y)$, and implicit differentiation shows that $\phi_{yz}=\frac{\mathcal W_x}{f_x^3}$. The derivative of the map $(x',y)\mapsto(h_{x',y}(z),\partial_zh_{x',y}(z))$ is
$$\begin{pmatrix}-1&\phi_y\\0&\phi_{yz}\end{pmatrix},$$
whose determinant $-\phi_{yz}=-\frac{\mathcal W_x}{f_x^3}$ is bounded away from zero on the chart. The small variation of the derivative makes these value--slope maps uniformly bi-Lipschitz, and smoothness controls the $C^2$ differences of the curves by the distance of the parameters. After fixed affine normalizations, these are the geometric hypotheses of Theorem~\ref{thm:rectangular-incidence}. Lemma~\ref{lem:conditional-jet-transport}, applied to the actual parameter law of $(X',Y)$, provides the required bound for the mass of preimages of rectangles under each value--slope map. Conditionally on $(X',Y)$, the input has the law of $X\mid Y$, and the bi-Lipschitz map $x\mapsto f(x,y)$ transfers its interval bound to the graph coordinate.

We now apply the incidence theorem, keeping track of the conditional input entropy. Before counting, we replace parameters and inputs by representatives of their cells at level $\ell$. The geometry of the chart allows only a bounded ambiguity in the quantization of the outputs and finitely many separation colors. Let $u$, $c$ and $b$ be the integer parts of the negative logarithms of the probabilities of the parameter cell $(D_\ell X',D_\ell Y)$, of the conditional input cell $D_\ell X\mid D_\ell Y$ and of the observed output square of $(F,X-X')$, respectively. These probabilities are computed under the copied law before the mass and separation labels are attached, and we make no claim of conditional independence after restricting to a class. On a class of probability $w$, each copied triple has mass comparable to $2^{-u-c}$, and since the graph coordinate is bi-Lipschitz, only boundedly many input cells give the same incidence between a curve and a square. Consequently, the class contains at least $c_fw2^{u+c}$ incidences, at most $C_f2^u$ curves and at most $C_f2^b$ output squares.

The bounds for parameters and fibers provide the counting constants
$$C_1\le C_fC^22^u\delta^{2s}\text{ and }C_2\le C_fC2^c\delta^{t_s}.$$
For the first, we divide the bound for the parameter mass of rectangles by the mass of a selected parameter cell, which is comparable to $2^{-u}$, and for the second, we divide the conditional interval bound by the mass of a selected input cell, which is comparable to $2^{-c}$. These are bounds for selected subsets of the original reference laws, so no uniform Frostman bound for the normalized class is needed. In the second bound, the exponent $\sigma=\min\{2s,2-2s\}$ for the number of points comes from
$$h^s\lesssim\delta^{t_s}\left(\frac{h}{\delta}\right)^\sigma\text{ for }\delta\le h\le O(1),$$
and a selected cell together with its reference bound allows us to enlarge both constants to at least one. We apply Theorem~\ref{thm:rectangular-incidence} with rate loss $\eta>0$. Taking logarithms and substituting the counting constants, we get
$$\log w+u+c\le u+\frac{1}{3}c+\frac{2}{3}b-\frac{2s+t_s}{3}\ell+\eta\ell+O(1+\log C),$$
in which the level $u$ of the parameter masses cancels. Since $\beta_s=s+\frac{t_s}{2}$, choosing $\frac{3\eta}{2}\le\varepsilon$ yields
$$b\ge c+\beta_s\ell+\frac{3}{2}\log w-\varepsilon\ell-O(1+\log C).$$
The incidence threshold does not depend on the counting constants, so the finitely many coarser scales are covered by increasing the fixed additive error.

Let $L$ be the complete class label, including the separation colors. The information levels have mean $O(\ell)$, so $H(L)=O(\log(\ell+2))$. All expectations below are taken under the unrestricted copied law. By the definition of the levels, $\mathbb Ec=H(D_\ell X\mid D_\ell Y)+O(1)$ and $\mathbb Eb=H_\ell(F,X-X')+O(1)$, where the second identity also uses the bounded quantization ambiguity, while the average of $\log w$ is $-H(L)$. Averaging the class inequality and absorbing the logarithmic cost of the label into an arbitrarily small additional rate loss, we obtain, after renaming the tolerance,
$$H_\ell(F,X-X')\ge H(D_\ell X\mid D_\ell Y)+(\beta_s-\varepsilon)\ell-O_\varepsilon(1+\log C).$$
To return to level $j$, we use
$$H(D_\ell X\mid D_\ell Y)\ge H(D_jX\mid D_\ell Y)\ge H(D_jX\mid D_jY)-r,$$
and every other refinement costs at most $r$ times the number of coordinates. Since $r$ is fixed by the chart, this proves
\begin{equation}\label{eq:inverse-copy-retained-fibre}
H_j(F,X-X')\ge H(D_jX\mid D_jY)+(\beta_s-\varepsilon)j-O_\varepsilon(1+\log C).
\end{equation}

It remains to eliminate the copy. We apply Lemma~\ref{lem:replica-ruzsa} to $D_\ell X$, $D_\ell X'$ and $D_\ell Y$, which have the required conditional independence. Sums and differences of these cells determine the corresponding real sums and differences up to a bounded dyadic ambiguity, and the input entropies and their mutual information change by $O(r)$ between the levels $\ell$ and $j$. Returning from $\ell$ to $j$, we therefore get
$$H_j(F,X-X')\le H_j(F)+2H_j(S)-H_j(Y)+2{\mathcal{I}_j}(X;Y)+O(1).$$
We combine this with \eqref{eq:inverse-copy-retained-fibre} and with $H(D_jX\mid D_jY)\ge sj-O(1+\log C)$, and since
$$H_j(F)+2H_j(S)-H_j(Y)+2{\mathcal{I}_j}(X;Y)=H_j(F)+2H_j(S)+2H_j(X)+H_j(Y)-2H_j(X,Y),$$
we obtain \eqref{eq:retained-copy-hereditary-new} for the filled law, which the initial comparison with bounded ambiguity returns to the original law on the chart. For an arbitrary event $A$, we apply Lemma~\ref{lem:conditional-cores} and repeat the construction under the normalized law of each core, with its own exact conditional constant. On each core, the copy is eliminated before the entropies of the fixed source observables are averaged, and the bounded entropy of the core label together with the bound for the mean logarithm of the conditional constants produces the error $O_{\varepsilon,\mathrm{chart}}(1+\log\frac{C}{a})$. Interchanging $X$ and $Y$, we obtain the version for $f_y\mathcal W_y$.
\end{proof}

\begin{lemma}\label{lem:conditional-posterior-new}
On a fixed buffered chart regular for $\mathcal W_x$, every event $A$ contained in the chart and of original probability $a>0$ satisfies
$$H_j(S\mid A)+H_j(F\mid A)+\frac{3}{2}H_j(Y\mid A)\ge\left(\frac{5s}{2}+\beta_s-\varepsilon\right)j-O_{\varepsilon,\mathrm{chart}}\left(1+\log\frac{C}{a}\right).$$
The same holds with $X$ and $Y$ interchanged on a chart regular for $\mathcal W_y$.
\end{lemma}

\begin{proof}
Consider first a law on the fixed chart whose two exact conditional constants are bounded by $C$. We apply Lemma~\ref{lem:copied-output-incidence} with $U$ and $V$ the dyadic representatives of $X$ and $Y$, with $V'$ a copy of $V$ that is independent of $V$ conditionally on $U$, and with
$$\Phi(x,z)=x+z\text{ and }\Psi(x,y)=f(x,y).$$
Since $(U,V')$ and $(U,V)$ have the same joint law, the first copied output has the law of the sum of the representatives, while the second is the value of $f$ at the representatives, so their entropies agree with those of $S$ and $F$ up to bounded quantization errors. The calculation of the determinant and the rescaling in the proof of Proposition~\ref{prop:joint-depth-seed} show that these curves are transversal on a buffered chart regular for $\mathcal W_x$. The marginal reference constant and the constants for the laws of $X\mid Y$ are $O(C)$, because averaging over the conditioning cells preserves these bounds at the mesh scales. Discarding the nonnegative mutual information of the copies, Lemma~\ref{lem:copied-output-incidence} yields
$$H_j(S)+H_j(F)\ge H_j(X)-\frac{5}{2}{\mathcal{I}_j}(X;Y)+\beta_sj-\varepsilon j-O_\varepsilon(1+\log C),$$
where the smoothness of $f$ and the fixed buffer of the chart account for the quantization errors.

The two conditional entropy bounds replace the remaining dependence term by the marginal entropy appearing in the conclusion. The first shows that $H_j(X)-{\mathcal{I}_j}(X;Y)=H(D_jX\mid D_jY)\ge sj-\log C$, and the second that ${\mathcal{I}_j}(X;Y)\le H_j(Y)-sj+\log C$, so that
$$\begin{aligned}
H_j(X)-\frac{5}{2}{\mathcal{I}_j}(X;Y)&=(H_j(X)-{\mathcal{I}_j}(X;Y))-\frac{3}{2}{\mathcal{I}_j}(X;Y)\\
&\ge sj-\frac{3}{2}(H_j(Y)-sj)-O(1+\log C).
\end{aligned}$$
Substituting this bound proves the desired inequality, which has positive coefficients, for the law on the chart. Lemma~\ref{lem:conditional-cores}, applied to its fixed source observables, extends it to every event with the stated logarithmic cost in the probability, and interchanging $X$ and $Y$ covers the orientation of $\mathcal W_y$.
\end{proof}

\subsubsection{Expansion in inverse coordinates}

In the following local estimate, we assume that the joint Frostman constant $C_{\mathrm{joint}}$ is at least one, as we may by enlarging it.

\begin{lemma}\label{lem:inverse-output-seed}
Let $\Phi(u,v)$ and $\Psi(u,v)$ be smooth on a fixed buffered product chart, with $\Phi_u$ bounded away from zero. Suppose that the curves $g_{v,z}(t)=\Psi(u(t,z),v)$, where $\Phi(u(t,z),z)=t$, have a common buffered domain, uniformly bi-Lipschitz value--slope maps $(v,z)\mapsto(g_{v,z}(t),g'_{v,z}(t))$, and $C^2$ differences bounded by the distance of the parameters. If $(U,V)$ is a jointly $(s,s,C_{\mathrm{joint}})$-Frostman pair on the chart and $0<s\le\frac{2}{3}$, then
\begin{equation}\label{eq:inverse-output-seed}
\begin{aligned}
H_n(\Phi(U,V))+H_n(\Psi(U,V))\ge{}&H_n(U)-\frac{1}{2}H_n(V)+2sn-\frac{5}{2}{\mathcal{I}_n}(U;V)\\
&-\eta n-O_{\eta,\mathrm{chart}}(1+\log C_{\mathrm{joint}}).
\end{aligned}
\end{equation}
\end{lemma}

\begin{proof}
Lemma~\ref{lem:copied-output-incidence} already controls the loss caused by dependence, so we only need to derive its conditional reference constants from the joint rectangle bound. Let $\delta=2^{-n}$, and let $\widetilde U$ and $\widetilde V$ be representatives of the dyadic intervals containing $U$ and $V$, so that $H(\widetilde U)=H_n(U)$, $H(\widetilde V)=H_n(V)$ and ${\mathcal{I}}(\widetilde U;\widetilde V)={\mathcal{I}_n}(U;V)$. The fixed smoothness bounds compare the original outputs with those of the representatives up to a bounded error in entropy. The buffer contains all representatives at sufficiently fine scales, and the finitely many remaining scales are absorbed into the additive constant.

Write $q_j=\mathbb P(\widetilde V=v_j)$, keeping only the indices with $q_j>0$. For an interval $I$ of length at least $\delta$, the joint rectangle bound implies
$$\mathbb P(\widetilde U\in I\mid\widetilde V=v_j)\le\frac{KC_{\mathrm{joint}}\delta^s}{q_j}|I|^s,$$
where $K\ge1$ depends only on the chart and includes fixed enlargements of the cells. Covering the bounded range of $U$ by finitely many intervals of length at most one also provides an $s$-Frostman bound for the marginal of $V$ with constant $KC_{\mathrm{joint}}$, and we increase $K$ once so that it serves for both bounds. Since $q_j\le KC_{\mathrm{joint}}\delta^s$, we may choose $K_{v_j}=\frac{KC_{\mathrm{joint}}\delta^s}{q_j}\ge1$, and then
$$\mathbb E\log K_{\widetilde V}\le H_n(V)-sn+\log C_{\mathrm{joint}}+O(1),$$
where the expectation is taken with respect to the probabilities $q_j$. The individual constants $K_{v_j}$ may be large, but only their mean logarithm enters the estimate of Lemma~\ref{lem:copied-output-incidence}.

We now let $\widetilde V$ and $\widetilde V'$ be conditionally independent copies given $\widetilde U$, and we use $KC_{\mathrm{joint}}$ as the marginal reference constant and the $K_{v_j}$ as the conditional reference constants in Lemma~\ref{lem:copied-output-incidence}. The lower bound for $|\Phi_u|$ and the assumptions on the curves ensure its geometric hypotheses for
$$W=(\Phi(\widetilde U,\widetilde V'),\Psi(\widetilde U,\widetilde V)).$$
Since $s\le\frac{2}{3}$, the coefficient in that lemma is $\beta_s=\frac{3s}{2}$. Substituting the bound for the mean logarithm of the constants and discarding the nonnegative mutual information of the copies, we get
$$H_n(W)\ge H_n(U)-\frac{1}{2}H_n(V)+2sn-\frac{5}{2}{\mathcal{I}_n}(U;V)-\eta n-O_{\eta,\mathrm{chart}}(1+\log C_{\mathrm{joint}}).$$
Because the copies preserve the joint law with $\widetilde U$, the first coordinate of $W$ has the law of $\Phi(\widetilde U,\widetilde V)$ and the second has the law of $\Psi(\widetilde U,\widetilde V)$. By subadditivity, $H_n(W)$ is at most the sum of these two output entropies, and comparing them with the original outputs proves \eqref{eq:inverse-output-seed}. All constants are uniform in the law, since the marginal reference constant enters only through $\log(KC_{\mathrm{joint}})=\log C_{\mathrm{joint}}+O(1)$.
\end{proof}

\begin{proposition}\label{prop:conditional-inverse-output}
For every admissible polynomial, every $0<s\le\frac{2}{3}$ and every $\varepsilon>0$,
\begin{equation}\label{eq:conditional-inverse-output-weighted}
11H_n(S)+20H_n(F)\ge(35s-\varepsilon)n-O_\varepsilon(1).
\end{equation}
Consequently, $M_n\ge(\frac{35s}{31}-\varepsilon)n-O_\varepsilon(1)$, uniformly over the prescribed conditional source laws.
\end{proposition}

\begin{proof}
We apply Lemma~\ref{lem:inverse-output-seed} to the inverse of the map $(x,y)\mapsto(S,F)$. This produces an inequality containing ${\mathcal{I}_n}(S;F)$, which is then cancelled by means of Lemma~\ref{lem:conditional-posterior-new}. The resulting combination has positive coefficients for the outputs and can therefore be averaged over regular charts, and to make the cancellation possible, we keep a single normalized source event fixed until it has taken place.

To construct the inverse charts, we first check that the polynomials controlling invertibility and curvature are nonzero. Set
$$D_f=f_x-f_y\text{ and }N_f=f_y^2f_{xx}-2f_xf_yf_{xy}+f_x^2f_{yy}.$$
None of $f_x$, $f_y$ and $D_f$ vanishes identically, since otherwise $f$ would depend only on $y$, on $x$ or on $x+y$, respectively. Also $N_f\not\equiv0$, because the absolute curvature of a regular level curve of $f$ is
$$\frac{|N_f|}{(f_x^2+f_y^2)^{\frac{3}{2}}}.$$
Indeed, if $N_f$ vanished identically, every regular level curve would contain a straight segment, and restricting the polynomial to the line containing such a segment would show that $f$ is constant on the whole line. Fix such a line, on which $f$ takes the value $c_0$. A regular point where $f$ takes a different value lies on another line on which $f$ is constant, and this line must be parallel to the first one, because two nonparallel complete lines intersect and cannot carry different values. Hence, the derivative of $f$ in this fixed direction vanishes on the open set of regular points where $f\ne c_0$, which is nonempty because a nonconstant polynomial is not constant near any regular point. By polynomial continuation, this derivative vanishes identically, and $f$ has the excluded form, a contradiction.

On an inverse branch, we write $(x,y)=(\Phi(u,v),\Psi(u,v))$ with $(u,v)=(S,F)$. By implicit differentiation,
$$\begin{aligned}
\Phi_u&=-\frac{f_y}{D_f}, & \Phi_v&=\frac{1}{D_f},\\
\Psi_u&=\frac{f_x}{D_f}, & \Psi_v&=-\frac{1}{D_f},
\end{aligned}$$
where the derivatives of $f$ are evaluated at $(\Phi(u,v),\Psi(u,v))$. In particular, $\Phi_u\ne0$ where $f_yD_f\ne0$, as Lemma~\ref{lem:inverse-output-seed} requires. For that lemma, we set $g_{v,z}(t)=\Psi(u(t,z),v)$, where $\Phi(u(t,z),z)=t$, and we have to show that the value--slope maps have nonzero determinant on sufficiently small boxes of parameters. At diagonal parameters $v=z$, we have $t=x$ and $f(t,y)=z$, so that
$$\begin{aligned}
g_v&=-\frac{1}{D_f}, & g_z&=\frac{f_x}{D_ff_y},\\
\frac{g_z}{g_v}&=-\frac{f_x}{f_y}, & \partial_t&=\partial_x-\frac{f_x}{f_y}\partial_y.
\end{aligned}$$
Hence,
\begin{equation}\label{eq:inverse-output-jet}
\det\begin{pmatrix}g_v&g_z\\\partial_tg_v&\partial_tg_z\end{pmatrix}=g_v^2\partial_t\left(\frac{g_z}{g_v}\right)=-\frac{N_f}{D_f^2f_y^3}.
\end{equation}
At a point where $f_xf_yD_fN_f\ne0$, we choose an inverse product rectangle with a larger buffer and then shrink its interval of parameters, so that all sampled values $\Phi(u,z)$ lie in one common interval of values of $t$ on which $u(t,z)$ is defined for every $z$ in the buffer. By continuity, the determinant stays away from zero for all nearby pairs of parameters, including those with $v\ne z$. After a further fixed subdivision, the derivatives of the value--slope maps are uniformly close to a single invertible matrix, integration along segments of parameters yields the bi-Lipschitz bounds, and smoothness controls the $C^2$ differences by the distance of the parameters.

Using Lemma~\ref{lem:joint-critical-multiplicity}, we choose a nonzero orientation $\mathcal W_x$ or $\mathcal W_y$. Outside its zero set and the zero set of $f_xf_yD_fN_f$, we choose common buffered charts for the estimate in inverse coordinates and for Lemma~\ref{lem:conditional-posterior-new}. On each chart, we apply Lemma~\ref{lem:conditional-jet-transport} to the map $(x,y)\mapsto(S,F)$ under the original law. Restricting to a source event $A$ of probability $a>0$ can only decrease the unnormalized mass, and after dividing by $a$, the output law satisfies a joint rectangle bound with constant $\frac{K_{\mathrm{chart}}C^2}{a}$. No exact conditional bound in the output coordinates is needed.

Now fix an event $A$ of original probability $a>0$ inside one of these charts. {In the rest of the local calculation, all entropies are taken under the normalized law on $A$, and we suppress this conditioning in their notation. We keep the output mutual information $\mathcal{I}_n(S;F\mid A)$ until the final cancellation.} Since the source can be recovered from $(S,F)$ by a bi-Lipschitz map, $H_n(X,Y)=H_n(S,F)+O(1)$, and the two bounds in Lemma~\ref{lem:event-conditional-entropy} yield
$$H_n(X)+H_n(Y)\le2H_n(S,F)-2sn+O\left(1+\log\frac{C}{a}\right).$$
We apply Lemma~\ref{lem:inverse-output-seed} with inputs $(S,F)$, outputs $(X,Y)$ and reference constant $\frac{K_{\mathrm{chart}}C^2}{a}$. Its lower bound for $H_n(X)+H_n(Y)$, combined with this upper bound and with $H_n(S,F)=H_n(S)+H_n(F)-{\mathcal{I}_n(S;F\mid A)}$, leads to
$$2H_n(S,F)-2sn\ge H_n(S)-\frac{1}{2}H_n(F)+2sn-\frac{5}{2}{\mathcal{I}_n(S;F\mid A)}-\eta n-O_{\eta,\mathrm{chart}}\left(1+\log\frac{C}{a}\right),$$
and substituting for $H_n(S,F)$ and collecting the output terms, we obtain
\begin{equation}\label{eq:inverse-output-retained}
H_n(S)+\frac{5}{2}H_n(F)+\frac{1}{2}{\mathcal{I}_n(S;F\mid A)}\ge(4s-\eta)n-O_{\eta,\mathrm{chart}}\left(1+\log\frac{C}{a}\right).
\end{equation}
In the orientation of $\mathcal W_x$, Lemma~\ref{lem:conditional-posterior-new} on the same event yields
$$H_n(S)+H_n(F)+\frac{3}{2}H_n(Y)\ge(4s-\eta)n-O_{\eta,\mathrm{chart}}\left(1+\log\frac{C}{a}\right),$$
because $\beta_s=\frac{3s}{2}$ in the present range. Substituting $H_n(Y)\le H_n(S,F)-sn+O(1+\log\frac{C}{a})$ and multiplying by two, we find
$$5(H_n(S)+H_n(F))-3{\mathcal{I}_n(S;F\mid A)}\ge(11s-\eta)n-O_{\eta,\mathrm{chart}}\left(1+\log\frac{C}{a}\right),$$
where the tolerance has been renamed after the multiplication. In the orientation of $\mathcal W_y$, we use $X$ in the intermediate marginal bound, and the resulting inequality is the same.

Adding this inequality to six times \eqref{eq:inverse-output-retained}, the term ${\mathcal{I}_n(S;F\mid A)}$ of the same event cancels, and after renaming the tolerance, we obtain $11H_n(S)+20H_n(F)\ge(35s-\eta)n$ with the stated logarithmic cost of the event. Both output coefficients are positive, so Lemma~\ref{lem:conditional-chart-averaging}, applied to the finite list of nonzero critical polynomials found above, proves \eqref{eq:conditional-inverse-output-weighted} for the full law. Its left-hand side is at most $31M_n$, and dividing by $31$ and renaming the tolerance proves the final assertion.
\end{proof}

\section{Radon estimates and entropy inequalities}
\label{sec:fourier-comparisons}

The estimates of this section compare a prescribed output law with the output of a product reference law. We keep track of the Frostman constants and of the costs of restriction, since both the classes of marginal masses and the entropy iteration require them.

\subsection{Radon estimates for product measures}
\label{sec:fourier}

We follow the paired Radon argument of \cite[Section~5.4]{LPPS}, extending its Fourier shell estimate from equal to unequal exponents and keeping track of the symmetric and asymmetric powers of the Frostman constants. The analytic input is the standard estimate for nondegenerate oscillatory integral operators, in the formulation of Greenleaf and Seeger \cite{GS02}. As before, $\mathcal K_g=\partial_x\partial_y\ln\left|\frac{g_x}{g_y}\right|$ denotes the curvature expression of $g$ on the set where $g_xg_y\ne0$.

\begin{proposition}\label{prop:entropy-radon-shell}
Let $g$ be real analytic on a buffered enlargement of a compact product rectangle $I\times J$, with $g_xg_y\mathcal K_g\ne0$ there. Let $\nu$ and $\mu$ be positive measures on $I$ and $J$ of mass at most one, with $\nu(U)\le C_1|U|^a$ and $\mu(U)\le C_2|U|^b$ for every interval $U$, where $0<a,b\le1$ and $C_1,C_2\ge1$. Let $\omega=g_*(\nu\otimes\mu)$, let $\Xi\ge2$, and let $\vartheta\in C_c^\infty(\mathbb R\setminus\{0\})$ be nonnegative. Then
\begin{equation}\label{eq:entropy-radon-shell}
\int\vartheta\left(\frac{t}{\Xi}\right)|\widehat\omega(t)|^2dt\lesssim(C_1C_2)^{\frac{3}{2}}\Xi^{\frac{3}{2}-a-b},
\end{equation}
and if $b<1$, then also
\begin{equation}\label{eq:entropy-radon-shell-linear}
\int\vartheta\left(\frac{t}{\Xi}\right)|\widehat\omega(t)|^2dt\lesssim C_1^{\frac{3}{2}}C_2\Xi^{\frac{3}{2}-a-b}.
\end{equation}
Interchanging the coordinates, the factor $C_1C_2^{\frac{3}{2}}$ is admissible when $a<1$. The implicit constants depend on the chart and on $g$, $\vartheta$, $a$ and $b$. For charts in a fixed bounded box, they are at most $C C_{\mathrm{ch}}^N$, provided that the finitely many bounds for the derivatives of the phase and the cutoffs, for the reciprocals of the nondegeneracy quantities and for the reciprocals of the widths of the buffers used below are at most $C_{\mathrm{ch}}\ge2$. Here $C$ and $N$ depend only on the box, the exponents and the cutoff.
\end{proposition}

\begin{proof}
We first prove the estimates on a small product rectangle and then sum over a finite subdivision of the given chart. On a small rectangle, we split the fourfold integral expressing the Fourier energy according to the sizes of the two coordinate gaps. Where the second gap is at most a fixed multiple of the first, a paired oscillatory integral operator provides the principal estimate, and we bound the complementary region and the gaps below $\Xi^{-1}$ by $O(C_1C_2\Xi^{1-a-b})$, which is smaller than either asserted bound. For these error terms, we use the bounds for the total mass that depend on the gaps, while the normalization $\mu(J)\le1$ only improves the principal $L^2$ estimate, where it produces the linear power of $C_2$.

We use $\widehat\omega(t)=\int e^{-itv}d\omega(v)$ and set $\Phi=g(x,y)-g(x',y')$. The left-hand side is the integral of $K_\Xi(\Phi)$ with respect to $\nu^{\otimes2}\otimes\mu^{\otimes2}$, where for every integer $m\ge1$,
\begin{equation}\label{eq:entropy-radon-kernel}
K_\Xi(v)=\int e^{-itv}\vartheta\left(\frac{t}{\Xi}\right)dt\text{ satisfies }|K_\Xi(v)|\lesssim_m\Xi(1+\Xi|v|)^{-m}.
\end{equation}
On the small rectangle, we choose a smooth dyadic partition in $x'-x\ne0$ and treat its two signs separately. On a piece where $|x'-x|\asymp r$, we write $x'=x+rq$ and $y'=y+rz$, so that $|q|\asymp1$, and we set
$$\phi_r=\frac{\Phi}{r}=-\int_0^1(qg_x+zg_y)(x+trq,y+trz)dt.$$
The integral formula extends the phase smoothly to $r=0$. Since $\partial_z\phi_r=-g_y(x+rq,y+rz)$ is bounded away from zero, the implicit function theorem represents the zero set of $\phi_r$ as $z={\mathcal{H}}_r(x,q;y)$, and at $r=0$,
$${\mathcal{H}}_0=-\frac{qg_x}{g_y}\text{ and }\det D_{(x,q)}({\mathcal{H}}_0,\partial_y{\mathcal{H}}_0)=-q\left(\frac{g_x}{g_y}\right)^2\mathcal K_g.$$
This determinant is the nondegeneracy quantity for the paired operator below, and the condition $\mathcal K_g\ne0$ makes it nonzero at $r=0$. Compactness and continuity therefore yield a uniform lower bound, after a fixed localization, for all $0\le r\le r_0$ with $r_0>0$ fixed. We choose $C_0>\frac{4\sup|g_x|}{\inf|g_y|}$, increased if necessary to account for the fixed range of $q$, so that the roots ${\mathcal{H}}_0$ lie strictly inside $|z|<C_0$. After decreasing $r_0$, the same is true for ${\mathcal{H}}_r$, and all translated points with $|z|\le2C_0$ remain in the buffered chart. Finally, we take the first support interval shorter than $r_0$. These choices ensure that the operator below covers every piece of the small rectangle with $|y'-y|\le2C_0r$.

For $\Lambda=\Xi r\ge1$, consider
$$T_{\Lambda,r}h(x,q)=\iiint e^{-i\tau\phi_r}\vartheta\left(\frac{\tau}{\Lambda}\right)\chi(x,q;y,z)h(y,z)d\tau dydz,$$
where $\chi$ is a smooth spatial cutoff to the region just described, equal to one on the paired supports below. We factor $\phi_r=A_r(z-{\mathcal{H}}_r)$, where $A_r$ is the average of $\partial_z\phi_r$ along the segment from ${\mathcal{H}}_r$ to $z$, so that $|A_r|$ is bounded below. The substitution $\tau=\frac{\Lambda\sigma}{A_r}$ puts the oscillation in the standard form with large parameter $\Lambda$, with the prefactor $\Lambda$, the phase $\sigma(z-{\mathcal{H}}_r)$ and the amplitude $|A_r|^{-1}\vartheta(\frac{\sigma}{A_r})\chi$. Since $|A_r|$ is bounded above and below, the amplitude is smooth and supported in a fixed annulus in $\sigma$, and the phase is nondegenerate because $\partial_z\partial_\sigma[\sigma(z-{\mathcal{H}}_r)]=1$. Its augmented mixed Hessian
$$\begin{pmatrix}
-\sigma {\mathcal{H}}_{r,xy}&0&-{\mathcal{H}}_{r,x}\\
-\sigma {\mathcal{H}}_{r,qy}&0&-{\mathcal{H}}_{r,q}\\
-{\mathcal{H}}_{r,y}&1&0
\end{pmatrix}$$
has determinant $\sigma\det D_{(x,q)}({\mathcal{H}}_r,\partial_y{\mathcal{H}}_r)$, which is bounded away from zero on the support of the amplitude. Hence, the localized phase satisfies the hypotheses of \cite[Section~1.1, (1.5)--(1.6)]{GS02}, with two spatial variables on each side and one frequency variable, and the resulting norm is $O(\Lambda^{-\frac{3}{2}})$ before the prefactor $\Lambda$ is restored. Restoring it, we get
\begin{equation}\label{eq:entropy-radon-operator}
\|T_{\Lambda,r}\|_{2\to2}\lesssim\Lambda^{-\frac{1}{2}}.
\end{equation}

To pair this operator with the original measures, we push the first dyadic piece of $\nu\otimes\nu$ forward under $(x,x')\mapsto(x,\frac{x'-x}{r})$ and call the result $m_r$. For the second pair, we multiply $\mu\otimes\mu$ by a cutoff with values in $[0,1]$ that equals one where $|y'-y|\le C_0r$ and vanishes where $|y'-y|\ge2C_0r$, push the product forward under $(y,y')\mapsto(y,\frac{y'-y}{r})$, and call the result $n_r$. Then, for $0<\rho\le1$,
$$\begin{aligned}
\|m_r\|&\lesssim C_1r^a, & m_r(B(u,\rho))&\lesssim C_1^2r^a\rho^{2a},\\
\|n_r\|&\le\min\{1,CC_2r^b\}, & n_r(B(v,\rho))&\lesssim C_2^2r^b\rho^{2b}.
\end{aligned}$$
Indeed, a ball of radius $\rho$ confines one of the two variables to an interval of length $O(\rho)$ and, once this variable is fixed, the other one to an interval of length $O(r\rho)$, so the two Frostman bounds imply the local estimates, while the restriction of the gaps controls the total masses. The additional inequality $\|n_r\|\le1$ comes from the normalization $\mu(J)\le1$.

Next, we bound the $L^2$ norms needed to pair the operator with $m_r$ and $n_r$, which need not have densities. Let $\Pi_\Lambda$ be a smooth annular Fourier multiplier at frequency $\Lambda$ with a real even symbol, whose annulus will be chosen to contain the spatial frequencies of the operator. If $\lambda$ is a finite positive measure on $\mathbb R^2$ with $\lambda(B(u,\rho))\le M\rho^\beta$ for some $0<\beta\le2$, we choose $m>2$ and integrate the kernel bound $C_m\Lambda^2(1+\Lambda|u-v|)^{-m}$ over dyadic annuli to obtain
$$\|\Pi_\Lambda\lambda\|_2^2\lesssim\|\lambda\|M\Lambda^{2-\beta}.$$
Applying this with the preceding bounds for balls and total masses, we find
\begin{equation}\label{eq:entropy-radon-paired-L2}
\|\Pi_\Lambda m_r\|_2\lesssim C_1^{\frac{3}{2}}r^a\Lambda^{1-a}\text{ and }\|\Pi_\Lambda n_r\|_2\lesssim C_2^{\frac{3}{2}}r^b\Lambda^{1-b},
\end{equation}
and, combining $\|n_r\|\le1$ with the ball bound for $n_r$ instead of the total mass depending on the gap,
\begin{equation}\label{eq:entropy-radon-paired-linear}
\|\Pi_\Lambda n_r\|_2\lesssim C_2r^{\frac{b}{2}}\Lambda^{1-b}.
\end{equation}
The last bound, in which $C_2r^{\frac{b}{2}}$ replaces $C_2^{\frac{3}{2}}r^b$, is the source of both the better power of the constant and the different exponent of the gap below.

To use the operator norm on these singular paired measures, we first localize the spatial frequencies of the operator. Since $\partial_q\phi_r=-g_x(x+rq,y+rz)$ and $\partial_z\phi_r=-g_y(x+rq,y+rz)$, both spatial frequencies have magnitude comparable to $\Lambda$, the nonvanishing derivatives in $q$ and $z$ giving the lower bounds and the remaining bounded derivatives of the phase giving the upper bounds. We choose the symbol of $\Pi_\Lambda$ equal to one on both ranges. Outside the chosen annuli, integration by parts in the spatial variables bounds the Fourier transform of the kernel by $O_m(\Lambda(\Lambda+|\xi|+|\eta|)^{-m})$, and integrating over the four Fourier variables, with $m$ increased as needed, we get
\begin{equation}\label{eq:entropy-radon-frequency}
T_{\Lambda,r}=\Pi_\Lambda T_{\Lambda,r}\Pi_\Lambda+E_{\Lambda,r},\text{ where }\|K_{E_{\Lambda,r}}\|_\infty\lesssim_m\Lambda^{-m}.
\end{equation}
This uniform bound for the kernel of the remainder allows us to pair the remainder with $m_r$ and $n_r$. For the principal term, we combine \eqref{eq:entropy-radon-operator} with \eqref{eq:entropy-radon-paired-L2}, or with \eqref{eq:entropy-radon-paired-linear} for the second measure. Since $K_\Xi(\Phi)=r^{-1}K_\Lambda(\phi_r)$, the principal contribution of the gap $r$ is bounded in the two cases by
$$\begin{aligned}
(C_1C_2)^{\frac{3}{2}}r^{-1}\Lambda^{-\frac{1}{2}}r^{a+b}\Lambda^{2-a-b}&=(C_1C_2)^{\frac{3}{2}}\Xi^{\frac{3}{2}-a-b}r^{\frac{1}{2}},\\
C_1^{\frac{3}{2}}C_2r^{-1}\Lambda^{-\frac{1}{2}}r^{a+\frac{b}{2}}\Lambda^{2-a-b}&=C_1^{\frac{3}{2}}C_2\Xi^{\frac{3}{2}-a-b}r^{\frac{1-b}{2}}.
\end{aligned}$$
Both expressions are summed over dyadic scales of the gap. The first sum is bounded uniformly in $\Xi$ because the exponent of $r$ is $\frac{1}{2}>0$, and the second because $\frac{1-b}{2}>0$. At $b=1$, the second sum would produce a logarithmic loss, which explains the restriction in \eqref{eq:entropy-radon-shell-linear}. This proves the two asserted bounds for the principal terms on the small rectangle.

We now bound the remaining terms. For the frequency remainders, we use the bounds $\|m_r\|\lesssim C_1r^a$ and $\|n_r\|\lesssim C_2r^b$, and taking $m>a+b$, we obtain
$$C_mC_1C_2\sum_{r\ge\Xi^{-1}}r^{a+b-1}(\Xi r)^{-m}\lesssim C_1C_2\Xi^{1-a-b}.$$
Outside the cutoff of the second pair, $|y'-y|$ is at least $C_0|x'-x|$ up to fixed dyadic constants, and the fixed sign of $g_y$ and the choice of $C_0$ ensure that $|\Phi|\gtrsim|y'-y|$. For $|y'-y|\asymp d\ge\Xi^{-1}$, the first gap is $O(d)$, so \eqref{eq:entropy-radon-kernel} leads to the same error,
$$C_mC_1C_2\sum_{d\ge\Xi^{-1}}\Xi(\Xi d)^{-m}d^{a+b}\lesssim C_1C_2\Xi^{1-a-b},$$
and the bounded overlap of the partition of the first gap prevents an additional sum over the first gap. If both gaps are $O(\Xi^{-1})$, their mass is $O(C_1C_2\Xi^{-a-b})$, and since $|K_\Xi|\lesssim\Xi$, we obtain the same bound. If the first gap is smaller than $\Xi^{-1}$ but the second is larger than a fixed multiple of $\Xi^{-1}$, the preceding estimate for the complementary region applies. These cases exhaust the gaps not treated by the operator, and the exact diagonal in the first coordinate has measure zero because $a>0$. Since $\Xi\ge2$ and $C_1,C_2\ge1$, the error $O(C_1C_2\Xi^{1-a-b})$ is bounded by either of the asserted right-hand sides.

For the given chart, let $\mathcal R$ be a finite subdivision into small product rectangles of the kind just considered, and apply the estimates to the unnormalized restrictions of the marginals, whose Frostman constants are still $C_1$ and $C_2$ and whose masses are at most one. If $\omega_B$ denotes the image measure of the piece over $B\in\mathcal R$, then
$$\left|\sum_{B\in\mathcal R}\widehat\omega_B\right|^2\le|\mathcal R|\sum_{B\in\mathcal R}|\widehat\omega_B|^2,$$
and integrating against the nonnegative cutoff proves both bounds on the original chart, with the same powers of $C_1$ and $C_2$.

It remains to check that the constants depend polynomially on the bounds for the chart, including the bound for the operator norm. All the bounds from inverse functions and from integration by parts used above are fixed powers of $C_{\mathrm{ch}}$. The derivative bounds and the lower bound for the determinant show that the smallest singular value of the augmented mixed Hessian is at least $C_{\mathrm{ch}}^{-O(1)}$, so after a subdivision into $C_{\mathrm{ch}}^{O(1)}$ patches of diameter $C_{\mathrm{ch}}^{-O(1)}$, the linear separation of the phase dominates its Taylor remainder. For each localized piece, we keep the notation $T_{\Lambda,r}$ and remove the prefactor by setting $\widetilde T_{\Lambda,r}=\Lambda^{-1}T_{\Lambda,r}$. The argument in \cite[Section~1.1, proof following (1.6)]{GS02}, applied to the adjoint, shows that
$$|K_{\widetilde T_{\Lambda,r}\widetilde T_{\Lambda,r}^*}(u,u')|\le CC_{\mathrm{ch}}^N\Lambda^{-1}(1+\Lambda|u-u'|)^{-3}\text{ for }u,u'\in\mathbb R^2.$$
The integral of $(1+\Lambda|u-u'|)^{-3}$ in either variable is $O(\Lambda^{-2})$, so Schur's test yields $\|\widetilde T_{\Lambda,r}\widetilde T_{\Lambda,r}^*\|\le CC_{\mathrm{ch}}^N\Lambda^{-3}$. Taking square roots and restoring the prefactor, we obtain $\|T_{\Lambda,r}\|_{2\to2}\le CC_{\mathrm{ch}}^N\Lambda^{-\frac{1}{2}}$ with the required polynomial dependence on the chart, and summing the operator norms over the $C_{\mathrm{ch}}^{O(1)}$ pieces preserves this dependence. The cutoffs, the subdivisions and the sums over the gaps contribute only fixed powers of $C_{\mathrm{ch}}$, and interchanging the coordinates proves the final assertion.
\end{proof}

\subsection{Entropy inequalities from output Frostman bounds}
\label{sec:entropy-profiles-latest}

A Frostman bound for the output improves the entropy lower bound obtained from the product reference measure. We prove this improvement first, and then we average over quantitative charts, keeping track of the bounded and logarithmic errors.

\begin{lemma}\label{lem:output-profile-transfer}
Let $\rho$ and $\omega$ be probability measures on a fixed bounded interval, where $\rho$ is $\zeta$-Frostman with $0<\zeta<1$ and
\begin{equation}\label{eq:profile-reference-shell}
\int_{\Xi\le|t|\le2\Xi}|\widehat\omega(t)|^2dt\le C\Xi^{-\epsilon}\text{ for all }\Xi\ge1.
\end{equation}
Let $\rho_n$ and $\omega_n$ be their dyadic probability vectors at level $n$. If $0<\epsilon<1-\zeta$ and $\lambda_*=\frac{2(1-\zeta)}{1-\zeta+\epsilon}$, then
\begin{equation}\label{eq:profile-transfer}
H(\rho_n)+\lambda_*D(\rho_n\Vert\omega_n)\ge n-\frac{\lambda_*}{2}\log(D(\rho_n\Vert\omega_n)+2)-O(1),
\end{equation}
and for every $\lambda>\lambda_*$, we have $H(\rho_n)+\lambda D(\rho_n\Vert\omega_n)\ge n-O(1)$. If $\epsilon=1-\zeta$, then
\begin{equation}\label{eq:profile-transfer-critical}
H(\rho_n)+D(\rho_n\Vert\omega_n)\ge n-\log(n+1)-O(1),
\end{equation}
and if $\epsilon>1-\zeta$, the logarithmic term can be omitted. For fixed exponents and a fixed support interval, the additive constants depend at most logarithmically on the Frostman and Fourier constants.
\end{lemma}

\begin{proof}
Write $\rho_n=(p_i)$ and $\omega_n=(q_i)$. The conclusion is trivial if $D(\rho_n\Vert\omega_n)=\infty$, so we assume that it is finite. We bound the cross entropy $-\sum_ip_i\log q_i=H(\rho_n)+D(\rho_n\Vert\omega_n)$ from below through an estimate for the probability that independent samples from the two laws fall in the same dyadic interval. Grouping the bins of comparable mass under $\rho$, we obtain two energy bounds, and splitting the frequency sum where they agree determines $\lambda_*$. Finally, we control the entropy of the group label by the entropy deficit $n-H(\rho_n)$, so that the grouping does not cause an uncontrolled loss.

Let $\delta=2^{-n}$ and $\beta=1-\zeta$. For each $k\ge0$, we group the bins of positive mass with $2^{-k-1}<p_i\le2^{-k}$, and we let $w_k$ be the mass of the group, $\rho^{(k)}$ the normalized restriction of $\rho$ to these bins and $K$ the label of the group, omitting empty groups. For intervals with $|I|\gtrsim\delta$, we have
$$\rho^{(k)}(I)\lesssim w_k^{-1}\min\{|I|^\zeta,2^{n-k}|I|\},$$
where the first bound comes from the Frostman estimate for $\rho$ and the second from counting the $O(2^n|I|)$ bins that meet $I$. To turn these interval bounds into energy bounds, we majorize the indicator of the frequency interval by a Gaussian. For an integer $N>2$, Fourier inversion produces a kernel bounded by $C_N\Xi(1+\Xi|x-y|)^{-N}$, and integrating over dyadic annuli in the distance, we get
\begin{equation}\label{eq:profile-two-energy-bounds}
\int_{|t|\le\Xi}|\widehat{\rho^{(k)}}(t)|^2dt\lesssim w_k^{-1}\min\{\Xi^\beta,2^{n-k}\}\text{ for }1\le\Xi\lesssim\delta^{-1}.
\end{equation}
These are estimates for the original measures, and no bin masses have been moved.

Points in the same bin are at distance at most $\delta$ from each other, and we majorize this collision event by a nonnegative Schwartz function at scale $\delta$ whose Fourier transform has compact support. Applying Fourier inversion and then the Cauchy--Schwarz inequality on each dyadic shell, with \eqref{eq:profile-two-energy-bounds} for $\rho^{(k)}$ and \eqref{eq:profile-reference-shell} for $\omega$, we obtain
$$\sum_{i\text{ in group }k}\frac{p_i}{w_k}q_i\lesssim\delta+\delta w_k^{-\frac{1}{2}}\sum_{1\le2^j\lesssim\delta^{-1}}\min\{2^{\frac{\beta j}{2}},2^{\frac{n-k}{2}}\}2^{-\frac{\epsilon j}{2}}.$$
Suppose first that $0<\epsilon<\beta$, and set $c=\frac{\beta-\epsilon}{2\beta}\in(0,\frac{1}{2})$. For $k\le n$, we split the sum at $2^{\beta j}=2^{n-k}$, where the two energy bounds agree. Below this threshold, the terms increase geometrically, and above it, they decrease geometrically, so the sum is $O(2^{c(n-k)})$, and if the threshold lies beyond the last shell, extending the first geometric sum up to the threshold only increases the bound. For $k>n$, we use the second energy bound throughout, and the sum is $O(2^{\frac{n-k}{2}})$, which is $O(2^{c(n-k)})$ because $c<\frac{1}{2}$. Since there are $O(2^n)$ bins, $w_k\lesssim2^{n-k}$, so $w_k^{-\frac{1}{2}}2^{c(n-k)}$ is bounded below when $k>n$, and the same holds for $k\le n$ because $w_k\le1$. In both cases, this bound absorbs the initial low-frequency term $\delta$, and we arrive at
\begin{equation}\label{eq:profile-cross-collision}
\sum_{i\text{ in group }k}\frac{p_i}{w_k}q_i\lesssim w_k^{-\frac{1}{2}}2^{-n+c(n-k)}.
\end{equation}
Jensen's inequality for $-\log$ now yields the cross-entropy bound for each group,
$$-\sum_{i\text{ in group }k}\frac{p_i}{w_k}\log q_i\ge(1-c)n+ck+\frac{1}{2}\log w_k-O(1).$$
On the group $k$, the information $-\log p_i$ lies in $[k,k+1)$, so that $\mathbb EK\le H(\rho_n)<\mathbb EK+1$. Averaging the preceding inequality with the weights $w_k$ and using $\sum_kw_k\log w_k=-H(K)$, we obtain
\begin{equation}\label{eq:profile-label-entropy}
(1-c)H(\rho_n)+D(\rho_n\Vert\omega_n)\ge(1-c)n-\frac{1}{2}H(K)-O(1).
\end{equation}

We now control the cost of the label by the amount by which $H(\rho_n)$ falls below $n$, which allows us to keep every bin of positive mass, however small. Since there are $O(2^n)$ bins,
$$\begin{aligned}
\mathbb P(K\ge n+j)&\lesssim2^{-j},\\
\mathbb E(K-n)_+&=O(1),\\
\mathbb E|K-n|&\le n-H(\rho_n)+O(1),
\end{aligned}$$
where $t_+=\max\{t,0\}$ denotes the positive part of $t$. Comparing the law of $K$ with a two-sided geometric distribution of scale $1+\mathbb E|K-n|$, we find $H(K)\le\log(n-H(\rho_n)+C_0)+O(1)$, where $C_0$ is fixed so that the logarithm is positive. With $J=n-H(\rho_n)$, the inequality \eqref{eq:profile-label-entropy} becomes
$$J\le\lambda_*D(\rho_n\Vert\omega_n)+\frac{1}{2}\lambda_*\log(J+C_0)+O(1),\text{ where }\lambda_*=(1-c)^{-1}=\frac{2\beta}{\beta+\epsilon}$$
is the positive coefficient in the statement. The bound on the support shows that $J$ is bounded below, and the logarithmic term is at most $\frac{J}{2}+O(1)$, so moving it to the left first shows that $J\lesssim D(\rho_n\Vert\omega_n)+1$, and substituting this bound back into the logarithm proves \eqref{eq:profile-transfer}. For $\lambda>\lambda_*$, we use the inequality $\frac{1}{2}\lambda_*\log(t+2)\le(\lambda-\lambda_*)t+O(1)$ for $t\ge0$, with $t=D(\rho_n\Vert\omega_n)$, to obtain the bound with bounded error.

It remains to treat the case $\epsilon\ge\beta$, in which no grouping of masses is needed and we apply the Frostman energy bound directly to $\rho$. The same collision majorant shows that
$$\sum_ip_iq_i\lesssim\delta\left(1+\sum_{1\le2^j\lesssim\delta^{-1}}2^{\frac{(\beta-\epsilon)j}{2}}\right).$$
At $\epsilon=\beta$, the sum has $O(n+1)$ terms, and for $\epsilon>\beta$, it is bounded, so Jensen's inequality proves \eqref{eq:profile-transfer-critical} and its version with bounded error. Before logarithms are taken, all constants are fixed powers of the given constants, which proves the last assertion.
\end{proof}

\begin{remark}\label{rem:profile-quantization}
The entropy estimates below rely on two different facts. Data processing controls the relative entropy when both inputs are sent through the same channel, while a small displacement allows us to compare the entropy of the actual channel output with the entropy of the output we are interested in. We record both facts here.

Fix an event $A$ of positive probability and an input resolution $m\ge n$. Let $\mu$ be either the law of $(D_mU,D_mV)$ under the normalized restriction to $A$, or the law obtained by spreading the mass of each occupied input cell uniformly over the cell, and write $\mu^\otimes$ for the product of the marginals of $\mu$. We send both laws through the same, possibly randomized, channel with real output, and quantize at level $n$, obtaining probability vectors $p$ and $q$. Then
$$D(p\Vert q)\le D(\mu\Vert\mu^\otimes)={\mathcal{I}_m}(U;V\mid A).$$
If the actual output of the channel can be coupled with $W$ under the restricted law so that the displacement is at most $c2^{-n}$, then
$$|H(p)-H_n(W\mid A)|\le\log(2\lceil c\rceil+3).$$
Consequently, for every $\lambda\ge0$,
$$H_n(W\mid A)+\lambda {\mathcal{I}_m}(U;V\mid A)\ge H(p)+\lambda D(p\Vert q)-O_c(1),$$
and in particular Jensen's inequality implies the following lower bound in terms of the collision probability,
$$H_n(W\mid A)+{\mathcal{I}_m}(U;V\mid A)\ge-\log\sum_ip_iq_i-O_c(1).$$
To prove the identity for the relative entropy, let $(a_{ij})$ be the conditional table of the inputs, with row masses $a_i$ and column masses $b_j$. On every occupied cell, the density of $\mu$ with respect to $\mu^\otimes$ equals $\frac{a_{ij}}{a_ib_j}$, both for the discrete table and for its uniform filling, so integrating the logarithm of this ratio produces exactly ${\mathcal{I}_m}(U;V\mid A)$, and data processing yields the inequality for the outputs. For the comparison of entropies, the displacement coupling leaves at most $2\lceil c\rceil+3$ possible bins in either direction, which bounds both conditional entropies as stated. Finally, $H(p)+D(p\Vert q)=-\sum_ip_i\log q_i$, and Jensen's inequality for $-\log$ leads to the bound in terms of the collision probability. Thus, the entropy is controlled by the bounded displacement, while the bound for the relative entropy comes from the use of a common channel and not from the displacement.

The global entropy estimate uses representatives of the input cells. Its discrete reference vector need not be the dyadic vector of the continuous reference $\omega=G_*(\mu_U^A\otimes\mu_V^A)$, where $\mu_U^A$ and $\mu_V^A$ are the restricted input marginals. However, these two reference outputs can be coupled with displacement $O(2^{-n})$ once $m-n$ is large enough in terms of the Lipschitz constant of the chart, so the collision argument uses the enlarged event $|x-y|\le C2^{-n}$ and the same Fourier majorant, and the shell estimate is required for the continuous measure $\omega$ and not for displaced atoms. The coupling of the actual output also bounds the mass of each bin by $O(\delta^\zeta)$, where $\delta=2^{-n}$, when the original output is $\zeta$-Frostman. Filling those bins uniformly extends this bound to smaller scales, because an interval of length $r<\delta$ meets at most two bins and each of them contributes $\frac{p_ir}{\delta}\lesssim\delta^{\zeta-1}r\le r^\zeta$. In the multiscale argument, on the other hand, $\mu$ is the filled law, both outputs receive the same independent noise, and the reference is $g_*\mu^\otimes$ convolved with this noise. The Fourier estimates apply to $g_*\mu^\otimes$, and the convolution can only decrease the modulus of its Fourier transform.
\end{remark}

\begin{lemma}\label{lem:quantitative-chart-partition}
Let $G$ be a nonspecial polynomial on a fixed bounded box. The complement of a fixed algebraic exceptional set has a countable partition into product rectangles $B_j$ with integer depths $k_j\ge0$ such that the following hold.
\begin{enumerate}[label=\textup{(\roman*)}]
\item Each $B_j$ has a buffered enlargement satisfying the hypotheses of Proposition~\ref{prop:entropy-radon-shell}.
\item At most $2^{C(k+1)}$ rectangles have depth $k$.
\item The chart constants of $B_j$ are at most $2^{C(k_j+1)}$, with the powers of the Frostman constants given in Proposition~\ref{prop:entropy-radon-shell}, and the input quantization shift of $B_j$, that is, the number of additional dyadic levels at which the inputs must be quantized so that the representatives stay in the buffer and the outputs move by $O(2^{-n})$, is at most $C(k_j+1)$.
\item If{, for some fixed $\theta>0$,} the law satisfies
$\mathbb P(\operatorname{dist}((U,V),Z)\le r)\lesssim r^\theta$
{for every $0<r\le1$ and} for the fixed algebraic curves
$Z$ involved, then the chart label $L$ satisfies
$\mathbb E k_L+H(L)<\infty$, uniformly for fixed constants
in this tube bound.
\end{enumerate}
The exceptional set may include the zero sets of any fixed finite list of nonzero polynomials, and the constants depend only on the fixed polynomial data and on the exponents. The statement also holds after either of the two fixed sum shears.
\end{lemma}

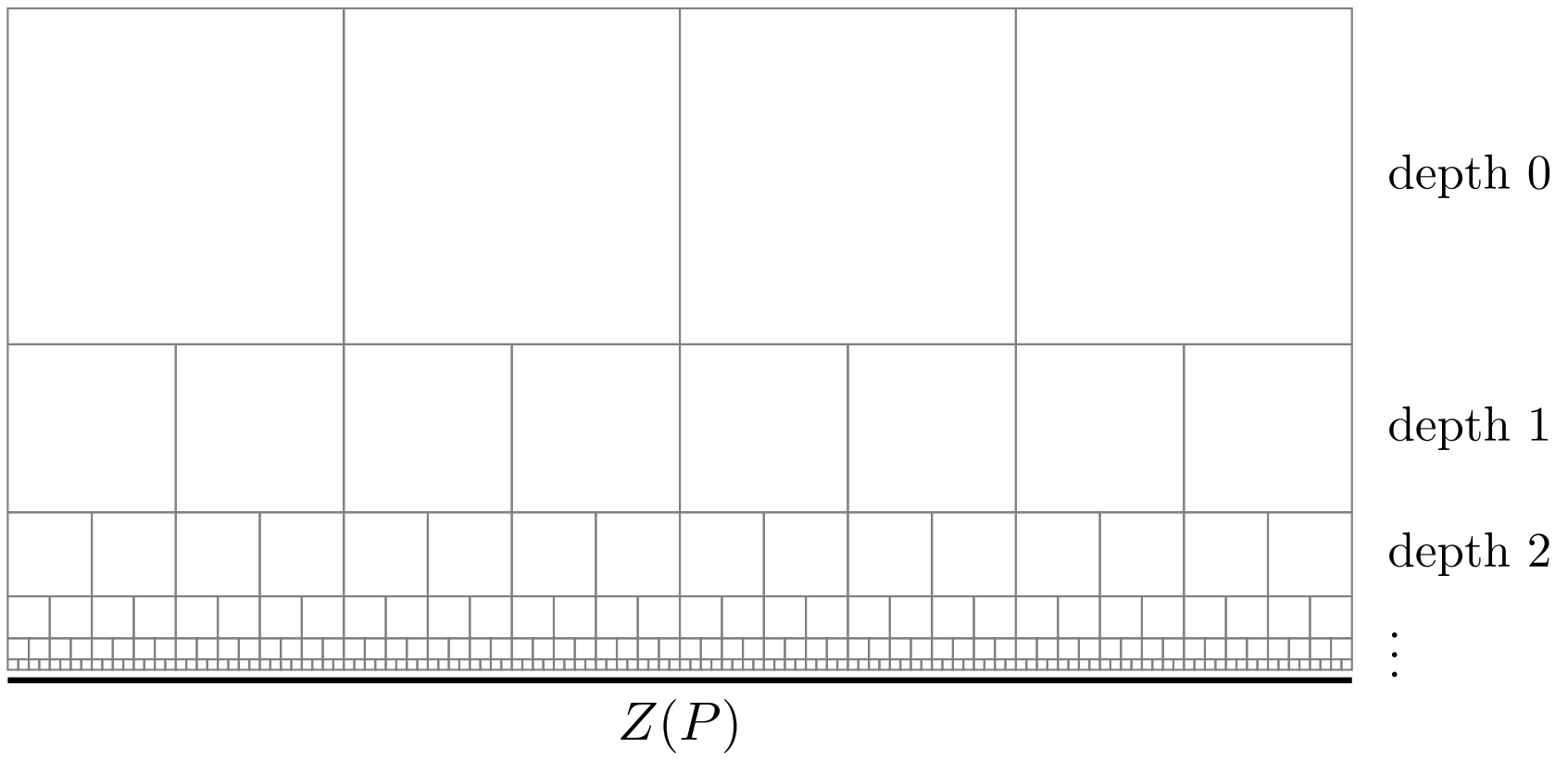
\begin{figure}[t]
\centering
\begin{tikzpicture}[scale=2.4]
\foreach \k in {0,1,2,3,4,5} {
  \pgfmathsetmacro{\side}{2^(-\k)}
  \pgfmathsetmacro{\ylow}{2^(-\k)}
  \pgfmathtruncatemacro{\num}{4*2^(\k)-1}
  \foreach \i in {0,...,\num} {
    \draw[gray] (\i*\side,\ylow) rectangle (\i*\side+\side,\ylow+\side);
  }
}
\draw[very thick] (0,0) -- (4,0);
\node[below] at (2,0) {$Z(P)$};
\node[right] at (4.05,1.5) {\small depth $0$};
\node[right] at (4.05,0.75) {\small depth $1$};
\node[right] at (4.05,0.375) {\small depth $2$};
\node[right] at (4.05,0.12) {\small $\vdots$};
\end{tikzpicture}
\caption{The Whitney decomposition in the proof of Lemma~\ref{lem:quantitative-chart-partition}, drawn near a smooth piece of the exceptional set $Z(P)$, here a horizontal segment. A square of depth $k$ has side $2^{-k}$ and lies at distance comparable to $2^{-k}$ from $Z(P)$. The constants of the charts it contains grow at most exponentially in $k$, while the tube bound makes the probability of depth at least $k$ decay like $2^{-\theta k}$.}
\label{fig:whitney}
\end{figure}

\begin{proof}
We choose the size of each chart according to its distance from the exceptional set. This bounds the analytic constants by an exponential in the depth, while the tube bound makes the probability of large depth decay exponentially. Let $P$ be the product of $G_u$, $G_v$, the numerator obtained by clearing denominators in $\mathcal K_G$, and the additional exceptional polynomials. All these factors are nonzero because $G$ is nonspecial. On the fixed neighborhood of the support, the polynomial {\L}ojasiewicz inequality \cite[Corollary~2.6.7]{BCR} shows that
$$|P(u,v)|\ge c\min\{1,\operatorname{dist}((u,v),Z(P))\}^N.$$
We take Whitney squares in the complement of $Z(P)$, with disjoint half-open representatives, so that each square has distance from $Z(P)$ comparable to its side length, as in Figure~\ref{fig:whitney}. A square of side $2^{-k}$ receives depth $k$, and the finitely many larger scales in the fixed box receive depth zero. A region separated from $Z(P)$ is covered by finitely many rectangles of depth zero, which also takes care of the case where the exceptional set is empty.

We subdivide each square of depth $k$ into product squares of side comparable to $2^{-N'k}$, where the fixed integer $N'$ is large enough that the buffered enlargements of these squares remain in the regular set. The lower bound above controls the reciprocal of $P$, and since all factors of $P$ are bounded above on the fixed neighborhood, it also controls the reciprocal of each factor by a fixed power of $2^k$. Together with the bounds for the derivatives of the polynomials, this bounds every inverse derivative and every reciprocal determinant used in the Radon estimate by $2^{C(k+1)}$, and the reciprocal widths of the buffers and the derivatives of the cutoffs satisfy bounds of the same form. After increasing $C$, there are at most $2^{C(k+1)}$ pieces of depth $k$. Proposition~\ref{prop:entropy-radon-shell} now yields (i)--(iii), since its constants depend polynomially on these quantities and the powers of the Frostman constants are unchanged, and an additional $O(k+1)$ input levels make the mesh smaller than the buffer and than the required tolerance for the output.

Each rectangle inherits the depth of its parent Whitney square, so every rectangle of depth at least $k$ lies within distance $O(2^{-k})$ of $Z(P)$. The tube bound first shows that $Z(P)$ has measure zero, so that $L$ is defined almost surely, and then yields $\mathbb P(k_L\ge k)\lesssim2^{-\theta k}$, which implies $\mathbb Ek_L+H(k_L)<\infty$. Given the depth, the chart label takes at most $2^{C(k_L+1)}$ values, so by the chain rule for entropy,
$$H(L)\le H(k_L)+\mathbb E\log(2^{C(k_L+1)})<\infty.$$
A fixed sum shear maps the exceptional sets to algebraic sets and changes distances by bounded factors, so applying the same construction in the sheared coordinates proves the last assertion.
\end{proof}

\begin{proposition}\label{prop:global-profile-correlation}
Let $F=G(U,V)$ be a nonspecial representation in the original or in one of the sum coordinate systems, where $U$ and $V$ satisfy Frostman bounds with exponents $a$ and $b$ such that $a+b>\frac{3}{2}$. Assume the tube bound of Lemma~\ref{lem:quantitative-chart-partition}, and assume that on each regular chart $B_j$ of probability $p_j$, the conditional law of the output satisfies a Frostman bound with exponent $\zeta\in(0,1)$ and constant at most $\frac{2^{C(k_j+1)}}{p_j}$. Set $\epsilon=a+b-\frac{3}{2}$. If $0<\epsilon<1-\zeta$, then
\begin{equation}\label{eq:global-profile-correlation}
H_n(F)+\frac{2(1-\zeta)}{1-\zeta+\epsilon}{\mathcal{I}_n}(U;V)\ge n-O(\log(n+1)),
\end{equation}
and if the coefficient of ${\mathcal{I}_n}(U;V)$ is replaced by any strictly larger number, the error becomes bounded. If $\epsilon=1-\zeta$, the inequality with coefficient one holds with error $O(\log(n+1))$, and if $\epsilon>1-\zeta$, it holds with coefficient one and a bounded error.
\end{proposition}

\begin{proof}
We apply the local inequality of Lemma~\ref{lem:output-profile-transfer} on each quantitative chart, with the same exponents $\epsilon$ and $\zeta$ on every chart. The constants of the charts and the probabilities of the restrictions vary, but they enter the local estimate only through logarithms, and after averaging, the errors of the charts, the refinements of the inputs and the costs of conditioning all have finite means.

Let $L$ be the label of Lemma~\ref{lem:quantitative-chart-partition} and write $p_j=\mathbb P(L=j)>0$. On the event $\{L=j\}$, we use the conditional law as the actual input and the product of its two marginals as the reference input. Restriction increases each marginal Frostman constant by a factor of at most $p_j^{-1}$. We apply Proposition~\ref{prop:entropy-radon-shell} to the two reference marginals, with exponents $a$ and $b$, so that the frequency exponent is $\frac{3}{2}-a-b=-\epsilon$ and the constant is at most $2^{C(k_j+1)}p_j^{-C}$, while the Frostman constant of the actual output is at most $\frac{2^{C(k_j+1)}}{p_j}$ by assumption. We choose a nonnegative integer $r_j\le C(k_j+1)$ so large that the input representatives stay in the buffered chart and that changing each input by at most $2^{-n-r_j}$ changes the output by $O(2^{-n})$, with a uniform implied constant, which the derivative bounds of the chart allow, as in Remark~\ref{rem:profile-quantization}. We then apply data processing to the tables of cells of the actual and of the product reference law under the same deterministic map
$$(D_{n+r_j}U,D_{n+r_j}V)\longmapsto D_nG(D_{n+r_j}U,D_{n+r_j}V).$$
Remark~\ref{rem:profile-quantization} provides the bound for the relative entropy of the outputs with the exact input cost ${\mathcal{I}_{n+r_j}}(U;V\mid L=j)$. The actual output differs from $D_nF$ by a bounded number of bins, so its entropy differs from $H_n(F\mid L=j)$ by $O(1)$, and the same remark allows these bounded displacements in the collision argument while the Fourier estimate is applied to the continuous reference measure. Lemma~\ref{lem:output-profile-transfer} therefore applies with the exponents $\epsilon$ and $\zeta$. The part of its additive error that comes from the constants of the chart and of the measures is $O(k_j+1+\log\frac{1}{p_j})$, because these constants enter only through logarithms, while the explicit loss in terms of the relative entropy or of the scale remains as in the lemma.

The bounds for the mean depth and for the entropy of the label in Lemma~\ref{lem:quantitative-chart-partition} make these errors summable. Refining each scalar input by $r_j$ levels increases the mutual information by at most $2r_j$, and conditioning costs at most $H(L)$, so that
$$\sum_jp_j{\mathcal{I}_{n+r_j}}(U;V\mid L=j)\le {\mathcal{I}_n}(U;V\mid L)+2\mathbb Er_L\le {\mathcal{I}_n}(U;V)+H(L)+2\mathbb Er_L.$$
By the concavity of entropy, $H_n(F)\ge\sum_jp_jH_n(F\mid L=j)$. Hence, the average of the left-hand sides of the local inequalities is at most the global left-hand side plus the finite costs of conditioning and refinement just estimated. In the subcritical case, Jensen's inequality moves the average over the charts inside the remaining logarithm of the relative entropy, and the preceding bound for the mutual information makes this term $O(\log({\mathcal{I}_n}(U;V)+2))=O(\log(n+1))$, since the inputs are bounded. This proves \eqref{eq:global-profile-correlation}. For a strictly larger coefficient, we use the part of the local lemma with bounded error before averaging, at $\epsilon=1-\zeta$, we use \eqref{eq:profile-transfer-critical}, and above this threshold, we use its version with bounded error. The same summable costs of the charts appear in all cases.
\end{proof}

\begin{corollary}\label{cor:profile-conditional}
If $f$ is admissible and $(X,Y)$ is a conditionally $s$-Frostman pair with $\frac{5}{6}\le s<1$, then
\begin{equation}\label{eq:profile-conditional-main}
M_n\ge\frac{1+2s}{3}n-O(\log(n+2)),
\end{equation}
and the error is $O(1)$ when $s>\frac{5}{6}$.
\end{corollary}

\begin{proof}
We apply Proposition~\ref{prop:global-profile-correlation} with $a=b=\zeta=s$ in a nonspecial representation, which exists by Corollary~\ref{cor:three-coordinate}. The original marginals are $s$-Frostman, and Lemma~\ref{lem:conditional-shear} transfers the same marginal bounds to both sum coordinate systems. The polynomial sublevel estimates provide the required tube bounds, also for the pullbacks of the exceptional sets under the fixed shears. To verify the hypothesis on the output, we include the relevant partial derivatives in the original coordinates in the exceptional set. Lemma~\ref{lem:regular-output} then yields the output exponent $s$ on each regular chart, by applying the original conditional bounds to the output event first. Intersecting with the chart can only reduce the unnormalized mass, so the normalization costs only the reciprocal of the probability of the chart, and this step does not require the restriction to a chart to preserve the exact conditional constants.

The shell exponent is $\epsilon=2s-\frac{3}{2}$, and $\epsilon\ge1-s$ exactly when $s\ge\frac{5}{6}$. Therefore, Proposition~\ref{prop:global-profile-correlation} shows that
$$H_n(F)+{\mathcal{I}_n}(U;V)\ge n-O(\log(n+2))$$
for each nonspecial representation, with error $O(1)$ when $s>\frac{5}{6}$. It remains to eliminate the input mutual information from this inequality, and how we do so depends on whether the nonspecial representation uses the coordinates $(X,Y)$ or a pair of sum coordinates. In the original coordinates, adding the two conditional entropy lower bounds, we obtain
$$H_n(X)+H_n(Y)\le2H_n(X,Y)-2sn+O(1).$$
Admissibility excludes $f=Q(x+y)$, so Theorem~\ref{thm:coordinate-entropy} applies to $(S,F)$ and yields $H_n(X,Y)\le H_n(S)+H_n(F)+O(1)\le2M_n+O(1)$, and hence
$${\mathcal{I}_n}(X;Y)\le2M_n-2sn+O(1).$$
In sum coordinates, the linear entropy identity and the conditional entropy bounds for the original pair lead instead to
$$\begin{aligned}
{\mathcal{I}_n}(X;S)&=H_n(S)-H(D_nY\mid D_nX)+O(1)\le M_n-sn+O(1),\\
{\mathcal{I}_n}(Y;S)&=H_n(S)-H(D_nX\mid D_nY)+O(1)\le M_n-sn+O(1).
\end{aligned}$$
Substituting the relevant bound for the mutual information and using $H_n(F)\le M_n$, we find $3M_n\ge(1+2s)n-O(\log(n+2))$ for a nonspecial original representation and $2M_n\ge(1+s)n-O(\log(n+2))$ for a nonspecial sum representation. The second bound is the stronger one, because $\frac{1+s}{2}-\frac{1+2s}{3}=\frac{1-s}{6}\ge0$, and by the classification, at least one of the two cases occurs. All errors are bounded when $s>\frac{5}{6}$, which proves the corollary together with its assertion at the endpoint.
\end{proof}

\begin{corollary}\label{cor:profile-joint}
If $(X,Y)$ is jointly $s$-Frostman with $\frac{9}{10}<s<1$ and a sum-coordinate representation of $f$ is nonspecial, then $M_{\mathrm j,n}\ge\frac{(1+3s)n}{4}-O(1)$.
\end{corollary}

\begin{proof}
Suppose first that $G_Y$ is nonspecial. We apply Proposition~\ref{prop:global-profile-correlation} in the coordinates $(Y,S)$ and then bound its mutual information term in terms of $M_{\mathrm j,n}$. The joint Frostman bound implies an $s$-Frostman bound for $Y$, and covering a strip $\{X+Y\in J\}$ by $O(|J|^{-1})$ squares of side $|J|$ produces the exponent $2s-1$ for $S$. On regular charts, Lemma~\ref{lem:regular-output}, applied in the original coordinates, provides the same exponent $\zeta=2s-1$ for $F$. The tube bound follows from the covering argument behind \eqref{eq:joint-tube}, since the sheared law retains the isotropic $2s$-Frostman bound that this argument requires. Thus, Proposition~\ref{prop:global-profile-correlation} applies with $a=s$, $b=2s-1$ and $\epsilon=3s-\frac{5}{2}$. Since $s>\frac{9}{10}$, we have $\epsilon>1-\zeta=2-2s$, and the case with bounded error shows that
$$H_n(F)+{\mathcal{I}_n}(Y;S)\ge n-O(1).$$
To bound the dependence term, note that $H_n(X)+H_n(Y)\le2M_{\mathrm j,n}$ by the definition of $M_{\mathrm j,n}$, while $H_n(X)\ge sn-O(1)$ by the marginal bound, so that $H_n(Y)\le2M_{\mathrm j,n}-sn+O(1)$. Also $H_n(Y,S)=H_n(X,Y)+O(1)\ge2sn-O(1)$, and consequently
$${\mathcal{I}_n}(Y;S)\le(2M_{\mathrm j,n}-sn)+M_{\mathrm j,n}-2sn+O(1)=3M_{\mathrm j,n}-3sn+O(1).$$
Substituting this bound and $H_n(F)\le M_{\mathrm j,n}$ into the preceding inequality, we obtain $4M_{\mathrm j,n}\ge(1+3s)n-O(1)$. If $G_X$ is nonspecial, we interchange $X$ and $Y$ throughout.
\end{proof}

\subsection{Entropy inequalities from marginal mass grouping}
\label{subsec:joint-mass-classes}

Grouping the bins of the original marginals according to their masses strengthens the output bounds available under joint Frostman regularity. The asymmetric estimate \eqref{eq:entropy-radon-shell-linear} keeps the resulting loss linear in sum coordinates and reduces it further in the original coordinates, and this loss can then be averaged by means of the original marginal entropy.

\begin{proposition}\label{prop:joint-mass-classes}
Let $(X,Y)$ be jointly $(s,s,C)$-Frostman with bounded support, where $\frac{1}{2}<s<1$, and set $S=X+Y$ and $F=f(X,Y)$. If the polynomial $G_Y(u,v)=f(v-u,u)$ is nonspecial, then
\begin{equation}\label{eq:joint-mass-class-comparison}
H_n(F)+{\mathcal{I}_n}(Y;S)+H_n(X)\ge\left(s+\min\left\{1,\frac{3s}{2}-\frac{1}{4}\right\}\right)n-O(\log(n+2)),
\end{equation}
and if $G_X$ is nonspecial, the same conclusion holds with $X$ and $Y$ interchanged. The error is uniform over laws with fixed support and Frostman constant. Independently of these cases, if $f$ itself is nonspecial, then
\begin{equation}\label{eq:joint-original-mass-class-comparison}
H_n(F)+{\mathcal{I}_n}(X;Y)+\frac{1}{2}H_n(X)\ge\left(\frac{s}{2}+\min\left\{1,\frac{3s}{2}-\frac{1}{4}\right\}\right)n-O(\log(n+2)),
\end{equation}
and the same estimate holds with $X$ and $Y$ interchanged.
\end{proposition}

\begin{proof}
We treat the sum coordinates and the original coordinates together, writing
$$(U,V,g,\theta)=
\begin{cases}
(Y,S,G_Y,1) & \text{if }G_Y\text{ is nonspecial},\\
(X,Y,f,\frac{1}{2}) & \text{if }f\text{ is nonspecial},
\end{cases}$$
so that $g(U,V)=F$ in both cases. We restrict the law to quantitative charts and group the dyadic intervals of the original variable $X$ according to their masses, and on each class, we compare the actual output with the output of a product reference law. The two choices of coordinates require different powers of the constant coming from the class of masses, and $\theta$ records this difference in the collision estimate. After averaging, it becomes the coefficient of $H_n(X)$. Constants denoted by $C_0$ depend only on the fixed data and may increase from line to line.

We apply Lemma~\ref{lem:quantitative-chart-partition} in the chosen coordinates. In the sum case, we add to the exceptional polynomial the polynomial $(G_Y)_u+(G_Y)_v=f_y(v-u,u)$, which is nonzero, since otherwise $f$ would depend only on $x$ and $G_Y$ would be special. The covering argument behind \eqref{eq:joint-tube} applies after the shear, because a fixed linear map preserves the isotropic $2s$-Frostman bound, so the hypotheses of Lemma~\ref{lem:quantitative-chart-partition} are satisfied. We obtain product charts $B_j$, a label $L$, and buffered enlargements on which the original derivative $f_y$ has a fixed sign and is bounded away from zero. Let $C_{\mathrm{ch},j}\ge2$ control the geometric constants of the chart $B_j$, including $|f_y|^{-1}$ and the factor of the chart in \eqref{eq:entropy-radon-shell-linear}. We choose shifts $r_j\ge0$ of the mesh such that the cells at level $n+r_j$ meeting $B_j$ stay inside its buffer, and such that filling these input cells uniformly and adding output noise of size $2^{-n-r_j}$ moves $F$ by at most $2^{-n}$. The bounds of Lemma~\ref{lem:quantitative-chart-partition} imply
\begin{equation}\label{eq:mass-class-chart-cost}
\log C_{\mathrm{ch},j}+r_j\le C_0(k_j+1)\text{ and }\mathbb E\log C_{\mathrm{ch},L}+\mathbb Er_L+H(L)\le C_0.
\end{equation}
The factor of the chart does not depend on the Frostman constants, whose powers in the asymmetric Fourier estimate are unchanged. The geometry of the charts therefore contributes logarithmic costs with bounded mean, and it remains to keep track of the probability of each class and of the bound for the masses of its intervals in $X$.

We group the dyadic intervals of $X$ at level $n$ according to their masses, before restricting to the charts. If $p_i=\mathbb P(D_nX=i)$, we define $K$ by $2^{-K-1}<p_{D_nX}\le2^{-K}$, and we set
$$\begin{aligned}
E_{j\ell}&=\{L=j,K=\ell\}, & w_{j\ell}&=\mathbb P(E_{j\ell}),\\
m_j&=n+r_j, & \eta_j&=2^{-m_j},
\end{aligned}$$
using only classes of positive probability. Since $\mathbb EK\le H_n(X)\le n+C_0$ and an integer-valued random variable with mean $t$ has entropy at most $\log(t+1)+\log e$, we have
\begin{equation}\label{eq:mass-class-label-entropy}
H(L,K)=O(\log(n+2)).
\end{equation}
The marginal Frostman bound implies $K\ge sn-C_0$, and since $(K-sn)_+\le K-sn+C_0$, it follows that
\begin{equation}\label{eq:mass-class-average-loss}
\mathbb E(K-sn)_+\le H_n(X)-sn+C_0.
\end{equation}
In view of \eqref{eq:mass-class-chart-cost}, \eqref{eq:mass-class-label-entropy} and \eqref{eq:mass-class-average-loss}, it suffices to prove on each class an entropy inequality whose loss is $\theta(\ell-sn)_+$ plus logarithmic costs for the chart and for the probability of the class, since these average to $\theta(H_n(X)-sn)$ and to $O(\log(n+2))$, respectively.

Fix a class $E_{j\ell}$ and write $w$ for its probability. A class of masses contains at most $2^{\ell+1}$ intervals of $X$ at level $n$, so this event meets at most $2^{\ell+r_j+1}$ intervals of $X$ at level $m_j$. For one such interval $I$ and an interval $J$ with $|J|\ge\eta_j$, we first bound the part of $\{X\in I,F\in J\}$ lying in the chart. In the original variables, the chart is either a product rectangle or the preimage of one under a linear shear, so it is convex, and on it $f_y$ has a fixed sign and $|f_y|\ge C_{\mathrm{ch},j}^{-1}$. The segment joining two points of this part stays in the chart, and integrating along it, the term involving the derivative in $y$ has no cancellation, while the term involving the derivative in $x$ can be moved to the other side, which leads to
$$C_{\mathrm{ch},j}^{-1}|y_1-y_2|\le|f(x_1,y_1)-f(x_2,y_2)|+\|f_x\|_\infty|x_1-x_2|\le|J|+\|f_x\|_\infty\eta_j.$$
Thus, this part lies in a rectangle with sides $\eta_j$ and $O(C_{\mathrm{ch},j}|J|)$, and in the sum case, the part where $S\in J$ has the same property, because then $|y_1-y_2|\le|J|+\eta_j$. Applying the original joint rectangle bound to these rectangles, summing over at most $2^{\ell+r_j+1}$ fine intervals of $X$, and dividing by $w$, we obtain
\begin{equation}\label{eq:mass-class-strip}
\mathbb P(Z\in J\mid E_{j\ell})\le\frac{C_{\mathrm{ch},j}^{C_0}2^{(1-s)r_j+(\ell-sn)_+}}{w}|J|^s,
\end{equation}
where $Z=F$, or $Z=S$ in sum coordinates. Here we used $2^{\ell+r_j}\eta_j^s=2^{\ell-sn+(1-s)r_j}$. After the restriction, the original marginals of $X$ and $Y$ keep $s$-Frostman constants $\frac{C_0}{w}$.

We now construct the actual and the reference laws used in the collision estimate. A continuous realization of the finite table of inputs preserves its exact relative entropy and makes the Frostman bounds available below the mesh of the inputs, so that the Fourier estimate and data processing can be applied to the same pair of laws. We select a pair of intervals at level $m_j$ according to the joint table of $(U,V)$ on the class, and conditionally on this pair, we choose independent uniform points in the two intervals. Let $\mu$ be the resulting joint law and $\mu^\otimes$ the product of its marginals. Both have densities that are constant on the input cells, with the same factors within the cells, so
\begin{equation}\label{eq:mass-class-input-kl}
D(\mu\Vert\mu^\otimes)={\mathcal{I}_{m_j}}(U;V\mid E_{j\ell}).
\end{equation}
The uniform filling extends each available marginal Frostman bound to all scales. An interval of length $t<\eta_j$ meets at most two cells, each contributing at most its mass times $\frac{t}{\eta_j}$, and since the masses of the cells are bounded by the original constant times $\eta_j^s$, the result is $O(t\eta_j^{s-1})\le O(t^s)$ times that constant. An interval of length $t\ge\eta_j$, we enlarge to include its two boundary cells, and we use \eqref{eq:mass-class-strip} for $S$ and the original marginal Frostman bounds for $X$ and $Y$. The two reference constants are therefore at most
\begin{equation}\label{eq:mass-class-reference-constants}
\begin{cases}
\frac{C_{\mathrm{ch},j}^{C_0}}{w}\text{ and }\frac{C_{\mathrm{ch},j}^{C_0}2^{(1-s)r_j+(\ell-sn)_+}}{w} & \text{if }\theta=1,\\
\frac{C_0}{w}\text{ and }\frac{C_0}{w} & \text{if }\theta=\frac{1}{2},
\end{cases}
\end{equation}
and the supports lie in the buffered product chart.

We apply $g$ and then independent uniform noise on $[-\eta_j,\eta_j]$ to both $\mu$ and $\mu^\otimes$, and we write $\rho$ and $\omega$ for the resulting laws and $a_i$ and $b_i$ for their masses on the dyadic intervals at level $n$. Data processing and the coupling with $F$ yield
\begin{equation}\label{eq:mass-class-output-kl}
D((a_i)\Vert(b_i))\le {\mathcal{I}_{m_j}}(U;V\mid E_{j\ell})
\end{equation}
and
\begin{equation}\label{eq:mass-class-output-entropy}
H((a_i))=H_n(F\mid E_{j\ell})+O(1).
\end{equation}
The bounds for the reference inputs do not by themselves give the regularity we need for the actual output $\rho$, and we deduce it from \eqref{eq:mass-class-strip}, using the displacement for intervals above the input mesh and the final convolution for smaller intervals. For intervals of length at least $\eta_j$, the coupling with the original $F$ enlarges the interval by $O(C_{\mathrm{ch},j}\eta_j)$, so \eqref{eq:mass-class-strip} applies after increasing the power of $C_{\mathrm{ch},j}$. For an interval $I$ of length $t<\eta_j$, the final uniform convolution shows that
$$\begin{aligned}
\rho(I)&\le\frac{t}{2\eta_j}\sup_z(g_*\mu)([z-\eta_j,z+\eta_j])\le\frac{C_{\mathrm{ch},j}^{C_0}2^{(1-s)r_j+(\ell-sn)_+}}{w}t\eta_j^{s-1}\\
&\le\frac{C_{\mathrm{ch},j}^{C_0}2^{(1-s)r_j+(\ell-sn)_+}}{w}t^s.
\end{aligned}$$
Thus, $\rho$ is $s$-Frostman with constant $\frac{C_{\mathrm{ch},j}^{C_0}2^{(1-s)r_j+(\ell-sn)_+}}{w}$, and a Gaussian majorant in frequency and a summation over annuli in the distance lead to
\begin{equation}\label{eq:mass-class-actual-energy}
\int_{|\xi|\le\Xi}|\widehat\rho(\xi)|^2d\xi\le\frac{C_{\mathrm{ch},j}^{C_0}2^{(1-s)r_j+(\ell-sn)_+}}{w}\Xi^{1-s}\text{ for }\Xi\ge2.
\end{equation}
For the reference energy, we apply \eqref{eq:entropy-radon-shell-linear} with $a=b=s$ to the two uniformized marginals with the constants in \eqref{eq:mass-class-reference-constants}, where the restriction $b<1$ holds by hypothesis. In sum coordinates, we assign the linear power of the constant to the marginal of $S$, and in the original coordinates, both constants equal $\frac{C_0}{w}$. The final convolution has a Fourier transform bounded by one, and therefore
\begin{equation}\label{eq:mass-class-reference-energy}
\int_{\Xi\le|\xi|\le2\Xi}|\widehat\omega(\xi)|^2d\xi\le\frac{C_{\mathrm{ch},j}^{C_0}2^{(2\theta-1)((1-s)r_j+(\ell-sn)_+)}}{w^{\frac{5}{2}}}\Xi^{\frac{3}{2}-2s}.
\end{equation}
In both cases, the normalization is $w^{-\frac{3}{2}}w^{-1}=w^{-\frac{5}{2}}$, and the factor coming from the class of masses is $2^{(\ell-sn)_+}$ in sum coordinates and absent in the original coordinates, exactly as recorded by $2^{(2\theta-1)(\ell-sn)_+}$.

It remains to turn the two energy bounds into an entropy inequality on the class. We choose a nonnegative Schwartz function $\psi$ with $\psi\ge1$ on $[-1,1]$ and with compactly supported Fourier transform, for instance a fixed multiple of the square of the inverse Fourier transform of a real, even, smooth function supported sufficiently close to zero. With $\delta=2^{-n}$, Fourier inversion and the Cauchy--Schwarz inequality on each shell yield
\begin{equation}\label{eq:mass-class-collision}
\begin{aligned}
\sum_ia_ib_i&\le\iint\psi\left(\frac{u-v}{\delta}\right)d\rho(u)d\omega(v)\\
&\le\frac{C_{\mathrm{ch},j}^{C_0}2^{\theta((1-s)r_j+(\ell-sn)_+)}}{w^{\frac{7}{4}}}2^{-n}\left(1+\sum_{\substack{\Xi\text{ dyadic}\\2\le\Xi\le C_02^n}}\Xi^{\frac{5}{4}-\frac{3s}{2}}\right).
\end{aligned}
\end{equation}
The contribution of the low frequencies is included in the initial term $1$, since the prefactor is at least one. On the dyadic shells, the geometric mean of the two energy constants is $C_{\mathrm{ch},j}^{C_0}2^{\theta((1-s)r_j+(\ell-sn)_+)}w^{-\frac{7}{4}}$, and the exponent of the shell is $\frac{1}{2}(1-s)+\frac{1}{2}(\frac{3}{2}-2s)=\frac{5}{4}-\frac{3s}{2}$. Hence, the sum is $O(2^{(\frac{5}{4}-\frac{3s}{2})n})$ for $s<\frac{5}{6}$, it is $O(n+2)$ at $s=\frac{5}{6}$, and it is bounded for $s>\frac{5}{6}$. By Jensen's inequality and \eqref{eq:mass-class-output-kl},
$$H((a_i))+{\mathcal{I}_{m_j}}(U;V\mid E_{j\ell})\ge\sum_ia_i\log\frac{1}{b_i}\ge-\log\sum_ia_ib_i,$$
and using \eqref{eq:mass-class-output-entropy}, we conclude that
\begin{equation}\label{eq:mass-class-local-comparison}
\begin{aligned}
H_n(F\mid E_{j\ell})+{\mathcal{I}_{m_j}}(U;V\mid E_{j\ell})\ge{}&\min\left\{1,\frac{3s}{2}-\frac{1}{4}\right\}n-\theta(\ell-sn)_+\\
&-C_0\left(1+\log C_{\mathrm{ch},j}+r_j+\log\frac{1}{w}\right)-O(\log(n+2)).
\end{aligned}
\end{equation}
This is the inequality we were looking for, since the only loss coming from the class of masses is $\theta(\ell-sn)_+$. The last logarithmic term comes from the number of shells at $s=\frac{5}{6}$ and is bounded when $s\ne\frac{5}{6}$.

Finally, we average \eqref{eq:mass-class-local-comparison}. By the concavity of entropy, the average of the conditional output entropies is at most $H_n(F)$. Refining both scalar inputs costs at most $2r_j$ in mutual information, and conditioning costs at most the entropy of the full label, so
$$\sum_{j,\ell}w_{j\ell}{\mathcal{I}_{n+r_j}}(U;V\mid E_{j\ell})\le {\mathcal{I}_n}(U;V)+H(L,K)+2\mathbb Er_L.$$
The cost of the probabilities of the classes averages exactly to $\sum_{j,\ell}w_{j\ell}\log\frac{1}{w_{j\ell}}=H(L,K)$, so \eqref{eq:mass-class-chart-cost} and \eqref{eq:mass-class-label-entropy} bound all the costs of the charts, of the refinements and of the normalizations by $O(\log(n+2))$. This averaging is legitimate also for countably many classes, because the costs of the charts and of the labels are integrable, while the bounded support bounds the local entropies by $O(n+r_j+1)$, whose mean is finite. Using \eqref{eq:mass-class-average-loss}, we obtain
$$H_n(F)+{\mathcal{I}_n}(U;V)\ge\min\left\{1,\frac{3s}{2}-\frac{1}{4}\right\}n-\theta(H_n(X)-sn)-O(\log(n+2)).$$
Moving $\theta H_n(X)$ to the left-hand side, the term $\theta sn$ on the right produces the asserted coefficient, with $\theta=1$ in sum coordinates and $\theta=\frac{1}{2}$ in the original coordinates. Interchanging $X$ and $Y$, including the choice of the original marginal that defines the classes of masses, we obtain the two symmetric versions.
\end{proof}

\section{Multiscale entropy inequalities}
\label{sec:multiscale-estimates}

In this section, we iterate the local entropy inequality, keeping track of the marginal entropies and of the masses of the bins at the final scale. The resulting bounds yield the remaining branches of Theorems~\ref{thm:joint-main} and~\ref{thm:conditional-main}.

\subsection{Local entropy inequality on regular charts}
\label{sec:multiscale-bootstrap}

In the Fourier argument, we now keep track of the masses of the bins at the scales of the frequencies and of the separations. Iterating the resulting entropy inequality improves the initial bound up to the coefficient needed in the main theorems.

\begin{theorem}\label{thm:local-multiscale-profile}
Let $1<\alpha\le2$, and let $g$ be analytic on a buffered product chart on which $g_ug_v\partial_u\partial_v\ln\left|\frac{g_u}{g_v}\right|\ne0$. Suppose that $(U,V)$ is supported in the core of the chart and satisfies $\mathbb P((U,V)\in B(z,r))\le Cr^\alpha$ for $0<r\le1$, where $C\ge1$, and set $W=g(U,V)$. For every $\varepsilon>0$,
\begin{equation}\label{eq:multiscale-local-comparison}
H_n(W)+{\mathcal{I}_n}(U;V)\ge\left(\min\left\{1,\alpha-\frac{1}{2}\right\}-\varepsilon\right)n-A_\varepsilon(1+\log C+\log C_{\mathrm{ch}}).
\end{equation}
Here $C_{\mathrm{ch}}\ge2$ bounds the derivatives, the reciprocals of the nondegeneracy quantities, the support and the reciprocal widths of the buffers required by the paired operator, and the constant $A_\varepsilon$ depends only on $\alpha$, $\varepsilon$ and the fixed orders of these bounds.
\end{theorem}

\subsubsection{Proof of Corollary~\ref{cor:sharp-three-images}}

\begin{proof}
Set $W=g(U,V)$. We obtain the lower bound by adding the input entropy to the local entropy inequality. Since $\alpha\ge\frac{3}{2}$, the coefficient in Theorem~\ref{thm:local-multiscale-profile} equals one, so for every $\eta>0$,
$$H_n(W)+{\mathcal{I}_n}(U;V)\ge(1-\eta)n-O_\eta(1+\log C),$$
where the fixed quantities of the chart are absorbed into the constant. The ball bound implies $H_n(U,V)\ge\alpha n-O(1+\log C)$, and adding the two inequalities and using the definition of ${\mathcal{I}_n}(U;V)$, we obtain
$$H_n(U)+H_n(V)+H_n(W)\ge(\alpha+1-\eta)n-O_\eta(1+\log C).$$
At least one of the three entropies is at least one third of the right-hand side, and taking $\eta=3\varepsilon$ proves \eqref{eq:sharp-three-image-lower}. Proposition~\ref{prop:three-image-fixed-sharpness} proves \eqref{eq:sharp-three-image-upper} with a single law for all scales.
\end{proof}

\subsubsection{One step of the entropy iteration}

Write
$$T_n=H_n(U)+H_n(V)+H_n(W)\text{ and }E_n=H_n(U,V).$$
By the ball bound, $E_n\ge\alpha n-O(1+\log C)$. The maps $(u,v)\mapsto(u,g(u,v))$ and $(u,v)\mapsto(v,g(u,v))$ are bi-Lipschitz with constants bounded by a power of $C_{\mathrm{ch}}$, as one sees by integrating the nonvanishing partial derivative and using the bound for the other derivative. Hence, the entropy of each pair among $U$, $V$ and $W$ is at least $\alpha n-O(1+\log C+\log C_{\mathrm{ch}})$. Adding the three bounds for the pairs and using subadditivity, which counts each scalar entropy twice, we get
\begin{equation}\label{eq:multiscale-initial}
T_n\ge\frac{3\alpha n}{2}-O(1+\log C+\log C_{\mathrm{ch}}).
\end{equation}
All these estimates survive the restriction to an event of mass $w$, with $C$ replaced by $\frac{C}{w}$.

For a number $\gamma$, consider the bound
\begin{equation}\label{eq:multiscale-inductive-property}
T_m\ge(\gamma-\eta)m-A_\eta(1+\log C+\log C_{\mathrm{ch}})\text{ for all }m\ge1,
\end{equation}
which we require for every $\eta>0$, with a constant $A_\eta$, uniformly over all laws as in Theorem~\ref{thm:local-multiscale-profile}. By \eqref{eq:multiscale-initial}, this bound holds with $\gamma=\frac{3\alpha}{2}$.

\begin{lemma}\label{lem:multiscale-iteration}
Let $1<\gamma\le3$ and assume \eqref{eq:multiscale-inductive-property} with coefficient $\gamma$. Then, for every $\varepsilon>0$,
\begin{equation}\label{eq:multiscale-one-step-cross}
T_n-E_n\ge(\kappa(\gamma)-\varepsilon)n-B_\varepsilon(1+\log C+\log C_{\mathrm{ch}}),\text{ where }\kappa(\gamma)=\min\left\{1,\frac{\gamma}{2}-\frac{1}{4}\right\}.
\end{equation}
Consequently, \eqref{eq:multiscale-inductive-property} holds with coefficient $\alpha+\kappa(\gamma)$.
\end{lemma}

\begin{proof}
Note that $T_n-E_n=H_n(W)+{\mathcal{I}_n}(U;V)$. The analytic step bounds the probability of a collision between the actual output and the output of the product of the input marginals. We carry out this estimate on classes on which the probabilities of the bins are controlled, so that both Fourier energies can be estimated, and data processing accounts for the replacement of the input law by the product of its marginals. Averaging over the classes must then cost only a logarithmic error.

Choose small $\tau,\eta>0$, to be fixed later in terms of $\varepsilon$. If $n$ is at most a fixed multiple of $1+\log C_{\mathrm{ch}}$, the claim follows from the nonnegativity of $H_n(W)+{\mathcal{I}_n}(U;V)$ after increasing $B_\varepsilon$. For larger $n$, we set $N=n+\lceil c_1(1+\log C_{\mathrm{ch}})\rceil$ and choose a grid of scales $0=j_0<\cdots<j_M=N$ with gaps at most $\tau n+2$, where $c_1$ is a fixed constant large enough to accommodate the margins in frequency used below. Then $N=O(n)$ and $M=O_\tau(1)$.

We first define the classes of masses and the reference laws. At each positive scale $j$ of the grid, an occupied dyadic interval of each $Z\in\{U,V,W\}$ receives the label $K_Z(j)$ determined by
$$2^{-K_Z(j)-1}<\mathbb P(D_jZ=i)\le2^{-K_Z(j)}.$$
We partition the original law according to all these labels and call each event of the partition a mass class. On a nonempty mass class $A$ of probability $w$, we set
$$\begin{aligned}
k_Z(j)&=\min\{K_Z(j),j\},\\
A_Z(j)&=2^{j-k_Z(j)},\\
k_Z(0)&=0.
\end{aligned}$$
There are at most $2^{K_Z(j)+1}$ intervals with the given mass and at most $C_{\mathrm{ch}}^{O(1)}2^j$ intervals in the support, and either count bounds the size of the support of the conditional law on $A$, so that
\begin{equation}\label{eq:multiscale-label-entropy-upper}
H_j(Z\mid A)\le k_Z(j)+O(1+\log C_{\mathrm{ch}}).
\end{equation}
The law on the class has ball constant $\frac{C}{w}$. Applying \eqref{eq:multiscale-inductive-property} to it with tolerance $\eta$, applying the entropy bound coming from the ball bound to its pair $(U,V)$, and combining these lower bounds with \eqref{eq:multiscale-label-entropy-upper}, we obtain
\begin{equation}\label{eq:multiscale-triple-profile}
k_U(j)+k_V(j)+k_W(j)\ge(\gamma-\eta)j-A_\eta\mathcal E_A
\end{equation}
and
\begin{equation}\label{eq:multiscale-pair-profile}
k_U(j)+k_V(j)\ge\alpha j-A_0\mathcal E_A,
\end{equation}
where $\mathcal E_A=1+\log C+\log C_{\mathrm{ch}}+\log\frac{1}{w}$. The bound \eqref{eq:multiscale-pair-profile} for the two input labels will control the mass of pairs at a given separation, while the bound \eqref{eq:multiscale-triple-profile} for all three labels enters later, when the energies of the actual and the reference outputs are combined. Keeping these two uses apart makes the sum over the gaps converge without any rectangle bound for the law on the class.

Recall that $n$ is the observation scale and $N$ the last scale of the grid, and let $n+r$ be the resolution of the inputs, where $r$ is chosen after the geometric margins defining $N$, with
$$n+r\ge N+c_4(1+\log C_{\mathrm{ch}})\text{ and }r=O(1+\log C_{\mathrm{ch}})$$
for a sufficiently large constant $c_4$. Then the filled input cells stay in the buffered chart, and the displacements of inputs and outputs, including the noise introduced below, are less than $2^{-N-10}$. Let $\mu$ be the smoothed law on the class and $\mu^\otimes$ the product of its marginals, both supported in the buffered product chart. Their densities are constant on the chosen input cells, so their relative entropy is exactly the mutual information of the table of cells,
\begin{equation}\label{eq:multiscale-input-relative-entropy}
D(\mu\Vert\mu^\otimes)={\mathcal{I}_{n+r}}(U;V\mid A).
\end{equation}
We send both laws through the channel that evaluates $g$ and adds independent uniform noise on $[-2^{-n-r},2^{-n-r}]$, and we denote the actual and the reference output laws by $\rho$ and $\omega$ and their dyadic vectors at level $n$ by $\rho_n$ and $\omega_n$. The argument with a common channel in Remark~\ref{rem:profile-quantization} yields
\begin{equation}\label{eq:multiscale-output-channel}
H(\rho_n)=H_n(W\mid A)+O(1)\text{ and }D(\rho_n\Vert\omega_n)\le {\mathcal{I}_{n+r}}(U;V\mid A).
\end{equation}
What remains to be shown on this class is a collision bound decaying like $2^{-\kappa(\gamma)n}$, up to small losses in the rate and a factor polynomial in $\frac{CC_{\mathrm{ch}}}{w}$. Indeed, by Jensen's inequality, the cross entropy is at least $-\log\sum_i\rho_n(i)\omega_n(i)$, while by \eqref{eq:multiscale-output-channel}, it is at most $H_n(W\mid A)+{\mathcal{I}_{n+r}}(U;V\mid A)+O(1)$, so that the polynomial factor turns into a logarithmic cost of the restriction when we return to entropy.

The mass labels also provide bounds for densities at the labeled scales. For the two input marginals $\nu_U$ and $\nu_V$ and for $\nu_W=\rho$, an interval $I$ of length at least $2^{-j}$ meets $O(2^j|I|)$ of the original dyadic intervals after enlargement by the displacement, and each of them has conditional mass at most $\frac{2^{-K_Z(j)}}{w}$. Since $k_Z(j)\le K_Z(j)$, this proves
\begin{equation}\label{eq:multiscale-truncated-density}
\nu_Z(I)\le D_AA_Z(j)|I|\text{ for }|I|\ge2^{-j}.
\end{equation}
Here and below, $D_A$ denotes a fixed power of $\frac{CC_{\mathrm{ch}}}{w}$ that may increase from line to line, whose exponent may depend on $\eta$ and on the fixed orders of the geometry, but not on $n$ or on the law.

We next estimate the Fourier energy of the output of the product reference law on each frequency shell, postponing the use of \eqref{eq:multiscale-triple-profile} until this estimate is paired with the energy of the actual output. Consider the shell $\Xi=2^q$ with $1\le q\le n+O(1)$, and choose the first scale $j$ of the grid with $j\ge q+c_2(1+\log C_{\mathrm{ch}})$, where $c_2$ is a sufficiently large fixed geometric constant and $c_1$ was chosen larger. Then $j-q\le\tau n+O(1+\log C_{\mathrm{ch}})$. We apply the paired decomposition from the proof of Proposition~\ref{prop:entropy-radon-shell} to the two smoothed input marginals and the phase $g$, whose partial derivatives and curvature expression do not vanish on the buffer. At a gap $r_0\asymp2^{-k}\gtrsim\Xi^{-1}$, the operator contributes the factor $r_0^{-1}\Lambda^{-\frac{1}{2}}$, where $\Lambda=\Xi r_0$. We denote the paired measures coming from the marginals of $U$ and $V$ by $m_{r_0}$ and $n_{r_0}$, and we need their total masses and their smoothed $L^2$ norms at this gap.

Let $\ell$ be the largest scale of the grid with $\ell\le k-c_3(1+\log C_{\mathrm{ch}})$, where $c_3$ is a sufficiently large fixed geometric constant. When such a positive scale exists, $k-\ell\le\tau n+O(1+\log C_{\mathrm{ch}})$. An interval of length $O_{C_{\mathrm{ch}}}(r_0)$ meets at most $C_{\mathrm{ch}}^{O(1)}$ dyadic intervals at level $\ell$, so the paired measures satisfy $\|m_{r_0}\|\le D_A2^{-k_U(\ell)}$ and $\|n_{r_0}\|\le D_A2^{-k_V(\ell)}$, and \eqref{eq:multiscale-pair-profile} implies
\begin{equation}\label{eq:multiscale-gap-mass-product}
\|m_{r_0}\|\|n_{r_0}\|\le D_A2^{O(\tau n)}r_0^\alpha.
\end{equation}
If there is no positive $\ell$, we take $\ell=0$ and use that the total masses are at most one, and since then $k\le\tau n+O(1+\log C_{\mathrm{ch}})$, the same bound holds with the displayed losses.

At a smoothing radius $h\gtrsim\Lambda^{-1}$, the paired coordinates confine the original variables to intervals of lengths $O(h)$ and $O(r_0h)$, where the second interval may depend on the first variable. The margins in frequency guarantee that these lengths are at least $C_{\mathrm{ch}}^{-O(1)}\Xi^{-1}$, so enlarging them to $2^{-j}$ costs only powers of $C_{\mathrm{ch}}$, and applying \eqref{eq:multiscale-truncated-density} twice, we obtain
$$m_{r_0}(B(z,h))\le D_AA_U(j)^2r_0h^2\text{ and }n_{r_0}(B(z,h))\le D_AA_V(j)^2r_0h^2.$$
Let $\Pi_\Lambda$ be the annular smoothing multiplier from \eqref{eq:entropy-radon-frequency}. Integrating its kernel estimate over dyadic annuli and combining the result with the preceding bounds for balls, we find
$$\|\Pi_\Lambda m_{r_0}\|_2\le D_AA_U(j)r_0^{\frac{1}{2}}\|m_{r_0}\|^{\frac{1}{2}}\text{ and }\|\Pi_\Lambda n_{r_0}\|_2\le D_AA_V(j)r_0^{\frac{1}{2}}\|n_{r_0}\|^{\frac{1}{2}}.$$
Together with the factor coming from the operator and \eqref{eq:multiscale-gap-mass-product}, this bounds the principal contribution of the gap $r_0$ by
$$D_A2^{O(\tau n)}A_U(j)A_V(j)\Xi^{-\frac{1}{2}}r_0^{\frac{\alpha-1}{2}}.$$
The exponent $\frac{\alpha-1}{2}$ is positive, so summing over the dyadic gaps costs only a constant depending on $\alpha$.

The remainders must be controlled with the same density factors, since these factors are what the iteration uses. For $r_0\ge\Xi^{-1}$, the paired masses are at most $D_AA_U(j)r_0$ and $D_AA_V(j)r_0$, so the frequency remainder of order $m\ge3$ sums to
$$D_AA_U(j)A_V(j)\sum_{r_0\ge\Xi^{-1}}r_0(\Xi r_0)^{-m}\lesssim D_AA_U(j)A_V(j)\Xi^{-1}.$$
Outside the cone of gaps treated by the operator, the absolute value of the phase is bounded below by a fixed multiple of the larger gap, and \eqref{eq:entropy-radon-kernel}, with a fixed sufficiently large order of decay, together with the same bounds for the masses of the gaps leads to $D_AA_U(j)A_V(j)\Xi^{-1}$ again. If both gaps are at most $\Xi^{-1}$, their fourfold mass is at most $D_AA_U(j)A_V(j)\Xi^{-2}$, while the kernel is $O(\Xi)$, and the uniformized marginals give no mass to the diagonals. Summing these contributions and observing that the final convolution of the output can only decrease the modulus of its Fourier transform, we arrive at
\begin{equation}\label{eq:multiscale-reference-shell}
\int_{\Xi\le|\xi|\le2\Xi}|\widehat\omega(\xi)|^2d\xi\le D_A2^{O(\tau n)}A_U(j)A_V(j)\Xi^{-\frac{1}{2}}.
\end{equation}
A Gaussian majorant together with \eqref{eq:multiscale-truncated-density}, applied only to intervals of length at least a fixed multiple of $\Xi^{-1}$, also implies
\begin{equation}\label{eq:multiscale-actual-energy}
\int_{|\xi|\le2\Xi}|\widehat\rho(\xi)|^2d\xi\le D_AA_W(j).
\end{equation}

We now combine the two energies and pass to entropy. The energy of the reference output carries the factor $A_U(j)A_V(j)$ and that of the actual output carries $A_W(j)$, so their geometric mean involves all three bin masses, and the bound \eqref{eq:multiscale-triple-profile} for the sum of the three labels controls the mixed collision. The geometric mean is at most
$$D_A2^{O(\tau n)}2^{-\frac{q}{4}+\frac{3j-k_U(j)-k_V(j)-k_W(j)}{2}},$$
and by \eqref{eq:multiscale-triple-profile} and the bound for $j-q$, it is at most
\begin{equation}\label{eq:multiscale-shell-profile-bound}
\left(\frac{CC_{\mathrm{ch}}}{w}\right)^{B_\eta}2^{O(\tau n)+\eta n}2^{(\frac{5}{4}-\frac{\gamma}{2})q}.
\end{equation}
The low frequencies satisfy the corresponding bound with $q=0$. We apply the positive Fourier majorant from the proof of Lemma~\ref{lem:output-profile-transfer} to the collision of $\rho$ and $\omega$ at level $n$ and use the Cauchy--Schwarz inequality on each shell. The exponent of the shell is $\frac{5}{4}-\frac{\gamma}{2}$, and the dyadic sum of these terms is bounded, of size $O(n+1)$, or of size $2^{(\frac{5}{4}-\frac{\gamma}{2})n}$, according as this exponent is negative, zero or positive. Since the collision majorant contributes the factor $2^{-n}$ of the observation scale, the three cases together show that
$$\sum_i\rho_n(i)\omega_n(i)\le\left(\frac{CC_{\mathrm{ch}}}{w}\right)^{B_\eta}(n+2)2^{-\kappa(\gamma)n+O(\tau n)+\eta n},$$
with $\kappa(\gamma)=\min\{1,\frac{\gamma}{2}-\frac{1}{4}\}$. By Jensen's inequality,
$$\sum_i\rho_n(i)\log\frac{1}{\omega_n(i)}\ge-\log\sum_i\rho_n(i)\omega_n(i),$$
where the left-hand side equals $H(\rho_n)+D(\rho_n\Vert\omega_n)$, and \eqref{eq:multiscale-output-channel} turns the collision bound into
$$H_n(W\mid A)+{\mathcal{I}_{n+r}}(U;V\mid A)\ge\kappa(\gamma)n-O(\tau n)-\eta n-O(\log(n+2))-B_\eta\mathcal E_A.$$

It remains to remove the conditioning without a linear loss. Each mass label has mean at most $j+O(1+\log C_{\mathrm{ch}})$, and hence entropy $O(\log(n+2)+\log C_{\mathrm{ch}})$ by comparison with a geometric distribution, so the full label has entropy $O_\tau(\log(n+2)+\log C_{\mathrm{ch}})$. By \eqref{eq:mixture-concavity}, the average of the output entropies is at most $H_n(W)$, by \eqref{eq:label-MI}, the average of the input mutual informations is at most ${\mathcal{I}_{n+r}}(U;V)$ plus the entropy of the label, and the average of $\log\frac{1}{w}$ is exactly this entropy. Finally, \eqref{eq:entropy-refinement} bounds the cost of refining the inputs by $2r=O(1+\log C_{\mathrm{ch}})$. Averaging the estimates for the classes, we therefore obtain
$$T_n-E_n\ge\kappa(\gamma)n-O(\tau n)-\eta n-O_{\tau,\eta}(\log(n+2)+1+\log C+\log C_{\mathrm{ch}}).$$
We first choose $\tau$ and $\eta$ so that their linear losses use at most half of the prescribed tolerance. With these choices fixed, the remaining multiple of $\log(n+2)$ is at most the other half of the tolerance times $n$, plus a fixed constant, which proves \eqref{eq:multiscale-one-step-cross} with the stated logarithmic dependence on $C$ and $C_{\mathrm{ch}}$. Adding $E_n\ge\alpha n-O(1+\log C)$, we arrive at \eqref{eq:multiscale-inductive-property} with coefficient $\alpha+\kappa(\gamma)$.
\end{proof}

\subsubsection{Proof of Theorem~\ref{thm:local-multiscale-profile}}

\begin{proof}
We apply Lemma~\ref{lem:multiscale-iteration} repeatedly to improve the coefficient for the sum of the three entropies, and then we use its conclusion \eqref{eq:multiscale-one-step-cross} one last time. By \eqref{eq:multiscale-initial}, we can start with $\gamma_0=\frac{3\alpha}{2}$, and we define
$$\gamma_{r+1}=\alpha+\min\left\{1,\frac{\gamma_r}{2}-\frac{1}{4}\right\}.$$
The initial value lies in $(1,3]$, and every later value lies between $\alpha$ and $\alpha+1\le3$, so each step satisfies the hypothesis of the lemma. The first update is at least $\gamma_0$, and since the update map is nondecreasing, the whole sequence is nondecreasing. As long as $\frac{\gamma_r}{2}-\frac{1}{4}<1$, the recurrence is affine with slope $\frac{1}{2}$ and fixed point $2\alpha-\frac{1}{2}$, and otherwise the next value is $\alpha+1$. These two alternatives show that the limit of the sequence is $\alpha+\min\{1,\alpha-\frac{1}{2}\}$, and consequently $\kappa(\gamma_r)$ tends to $\min\{1,\alpha-\frac{1}{2}\}$, also at the boundary value $\alpha=\frac{3}{2}$ between the two regimes.

For the prescribed $\varepsilon$, we choose finitely many steps so that the limiting coefficient is approached within $\frac{\varepsilon}{2}$, and we choose the tolerances of the steps so small that their accumulated loss in the rate is at most $\frac{\varepsilon}{2}$. Only these finitely many applications of the lemma are used, also when $\alpha=\frac{3}{2}$, so we never need an estimate at a limiting coefficient that is not attained. The constants may depend on the number of steps, but their dependence on $C$ and $C_{\mathrm{ch}}$ remains linear in $1+\log C+\log C_{\mathrm{ch}}$, and the last application of the lemma yields the asserted bound for $T_n-E_n=H_n(W)+{\mathcal{I}_n}(U;V)$.
\end{proof}

\subsubsection{Global entropy inequalities in the three coordinate systems}

\begin{corollary}\label{cor:global-multiscale-profile}
Let $\frac{1}{2}<s<1$, let $(X,Y)$ be jointly $(s,s,C)$-Frostman in a fixed bounded box, and set $\kappa_s=\min\{1,2s-\frac{1}{2}\}$. For each nonspecial representation $F=G(U,V)$ with $(U,V)$ equal to $(X,Y)$, $(X,S)$ or $(Y,S)$, and for every $\varepsilon>0$,
\begin{equation}\label{eq:global-multiscale-profile}
H_n(F)+{\mathcal{I}_n}(U;V)\ge(\kappa_s-\varepsilon)n-O_\varepsilon(1).
\end{equation}
\end{corollary}

\begin{proof}
Fix a nonspecial representation $F=G(U,V)$. We apply the local theorem on finitely many regular charts, making the discarded mass small enough to fit within the tolerance. The joint Frostman bound implies an isotropic $2s$-Frostman bound, which the fixed sum shears preserve, and since $2s>1$, the covering argument behind \eqref{eq:joint-tube} makes the mass near any fixed algebraic curve uniformly small. We apply it to the zero sets of the two partial derivatives of $G$ and of the numerator of its curvature expression, which are nonzero polynomials because $G$ is nonspecial.

For a fixed $\rho>0$, we remove a neighborhood of these zero sets of mass at most $\rho$. The remaining compact set has a finite cover by buffered regular product charts, and we partition it into disjoint events subordinate to the cover, the removed part receiving one more value of the label. On an event of probability $w$, we apply Theorem~\ref{thm:local-multiscale-profile} with $\alpha=2s$ and ball constant $\frac{C}{w}$, up to a fixed geometric factor, and its coefficient is $\min\{1,2s-\frac{1}{2}\}=\kappa_s$.

By the concavity of entropy and \eqref{eq:label-MI}, the average of the local left-hand sides is at most $H_n(F)+{\mathcal{I}_n}(U;V)$ plus the finite entropy of the label, and the errors coming from the normalizations also average to a bounded multiple of this entropy. On the removed part, the local left-hand side is nonnegative, so if the local tolerance is $\eta$, the retained events, of total mass at least $1-\rho$, contribute a coefficient of at least $\kappa_s-\kappa_s\rho-\eta$. We first choose $\rho$ with $\kappa_s\rho<\frac{\varepsilon}{2}$ and then $\eta<\frac{\varepsilon}{2}$. The tube estimate is uniform over the laws, so the neighborhood, the finite cover and the constants of the charts can all be fixed independently of the law and of $n$, and averaging proves \eqref{eq:global-multiscale-profile}.
\end{proof}

\subsubsection{Proof of Corollary~\ref{cor:independent-single-image}}

\begin{proof}
Set $\alpha=s_1+s_2$ and $\kappa=\min\{1,\alpha-\frac{1}{2}\}$. By independence and the marginal Frostman bounds,
$$\mathbb P((X,Y)\in B(z,r))\le2^\alpha C_1C_2r^\alpha\text{ for }z\in\mathbb R^2\text{ and }0<r\le1,$$
and ${\mathcal{I}_n}(X;Y)=0$. We apply Theorem~\ref{thm:local-multiscale-profile} after the localization used in the proof of Corollary~\ref{cor:global-multiscale-profile}, now with $\alpha=s_1+s_2$. We may assume that $0<\varepsilon<\kappa$, since otherwise the claimed lower bound follows from the nonnegativity of entropy.

Let $\mathcal Z$ be the union of the zero sets of $f_x$, of $f_y$ and of the numerator of $\partial_x\partial_y\ln\left|\frac{f_x}{f_y}\right|$, which are nonzero polynomials because $f$ is not analytic-special. The covering argument preceding \eqref{eq:joint-tube} covers the $r$-neighborhood of $\mathcal Z$ in the support box by $O_{f,R}(r^{-1})$ balls of radius $O_{f,R}(r)$, so the ball bound implies
$$\mathbb P(\operatorname{dist}((X,Y),\mathcal Z)\le r)\lesssim_{f,R,s_1,s_2,C_1,C_2}r^{\alpha-1}.$$
Since $\alpha>1$, we can choose $0<\rho<\frac{\varepsilon}{2\kappa}$ and remove a fixed neighborhood of $\mathcal Z$ of mass at most $\rho$. The remaining compact set has a finite cover by cores of buffered regular product charts, and we partition it into disjoint events subordinate to this cover, labeled by $L\in\{0,\ldots,m\}$, where $L=0$ denotes the removed part. The charts and $m$ depend only on the fixed data and on $\varepsilon$.

For each $j\ge1$ with $w_j=\mathbb P(L=j)>0$, the conditional law has ball constant $\frac{2^\alpha C_1C_2}{w_j}$, and Theorem~\ref{thm:local-multiscale-profile} with local tolerance $\frac{\varepsilon}{2}$ shows that
$$H_n(f(X,Y)\mid L=j)+{\mathcal{I}_n}(X;Y\mid L=j)\ge\left(\kappa-\frac{\varepsilon}{2}\right)n-O_\varepsilon\left(1+\log\frac{1}{w_j}\right),$$
while the corresponding left-hand side on the removed part is nonnegative. The conditional inputs need not be independent, but \eqref{eq:label-MI} and \eqref{eq:mixture-concavity} imply
$${\mathcal{I}_n}(X;Y\mid L)\le {\mathcal{I}_n}(X;Y)+H(L)=H(L)\text{ and }H_n(f(X,Y)\mid L)\le H_n(f(X,Y)).$$
Since $\sum_{j=1}^mw_j\log\frac{1}{w_j}\le H(L)\le\log(m+1)$, where terms of probability zero are omitted, averaging the local inequalities yields
$$H_n(f(X,Y))\ge(1-\rho)\left(\kappa-\frac{\varepsilon}{2}\right)n-O_\varepsilon(1)\ge(\kappa-\varepsilon)n-O_\varepsilon(1),$$
and all choices were uniform in the input laws and in $n$.

If $s_1+s_2\ge\frac{3}{2}$, then $\kappa=1$, while the bounded support of $f(X,Y)$ forces $H_n(f(X,Y))\le n+O_{f,R}(1)$. Dividing these bounds by $n$, letting $n\to\infty$ and then letting $\varepsilon\downarrow0$, we obtain the asserted limit.
\end{proof}

\subsection{Conditional entropy bounds}
\label{sec:conditional-profiles-new}

We now combine the inequalities for conditional copies with the multiscale entropy inequality. Throughout this subsection, the original pair satisfies the conditional hypotheses of Theorem~\ref{thm:conditional-main}, $C\ge1$ bounds both exact conditional Frostman constants, and as before, $\beta_s=\min\{\frac{3s}{2},1\}$. For a source representation $(U,V)$ on a fixed regular chart, we write
$$\begin{aligned}
T_j&=H_j(U)+H_j(V)+H_j(F), & E_j&=H_j(U,V),\\
C_j&=T_j-E_j=H_j(F)+{\mathcal{I}_j}(U;V){.}
\end{aligned}$$
{Local statements concern an event $A$ of original probability $a>0$ in the chart. Entropies and the quantities $T_j,E_j,C_j$ are computed under the normalized law on $A$, with this conditioning suppressed from their notation; the conditioning event is displayed explicitly in the output mutual information.} We set $L_a=1+\log\frac{C}{a}$ and denote the fixed constant of the chart by $C_{\mathrm{ch}}$. A bound is called hereditary if it holds on every such event with an error $O_{\eta,C_{\mathrm{ch}}}(L_a)$ that is uniform in $A$. All our charts also make $(S,F)$ a source pair related to $(X,Y)$ by a bi-Lipschitz map, so that ${\mathcal{I}_j(S;F\mid A)}=H_j(S)+H_j(F)-E_j+O_{C_{\mathrm{ch}}}(1)$. Inequalities with positive coefficients for the outputs will be averaged over the charts by means of Lemma~\ref{lem:conditional-chart-averaging}.

\subsubsection{Retaining output mutual information}
\label{sec:retained-information}

\begin{lemma}\label{lem:retained-information}
Fix $\frac{1}{2}<s<1$ and a common buffered chart as above, with source representation $(U,V)$ and $F=g(U,V)$ satisfying the hypotheses of Theorem~\ref{thm:local-multiscale-profile}, so that the source exponent is $\alpha=2s>1$. {For an event $A$ of original probability $a>0$ contained in the chart, write $T_j^A$ and $C_j^A$ for the corresponding quantities under the normalized restriction to $A$.} Fix $1<\gamma\le3$ and $b\ge0$, and suppose that for every $\eta>0$, there is a constant $K_\eta$ such that
$$T_j^A\ge(\gamma-\eta)j-b{\mathcal{I}_j(S;F\mid A)}-K_\eta L_a\text{ for all }j\ge1$$
and for every such event $A$, where the constants are uniform over the laws, events and scales under consideration and may depend on $s$, $\gamma$, $b$, $C_{\mathrm{ch}}$ and the other fixed data. Then, for every $\varepsilon>0$, there is a constant $K_\varepsilon'$, which may also depend on the constants $K_\eta$ but not on the law, the event or the scale, such that
\begin{equation}\label{eq:retained-information-iteration}
C_n^A+\frac{b}{2}{\mathcal{I}_n(S;F\mid A)}\ge\left(\min\left\{1,\frac{\gamma}{2}-\frac{1}{4}\right\}-\varepsilon\right)n-K_\varepsilon'L_a\text{ for all }n\ge1
\end{equation}
on every event $A$ in the chart.
\end{lemma}

\begin{proof}
We adapt the iteration of Lemma~\ref{lem:multiscale-iteration}, keeping the term $b{\mathcal{I}_j(S;F\mid A)}$ instead of discarding it. On each class of masses, data processing bounds this information by its value at the observation scale, and taking the geometric mean of the two Fourier energies then produces the coefficient $\frac{b}{2}$ in the conclusion.

Fix an event $A$ in the chart and small $\tau,\eta>0$, to be chosen below. We repeat, for the normalized law on $A$, the construction in the proof of Lemma~\ref{lem:multiscale-iteration}, with the finite grid of frequencies and gaps of relative spacing $\tau$, the mass labels and the common channel. Let $N=n+O_{C_{\mathrm{ch}}}(1)$ be the last scale of the grid, and let $B$ be a class of masses of probability $w$ under this law, so that $B\subset A$ has original probability $aw$. {Throughout, $\mathcal{I}_j(S;F\mid B)$ is computed under the normalized restriction of the original law to $B$, and the smoothing is used only in the collision estimate.} The isotropic ball bound on $B$ has constant at most $\frac{C_{\mathrm{ch}}^{O(1)}C^2}{aw}$, whose logarithm is $O_{C_{\mathrm{ch}}}(L_{aw})$. This only uses the restriction of the joint square bound, and no exact conditional Frostman bound for the law on $B$ is needed.

The bounds \eqref{eq:multiscale-label-entropy-upper} for the supports of the labels, together with the hereditary hypothesis applied on $B$, imply that at every positive scale of the grid
$$k_U(j)+k_V(j)+k_F(j)\ge(\gamma-\eta)j-b{\mathcal{I}_j(S;F\mid B)}-O_{\eta,s,\gamma,b,C_{\mathrm{ch}}}(L_{aw}).$$
Since $b\ge0$, any upper bound for ${\mathcal{I}_j(S;F\mid B)}$ yields a valid lower bound for the sum of the three labels. For $j\le n$, data processing shows that ${\mathcal{I}_j(S;F\mid B)}\le {\mathcal{I}_n(S;F\mid B)}$, while for $n<j\le N$, refinement implies ${\mathcal{I}_j(S;F\mid B)}\le {\mathcal{I}_n(S;F\mid B)}+2(j-n)={\mathcal{I}_n(S;F\mid B)}+O_{C_{\mathrm{ch}}}(1)$, so the replacement of ${\mathcal{I}_j(S;F\mid B)}$ by ${\mathcal{I}_n(S;F\mid B)}$ costs only $O_{C_{\mathrm{ch}}}(b)$, which is absorbed into the error. Using this bound in the geometric mean of \eqref{eq:multiscale-reference-shell} and \eqref{eq:multiscale-actual-energy}, we find that on a shell $\Xi=2^h$, the geometric mean is at most
$$2^{O_{\eta,s,\gamma,b,C_{\mathrm{ch}}}(L_{aw})}2^{O_{s,\gamma,b}(\tau n+\eta n)}2^{(\frac{5}{4}-\frac{\gamma}{2})h+\frac{b}{2}{\mathcal{I}_n(S;F\mid B)}}.$$
The sum over the gaps converges because its exponent is $\frac{2s-1}{2}>0$, and since $b{\mathcal{I}_n(S;F\mid B)}\ge0$, the low frequencies satisfy the same bound with $h=0$ after the constant is increased.

Let $\kappa=\min\{1,\frac{\gamma}{2}-\frac{1}{4}\}$. Summing over the frequency shells and using the positive collision majorant, we get
$$\sum_i\rho_n(i)\omega_n(i)\le2^{O_{\eta,s,\gamma,b,C_{\mathrm{ch}}}(L_{aw})}(n+2)2^{-\kappa n+O_{s,\gamma,b}(\tau n+\eta n)+\frac{b}{2}{\mathcal{I}_n(S;F\mid B)}},$$
where $\rho_n$ and $\omega_n$ are the actual and the reference output vectors for the class, obtained from the same channel. Jensen's inequality, data processing as in \eqref{eq:multiscale-output-channel}, and the cost $O_{C_{\mathrm{ch}}}(1)$ of refining the inputs therefore lead to
$$C_n^B+\frac{b}{2}{\mathcal{I}_n(S;F\mid B)}\ge\kappa n-O_{s,\gamma,b}(\tau n+\eta n)-O(\log(n+2))-O_{\eta,s,\gamma,b,C_{\mathrm{ch}}}(L_{aw}).$$

Let $L$ denote the full mass label under the law on $A$. By the concavity of entropy and \eqref{eq:label-MI},
$$\sum_Bw\left(C_n^B+\frac{b}{2}{\mathcal{I}_n(S;F\mid B)}\right)\le C_n^A+\frac{b}{2}{\mathcal{I}_n(S;F\mid A)}+\left(1+\frac{b}{2}\right)H(L),$$
because the average entropy of $F$ over the classes is at most its entropy on $A$, while averaging the two mutual informations costs $H(L)$ and $\frac{b}{2}H(L)$, respectively. Moreover, $\sum_BwL_{aw}=L_a+H(L)$, and since the grid is finite, $H(L)=O_{\tau,C_{\mathrm{ch}}}(\log(n+2))$. These two averaging bounds take care separately of the loss of information and of the logarithmic cost of the normalization. We choose $\tau$ and $\eta$ so that the linear losses use at most half of $\varepsilon n$, and after these choices, the remaining logarithmic terms are bounded by the other half plus a fixed constant times $L_a$, since $L_a\ge1$. In particular, these choices do not depend on the mass $a$ of the original event. The finitely many scales below the geometric cutoff are covered by increasing this constant, because $C_n^A,{\mathcal{I}_n(S;F\mid A)}\ge0$ and $L_a\ge1$, and this proves \eqref{eq:retained-information-iteration} uniformly on every event $A$.
\end{proof}

\subsubsection{Two-copy entropy inequality in the original coordinates}\label{sec:two-copy-original}

\begin{proposition}\label{prop:two-copy-original}
If $f$ is nonspecial in its original coordinates and $\frac{1}{2}<s<\frac{3}{4}$, then
$$H_n(S)+2H_n(F)\ge\left(2s+\min\left\{1,2s+\frac{\beta_s}{3}-\frac{2}{3}\right\}-\varepsilon\right)n-O_\varepsilon(1).$$
\end{proposition}

\begin{proof}
We first prove an entropy inequality on every event in a common regular chart. The final inequality involves only a positive combination of the two output entropies, so Lemma~\ref{lem:conditional-chart-averaging} then makes it global. Both $\mathcal W_x$ and $\mathcal W_y$ are nonzero, because $\partial_x\ln\left|\frac{f_x}{f_y}\right|=\frac{\mathcal W_x}{f_xf_y}$ and $\partial_y\ln\left|\frac{f_x}{f_y}\right|=-\frac{\mathcal W_y}{f_xf_y}$, and the vanishing of either would make the curvature expression of $f$ vanish, contrary to nonspecialness. We use the original representation $(U,V)=(X,Y)$ on a buffered chart that is regular for both inequalities of Lemma~\ref{lem:retained-copy-hereditary} and for the recovery of the source from $(S,F)$, {together with the notation $T_j,E_j,C_j,L_a$ introduced above, for an event $A$ of original probability $a>0$ in the chart.}

We increase a coefficient $p$ while allowing a controlled coefficient $0\le q\le1$ for the output mutual information, in the hereditary bound
\begin{equation}\label{eq:two-copy-hereditary}
C_j+q{\mathcal{I}_j(S;F\mid A)}\ge(p-\eta)j-O_{\eta,C_{\mathrm{ch}}}(L_a).
\end{equation}
It is enough to reach values of $p$ arbitrarily close to $\min\{1,2s+\frac{\beta_s}{3}-\frac{2}{3}\}$. Indeed, on the same event,
$$\begin{aligned}&H_n(S)+2H_n(F)-(C_n+q{\mathcal{I}_n(S;F\mid A)})\\&\quad=2E_n-H_n(X)-H_n(Y)+(1-q){\mathcal{I}_n(S;F\mid A)}+O_{C_{\mathrm{ch}}}(1)\\&\quad\ge2sn-O_{C_{\mathrm{ch}}}(L_a),\end{aligned}$$
where the last information term was discarded because $1-q\ge0$ and ${\mathcal{I}_n(S;F\mid A)}\ge0$, and the remaining lower bound is the sum of the two conditional entropy bounds on the event. Thus, raising $p$ in \eqref{eq:two-copy-hereditary} raises the positive coefficient of the outputs to $2s+p$, which is the quantity we need for averaging over the charts.

We now establish and iterate \eqref{eq:two-copy-hereditary}. Theorem~\ref{thm:local-multiscale-profile}, applied on the original chart with $\alpha=2s<\frac{3}{2}$ and the ball constant of the event, provides the starting point $(p_0,q_0)=(2s-\frac{1}{2},0)$. To apply Lemma~\ref{lem:retained-information}, we need a lower bound for $T_j$ in which only ${\mathcal{I}_j(S;F\mid A)}$ remains on the right, and we obtain it by combining two inequalities with opposite coefficients of $H_j(F)$. First, we average the two orientations of Lemma~\ref{lem:retained-copy-hereditary} on the same event. The averaged left-hand side is
$$H_j(F)+2H_j(S)+\frac{3}{2}(H_j(X)+H_j(Y))-2E_j=\frac{3}{2}T_j-\frac{5}{2}H_j(F)+2{\mathcal{I}_j(S;F\mid A)}+O_{C_{\mathrm{ch}}}(1),$$
and after rearranging, we arrive at
\begin{equation}\label{eq:two-copy-lower-profile}
T_j\ge\frac{2}{3}(s+\beta_s-\eta)j+\frac{5}{3}H_j(F)-\frac{4}{3}{\mathcal{I}_j(S;F\mid A)}-O_{\eta,C_{\mathrm{ch}}}(L_a).
\end{equation}
A second lower bound follows from the current hereditary inequality. By the conditional entropy bounds on the event, $H_j(X)+H_j(Y)\le2E_j-2sj+O(L_a)$, and substituting $E_j=T_j-C_j$ and then \eqref{eq:two-copy-hereditary}, we obtain
\begin{equation}\label{eq:two-copy-upper-profile}
T_j\ge2(p+s-\eta)j-H_j(F)-2q{\mathcal{I}_j(S;F\mid A)}-O_{\eta,C_{\mathrm{ch}}}(L_a).
\end{equation}
We give \eqref{eq:two-copy-upper-profile} the weight $\frac{5}{8}$ and \eqref{eq:two-copy-lower-profile} the weight $\frac{3}{8}$, so that the coefficient of $H_j(F)$ becomes $-\frac{5}{8}+\frac{3}{8}\cdot\frac{5}{3}=0$, and we arrive at
$$T_j\ge\left(\frac{5}{4}p+\frac{3}{2}s+\frac{1}{4}\beta_s-O(\eta)\right)j-\left(\frac{5}{4}q+\frac{1}{2}\right){\mathcal{I}_j(S;F\mid A)}-O_{\eta,C_{\mathrm{ch}}}(L_a).$$

We apply Lemma~\ref{lem:retained-information} with
$$\gamma=\frac{5}{4}p+\frac{3}{2}s+\frac{1}{4}\beta_s\text{ and }b=\frac{5}{4}q+\frac{1}{2}.$$
For the iterates considered below, we have $p_0=2s-\frac{1}{2}\le p\le1$ and $0\le q\le\frac{2}{3}$, so $1<\gamma\le\frac{21}{8}<3$ and $\frac{1}{2}\le b\le\frac{4}{3}$, and the lemma applies. It replaces $(p,q)$ in the hereditary inequality by
\begin{equation}\label{eq:two-copy-recurrence}
p'=\min\left\{1,\frac{5}{8}p+\frac{3}{4}s+\frac{1}{8}\beta_s-\frac{1}{4}\right\}\text{ and }q'=\frac{5}{8}q+\frac{1}{4}.
\end{equation}
The update of $q$ maps $[0,\frac{2}{3}]$ into itself, and the update of $p$ is nondecreasing, capped at one, and sends $p_0$ to a value at least $p_0$, as the computation of the fixed point below shows, so the bounds on the parameters hold inductively. The tolerance in the convex combination can be chosen arbitrarily small, so each application preserves the hereditary error $O_{\eta,C_{\mathrm{ch}}}(L_a)$, with constants depending also on $s$ and on the fixed parameters of the iteration.

Solving the recurrence, we find
$$q_r=\frac{2}{3}\left(1-\left(\frac{5}{8}\right)^r\right)\text{ and }p_r\uparrow p_*=\min\left\{1,2s+\frac{\beta_s}{3}-\frac{2}{3}\right\},$$
and in particular $0\le q_r<1$ at every step. Indeed, the affine part of the recurrence for $p$ contracts by the factor $\frac{5}{8}$, and its fixed point exceeds $p_0$ by $\frac{\beta_s}{3}-\frac{1}{6}>0$. Given the target tolerance, we choose finitely many iterations that bring $p$ within this tolerance of $p_*$, and then we choose the finitely many intermediate tolerances sufficiently small. After these choices, all constants are fixed and uniform in the law, the event and the observation scale. The initial reduction now produces an inequality with positive coefficients for the outputs and with coefficient arbitrarily close to $2s+p_*$, and applying Lemma~\ref{lem:conditional-chart-averaging} to it, with the local tolerances allotted within $\varepsilon$, proves the proposition.
\end{proof}

\subsubsection{Conditional entropy bounds using bin masses at the final scale}\label{sec:terminal-total-refinement}

The next refinement uses the masses of the bins at the final scale in addition to the entropy bounds at the intermediate scales. Two different ways of eliminating the same information term then lead to the two coefficients in the middle range.

\begin{proposition}\label{prop:conditional-terminal-total-new}
If $f$ is nonspecial in the original coordinates, then for $\frac{5}{9}\le s<\frac{2}{3}$ and every $\varepsilon>0$,
\begin{equation}\label{eq:terminal-total-output}
M_n\ge\left(\frac{(5s-1)(7-6s)}{19-15s}-\varepsilon\right)n-O_\varepsilon(1)
\end{equation}
and
\begin{equation}\label{eq:terminal-inverse-output}
M_n\ge\left(\frac{49s-30s^2-7}{26-15s}-\varepsilon\right)n-O_\varepsilon(1).
\end{equation}
At $s=\frac{2}{3}$, the first coefficient follows from Proposition~\ref{prop:two-copy-original}.
\end{proposition}

\begin{proof}
{We reduce the proposition to a weighted entropy inequality for $C_n$, $T_n$ and $\mathcal{I}_n(S;F\mid A)$ on every event $A$ in the chart.} The conditional entropy bounds eliminate $T_n$, and two different treatments of the remaining information term yield the two inequalities with positive coefficients for the outputs. We first carry out this reduction, and then we prove the weighted inequality by keeping the mass labels at the final scale in the Fourier argument. Since the coefficients of the outputs are positive, the inequalities can be averaged over the charts before we pass to the maximum.

We use the original representation $(U,V)=(X,Y)$ and choose the common chart so that it is also regular for \eqref{eq:inverse-output-retained}, which is possible because the additional polynomial $N_f$ is nonzero. For $s<\frac{2}{3}$, the recurrence \eqref{eq:two-copy-recurrence} converges to $p=\frac{5s}{2}-\frac{2}{3}<1$ and $q=\frac{2}{3}$. For a given $\eta>0$, we choose a finite iterate whose coefficient $p$ is sufficiently close to this limit. Its coefficient $q$ is at most $\frac{2}{3}$, and since ${\mathcal{I}_j(S;F\mid A)}\ge0$, we may weaken the bound by replacing $q$ with $\frac{2}{3}$. Substituting this bound into the same convex combination as in the proof of Proposition~\ref{prop:two-copy-original}, we obtain the hereditary bound
\begin{equation}\label{eq:terminal-total-coarse}
T_j\ge\left(5s-\frac{5}{6}-\eta\right)j-\frac{4}{3}{\mathcal{I}_j(S;F\mid A)}-O_{\eta,C_{\mathrm{ch}}}(L_a),
\end{equation}
with constants depending on $\eta$ but uniform over events and scales. In particular, this bound uses no estimate at a limiting coefficient that is not attained.

Set
$$c_s=\frac{10-15s}{23-30s}\text{ and }q_s^{\mathrm{tot}}=\frac{2}{3}(1-2c_s).$$
We claim that on every event $A$ of original probability $a>0$ in the common chart, and for every $\varepsilon>0$,
\begin{equation}\label{eq:terminal-total-weighted}
C_n-c_sT_n+q_s^{\mathrm{tot}}{\mathcal{I}_n(S;F\mid A)}\ge(1-3c_s-\varepsilon)n-O_{\varepsilon,C_{\mathrm{ch}}}(L_a),
\end{equation}
where all quantities refer to the normalized law on $A$.

Let us assume this weighted inequality for the moment and derive both output bounds from it. The conditional entropy bounds on the event show that $T_n\le2E_n-2sn+H_n(F)+O(L_a)$, and since $C_n-c_sT_n=(1-c_s)T_n-E_n$ with $1-c_s>0$, this bounds the left-hand side of \eqref{eq:terminal-total-weighted} from above. Substituting also $E_n=H_n(S)+H_n(F)-{\mathcal{I}_n(S;F\mid A)}+O_{C_{\mathrm{ch}}}(1)$, we obtain
\begin{equation}\label{eq:terminal-output-retained}
\begin{aligned}
&(1-2c_s)H_n(S)+(2-3c_s)H_n(F)-\frac{1-2c_s}{3}{\mathcal{I}_n(S;F\mid A)}\\
&\ge(1-3c_s+2(1-c_s)s-\varepsilon)n-O_{\varepsilon,C_{\mathrm{ch}}}(L_a).
\end{aligned}
\end{equation}
For the first conclusion, we simply discard the nonpositive information term. For the second, we add $\frac{2(1-2c_s)}{3}$ times \eqref{eq:inverse-output-retained}, on the same event, and multiply by three, so that the information terms cancel and
$$5(1-2c_s)H_n(S)+(11-19c_s)H_n(F)\ge(3-9c_s+(14-22c_s)s-\varepsilon)n-O_{\varepsilon,C_{\mathrm{ch}}}(L_a).$$
For $\frac{5}{9}\le s<\frac{2}{3}$, we have $0<c_s\le\frac{5}{19}$, so all coefficients of the outputs in both inequalities are positive. We apply Lemma~\ref{lem:conditional-chart-averaging} and only then bound the left-hand sides by $(3-5c_s)M_n$ and $(16-29c_s)M_n$, respectively. Both divisors are positive, and a direct simplification shows that
$$\frac{1-3c_s+2(1-c_s)s}{3-5c_s}=\frac{(5s-1)(7-6s)}{19-15s}\text{ and }\frac{3-9c_s+(14-22c_s)s}{16-29c_s}=\frac{49s-30s^2-7}{26-15s}.$$
Once the weighted inequality is established, renaming the tolerance proves both bounds on the half-open interval $[\frac{5}{9},\frac{2}{3})$. At $s=\frac{2}{3}$, Proposition~\ref{prop:two-copy-original} yields $H_n(S)+2H_n(F)\ge(\frac{7}{3}-\varepsilon)n-O_\varepsilon(1)$, and dividing by three, we obtain the first coefficient $\frac{7}{9}$ at this endpoint. It remains to prove \eqref{eq:terminal-total-weighted} for $\frac{5}{9}\le s<\frac{2}{3}$.

Fix an event $A$ of original mass $a$ in the chart. We use the finite grid of frequencies and gaps from the proof of Lemma~\ref{lem:multiscale-iteration}, with relative spacing $\tau$ and final scale $N=n+O_{C_{\mathrm{ch}}}(1)\ge n$, which we include among the scales of the grid. Before restricting to a class of masses $B$ of conditional probability $w$ under $A$, we define
$$\begin{aligned}
K_Z(N)&=\lfloor-\log\mathbb P(D_NZ\mid A)\rfloor,\\
k_Z(N)&=\min\{K_Z(N),N\},\\
b_B&=k_X(N)+k_Y(N)+k_F(N),
\end{aligned}$$
where the probability is the mass of the realized dyadic interval of $D_NZ$ and $Z$ ranges over $X$, $Y$ and $F$. Thus, $b_B$ is constant on $B$ and $0\le b_B\le3N$. {We write $C_n^B$ for the quantity $C_n$ computed under the normalized law on $B$, an event of original probability $aw$.} We shall see at the averaging stage that the mean of $b_B$ is at least $T_n-O_{C_{\mathrm{ch}}}(1)$, so keeping this label produces the term $-c_sT_n$ required in \eqref{eq:terminal-total-weighted}.

Write $\Xi=2^h$ for the frequency shell. We need two bounds for the mixed Fourier energy on this shell. The entropy bound at the intermediate scales contributes a factor that increases with $\Xi$, whereas the labels at the final scale contribute a decreasing factor, and taking the minimum retains the information from both. For the first bound, we use \eqref{eq:terminal-total-coarse} on $B$ in \eqref{eq:multiscale-shell-profile-bound}. For $j\le n$, we have ${\mathcal{I}_j(S;F\mid B)}\le {\mathcal{I}_n(S;F\mid B)}$ by data processing, and the additional scales up to $N$ contribute at most $2(N-n)=O_{C_{\mathrm{ch}}}(1)$ by refinement. The resulting bound is
$$D_B2^{O(\tau n+\eta n)}\Xi^{a_s}2^{\frac{2{\mathcal{I}_n(S;F\mid B)}}{3}},\text{ where }a_s=\frac{5}{3}-\frac{5s}{2}>0\text{ and }\log D_B=O_{\eta,C_{\mathrm{ch}}}(L_{aw}).$$
For the second bound, the labels at the final scale provide density bounds $D_B2^{N-k_Z(N)}$ for intervals of length at least $2^{-N}$. Applying \eqref{eq:multiscale-reference-shell} and \eqref{eq:multiscale-actual-energy} with these density bounds, we obtain the energy bounds
$$D_B2^{2N-k_X(N)-k_Y(N)}\Xi^{-\frac{1}{2}}\text{ and }D_B2^{N-k_F(N)},$$
whose geometric mean is $D_B2^{O(\tau n)}2^{\frac{3N-b_B}{2}}\Xi^{-\frac{1}{4}}$. The sum over the gaps still converges, because the law on the event satisfies an isotropic $2s$-Frostman bound and $2s>1$. Note that these density estimates use only the masses of the bins and require no conditional Frostman bound for the restricted law.

On each shell, we take the smaller of the two bounds. The increasing and the decreasing factors cross at
$$h_*=\frac{\frac{3N-b_B}{2}-\frac{2{\mathcal{I}_n(S;F\mid B)}}{3}}{a_s+\frac{1}{4}},$$
and splitting the dyadic sum at this value, or at zero if it is negative, and summing the two geometric series, we find
$$\sum_{h\ge0}\min\{2^{a_sh+\frac{2{\mathcal{I}_n(S;F\mid B)}}{3}},2^{-\frac{h}{4}+\frac{3N-b_B}{2}}\}\lesssim_s2^{\frac{2}{3}(1-2c_s){\mathcal{I}_n(S;F\mid B)}+c_s(3N-b_B)},$$
with $c_s=\frac{a_s}{2(a_s+\frac{1}{4})}$, which is the number defined above. If $h_*<0$, the decreasing bound is already the smaller one at $h=0$ and remains smaller for all $h\ge0$, and its value at zero is at most the displayed weighted geometric mean because $0<2c_s<1$, so the same bound covers this case. The positive collision majorant and Jensen's inequality from the proof of Lemma~\ref{lem:multiscale-iteration} therefore apply. Choosing the spacing of the grid and the tolerance in the intermediate entropy bound small enough that their combined linear loss is at most $\eta n$, we obtain the estimate
\begin{equation}\label{eq:terminal-total-class}
C_n^B+q_s^{\mathrm{tot}}{\mathcal{I}_n(S;F\mid B)}\ge(1-3c_s-\eta)n+c_sb_B-O_{\eta,C_{\mathrm{ch}}}(L_{aw})-O_{\tau,C_{\mathrm{ch}}}(\log(n+2))
\end{equation}
for each class, where we used $N-n=O_{C_{\mathrm{ch}}}(1)$ and absorbed all fixed changes of the mesh into the error of the chart.

Finally, we average over the labels at the final scale. Their expectations satisfy
$$\mathbb EK_Z(N)\ge H_N(Z\mid A)-1\text{ and }\mathbb P(K_Z(N)\ge N+r\mid A)\le\min\{1,C_{\mathrm{ch}}^{O(1)}2^{-r}\}\text{ for }r\ge0,$$
where the second inequality follows by multiplying the largest mass permitted by the label by the number of bins in the fixed support. Summing this tail bound, we find that $\mathbb E(K_Z(N)-N)_+=O_{C_{\mathrm{ch}}}(1)$, so the truncation of $K_Z(N)$ costs only $O_{C_{\mathrm{ch}}}(1)$ in mean. Since $N\ge n$, refinement implies $H_N(Z\mid A)\ge H_n(Z\mid A)$, and summing over $X$, $Y$ and $F$ yields $\sum_Bwb_B\ge T_n-O_{C_{\mathrm{ch}}}(1)$, as announced. On the other side, \eqref{eq:label-MI} and the concavity of entropy bound the average of $C_n^B+q_s^{\mathrm{tot}}{\mathcal{I}_n(S;F\mid B)}$ by $C_n+q_s^{\mathrm{tot}}{\mathcal{I}_n(S;F\mid A)}+(1+q_s^{\mathrm{tot}})H(L)$, where $L$ is the full label, and the finite grid makes both the entropy of the label and the mean logarithmic cost of the normalization logarithmic in $n$. Choosing $\tau$ and $\eta$ small enough and absorbing these costs into the tolerance proves \eqref{eq:terminal-total-weighted}, and with it both output bounds and the statement at the endpoint.
\end{proof}

\begin{corollary}\label{cor:conditional-original-high-new}
If $f$ is nonspecial in the original coordinates and $\frac{2}{3}\le s<1$, then $M_n\ge(\frac{1+2s}{3}-\varepsilon)n-O_\varepsilon(1)$.
\end{corollary}

\begin{proof}
It suffices to prove $H_n(S)+2H_n(F)\ge(1+2s-\varepsilon)n-O_\varepsilon(1)$, since the left-hand side is at most $3M_n$ and the tolerance may be renamed after division. For $\frac{2}{3}\le s<\frac{3}{4}$, we apply Proposition~\ref{prop:two-copy-original}. Here $\beta_s=1$ and $2s+\frac{\beta_s}{3}-\frac{2}{3}=2s-\frac{1}{3}\ge1$, so the minimum in its coefficient equals one, which is the required inequality.

For $\frac{3}{4}\le s<1$, the conditional law is also jointly Frostman, and Corollary~\ref{cor:global-multiscale-profile}, applied in the original nonspecial representation with coefficient $\min\{1,2s-\frac{1}{2}\}=1$, shows that
$$H_n(F)+{\mathcal{I}_n}(X;Y)\ge(1-\varepsilon)n-O_\varepsilon(1).$$
By the two conditional entropy bounds, $H_n(X)+H_n(Y)\le2H_n(X,Y)-2sn+O(1)$. Since $f$ is nonspecial, it is not of the form $P(x+y)$, so $f_y-f_x\not\equiv0$, and Theorem~\ref{thm:coordinate-entropy} applied to $(S,F)$ yields $H_n(X,Y)\le H_n(S)+H_n(F)+O(1)$. Hence,
$${\mathcal{I}_n}(X;Y)\le H_n(S)+H_n(F)-2sn+O(1),$$
and substituting this bound, we arrive at the same weighted inequality for the outputs. The two ranges cover the whole interval $[\frac{2}{3},1)$, including $s=\frac{2}{3}$ and $s=\frac{3}{4}$, and the initial reduction proves the corollary.
\end{proof}

\subsubsection{Conditional entropy bounds in sum coordinates}

When a sum-coordinate representation is nonspecial, Lemma~\ref{lem:conditional-posterior-new} provides a stronger entropy bound, because the marginal entropy cancels before the iteration.

\begin{proposition}\label{prop:conditional-sum-profile-new}
If a sum-coordinate representation is nonspecial and $\frac{1}{2}<s<1$, then $M_n\ge(\Psi(s)-\varepsilon)n-O_\varepsilon(1)$, where
$$\Psi(s)=\begin{cases}
\frac{28s-3}{20}, & \frac{1}{2}<s\le\frac{2}{3},\\
\frac{22s+1}{20}, & \frac{2}{3}\le s\le\frac{3}{4},\\
\frac{1+s}{2}, & \frac{3}{4}\le s<1.
\end{cases}$$
\end{proposition}

\begin{proof}
We prove the estimate for the representation $(U,V)=(Y,S)$ when $G_Y$ is nonspecial. Exchanging $X$ and $Y$ preserves the hypotheses, the outputs $S$ and $F$ and the maximum $M_n$, and therefore covers the case of $G_X$. The associated original orientation satisfies $\mathcal W_x\not\equiv0$. Indeed, if $\mathcal W_x\equiv0$, then on an inverse chart $x=\phi(y,z)$ of $z=f(x,y)$, we have $\phi_{yz}=\frac{\mathcal W_x}{f_x^3}=0$, so that $\phi(y,z)=A(z)+B(y)$ there. Since $\phi_z=\frac{1}{f_x}\ne0$, the function $A$ has a local inverse, and the identity $S=Y+\phi(Y,F)$ would give $F=A^{-1}(S-Y-B(Y))$, a special representation of $G_Y$, contrary to the hypothesis.

We choose common buffered charts that are regular for $G_Y$, for Lemma~\ref{lem:conditional-posterior-new} in the orientation of $\mathcal W_x$, and for the recovery of the source from $(S,F)$, and on each such chart, we improve a hereditary lower bound for $C_j$. This suffices, because the conditional entropy bound on the event and the sum shear imply
$$E_n-H_n(Y)\ge sn-O_{C_{\mathrm{ch}}}(L_a)\text{ and }C_n\le H_n(S)+H_n(F)-sn+O_{C_{\mathrm{ch}}}(L_a),$$
so that a hereditary lower coefficient $k$ for $C_n$ translates into the positive coefficient $s+k$ for the sum of the outputs. Once this local coefficient has been improved, Lemma~\ref{lem:conditional-chart-averaging} passes the inequality for the outputs to the whole law.

The starting coefficient is $k_0=s$. Indeed, by the recovery of coordinates, $H_j(S)+H_j(F)\ge E_j-O_{C_{\mathrm{ch}}}(1)$, and hence $C_j\ge H_j(Y)-O_{C_{\mathrm{ch}}}(1)\ge sj-O_{C_{\mathrm{ch}}}(L_a)$. Suppose now that $C_j\ge(k-\eta)j-O_{\eta,C_{\mathrm{ch}}}(L_a)$ on every event in the chart. The identity $T_j=C_j+E_j$ and the conditional entropy bound for $E_j-H_j(Y)$ on the event lead to the first of the inequalities
$$\begin{aligned}
T_j&\ge(k+s-\eta)j+H_j(Y)-O_{\eta,C_{\mathrm{ch}}}(L_a),\\
T_j&\ge\left(\frac{5s}{2}+\beta_s-\eta\right)j-\frac{1}{2}H_j(Y)-O_{\eta,C_{\mathrm{ch}}}(L_a),
\end{aligned}$$
while the second follows from Lemma~\ref{lem:conditional-posterior-new}, applied in the orientation of $\mathcal W_x$ on the same event. With the weights $\frac{1}{3}$ and $\frac{2}{3}$, the coefficient of $H_j(Y)$ becomes $\frac{1}{3}-\frac{2}{3}\cdot\frac{1}{2}=0$, and we are left with the coefficient $\frac{k}{3}+2s+\frac{2\beta_s}{3}$ for $T_j$ that the next step requires.

The sum shear preserves the isotropic $2s$-Frostman bound, and the combined inequality holds on every event in the chart, so we can apply Lemma~\ref{lem:retained-information} with
$$\gamma=\frac{k}{3}+2s+\frac{2\beta_s}{3}\text{ and }b=0.$$
For the iterates below, $s\le k\le1$, and therefore $\frac{5}{3}<\gamma<3$ for $\frac{1}{2}<s<1$. From the conclusion of the lemma, we obtain the recurrence
$$k_0=s\text{ and }k_{r+1}=\min\left\{1,\frac{k_r}{6}+s+\frac{\beta_s}{3}-\frac{1}{4}\right\},$$
whose terms increase to
$$\lim_rk_r=\min\left\{1,\frac{12s+4\beta_s-3}{10}\right\}.$$
The update is nondecreasing, sends $s$ to a value at least $s$ and is capped at one, so it preserves the range $s\le k\le1$ and the hypotheses of the lemma hold at every step. The affine part of the recurrence contracts by the factor $\frac{1}{6}$, its fixed point exceeds $s$ when $s>\frac{1}{2}$, and the limit reaches one at $s=\frac{3}{4}$. We fix the target tolerance and then choose finitely many iterations and their intermediate tolerances, and after these choices, the hereditary error remains $O_{\eta,C_{\mathrm{ch}}}(L_a)$, with constants independent of the law, the event and the observation scale.

The initial reduction therefore yields every coefficient strictly below $s+\lim_rk_r$ for the sum of the outputs. To identify it, we substitute $\beta_s=\frac{3s}{2}$ for $s\le\frac{2}{3}$ and $\beta_s=1$ for $\frac{2}{3}\le s\le\frac{3}{4}$, and we use the limit one for $s\ge\frac{3}{4}$. Dividing the resulting coefficients by two, we obtain the three branches of $\Psi(s)$, which agree at both junctions. We apply Lemma~\ref{lem:conditional-chart-averaging} to the positive sum of the outputs, divide by two and rename the tolerance, which proves the result for $G_Y$, and the exchange of $X$ and $Y$ proves it for $G_X$.
\end{proof}

\subsection{Proofs of Theorems~\ref{thm:joint-main} and~\ref{thm:conditional-main}}\label{sec:main-proofs}

We prove the joint theorem first and then the conditional one. In each case, we first obtain the stated coefficients up to an arbitrary loss $\varepsilon$, while the estimates without this loss come from separate entropy estimates and not from a passage to the limit $\varepsilon\downarrow0$.

For the joint theorem, Proposition~\ref{prop:joint-mass-classes} provides estimates with error $O(\log(n+2))$ for $\frac{5}{6}\le s<1$. For a nonspecial sum-coordinate representation, Corollary~\ref{cor:profile-joint} uses $(a,b,\zeta)=(s,2s-1,2s-1)$, and the Fourier estimate is critical at $s=\frac{9}{10}$, with a bounded error for $s>\frac{9}{10}$. In the remaining case, where only the original representation is nonspecial, the coefficient $\frac{2+5s}{7}$ obtained from the classes of masses exceeds the target $\frac{1+3s}{4}$ by $\frac{1-s}{28}$, and for each fixed $s<1$, this positive surplus absorbs the logarithmic error. For the conditional theorem, Corollary~\ref{cor:profile-conditional} uses $a=b=\zeta=s$. The decay exponent $2s-\frac{3}{2}$ of the reference equals $1-\zeta$ at $s=\frac{5}{6}$ and exceeds it for $s>\frac{5}{6}$, so Lemma~\ref{lem:output-profile-transfer} produces a logarithmic error at the endpoint and a bounded error above it.

\subsubsection{Proof of Theorem~\ref{thm:joint-main}}

\begin{proof}
We first derive the rates in \eqref{eq:joint-main} and the estimate \eqref{eq:joint-main-bounded} with bounded error from the incidence and the multiscale bounds, and then we prove \eqref{eq:joint-main-endpoint} and its version with bounded error by means of the direct entropy estimates.

Theorem~\ref{thm:joint-incidence-bound} proves \eqref{eq:joint-main-bounded} throughout $0<s<1$. For $0<s\le\frac{1}{2}$, we have $\min\{s,2-2s\}=s$, so this estimate also yields the gain $\frac{s^2}{64d}$ required in the first branch of \eqref{eq:joint-main-gain}. For $\frac{1}{2}<s\le\frac{2}{3}$, the same minimum equals $s$, and the gain $\frac{s}{16}$ follows from \eqref{eq:joint-regular-gain}. Corollary~\ref{cor:joint-multiscale-gain} accounts for the gain $\frac{2s-1}{8}$ for $\frac{2}{3}\le s\le\frac{3}{4}$ and $\frac{1-s}{4}$ for $\frac{3}{4}\le s<1$. Adjacent formulas agree at $\frac{2}{3}$ and $\frac{3}{4}$, the point $s=\frac{1}{2}$ belongs to the first range, and the intervals cover $0<s<1$. Since each estimate holds for every $\varepsilon>0$, this proves \eqref{eq:joint-main}, including its special case for $\frac{3}{4}\le s<1$.

Now let $\frac{5}{6}\le s<1$. We seek a lower bound for $4M_{\mathrm j,n}$ with coefficient $1+3s$ and logarithmic error. By Corollary~\ref{cor:three-coordinate}, one of $f$, $G_Y$ and $G_X$ is nonspecial. Suppose first that $G_Y$ is nonspecial. Since $\frac{3s}{2}-\frac{1}{4}\ge1$, the minimum in \eqref{eq:joint-mass-class-comparison} equals one, and Proposition~\ref{prop:joint-mass-classes} yields
$$H_n(F)+{\mathcal{I}_n}(Y;S)+H_n(X)\ge(1+s)n-O(\log(n+2)).$$
The shear $(X,Y)\mapsto(Y,S)$ preserves the joint dyadic entropy up to $O(1)$, so by the joint Frostman bound, $H_n(Y,S)=H_n(X,Y)+O(1)\ge2sn-O(1)$. Adding this lower bound to the preceding inequality cancels the negative term $-H_n(Y,S)$ contained in ${\mathcal{I}_n}(Y;S)$, and hence
$$4M_{\mathrm j,n}\ge H_n(X)+H_n(Y)+H_n(S)+H_n(F)\ge(1+3s)n-O(\log(n+2)),$$
where the first inequality follows from the definition of $M_{\mathrm j,n}$. Dividing by four proves \eqref{eq:joint-main-endpoint} in this case. If $G_X$ is nonspecial, we use the version of the same proposition with $X$ and $Y$ interchanged and the shear $(X,Y)\mapsto(X,S)$, and the sum of the four entropies and the resulting bound for $M_{\mathrm j,n}$ are the same.

It remains to treat the case where both sum-coordinate representations are special. Then $f$ is nonspecial by Corollary~\ref{cor:three-coordinate}, and averaging \eqref{eq:joint-original-mass-class-comparison} over its two orientations, in which the same minimum equals one, we get
$$H_n(F)+{\mathcal{I}_n}(X;Y)+\frac{1}{4}(H_n(X)+H_n(Y))\ge\left(1+\frac{s}{2}\right)n-O(\log(n+2)).$$
Expanding ${\mathcal{I}_n}(X;Y)$ and adding the lower bound $H_n(X,Y)\ge2sn-O(1)$, we find
$$H_n(F)+\frac{5}{4}(H_n(X)+H_n(Y))\ge\left(1+\frac{5s}{2}\right)n-O(\log(n+2)).$$
By definition, $H_n(F)\le M_{\mathrm j,n}$ and $H_n(X)+H_n(Y)\le2M_{\mathrm j,n}$, so the left-hand side is at most $\frac{7}{2}M_{\mathrm j,n}$, and after dividing by $\frac{7}{2}$, we obtain
$$M_{\mathrm j,n}\ge\frac{2+5s}{7}n-O(\log(n+2)).$$
Since $\frac{2+5s}{7}-\frac{1+3s}{4}=\frac{1-s}{28}>0$, this coefficient is larger than the one required in \eqref{eq:joint-main-endpoint}, which proves that estimate in the remaining case and hence for every admissible polynomial.

Finally, let $\frac{9}{10}<s<1$. If one of the sum-coordinate representations is nonspecial, Corollary~\ref{cor:profile-joint} shows directly that $M_{\mathrm j,n}\ge\frac{(1+3s)n}{4}-O(1)$. Otherwise, the estimate of the preceding paragraph has the fixed positive surplus $\frac{1-s}{28}$ in the coefficient, and for fixed $s<1$, a logarithmic loss minus $\frac{(1-s)n}{28}$ is bounded above uniformly in $n$, so the surplus absorbs the logarithmic error and the same estimate holds with a bounded error. All estimates are uniform over laws with the prescribed support and Frostman constant, and this proves the final assertion of the theorem.
\end{proof}

\subsubsection{Proof of Theorem~\ref{thm:conditional-main}}

\begin{proof}
We first prove \eqref{eq:conditional-main} by selecting the appropriate coefficient on each interval of exponents. The available estimates depend on whether the nonspecial representation uses the original or the sum coordinates, and we compare them so that Corollary~\ref{cor:three-coordinate} yields the required bound for every admissible polynomial. The stronger additive errors are obtained separately at the end.

Proposition~\ref{prop:conditional-inverse-output} applies to every admissible polynomial for $0<s\le\frac{2}{3}$. It establishes the coefficient $\frac{35s}{31}=s+\frac{4s}{31}$ throughout $0<s\le s_{\mathrm a}$, which proves the first branch of \eqref{eq:conditional-main-gain}.

For the two middle branches, we first compare the estimates in the original and in the sum coordinates on $\frac{5}{9}\le s<\frac{2}{3}$, an interval containing $[s_{\mathrm a},\frac{2}{3})$, since $\frac{5}{9}<s_{\mathrm a}<s_{\mathrm b}<\frac{2}{3}$. If $f$ is nonspecial in the original coordinates, Proposition~\ref{prop:conditional-terminal-total-new} provides both coefficients in \eqref{eq:terminal-total-output} and \eqref{eq:terminal-inverse-output}, and if a sum-coordinate representation is nonspecial, the first branch of Proposition~\ref{prop:conditional-sum-profile-new} leads to $\frac{28s-3}{20}$. The exact comparisons are
$$\begin{aligned}
\frac{(5s-1)(7-6s)}{19-15s}-\frac{28s-3}{20}&=\frac{180s^2-243s+83}{20(15s-19)}<0,\\
\frac{49s-30s^2-7}{26-15s}-\frac{28s-3}{20}&=\frac{180s^2-207s+62}{20(15s-26)}<0.
\end{aligned}$$
Indeed, the denominators are negative on this interval, while both numerators have positive leading coefficient and negative discriminant, equal to $-711$ and $-1791$, respectively, and are therefore positive. Thus, the coefficient in sum coordinates implies both bounds in the original coordinates with the same tolerance, and since Corollary~\ref{cor:three-coordinate} guarantees that at least one of the three representations is nonspecial, both bounds hold for every admissible polynomial on $\frac{5}{9}\le s<\frac{2}{3}$.

We now compare these bounds with the coefficient $\frac{35s}{31}$ in order to identify the two middle branches. Subtracting, we find
$$\begin{aligned}
\frac{49s-30s^2-7}{26-15s}-\frac{35s}{31}&=\frac{405s^2-609s+217}{31(15s-26)},\\
\frac{49s-30s^2-7}{26-15s}-\frac{(5s-1)(7-6s)}{19-15s}&=\frac{90s^2-135s+49}{(15s-26)(15s-19)}.
\end{aligned}$$
The smaller roots of the two numerators are $s_{\mathrm a}$ and $s_{\mathrm b}$, respectively, and their other roots lie above $\frac{2}{3}$. The first denominator is negative, so the coefficient in \eqref{eq:terminal-inverse-output} is at least $\frac{35s}{31}$ for $s_{\mathrm a}\le s<\frac{2}{3}$. The second denominator is positive, while its numerator is nonnegative up to $s_{\mathrm b}$ and nonpositive between $s_{\mathrm b}$ and $\frac{2}{3}$. We therefore use \eqref{eq:terminal-inverse-output} for $s_{\mathrm a}\le s\le s_{\mathrm b}$ and \eqref{eq:terminal-total-output} for $s_{\mathrm b}\le s<\frac{2}{3}$. Subtracting $s$ from these two coefficients, we obtain
$$\frac{23s-15s^2-7}{26-15s}\text{ and }\frac{(1-s)(15s-7)}{19-15s},$$
which are the second and the third branches of \eqref{eq:conditional-main-gain}, and the differences displayed above also show that adjacent coefficients agree at $s_{\mathrm a}$ and at $s_{\mathrm b}$.

For $\frac{2}{3}\le s<1$, Corollary~\ref{cor:conditional-original-high-new} gives the coefficient $\frac{1+2s}{3}$ when $f$ is nonspecial in the original coordinates, and when a sum-coordinate representation is nonspecial, Proposition~\ref{prop:conditional-sum-profile-new} yields the coefficient $\Psi(s)$, which exceeds $\frac{1+2s}{3}$ by
$$\frac{26s-17}{60}>0\text{ for }\frac{2}{3}\le s\le\frac{3}{4}\text{ and by }\frac{1-s}{6}>0\text{ for }\frac{3}{4}\le s<1.$$
By the classification, the coefficient $\frac{1+2s}{3}$ therefore holds for every admissible polynomial on this whole interval, and its gain over $s$ is $\frac{1-s}{3}$, as required in the fourth branch. At $s=\frac{2}{3}$, this agrees with the formula of the third branch, so the endpoint is covered without extending the preceding proposition beyond its range. Choosing the tolerances in the cited estimates no larger than the prescribed $\varepsilon$ proves \eqref{eq:conditional-main} throughout $0<s<1$.

It remains to prove the stronger error bounds for $\frac{5}{6}\le s<1$. Corollary~\ref{cor:profile-conditional} applies directly to the given admissible polynomial and conditionally Frostman law, and it provides the coefficient $\frac{1+2s}{3}$ with error $O(\log(n+2))$ at $s=\frac{5}{6}$ and $O(1)$ for $\frac{5}{6}<s<1$. These bounds come from the direct entropy estimate and require no passage to the limit $\varepsilon\downarrow0$ in \eqref{eq:conditional-main}, and the errors are uniform over the prescribed source laws and $n\ge1$. This proves \eqref{eq:conditional-main-endpoint} and its version with bounded error.
\end{proof}

\section{Proof of Theorem~\ref{thm:conditional-classification}}
\label{sec:sharpness}

We first establish the lower bounds. The matching laws are constructed in Section~\ref{sec:conditional-constructions}, and invoking those results completes the proof, including the assertions about uniformity.

\begin{proof}
Throughout, $S=X+Y$ and $F=f(X,Y)$. By Theorem~\ref{thm:two-transition}, both sum-coordinate representations are special exactly when $f\in\mathfrak E_S$. If $f\in\mathfrak E_S$, Theorem~\ref{thm:conditional-main} shows that
$$M_n\ge\left(\frac{1+2s}{3}-\varepsilon\right)n-O_\varepsilon(1)\text{ for }\frac{2}{3}\le s<1,$$
and the same theorem provides the coefficient $\frac{1+2s}{3}$ with error $O(\log(n+2))$ at $s=\frac{5}{6}$ and $O(1)$ for $\frac{5}{6}<s<1$. Thus, in this case, the lower bounds are inherited from Theorem~\ref{thm:conditional-main}.

If $f\notin\mathfrak E_S$, at least one sum-coordinate representation is nonspecial, and the branch $\frac{3}{4}\le s<1$ of Proposition~\ref{prop:conditional-sum-profile-new} gives
$$M_n\ge\left(\frac{1+s}{2}-\varepsilon\right)n-O_\varepsilon(1).$$
To obtain the exact coefficient for $\frac{5}{6}\le s<1$, we use the entropy inequality established in the proof of Corollary~\ref{cor:profile-conditional} in the nonspecial sum-coordinate representation. Suppose first that $G_Y$ is nonspecial, so that the source coordinates are $(Y,S)$. Then this inequality takes the form
$$H_n(F)+{\mathcal{I}_n}(Y;S)\ge n-O(\log(n+2)),$$
with error $O(1)$ when $\frac{5}{6}<s<1$.

We bound the mutual information in terms of the entropy of the sum by means of the exact conditional Frostman hypothesis. Given $Y=y$, the map from $X$ to $S$ is a translation, so the conditional law of $S$ given $Y=y$ is $(s,C_1)$-Frostman for almost every $y$, and every atom of $D_nS$ under this conditional law has mass at most $C_12^{-sn}$. Since $D_nY$ is a function of $Y$, conditioning on $D_nY$ leaves at least as much entropy as conditioning on $Y$, and hence
$$H(D_nS\mid D_nY)\ge H(D_nS\mid Y)\ge sn-\log C_1,$$
so that
$${\mathcal{I}_n}(Y;S)=H_n(S)-H(D_nS\mid D_nY)\le H_n(S)-sn+\log C_1.$$
Substituting this upper bound into the entropy inequality, we get
$$2M_n\ge H_n(S)+H_n(F)\ge(1+s)n-O(\log(n+2)),$$
with error $O(1)$ when $\frac{5}{6}<s<1$, and dividing by two, we obtain the coefficient $\frac{1+s}{2}$ with the corresponding error. If $G_X$ is nonspecial, we use the entropy inequality with ${\mathcal{I}_n}(X;S)$, condition on $X$ and apply the exact conditional bound for $Y$ with constant $C_2$, and the same calculation leads to the same bound for the outputs. This proves all the stated lower bounds.

For the matching examples, we apply the constructions in the open square $(-\frac{R}{2},\frac{R}{2})^2$. If $f\in\mathfrak E_S$, we write $f=h(ax+P(x+y))$. Admissibility forces $h$ to be nonconstant, $a\ne0$ and $\deg P\ge2$, since a constant $h$ makes $f$ constant, $a=0$ makes $f$ a function of $x+y$, and $\deg P\le1$ makes $f$ a function of one affine coordinate. We fix this representation before choosing $n$. Proposition~\ref{upper-modular:prop} applies throughout $\frac{2}{3}\le s<1$, since its construction is valid for $\frac{1}{4}\le s<1$, and for every sufficiently large $n$, it produces a conditionally Frostman law for which each output entropy equals $\frac{1+2s}{3}n+O(1)$, with constants independent of $n$. As explained after that proposition, a fixed uniform product law takes care of the finitely many remaining coarse scales, and enlarging the Frostman and entropy constants once covers all of these scales.

If $f\notin\mathfrak E_S$, we apply Proposition~\ref{prop:balanced-conditional-source} in the same open square. It provides a single fixed conditionally Frostman law for which each output entropy equals $\frac{1+s}{2}n+O(1)$ for every $n\ge1$, and since the proposition applies for every $0<s<1$, it covers the required range $\frac{3}{4}\le s<1$.

All these laws are supported in $[-R,R]^2$, their conditional Frostman constants are bounded by a single finite constant $C_0=C_0(s,f,R)$, and their entropy errors are bounded independently of $n$. Outside $\mathfrak E_S$, the law is fixed, whereas inside the class, it may vary with $n$. Together with the lower bounds, these examples prove every assertion of the theorem.
\end{proof}

\section{Joint output bounds and polynomial image dimensions}
\label{sec:joint-models}

We prove the lower bounds in Proposition~\ref{thm:joint-models}, and then we deduce Corollary~\ref{cor:polynomial-image-product} from the entropy estimate for independent inputs. The matching laws for the proposition are constructed in Section~\ref{sec:prescribed-output-laws}.

\subsection{Lower bounds in Proposition~\ref{thm:joint-models}}

For $s>\frac{1}{2}$, the lower bound follows from the recovery of polynomial coordinates whenever $f_y-f_x\not\equiv0$. For the quadratic models and smaller exponents, the degeneracy of the joint map of the two outputs is compensated by coarse information that the two outputs share.

\subsubsection{The quadratic entropy estimate}

\begin{lemma}\label{joint-models:radial-square}
Let $0<t\le2$ and $A\ge1$, and let $(U,V)$ be supported in $[-1,1]^2$, with every axis-parallel square of side $h$ having mass at most $Ah^t$. There is an absolute constant $K$ such that, for either choice of the sign,
\begin{equation}\label{joint-models:radial-square-bound}
H_n(U)+H_n(U^2\pm V^2)\ge tn-K(\log(A+2)+\log(n+2)).
\end{equation}
\end{lemma}

\begin{proof}
Fix the sign and write $W=U^2\pm V^2$. We first work on an annulus, where coarse information shared by the two outputs compensates for the possible loss in their joint entropy, and the dependence on $A$ will then allow us to sum over dyadic annuli around the origin.

Suppose first that $\frac{1}{2}<\max\{|U|,|V|\}\le1$. With $H_0$ defined by the same convention of quantization, we claim that
\begin{equation}\label{joint-models:annulus}
H_m(U)+H_m(W)\ge tm-2\log A-K\log(m+2)\text{ for all }m\ge0.
\end{equation}
For $0\le m\le5$, the claim follows from the nonnegativity of entropy and $t\le2$ after increasing the absolute constant $K$. For $m\ge6$, we set $\delta=2^{-m}$ and $k_*=\lfloor\frac{m}{2}\rfloor$, a threshold chosen because $2^{-2k_*}$ is comparable to the output resolution $\delta$. We partition the annulus into the region where $|V|>\frac{1}{4}$, the shells where $2^{-k-1}<|V|\le2^{-k}$ for $2\le k<k_*$, and the strip where $|V|\le2^{-k_*}$. We split the first region according to the sign of $V$, the intermediate shells according to the signs of both $U$ and $V$, and the last strip according to the sign of $U$, so that the inverse estimates below are single-valued on each piece. These disjoint pieces cover the annulus, and there are $O(m+1)$ of them. On a piece of probability $w>0$, the square bound holds with constant $\frac{A}{w}$, and until the averaging step, all entropies refer to the conditional law on the piece.

Consider first a piece on which $|V|>\frac{1}{4}$ and $V$ has a fixed sign. For two points $(U,V)$ and $(U',V')$ of this piece, with $W'=U'^2\pm V'^2$, the outputs satisfy
$$|V-V'||V+V'|=|W-W'-(U-U')(U+U')|,$$
and since $|V+V'|>\frac{1}{2}$, this provides a uniformly Lipschitz inverse for $(U,V)\mapsto(U,W)$. Each output square of side $\delta$ therefore lifts into an input square of side $O(\delta)$, and the square bound implies
$$H_m(U)+H_m(W)\ge H_m(U,W)\ge tm-\log\frac{A}{w}-O(1)\ge tm-2\log\frac{A}{w}-O(1),$$
where the last inequality uses $\frac{A}{w}\ge1$.

On an intermediate shell, we set $r=2^{-k}$. An output square at level $m$ determines $U$ up to $O(\delta)$ and $V^2=|W-U^2|$ up to $O(\delta)$, and since $V$ has a fixed sign and $|V|>\frac{r}{2}$, it determines $V$ up to $O(\frac{\delta}{r})$. Moreover, a dyadic interval of $U$ at level $2k$ together with the shell lies in a square of side $O(r)$. Applying the square bound at these two scales, we get
$$H_m(U,W)\ge t(m-k)-\log\frac{A}{w}-O(1)\text{ and }H_{2k}(U)\ge tk-\log\frac{A}{w}-O(1).$$
The first estimate loses $tk$ compared with the target $tm$, and we recover this amount from the second one by showing that the coarse information about $U$ is also determined, up to bounded ambiguity, by the observed value of $W$. On these shells, $|U|>\frac{1}{2}$ and $U$ has a fixed sign. Since $|W-U^2|=V^2\le r^2$ for either sign and $\delta\le r^2$, a dyadic interval of $W$ at level $m$ allows only $O(1)$ dyadic intervals of $U$ at level $2k$. Thus,
$$H(D_{2k}U\mid D_mW)=O(1)\text{ and }{\mathcal{I}}(D_mU;D_mW)\ge H_{2k}(U)-O(1),$$
where the second inequality uses that $D_{2k}U$ is a function of $D_mU$. Adding the bounds for the joint entropy and for the mutual information yields
$$H_m(U)+H_m(W)=H_m(U,W)+{\mathcal{I}}(D_mU;D_mW)\ge tm-2\log\frac{A}{w}-O(1).$$

On the last strip, we set $r=2^{-k_*}$, so that $\delta\le r^2\le2\delta$. A dyadic interval of $U$ at level $m$ intersected with this strip lies in a square of side $O(r)$, and therefore
$$H_m(U)\ge tk_*-\log\frac{A}{w}-O(1).$$
Here again, $|U|>\frac{1}{2}$ and $U$ has a fixed sign, so the relation $W=U^2+O(\delta)$ implies $H(D_mU\mid D_mW)=O(1)$, and hence
$$H_m(U)+H_m(W)\ge2H_m(U)-O(1)\ge2tk_*-2\log\frac{A}{w}-O(1)\ge tm-2\log\frac{A}{w}-O(1),$$
where the last step uses $m-2k_*\le1$ and $t\le2$. All the geometric errors above are absolute, because the fixed geometric constants appear with exponents at most two.

Let $L$ be the label of the piece. The concavity of entropy and the costs of the normalizations lead to
$$H_m(U)+H_m(W)\ge tm-2\log A-2H(L)-O(1),$$
and since $H(L)=O(\log(m+2))$, this proves \eqref{joint-models:annulus}.

We now remove the restriction to an annulus. Let $J=\lfloor\frac{n}{2}\rfloor$ and set
$$E_j=\{2^{-j-1}<\max\{|U|,|V|\}\le2^{-j}\}\text{ and }p_j=\mathbb P(E_j)\text{ for }0\le j\le J,$$
omitting annuli of probability zero, and let $p_\circ$ be the probability of the remaining central piece. This piece lies in a square of side $2^{-J}$, so
$$p_\circ\le\min\{1,A2^{-tJ}\}\le\min\{1,4A2^{-\frac{tn}{2}}\}.$$
On $E_j$, multiplying the inputs by $2^j$ produces the normalized annulus, with constant $\frac{4A2^{-tj}}{p_j}$ in the square bound, which is at least one because $p_j\le4A2^{-tj}$. For either sign, the homogeneity of the quadratic forms shows that $W$ is observed at the normalized precision $m=n-2j\ge0$, while $U$ is observed at the finer precision $n-j$. These identities between entropies are exact, because the rescaling is dyadic, and since entropy increases under refinement, we may lower the precision for $U$ from $n-j$ to $m=n-2j$. Applying \eqref{joint-models:annulus} at this common precision, we find
$$\begin{aligned}
H_n(U\mid E_j)+H_n(W\mid E_j)&\ge t(n-2j)-2\log\frac{4A2^{-tj}}{p_j}-K\log(n+2)\\
&=tn-2\log(4A)+2\log p_j-K\log(n+2).
\end{aligned}$$
The terms involving $j$ cancel, and this is why the estimate on annuli can be averaged without a loss in the leading coefficient.

The label of the annuli, including the central piece, takes at most $J+2$ values. Averaging by concavity and using the nonnegativity of entropy on the central piece, we obtain
$$H_n(U)+H_n(W)\ge tn-tnp_\circ-2\log(4A)-O(\log(n+2)).$$
It remains to bound the loss caused by the central mass uniformly in $t$. Let $z=\frac{tn}{2}$ and $b=\log(4A)\ge2$. If $z\le b$, then $z\min\{1,2^{b-z}\}\le b$, and if $z>b$, then
$$z\min\{1,2^{b-z}\}=b2^{-(z-b)}+(z-b)2^{-(z-b)}\le b+O(1).$$
Hence, $tnp_\circ=O(\log(A+2))$, and substituting this bound proves \eqref{joint-models:radial-square-bound}, uniformly for $0<t\le2$ and for both choices of the sign.
\end{proof}

\subsubsection{Proof of the lower bounds}

\begin{proof}[Proof of the lower bounds in Proposition~\ref{thm:joint-models}]
We prove the bound for the triangular polynomials by recovering the input cells from the two output cells, we reduce the quadratic models to Lemma~\ref{joint-models:radial-square} by a fixed linear change of coordinates, and we obtain the general bound for $s>\frac{1}{2}$ from the recovery of polynomial coordinates. Throughout, we write $\mathcal T(x,y)=(x+y,f(x,y))$.

Suppose first that $f=ax+P(x+y)$. Since $a\ne0$, the map $\mathcal T$ has the inverse
$$\mathcal T^{-1}(u,v)=\left(\frac{v-P(u)}{a},u-\frac{v-P(u)}{a}\right),$$
which is Lipschitz on a fixed bounded neighborhood of the range of the outputs. Hence, a dyadic output square meets the images of only boundedly many input squares of the same scale, so that
$$H(D_n(X,Y)\mid D_n(S,F))=O_{a,P,R}(1),$$
and by the chain rule, $H_n(S,F)\ge H_n(X,Y)-O_{a,P,R}(1)$. By the joint Frostman bound, each input square has mass at most $C2^{-2sn}$, and combining this atom bound with subadditivity, we get
$$2M_n\ge H_n(S)+H_n(F)\ge H_n(S,F)\ge H_n(X,Y)-O_{a,P,R}(1)\ge2sn-\log C-O_{a,P,R}(1).$$
Dividing by two proves \eqref{joint-models:triangular-lower}, with the stated dependence of the additive constant.

For the two quadratic models, we set
$$U=\frac{X+Y}{1+2R}\text{ and }V=\frac{X-Y}{1+2R},$$
so that $(U,V)$ is supported in $[-1,1]^2$. To check the square bound required by Lemma~\ref{joint-models:radial-square}, observe that a square of side $h$ in these coordinates is contained, in the original coordinates, in an axis-parallel square of side $(1+2R)h$. When $(1+2R)h\le1$, the joint Frostman bound assigns mass at most $C(1+2R)^{2s}h^{2s}$, and when $(1+2R)h>1$, the trivial bound by one is at most $(1+2R)^{2s}h^{2s}$. Hence, the required constant is
$$A=\max\{C,1\}(1+2R)^{2s}\ge1,$$
and the lemma with $t=2s\in(0,2]$ yields
$$H_n(U)+H_n(U^2\pm V^2)\ge2sn-K(\log(A+2)+\log(n+2)).$$
The original outputs are $S=(1+2R)U$ and
$$F=\begin{cases}
\frac{(1+2R)^2(U^2+V^2)}{2} & \text{if }f=x^2+y^2,\\
\frac{(1+2R)^2(U^2-V^2)}{4} & \text{if }f=xy,
\end{cases}$$
so we use the plus sign for $x^2+y^2$ and the minus sign for $xy$. Each fixed rescaling changes dyadic entropies by $O_R(1)$, and $\log(A+2)\le\log(C+2)+O_R(1)$ uniformly for $0<s\le1$, so substituting proves \eqref{joint-models:radial-lower}, including the endpoint $s=1$.

For the general recovery assertion, let $\frac{1}{2}<s\le1$ and $f_y-f_x\not\equiv0$. We apply the isotropic part of Theorem~\ref{thm:coordinate-entropy} with $\mathsf L=x+y$, $P=f$ and $t=2s\in(1,2]$. Its Jacobian condition is exactly $J_{\mathsf L,P}=f_y-f_x\not\equiv0$, and its square bound follows from the joint Frostman hypothesis, so the theorem shows that
$$H_n(S,F)=H_n(X,Y)+O_{s,f,R,C}(1).$$
The square bound also implies $H_n(X,Y)\ge2sn-\log C$, which proves \eqref{joint-models:general-recovery}. For an admissible polynomial, the Jacobian condition holds automatically, since otherwise the last assertion of Theorem~\ref{thm:coordinate-entropy} would give $f(x,y)=q(x+y)$ for a real polynomial $q$ of one variable, contrary to admissibility. This proves all the lower bounds in the proposition, and the matching laws are constructed in Section~\ref{sec:prescribed-output-laws}.
\end{proof}

\subsection{Proof of Corollary~\ref{cor:polynomial-image-product}}
\label{sec:polynomial-image-products}

\begin{proof}[Proof of Corollary~\ref{cor:polynomial-image-product}]
Fix $0<t<s$ and $\varepsilon>0$. By Frostman's lemma, some compact subset $K\subset A$ supports a probability measure $\mu$ with $\mu(I)\le C_t|I|^t$ for every interval $I$ of length at most one. Let $X$ and $Y$ be independent with law $\mu$, and choose a fixed interval $[-R,R]$ containing their supports. Since $0<t<1$ and $P$ and $Q$ satisfy the hypotheses of Theorem~\ref{thm:independent-polynomial-pair}, we may apply that theorem with exponent $t$ and both Frostman constants equal to $C_t$, and we obtain for every $n\ge1$
$$H_n(P(X,Y))+H_n(Q(X,Y))\ge\left(\min\left\{\frac{5t}{2},1+t\right\}-\varepsilon\right)n-O_{t,\varepsilon}(1),$$
with an additive constant independent of $n$.

For a bounded set $E$ in a Euclidean space, let $N_n(E)$ denote the number of half-open dyadic cubes of side $2^{-n}$ that meet $E$. The two output laws are supported in $P(A,A)$ and $Q(A,A)$, so their entropies are at most the logarithms of the corresponding covering numbers, and the Cartesian structure of the dyadic grid leads to the exact identity
$$N_n(P(A,A)\times Q(A,A))=N_n(P(A,A))N_n(Q(A,A)).$$
Consequently,
$$\begin{aligned}
\log N_n(P(A,A)\times Q(A,A))&\ge H_n(P(X,Y))+H_n(Q(X,Y))\\
&\ge\left(\min\left\{\frac{5t}{2},1+t\right\}-\varepsilon\right)n-O_{t,\varepsilon}(1).
\end{aligned}$$
Dividing by $n$ and taking the lower limit as $n\to\infty$, the additive error disappears, and the definition of the lower box dimension through dyadic coverings yields
$$\underline{\dim}_{\mathrm B}(P(A,A)\times Q(A,A))\ge\min\left\{\frac{5t}{2},1+t\right\}-\varepsilon.$$
Letting $\varepsilon\downarrow0$ and then $t\uparrow s$, and using the continuity of the coefficient, we obtain \eqref{eq:polynomial-image-product}, including the case $s=1$.
\end{proof}

\section{Constructions and counterexamples}
\label{sec:sharpness-counterexamples}

In this section, we collect the constructions used for the matching assertions and the upper bounds throughout the paper. After some common scalar sources and estimates for mixtures, the joint and the conditional examples follow the order of Theorems~\ref{thm:joint-main}--\ref{thm:conditional-classification}. We then prescribe the two output laws for Proposition~\ref{thm:joint-models}, construct the fixed law of Corollary~\ref{cor:sharp-three-images}, and give the arithmetic counterexamples. Each subsection states which coefficients are attained and which further improvements the examples rule out. For each fixed choice of the parameters, the Frostman constants and the entropy errors of the constructions do not depend on the observation scale.

\subsection{Scalar laws, progressions, and mixtures}
\label{sec:construction-preliminaries}

\subsubsection{Scalar laws with prescribed entropy}

\begin{lemma}\label{joint-models:scalar-source}
For $0<s\le1$ and every compact nondegenerate interval $J$, there is a probability measure $\mu$ on $J$ such that
$$\mu(I)\le C_{s,J}|I|^s\text{ for every interval }I\text{ and }H_n(\mu)=sn+O_{s,J}(1).$$
Its support meets $O_{s,J}(2^{sn})$ dyadic intervals at level $n$, and it is a Lebesgue null set if $s<1$. For $s=1$, one may take the uniform law on $J$.
\end{lemma}

\begin{proof}
For $s=1$, normalized Lebesgue measure on $J$ assigns to intervals mass at most $C_J|I|$ and meets $O_J(2^n)$ dyadic intervals at level $n$, and together these two bounds show that $H_n(\mu)=n+O_J(1)$.

For $0<s<1$, we use a Cantor law in which the masses of the construction intervals are the $s$th powers of their lengths. Let $\lambda=2^{-\frac{1}{s}}<\frac{1}{2}$, and take the self-similar law on $[0,1]$ with equal weights generated by the maps $x\mapsto\lambda x$ and $x\mapsto1-\lambda+\lambda x$. Its construction intervals at step $k$ have length $\lambda^k$ and mass $2^{-k}=\lambda^{ks}$, and they are separated by at least $(1-2\lambda)\lambda^{k-1}$. Given $0<r<1$, we choose $k$ with $\lambda^k\le r<\lambda^{k-1}$. An interval of length $r$ meets only $O_s(1)$ construction intervals of step $k$, so it has mass $O_s(\lambda^{ks})=O_s(r^s)$, and intervals of length at least one satisfy the same bound after enlarging the constant.

At $r=2^{-n}$, the choice of $k$ ensures that $k<sn+1$. Each construction interval of step $k$ has length at most $2^{-n}$, so it meets at most two dyadic intervals at level $n$, while there are $2^k=O_s(2^{sn})$ of them. Thus, the support meets $O_s(2^{sn})$ dyadic intervals at level $n$. The interval bound implies $H_n(\mu)\ge sn-O_s(1)$, and this count implies $H_n(\mu)\le sn+O_s(1)$. The covering intervals have total length $O_s(2^{-(1-s)n})$, which tends to zero, so the support is a Lebesgue null set. An affine map from $[0,1]$ onto $J$ preserves all these conclusions and changes only the constants.
\end{proof}

\subsubsection{Progressions and finite mixtures}

Let $0<s<1$, $\delta=2^{-n}$ and $N\asymp\delta^{-s}$. By a smoothed progression of exponent $s$ at level $n$, we mean the uniform law on the centers $a+\frac{\ell_0j}{N}$, $0\le j<N$, with fixed $a\in\mathbb R$ and $\ell_0>0$, convolved with the uniform law on an interval of length $c\delta$ for a fixed $c>0$. Two smoothed progressions are aligned if they have the same step $\frac{\ell_0}{N}$.

\begin{lemma}\label{lem:progressions-mixtures}
A smoothed progression of exponent $s$ at level $n$ is $s$-Frostman with a constant independent of $n$, and its dyadic entropy at level $n$ is $sn+O(1)$. Fixed bi-Lipschitz changes of variables preserve these conclusions. The sum and the difference of two independent aligned smoothed progressions have entropy at most $sn+O(1)$, and a fixed Lipschitz scalar function of two such variables has entropy at most $\min\{2s,1\}n+O(1)$.

Finite mixtures preserve common joint and common exact conditional Frostman bounds. If $L$ is the label of the component, then, for every scalar or vector observable $W$,
\begin{equation}\label{eq:finite-mixture-upper}
\sum_\ell\mathbb P(L=\ell)H_n(W\mid L=\ell)\le H_n(W)\le H(L)+\sum_\ell\mathbb P(L=\ell)H_n(W\mid L=\ell).
\end{equation}
\end{lemma}

\begin{proof}
Let $Z$ be a smoothed progression. We first prove an interval bound uniform in $n$ and then use the count of the support to obtain both entropy bounds, with constants depending on the fixed spacing and noise parameters. For $\delta\le r\le1$, at most $O(Nr+1)$ centers contribute to an interval of length $r$, so
$$\mathbb P(Z\in I)\le C(r+N^{-1})\le C'r^s,$$
because $N^{-1}\asymp\delta^s\le r^s$ and $r\le r^s$. Below the smoothing scale, the number of centers alone does not give the required bound, and we use the density instead. The ratio $\frac{cN\delta}{\ell_0}$ of the length of the noise intervals to the spacing of the centers is $O(\delta^{1-s})=O(1)$, so the noise intervals have bounded overlap and the density of $Z$ is at most $\frac{C}{N\delta}$. Hence, for $0<r<\delta$,
$$\mathbb P(Z\in I)\le Cr\delta^{s-1}=Cr^s\left(\frac{r}{\delta}\right)^{1-s}\le Cr^s,$$
while intervals longer than one are controlled by the total mass. The support meets $O(N)$ dyadic intervals at level $n$, and each of them has mass at most $C\delta^s$, so
$$sn-\log C\le H_n(Z)\le\log O(N)=sn+O(1).$$
A fixed bi-Lipschitz change of variables preserves the interval bound and changes dyadic entropies by $O(1)$, which proves the corresponding assertions.

For the bounds on images, recall that aligned progressions have the same step. Their sums and differences therefore have $O(N)$ centers, perturbed by $O(\delta)$. For a fixed Lipschitz function of two inputs, each pair of noise intervals is mapped into an interval of length $O(\delta)$, so the image meets $O(N^2)$ output intervals at level $n$, and this remains true after fixed bi-Lipschitz changes of the inputs. Since the support is bounded, the image also meets only $O(\delta^{-1})$ intervals, and taking the logarithm of the smaller count, we obtain the coefficient $\min\{2s,1\}$.

Joint Frostman bounds pass to mixtures by averaging. For the conditional bounds, we use the disintegration
$$\mathcal L(X\mid Y=y)=\sum_\ell\mathbb P(L=\ell\mid Y=y)\mathcal L(X\mid Y=y,L=\ell),$$
valid for almost every $y$. For each component, the exceptional conditioning values have measure zero under the marginal of $Y$ in that component, and therefore they have posterior weight zero almost surely. Since there are finitely many components, on a common set of full measure, every component with positive posterior weight satisfies the common Frostman bound, and so does their convex combination. Interchanging $X$ and $Y$, we obtain the other conditional bound. Finally, \eqref{eq:finite-mixture-upper} is the chain of inequalities $H(D_nW\mid L)\le H_n(W)\le H(D_nW,L)$.
\end{proof}

\subsection{Joint Frostman constructions}
\label{sec:joint-constructions}

The examples in this subsection concern the quantity $M_{\mathrm j,n}$ of Theorem~\ref{thm:joint-main}. They attain its coefficient $\frac{1+3s}{4}$ for $\frac{3}{4}\le s<1$ and limit the possible gains at smaller exponents. In these examples, the marginal term is balanced against the two output entropies. We first obtain the sharp coefficient for large $s$ in every degree, and then a stronger upper bound for a quadratic polynomial when $s\le\frac{1}{2}$.

\begin{proposition}\label{prop:joint-high-mixture-upper}\label{prop:joint-quadratic-mixture-upper}
The following families of laws have a fixed bounded support and are jointly $(s,s,C)$-Frostman with a constant $C$ independent of $n$.
\begin{enumerate}[label=\textup{(\roman*)}]
\item For $0<s<1$ and $d\ge2$, the polynomial $f=x^d+y$ admits such laws with
\begin{equation}\label{eq:joint-high-mixture-upper}
\max\left\{\frac{H_n(X)+H_n(Y)}{2},H_n(X+Y),H_n(f(X,Y))\right\}\le\min\left\{\frac{5s}{4},\frac{1+3s}{4}\right\}n+O_{s,d}(1).
\end{equation}
\item For $0<s\le\frac{1}{2}$, the polynomial $f=x^2-y^2$ admits such laws with
\begin{equation}\label{eq:joint-quadratic-mixture-upper}
\max\left\{\frac{H_n(X)+H_n(Y)}{2},H_n(X+Y),H_n(f(X,Y))\right\}\le\frac{6s}{5}n+O_s(1).
\end{equation}
\end{enumerate}
Both polynomials are admissible. In each part, a law may be chosen for every integer $n\ge1$, and while the law itself may vary with $n$, its bounded support, its Frostman constant and the additive entropy error can be chosen independently of $n$ for fixed parameters.
\end{proposition}

\begin{proof}
In each part, we list upper bounds for the coefficients in the order of the three terms in the maximum, and every estimate holds with a bounded error uniform in $n$. Let $t=\min\{2s,1\}$. For (i), we construct the following three laws and then balance their coefficients.
$$\begin{array}{c|ccc}
&\frac{H_n(X)+H_n(Y)}{2}&H_n(X+Y)&H_n(X^d+Y)\\\hline
\mathrm A&s&s&t\\
\mathrm B&s&t&s\\
\mathrm C&\frac{s+t}{2}&s&s
\end{array}$$
All intervals and inverse branches in the three constructions are fixed before $n$ is chosen, while the lengths of the progressions and the scales of the noise may depend on $n$. For $\mathrm A$, we take independent aligned smoothed progressions as $X$ and $Y$, and Lemma~\ref{lem:progressions-mixtures} provides the first row. For $\mathrm B$, we take independent aligned smoothed progressions as $X^d$ and $Y$, with the support of $X^d$ in a fixed positive interval and with a fixed margin for the noise. The positive root is bi-Lipschitz there, so the marginal coefficient is $s$. The aligned sum $X^d+Y$ has coefficient $s$, and the other output is a Lipschitz image of two inputs, so its coefficient is at most $t$.

For $\mathrm C$, we choose a fixed interval on which the derivative of $x^d-x$ is bounded away from zero. We take independent aligned smoothed progressions $T$ and $S$ on fixed intervals of the same length, choosing the interval and the margin of the noise for $T$ so that the whole support of $T$ lies in the image of this branch. We then invert $T=X^d-X$ on the branch and set $Y=S-X$, so that $X^d+Y=T+S$. The inverse branch is bi-Lipschitz, so the same lemma shows that $H_n(X)=sn+O(1)$. The variable $Y$ is a fixed Lipschitz function of $T$ and $S$, while the outputs $X+Y=S$ and $X^d+Y=T+S$ are a single progression and an aligned sum. Hence, $H_n(Y)\le tn+O(1)$, and both output entropies are at most $sn+O(1)$.

The first two laws are products of $s$-Frostman marginals with constants independent of $n$. For the third, $X$ and $S$ are independent and $s$-Frostman with such constants, so that
$$\mathbb P(X\in I,Y\in J)=\int_I\mathbb P(S\in x+J)d\mathcal L(X)(x)\le C|I|^s|J|^s.$$
Thus, all three laws satisfy a common joint Frostman bound. Their mixture with weights $\frac{1}{4}$, $\frac{1}{4}$ and $\frac{1}{2}$ has each of the three coefficients bounded by
$$\frac{3s+t}{4}=\min\left\{\frac{5s}{4},\frac{1+3s}{4}\right\},$$
and applying \eqref{eq:finite-mixture-upper} to both outputs and separately to the two marginals, we obtain (i), including the junction $s=\frac{1}{2}$.

For (ii), independent aligned progressions as $X$ and $Y$ have the coefficients $(s,s,2s)$, since $2s\le1$. For the second component, we take independent aligned progressions as $X^2$ and $Y^2$, on fixed intervals bounded away from zero, and take the positive square roots. The square root is bi-Lipschitz there, so the marginal coefficient is $s$, while the difference of the squares is an aligned difference and the sum is a Lipschitz image of the inputs, so the coefficients are $(s,2s,s)$. Both laws satisfy joint Frostman bounds with constants independent of $n$. For a third law, we take a fixed $2s$-Frostman variable $Z$ with $H_n(Z)=2sn+O(1)$ from Lemma~\ref{joint-models:scalar-source} and set $(X,Y)=(Z,-Z)$. Both outputs vanish, and
$$\mathbb P(Z\in I,-Z\in J)\le C\min\{|I|,|J|\}^{2s}\le C|I|^s|J|^s,$$
so the vector of coefficients is $(2s,0,0)$. At $s=\frac{1}{2}$, the scalar lemma uses the uniform law on an interval, so this component also covers the endpoint. The mixture with weights $\frac{2}{5}$, $\frac{2}{5}$ and $\frac{1}{5}$ has the coefficients
$$\frac{2}{5}(s,s,2s)+\frac{2}{5}(s,2s,s)+\frac{1}{5}(2s,0,0)=\left(\frac{6s}{5},\frac{6s}{5},\frac{6s}{5}\right),$$
and the mixture lemma proves (ii), with a common bounded support and a single joint Frostman constant. Each component lies in a fixed bounded rectangle, so we can enclose the finitely many supports in one such rectangle and use a uniform product law inside it for the finitely many scales not covered by the constructions with progressions. Enlarging the common Frostman constant and the bounded entropy errors, we obtain the stated families for every $n\ge1$.

Finally, every polynomial $h(ax+by+c)$ is annihilated by a nonzero constant vector field $u\partial_x+v\partial_y$, since we may take $(u,v)=(b,-a)$ unless $a=b=0$, and in that case every such field annihilates it. For our two polynomials, $u\partial_xf+v\partial_yf$ equals $udx^{d-1}+v$ and $2ux-2vy$, respectively, and each of these vanishes identically only when $u=v=0$. Hence, both polynomials are admissible.
\end{proof}

The first bound matches Theorem~\ref{thm:joint-main} when $s\ge\frac{3}{4}$, and the second provides the upper coefficient conjectured in Section~\ref{sec:conjectures-support} for quadratic polynomials at small exponents. The next example shows that the possible gain in degree $d$ is at most a constant multiple of $\frac{s}{d}$.

\begin{proposition}\label{prop:degree-upper-latest}
For $0<s\le\frac{1}{2}$ and every integer $d\ge2$, there are a fixed nonspecial polynomial $f_d$ of degree $d$ and constants $C_{s,d},K_{s,d}<\infty$ such that, for every integer $n\ge1$, one can choose a jointly $(s,s,C_{s,d})$-Frostman law $(X,Y)$ in $[0,1]^2$ with
\begin{equation}\label{eq:degree-upper-all-degrees}
\max\left\{\frac{H_n(X)+H_n(Y)}{2},H_n(X+Y),H_n(f_d(X,Y))\right\}\le\left(s+\frac{s}{2d-1}\right)n+K_{s,d}.
\end{equation}
Only the law varies with $n$, while the polynomial and both constants remain fixed.
\end{proposition}

\begin{proof}
Let $r=d-1$ and $f_d(x,y)=(1+y)(x+y-1)^r$, a polynomial of degree $d$ vanishing to order $r$ on the line $x+y=1$. We place one component of the law near this line, with a width chosen so that the polynomial image lies in a single dyadic interval at level $n$. The price is a larger marginal entropy, which we balance against the polynomial entropy by mixing with a product law built from progressions. To see that $f_d$ is nonspecial, note that
$$\frac{(f_d)_x}{(f_d)_y}=\frac{r(1+y)}{x+(r+1)y+r-1}\text{ and }\partial_x\partial_y\ln\left|\frac{(f_d)_x}{(f_d)_y}\right|=\frac{r+1}{(x+(r+1)y+r-1)^2}\not\equiv0,$$
whereas by \eqref{eq:separated-derivative-ratio}, this mixed derivative would vanish on every regular rectangle if $f_d$ were special.

We choose a fixed $2s$-Frostman variable $T$ on $[\frac{1}{4},\frac{3}{8}]$ from Lemma~\ref{joint-models:scalar-source}, which includes the case $2s=1$. For integers $m\ge2$, let $h=2^{-m}$. The atoms of ${D_mT}$ are $h$-separated and have mass at most $Ch^{2s}$, and pulling an interval back under the quantization enlarges its length by at most $h$, so that
\begin{equation}\label{eq:degree-upper-coarse}
\mathbb P({D_mT}\in I)\le C(|I|+h)^{2s}\text{ and }H({D_mT})=H_m(T)=2sm+O(1).
\end{equation}
The variable $T$ is fixed once $s$ is fixed, and only its level of quantization varies with $n$. This coarse law spreads the mass along the line on which $f_d$ vanishes. For the finer scales, let $U_q$ and $V_q$ be independent aligned smoothed progressions of exponent $s$ at level $q$, rescaled into $[0,\frac{1}{32}]$. Lemma~\ref{lem:progressions-mixtures} provides Frostman bounds with constants independent of $q$ and
\begin{equation}\label{eq:degree-upper-fine}
\max\{H_q(U_q),H_q(V_q),H_q(U_q+V_q),H_q(V_q-U_q)\}\le sq+O(1).
\end{equation}
For $q=0$, we use independent copies of a fixed uniform law on a smaller interval, with the same quantization formula $D_0z=\lfloor z\rfloor$, and the displayed entropy bound and the Frostman bound remain valid after enlarging the constant once. In the rest of the proof, vectors of coefficients list the average marginal entropy, the entropy of the sum and the entropy of the polynomial, in this order. The product law $(X,Y)=(U_n,V_n)$ is jointly $s$-Frostman, its aligned sum has coefficient $s$, and the fixed Lipschitz polynomial image has coefficient at most $2s$, so this law has the vector $(s,s,2s)$.

For the component near the line, we take ${D_mT}$, $U_q$ and $V_q$ independent and set
\begin{equation}\label{eq:degree-upper-diagonal}
\begin{aligned}
m&=\left\lceil\frac{n}{r}\right\rceil+2, & q&=\max\{n-m,0\},\\
Y&={D_mT}+hU_q, & X&=1-{D_mT}+h(V_q-U_q).
\end{aligned}
\end{equation}
Since $m\ge3$, both endpoints of $[\frac{1}{4},\frac{3}{8}]$ belong to the dyadic grid at level $m$, so the quantization keeps ${D_mT}$ in $[\frac{1}{4},\frac{3}{8}]$, and both coordinates lie in $[0,1]$. Their sum is $1+hV_q$, and the choice of $m$ ensures $rm\ge n+2r$, which will force the polynomial output into a single dyadic interval at level $n$.

To check the rectangle bound at every scale, let $I$ and $J$ be intervals of lengths $a,b\le1$. The threshold $h$ separates the intervals that resolve individual coarse atoms from those that contain several of them. If $a,b\le h$, the condition on $Y$ confines ${D_mT}$ to an interval of length $b+O(h)$, which contains only boundedly many atoms. For each such atom, we condition further on $U_q$. The condition on $X$ then confines the independent variable $V_q$ to an interval of length $\frac{a}{h}$, while the condition on $Y$ confines $U_q$ to an interval of length $\frac{b}{h}$, and the two Frostman bounds for the fine variables imply
$$\mathbb P(X\in I,Y\in J\mid {D_mT}=t)\le C\left(\frac{a}{h}\right)^s\mathbb P\left(U_q\in\frac{J-t}{h}\right)\le Ch^{-2s}a^sb^s.$$
Multiplying by the masses $O(h^{2s})$ of the boundedly many contributing atoms and summing, we obtain the required bound. If $a\le h\le b$, the condition on $X$ confines ${D_mT}$ to an interval of length $a+O(h)$, and conditionally on each contributing atom and on $U_q$, the interval bound for $V_q$ contributes $C(\frac{a}{h})^s$, while the condition on $Y$ is discarded. If $b\le h\le a$, we use the condition on $Y$ and the interval bound for $U_q$ instead. When $a,b\ge h$, the coarse variable lies in an interval of length at most $\min\{a,b\}+O(h)$, and we use \eqref{eq:degree-upper-coarse}. The four bounds are
$$\mathbb P(X\in I,Y\in J)\lesssim
\begin{cases}
a^sb^s, & a,b\le h,\\
h^{2s}\left(\frac{a}{h}\right)^s, & a\le h\le b,\\
h^{2s}\left(\frac{b}{h}\right)^s, & b\le h\le a,\\
\min\{a,b\}^{2s}, & a,b\ge h.
\end{cases}$$
In the mixed cases, we use $h\le b$ or $h\le a$, respectively, and in the last case, $\min\{a,b\}^{2s}\le a^sb^s$, so every entry is at most a constant times $a^sb^s$. These cases include the equalities at the thresholds, and they prove a rectangle bound uniform in $m$, $q$ and $n$.

For the entropy bounds, we first condition on ${D_mT}$. If $m\le n$, scaling by $h$ turns the resolution at level $n$ into resolution at level $q=n-m$, and the translation changes the entropy by at most $O(1)$, so the chain rule and \eqref{eq:degree-upper-fine} lead to
$$\begin{aligned}
\frac{H_n(X)+H_n(Y)}{2}&\le2sm+sq+O(1)=\left(s+\frac{s}{r}\right)n+O_{s,d}(1),\\
H_n(X+Y)&\le sq+O(1)=\left(s-\frac{s}{r}\right)n+O_{s,d}(1).
\end{aligned}$$
If $m>n$, the fine perturbations fit into boundedly many cells at level $n$ and $0<m-n\le3$, and since $m=\frac{n}{r}+O(1)$, this case forces $(1-\frac{1}{r})n=O(1)$. The average marginal entropy is then at most $2sm+O(1)$ and the entropy of the sum is $O(1)$, so the same coefficients hold with bounded errors, and this also covers the case $r=1$, in which $m=n+2$ for every $n$. Moreover,
$$0\le f_d(X,Y)=(1+Y)(hV_q)^r\le2\left(\frac{h}{32}\right)^r<2^{-n},$$
so the polynomial entropy is zero, and this component has the vector of coefficients $(s+\frac{s}{r},s-\frac{s}{r},0)$.

We give this component the weight $p$ and the product component the weight $1-p$. The coefficients of the marginal term and of the polynomial are then $s+\frac{ps}{r}$ and $2s(1-p)$, which agree for $p=\frac{r}{2r+1}\in(0,1)$, when both equal $s+\frac{s}{2r+1}$, while the coefficient of the sum is the smaller number $s-\frac{ps}{r}$. The weight depends only on $d$. Applying the mixture lemma to each marginal and each output preserves the common joint Frostman bound and increases each average entropy by at most one bit. Since $2r+1=2d-1$, this proves \eqref{eq:degree-upper-all-degrees}, with all the stated properties of the support and of uniformity.
\end{proof}

\subsection{Conditional Frostman constructions}
\label{sec:conditional-constructions}

The modular construction provides matching laws for every admissible $f\in\mathfrak E_S$ in the range $\frac{2}{3}\le s<1$, as required in Theorem~\ref{thm:conditional-classification}. Taking $f(x,y)=x+(x+y)^2$, it shows that the coefficient in Theorem~\ref{thm:conditional-main} cannot be increased in a bound valid for every admissible polynomial in that range. The construction with level sets supplies the fixed matching laws outside $\mathfrak E_S$ for $\frac{3}{4}\le s<1$, and for $0<s\le\frac{1}{2}$, the final mixtures lead to upper bounds for the possible entropy gains.

\subsubsection{Modular laws}
\label{upper-modular:subsection}

For admissible $f\in\mathfrak E_S$ and $\frac{1}{4}\le s<1$, the following construction produces the coefficient $\frac{1+2s}{3}$ for both output entropies. We first prove the discrepancy estimate on intervals of consecutive indices that the construction requires.

\begin{lemma}\label{upper-modular:discrepancy-lemma}
Let $P$ be a real polynomial such that $P''$ has constant sign and is bounded above and away from zero on a neighborhood of a compact interval $U$. Write $\{t\}=t-\lfloor t\rfloor$. For every sufficiently large integer $q$, every interval $J$ of consecutive integers with $q^{-1}J\subset U$, and every interval $W\subset\mathbb R/\mathbb Z$, the number $N_W(J)$ of indices $j\in J$ with $\{qP(\frac{j}{q})\}\in W$ satisfies
\begin{equation}\label{upper-modular:discrepancy}
\left|N_W(J)-|W||J|\right|\le K_P(|J|q^{-\frac{1}{3}}+q^{\frac{1}{2}}).
\end{equation}
\end{lemma}

\begin{proof}
The Erd\H{o}s--Tur\'an inequality reduces the count to exponential sums with phases $vqP(\frac{t}{q})$. We first prove an estimate for such sums based on the second derivative, which retains the length of the interval of indices, and then the choice of the Fourier cutoff produces the two error terms.

Let $\varphi$ be a twice differentiable phase on an interval containing $J$ together with a buffer of length one on each side, and suppose that $\lambda\le|\varphi''|\le K\lambda$ with constant sign, where $0<\lambda\le\frac{1}{16}$. We claim that
$$\left|\sum_{j\in J}e^{2\pi i\varphi(j)}\right|\lesssim_K|J|\sqrt\lambda+\lambda^{-\frac{1}{2}}.$$
The sign condition makes $\varphi'$ monotone, and its range crosses $O_K(|J|\lambda+1)$ integers. For $0<\alpha\le\frac{1}{4}$, we split the domain into the pieces where $\varphi'$ lies within $\alpha$ of an integer and the complementary pieces. Each piece of the first kind has length $O(\frac{\alpha}{\lambda})$, so the integer points in these pieces contribute $O_K(|J|\alpha+\frac{\alpha}{\lambda}+|J|\lambda+1)$. On a piece of the second kind, we discard at most two integer points at its ends, so that each remaining interval $[j,j+1]$ lies in the piece. The increments $d_j=\varphi(j+1)-\varphi(j)$ are then monotone and lie in $[m+\alpha,m+1-\alpha]$ for a single integer $m$, and summation by parts in
$$e^{2\pi i\varphi(j)}=\frac{e^{2\pi i\varphi(j)}-e^{2\pi i\varphi(j+1)}}{1-e^{2\pi id_j}}$$
costs $O(\frac{1}{\alpha})$ on the piece, because the function $(1-e^{2\pi it})^{-1}=\frac{1}{2}+\frac{i}{2}\cot(\pi t)$ has supremum and total variation $O(\frac{1}{\alpha})$ on the range of the increments. Including the discarded endpoints, the pieces of the second kind contribute $O_K(\frac{|J|\lambda+1}{\alpha}+|J|\lambda+1)$. We take $\alpha=\sqrt\lambda$, which is at most $\frac{1}{4}$ by the assumption on $\lambda$, and then both contributions are $O_K(|J|\sqrt\lambda+\lambda^{-\frac{1}{2}})$, which proves the claim.

Let $M=\lfloor q^{\frac{1}{3}}\rfloor$, a choice that balances the error $\frac{|J|}{M}$ coming from the cutoff against the contribution $|J|\sqrt{\frac{M}{q}}$ of the estimate above. We apply the claim to $\varphi(t)=vqP(\frac{t}{q})$ for $1\le v\le M$, whose second derivative is comparable to $\frac{v}{q}$, so that the condition $\lambda\le\frac{1}{16}$ holds for all sufficiently large $q$, while the fixed neighborhood of $U$ provides the buffers at the ends. Consequently,
$$\left|\sum_{j\in J}e^{2\pi ivqP(\frac{j}{q})}\right|\le K_P\left(|J|\sqrt{\frac{v}{q}}+\sqrt{\frac{q}{v}}\right).$$
The Erd\H{o}s--Tur\'an inequality \cite[Part~I, Theorem~III, p.~374]{ErdosTuran}, applied to the fractional parts $\{qP(\frac{j}{q})\}$ with cutoff $M$, now bounds the discrepancy by
$$K_P\left(\frac{|J|}{M+1}+\frac{|J|}{\sqrt q}\sum_{v=1}^Mv^{-\frac{1}{2}}+\sqrt q\sum_{v=1}^Mv^{-\frac{3}{2}}\right)\le K_P'(|J|q^{-\frac{1}{3}}+q^{\frac{1}{2}}),$$
which is the required estimate, uniformly over the interval of indices and the target interval $W$.
\end{proof}

\begin{proposition}\label{upper-modular:prop}
Let $f(x,y)=h(ax+P(x+y))$ with real coefficients, where $a\ne0$, $\deg P\ge2$ and $h$ is nonconstant. Fix a nonempty bounded open rectangle $B$ and $\frac{1}{4}\le s<1$. There are constants $A,C>0$, independent of $n$, such that, for every sufficiently large integer $n$, there is a conditionally $(s,C;s,C)$-Frostman law supported in $B$ with
\begin{equation}\label{upper-modular:entropies}
\begin{aligned}
H_n(X+Y)&=\frac{1+2s}{3}n+O(1),\\
H_n(f(X,Y))&=\frac{1+2s}{3}n+O(1).
\end{aligned}
\end{equation}
With $\delta=2^{-n}$, each such law lies within coordinate distance $\frac{\delta}{100}$ of a grid $q^{-1}\mathbb Z^2$ for a prime $q$, where
$$\begin{aligned}
A\delta^{-\frac{2+s}{3}}&\le q\le2A\delta^{-\frac{2+s}{3}},\\
H_n(X)=H_n(Y)&=\log q+O(1),\\
H_n(X,Y)&=2\log(q^2\delta)+O(1).
\end{aligned}$$
Any fixed finite set of primes may be excluded. All constants are independent of $n$, and $f$ is admissible.
\end{proposition}

\begin{proof}
We regard a point $(\frac{k}{q},\frac{l}{q})$ of the grid as an edge joining the row index $k$ to the column index $l$. We first select candidate sums and then a second set of labels that constrains the other output. The exact conditional bounds will follow from two properties of the retained edges, namely that all rows and all columns have comparable degrees, and that the number of neighbors in every interval of consecutive indices is uniformly bounded. Once these properties are established, we smooth the law on the grid and compute its entropies.

Let $G=x+\frac{P(x+y)}{a}$, so that $f=h(aG)$. We choose a compact rectangle $I_0\times J_0$ inside $B$ and a slightly larger rectangle, still inside $B$, on which the four quantities $|\frac{P'(x+y)}{a}|$, $|\frac{P''(x+y)}{a}|$, $|1+\frac{P'(x+y)}{a}|$ and $|h'(aG(x,y))|$ are bounded away from zero. The first three conditions exclude finitely many values of the sum, and the last one excludes finitely many level sets of $G$, which cannot cover an open set because $G_x-G_y=1$. We shrink the rectangles so that $|I_0+J_0|<1$ and so that $t\mapsto h(at)$ is bi-Lipschitz on the range of $G$ over the larger rectangle, which will contain all later perturbations.

Let $\delta=2^{-n}$ and fix $A\ge1$, whose size will be determined by the estimate at the endpoint $s=\frac{1}{4}$ below. By Bertrand's postulate, there is a prime $q$ with $A\delta^{-\frac{2+s}{3}}\le q\le2A\delta^{-\frac{2+s}{3}}$ for all sufficiently large $n$, and such primes eventually avoid any fixed finite set. We define
\begin{equation}\label{upper-modular:parameters}
\begin{aligned}
p&=q\delta, & N&=q^2\delta,\\
\theta&=\frac{1+2s}{2+s}, & b&=\frac{1+s}{2+s}.
\end{aligned}
\end{equation}
Then $p\to0$, $pN\asymp_A\delta^{-s}$ and $N=\kappa_qq^\theta$, where $A^{\frac{3}{2+s}}\le\kappa_q\le(2A)^{\frac{3}{2+s}}$. These parameters play different roles. The number $N$ bounds the number of centers of the outputs, while $pN$ will be the number of neighbors in each row and each column, and the relation $pN\asymp_A\delta^{-s}$ then assigns to each neighbor conditional mass $O(\delta^s)$.

We select the candidate sum indices
$$\mathcal B_q=\left\{i\in q(I_0+J_0)\cap\mathbb Z:\left\{\frac{qP(\frac{i}{q})}{a}\right\}\in[0,p)\right\},\text{ and we write }j_i=\left\lfloor\frac{qP(\frac{i}{q})}{a}\right\rfloor.$$
We claim that for every nonempty interval $T$ of consecutive integers in this range of sums,
\begin{equation}\label{upper-modular:residue-bound}
|\mathcal B_q\cap T|\le KN\left(\frac{|T|}{q}\right)^s.
\end{equation}
This bound is needed for short intervals as well as for full rows, because the conditional Frostman estimate must hold at every scale. The threshold $q^b$ is chosen so that $q^{b(1-s)}=q^{\theta-s}$, and therefore the trivial count $|T|$ proves the claim when $|T|\le q^b$. For $q^b\le|T|\le q$, we apply Lemma~\ref{upper-modular:discrepancy-lemma} to $\frac{P}{a}$ with $W=[0,p)$. Its main term is $p|T|\le N(\frac{|T|}{q})^s$, and its errors satisfy
$$|T|q^{-\frac{1}{3}}\le q^{\theta-s}|T|^s\text{ and }q^{\frac{1}{2}}\le q^b\le q^{\theta-s}|T|^s.$$
For the first inequality, we divide by $|T|^s$ and use $|T|\le q$, and the resulting comparison is exactly $\theta\ge\frac{2}{3}$, which is equivalent to $s\ge\frac{1}{4}$. For the second, we use $b>\frac{1}{2}$, the identity $\theta-s+bs=b$ and $|T|\ge q^b$. Since $N\ge q^\theta$, the claim follows. For large $q$, the range of sums contains fewer than $q$ integers, so these cases cover every $T$, and they also show that $|\mathcal B_q|\le KN$.

To normalize the exact conditional laws, we also need a lower bound for the counts. A full row have $|J_0|q+O(1)$ possible sums and a full column has $|I_0|q+O(1)$, and the main terms in the discrepancy estimate are fixed positive multiples of $N$. For $s>\frac{1}{4}$, they dominate the errors, because $Nq^{-\frac{2}{3}}\to\infty$. At $s=\frac{1}{4}$, we choose $A$ once and for all so large that the lower bound for $\kappa_q$ dominates the discrepancy constants divided by $\min\{|I_0|,|J_0|\}$, and we fix this choice before any scale $n$ is chosen. Thus, every full row and every full column contains between $K^{-1}N$ and $KN$ candidate edges, also at the endpoint.

All labels $k+j_{k+l}$ arising from the rectangle lie in a fixed interval of $O(q)$ integers. We select each label independently with probability $p$, call the set of selected labels $\mathcal C_q$, and retain the edges
$$\mathcal E_q=\{(k,l)\in(qI_0\cap\mathbb Z)\times(qJ_0\cap\mathbb Z):k+l\in\mathcal B_q\text{ and }k+j_{k+l}\in\mathcal C_q\}.$$
A fixed row uses each label at most $M_0$ times, because the monotone function $\frac{qP(\frac{i}{q})}{a}$ has derivative bounded away from zero, so that a unit interval has only boundedly many integer preimages under it. In a column, we apply the same argument to $i+j_i-l=\lfloor i+\frac{qP(\frac{i}{q})}{a}\rfloor-l$, whose underlying function has derivative $1+\frac{P'}{a}$, and we choose a single integer $M_0$ that works in both directions.

Because labels repeat, the edges at a given vertex are not selected independently, and we handle this by grouping the edges with the same label and using their multiplicities as bounded weights. Let $Z=\sum_cw_c\xi_c$, where the $\xi_c$ are independent Bernoulli variables with parameter $p$ and $0\le w_c\le M_0$. By convexity, for every real $t$,
$$\mathbb Ee^{tZ}\le\exp\left(\frac{p\sum_cw_c}{M_0}(e^{tM_0}-1)\right),$$
and since $\mathbb EZ=p\sum_cw_c$, the exponential Markov inequality with positive and negative $t$ yields
$$\mathbb P\left(|Z-\mathbb EZ|\ge\frac{\mathbb EZ}{2}\right)\le2e^{-\frac{\mathbb EZ}{12M_0}}\text{ and }\mathbb P(Z>K_{M_0}(\mathbb EZ+u))\le e^{-u}\text{ for }u\ge1.$$
The expected full degrees are comparable to $pN$, a positive power of $q$, and for the counts in intervals, we take $u$ equal to a sufficiently large fixed multiple of $\log q$. There are $O(q)$ full degrees and $O(q^3)$ choices of a row or column and an interval, so a union bound shows that with probability tending to one
\begin{equation}\label{upper-modular:neighbor-bound}
d_k^{\mathrm r},d_l^{\mathrm c}\asymp pN\text{ and }|\{\text{neighbors with sum indices in }T\}|\le KpN\left(\frac{|T|}{q}\right)^s,
\end{equation}
where $d_k^{\mathrm r}$ and $d_l^{\mathrm c}$ are the degrees of the row $k$ and of the column $l$. Indeed, for a fixed row or column and an interval $T$ inside its range of sums, the expected number of neighbors with sum index in $T$ is $p|\mathcal B_q\cap T|$, since each candidate edge is retained with probability $p$, even when labels repeat. By \eqref{upper-modular:residue-bound}, this expectation is at most $KpN(\frac{|T|}{q})^s$, and the additive term $u$ in the tail bound is absorbed because
$$pNq^{-s}\asymp_A\delta^{-\frac{s(1-s)}{3}}\gg\log q.$$
The right-hand side of the desired bound for the neighbors is smallest when $|T|=1$, so this growth absorbs the additive term simultaneously for all nonempty intervals, for each fixed $s<1$. A further binomial tail bound shows that $|\mathcal C_q|\le KN$. For all sufficiently large $n$, all these events occur simultaneously with positive probability, and we fix a selection for which they hold. The number of rows is comparable to $q$, with constants depending on $I_0$, so summing the degrees, we find that $|\mathcal E_q|\asymp qpN=N^2$.

This completes the combinatorial construction. To obtain a law with exact conditional bounds below the spacing of the grid, we choose an edge uniformly at random and, independently, random variables $U$ and $V$ uniform on $[0,\frac{\delta}{100}]$, and we set $X=\frac{k}{q}+U$ and $Y=\frac{l}{q}+V$. For large $n$, the support lies in the buffered rectangle, and distinct noise intervals in each coordinate are disjoint because $q\delta=p\to0$. For almost every value $Y=y$, the column $l$ is determined, and since the edges have equal weights and the noise added to $Y$ has a common density, the conditional row index is uniform among the $d_l^{\mathrm c}$ neighbors of that column. The exact conditional density is therefore
$$\frac{d\mathcal L(X\mid Y=y)}{dx}(x)=\frac{100}{\delta d_l^{\mathrm c}}\sum_{k:(k,l)\in\mathcal E_q}1_{[\frac{k}{q},\frac{k}{q}+\frac{\delta}{100}]}(x),$$
because the independent noise added to $X$ is not affected by this conditioning. Interchanging rows and columns, we obtain the formula for $Y\mid X$.

We check the conditional Frostman bounds in three ranges of lengths, separated by the noise scale $\delta$ and the spacing $\frac{1}{q}$ of the grid, which are ordered because $q\delta\to0$. For a test interval of length $\frac{1}{q}\le r\le1$, at most $qr+O(1)$ consecutive centers contribute, and dividing \eqref{upper-modular:neighbor-bound} by the lower bound for the degree, we obtain conditional mass $O(r^s)$. For the remaining scales,
$$\begin{cases}
O((pN)^{-1})=O(\delta^s)\le O(r^s) & \text{if }\delta\le r<\frac{1}{q},\\
O\left(\frac{r}{\delta pN}\right)=O(r\delta^{s-1})\le O(r^s) & \text{if }0<r<\delta,
\end{cases}$$
where the first line counts boundedly many noise intervals and the second uses the conditional density. Longer intervals are controlled by the total mass one. Thus, both exact conditional bounds hold at every scale with constants independent of $n$.

It remains to compute the entropies. The probability of a row is its degree divided by $|\mathcal E_q|$, and the same holds for a column, so the counts of degrees and edges make the probabilities of all rows and columns comparable to $\frac{1}{q}$. A noise interval meets boundedly many cells at level $n$, while for large $n$, a cell determines its row or column, and since the edge is uniform,
$$H_n(X)=H_n(Y)=\log q+O(1)\text{ and }H_n(X,Y)=\log|\mathcal E_q|+O(1)=2\log N+O(1).$$
For an edge with $i=k+l$, the sum is $\frac{i}{q}+O(\delta)$ with $i\in\mathcal B_q$, and
$$G\left(\frac{k}{q},\frac{l}{q}\right)=\frac{k+j_i}{q}+\frac{1}{q}\left\{\frac{qP(\frac{i}{q})}{a}\right\}=\frac{c}{q}+O(\delta)\text{ with }c\in\mathcal C_q,$$
where the error is at most $\frac{p}{q}=\delta$. The noise changes $G$ by another $O(\delta)$, so both output entropies are at most $\log N+O(1)$. The map $(x,y)\mapsto(S,G)$, with $S=x+y$, has the Lipschitz inverse $x=G-\frac{P(S)}{a}$, $y=S-x$ on the fixed ranges, so the joint entropy of $(S,G)$ is $2\log N+O(1)$. By subadditivity, this joint entropy is at most the sum of the two scalar entropies, so subtracting either scalar upper bound yields the matching lower bound for the other. The bi-Lipschitz map $t\mapsto h(at)$ preserves entropy up to $O(1)$, and since $\log N=\frac{(1+2s)n}{3}+O(1)$, all the claimed identities follow.

Finally, if a constant vector field $u\partial_x+v\partial_y$ annihilates $f$, then $h'(ax+P(x+y))(au+(u+v)P'(x+y))\equiv0$. The first factor is a nonzero polynomial, so $au+(u+v)P'\equiv0$, and since $P'$ is nonconstant, this forces $u+v=0$ and then $u=0$, because $a\ne0$. As we observed in the proof of Proposition~\ref{prop:joint-high-mixture-upper}, a polynomial in one affine coordinate is annihilated by a nonzero constant vector field, so $f$ is admissible. All constants, as well as the choice of $A$ at the endpoint, were fixed before $n$.
\end{proof}

For applications requiring every $n\ge1$, we fill the finitely many remaining coarse scales with a fixed uniform product law on a smaller rectangle in the box. At each such scale, we choose a prime $q\ge50\delta^{-1}$ outside the prescribed finite set, so that every point lies within coordinate distance $\frac{1}{2q}\le\frac{\delta}{100}$ of the grid $q^{-1}\mathbb Z^2$. Since only finitely many scales are added, a single enlargement of the constants covers their entropy errors and Frostman bounds and preserves the relation $q\asymp\delta^{-\frac{2+s}{3}}$.

\subsubsection{Laws built from level sets}
\label{sec:conditional-sharpness}

We begin with a fixed law that has a prescribed entropy for one output and full entropy for every output transverse to it. The building blocks are the scalar measures of Lemma~\ref{joint-models:scalar-source}, whose supports for $s<1$ are null sets meeting $O(2^{sn})$ dyadic intervals at level $n$.

\begin{lemma}\label{lem:cantor-level-source}
Let $0<s<1$, let $I$ and $J$ be compact nondegenerate intervals, and suppose that $f$ is continuously differentiable near $I\times J$ and that $f_xf_y$ does not vanish there. There is a fixed law on $I\times J$ with exact conditional $s$-Frostman bounds in both directions such that
\begin{equation}\label{eq:cantor-level-entropies}
\begin{aligned}
H_n(f(X,Y))&=sn+O(1),\\
H_n(X)=H_n(Y)&=n+O(1),\\
H_n(X,Y)&=(1+s)n+O(1).
\end{aligned}
\end{equation}
Its marginal densities are bounded above and below on $I$ and $J$, respectively, and its exact conditional laws are almost surely singular. If $g$ is continuously differentiable and $|f_xg_y-f_yg_x|$ is bounded away from zero on the rectangle, then $g(X,Y)$ has a bounded density and $H_n(g(X,Y))=n+O(1)$.
\end{lemma}

\begin{proof}
The idea is to choose a measure on the values of $f$ and to spread its mass along the corresponding level sets. The weight along each level set is chosen so that the law has symmetric descriptions by rows and by columns, and this symmetry allows us to check both exact conditional laws in the same way.

We fix constants $0<\lambda\le\Lambda$ with $\lambda\le|f_x|\le\Lambda$ and $\lambda\le|f_y|\le\Lambda$ on $I\times J$. Each partial derivative has a constant sign there, so the values of $f$ along any row or column of the rectangle form an interval of length at least $\ell=\lambda\min\{|I|,|J|\}$. We choose finitely many nondegenerate closed intervals $K_1,\dots,K_m$ in the interior of $f(I\times J)$ such that every subinterval of $f(I\times J)$ of length $\ell$ contains one of them. Such a choice is possible by compactness. Indeed, the subintervals of length $\ell$ are parametrized by their left endpoints, which form a compact interval, and a closed interval lying strictly inside one subinterval remains inside all sufficiently small translates of it, so finitely many of these neighborhoods of left endpoints cover the parameter interval.

We apply Lemma~\ref{joint-models:scalar-source} with exponent $s$ on each $K_j$ and let $\alpha$ be the equal mixture of the resulting probability measures. Then $\alpha$ is $s$-Frostman, its support is a null set covered by $O(2^{sn})$ dyadic intervals at level $n$, and
$$\alpha(f(x,J))\ge\frac{1}{m}\text{ and }\alpha(f(I,y))\ge\frac{1}{m}\text{ for all }x\in I\text{ and }y\in J.$$
These lower bounds will keep all conditional normalizations away from zero.

Since $f$ is strictly monotone along each row and each column, we have inverse functions in both variables. Let $\phi(x,t)$ be the solution of $f(x,\phi(x,t))=t$ and let $\psi(y,t)$ be the solution of $f(\psi(y,t),y)=t$, whenever $t$ belongs to the corresponding range. We define a measure $\mu_0$ on $I\times J$ by setting, for bounded Borel functions $\Phi$,
$$\begin{aligned}
\int\Phi d\mu_0&=\int\left(\int_{\{x\in I:t\in f(x,J)\}}\frac{\Phi(x,\phi(x,t))}{|f_y(x,\phi(x,t))|}dx\right)d\alpha(t)\\
&=\int\left(\int_{\{y\in J:t\in f(I,y)\}}\frac{\Phi(\psi(y,t),y)}{|f_x(\psi(y,t),y)|}dy\right)d\alpha(t).
\end{aligned}$$
The weights $\frac{1}{|f_y|}$ and $\frac{1}{|f_x|}$ are chosen precisely so that these two descriptions agree. To see this, we fix $t$. The values of $f$ at the two ends of a row are monotone in $x$ in the same direction, so the set of $x$ with $t\in f(x,J)$ is an interval, and the points $(x,\phi(x,t))$ form a connected piece of the level set $\{f=t\}$. Along this piece, $\frac{dy}{dx}=-\frac{f_x}{f_y}$ has a constant nonzero sign and $\frac{dx}{|f_y|}=\frac{|dy|}{|f_x|}$, so the change of variables from $x$ to $y$ turns the first formula into the second. We set $Z=\mu_0(I\times J)$ and $\mu=\frac{\mu_0}{Z}$. The bounds for the slices show that
$$\frac{|I|}{m\Lambda}\le Z\le\frac{|I|}{\lambda},$$
so $\mu$ is a probability law with a fixed normalization.

Next, we check the regularity of the marginal and conditional laws. The unnormalized density of $X$ is
$$w_X(x)=\int_{f(x,J)}\frac{d\alpha(t)}{|f_y(x,\phi(x,t))|},\text{ which satisfies }\frac{1}{m\Lambda}\le w_X(x)\le\frac{1}{\lambda}.$$
The column formula yields the same bounds for the unnormalized density of $Y$, and dividing by $Z$, we obtain the upper and lower bounds for both marginal densities. Given $X=x$, the conditional law of $Y$ is the pullback of the restriction of $\alpha$ to $f(x,J)$ under $y\mapsto f(x,y)$, with weight $\frac{1}{|f_y|}$ and normalized by $w_X(x)$. This map sends an interval of length $r$ into an interval of length at most $\Lambda r$, so the conditional mass of such an interval is at most
$$\frac{C(\Lambda r)^s}{\lambda w_X(x)}\le C'r^s.$$
The column formula, with $|f_x|$ in place of $|f_y|$, proves the other exact conditional bound. Both conditional laws are singular, because the maps along rows and columns are bi-Lipschitz and the support of $\alpha$ is a null set.

We now turn to the entropies. The law of $f(X,Y)$ is absolutely continuous with respect to $\alpha$, with density
$$\frac{1}{Z}\int_{\{x\in I:t\in f(x,J)\}}\frac{dx}{|f_y(x,\phi(x,t))|}\le\frac{|I|}{Z\lambda}.$$
This upper bound transfers the $s$-Frostman estimate from $\alpha$ to the law of $f(X,Y)$, which shows that $H_n(f(X,Y))\ge sn-O(1)$. The law of $f(X,Y)$ is also supported in the support of $\alpha$, so the covering bound yields the reverse inequality. The bounded marginal densities show that $H_n(X)=H_n(Y)=n+O(1)$. For the lower bound on the joint entropy, we observe that conditioning on the exact value of $X$ can only decrease the entropy compared with conditioning on $D_nX$, and then the exact conditional bound on atoms implies that
$$H_n(X,Y)\ge H_n(X)+H(D_nY\mid X)\ge(1+s)n-O(1).$$
Conversely, the map $(x,y)\mapsto(x,f(x,y))$ has a Lipschitz inverse, since
$$|y-y'|\le\frac{|f(x,y)-f(x',y')|+\Lambda|x-x'|}{\lambda}.$$
Consequently, each dyadic square in the coordinates $(x,f)$ meets the images of only boundedly many dyadic squares of the same size, and $H_n(X,Y)\le H_n(X)+H_n(f(X,Y))+O(1)$, which completes the proof of \eqref{eq:cantor-level-entropies}.

It remains to treat a transverse function $g$, and here the connectedness of the pieces of level sets established above is what we need. Let $\kappa=\inf_{I\times J}|f_xg_y-f_yg_x|>0$. Along each such piece, parametrized by $x$,
$$\frac{dg}{dx}=-\frac{f_xg_y-f_yg_x}{f_y}$$
has a constant nonzero sign, so $g$ is injective along the piece, and the pushforward of $\frac{dx}{|f_y|}$ under $g$ has density $\frac{1}{|f_xg_y-f_yg_x|}\le\frac{1}{\kappa}$. Integrating against $\alpha$ and dividing by $Z$, we see that $g(X,Y)$ has density at most $\frac{1}{Z\kappa}$. Its support is bounded, and therefore $H_n(g(X,Y))=n+O(1)$. The construction and all the constants were fixed before $n$.
\end{proof}

Two applications of this lemma, with the roles of the sum and of $f$ interchanged, lead to the matching laws outside $\mathfrak E_S$.

\begin{proposition}\label{prop:balanced-conditional-source}
For every admissible $f$, every $0<s<1$ and every bounded rectangle $B$ with nonempty interior, there are a constant $C=C(s,f,B)<\infty$ and a fixed conditionally $(s,C;s,C)$-Frostman law supported in $B$ such that, for every integer $n\ge1$,
\begin{equation}\label{eq:balanced-conditional-outputs}
\begin{aligned}
H_n(X+Y)&=\frac{1+s}{2}n+O(1),\\
H_n(f(X,Y))&=\frac{1+s}{2}n+O(1),
\end{aligned}
\end{equation}
where the implied constants depend only on $s$, $f$ and $B$. Its marginal entropies are $n+O(1)$ and its joint entropy is $(1+s)n+O(1)$. Both marginal densities are bounded above and below on their support intervals, and both exact conditional laws are singular.
\end{proposition}

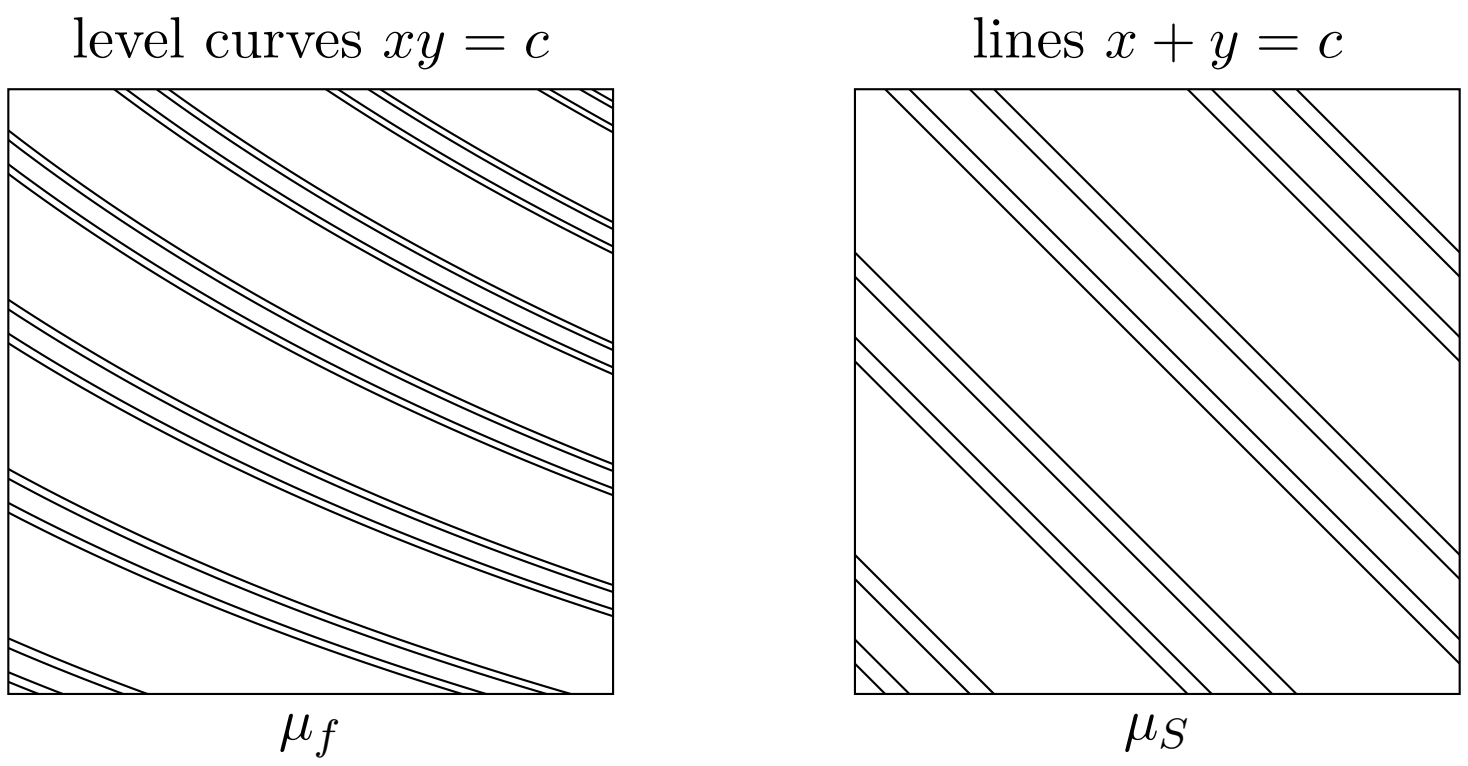
\begin{figure}[ht]
\centering
\begin{tikzpicture}[scale=4]
\begin{scope}[shift={(-2.5,-1)}]
\begin{scope}
\clip (2.5,1) rectangle (3.5,2);
\foreach \a in {2.55,3.25,3.95,4.65,5.35,6.05,6.75}{
\foreach \o in {0,0.04,0.14,0.18}{
\draw[thin,domain=2.5:3.5,samples=50] plot (\x,{(\a+\o)/\x});
}
}
\end{scope}
\draw (2.5,1) rectangle (3.5,2);
\node[below] at (3,1) {$\mu_f$};
\node[above] at (3,2) {level curves $xy=c$};
\end{scope}
\begin{scope}[shift={(-1.1,-1)}]
\begin{scope}
\clip (2.5,1) rectangle (3.5,2);
\foreach \a in {3.55,4.05,4.55,5.05}{
\foreach \o in {0,0.04,0.14,0.18}{
\draw[thin] (\a+\o-2.2,2.2) -- (\a+\o-0.8,0.8);
}
}
\end{scope}
\draw (2.5,1) rectangle (3.5,2);
\node[below] at (3,1) {$\mu_S$};
\node[above] at (3,2) {lines $x+y=c$};
\end{scope}
\end{tikzpicture}
\caption{A schematic picture of the two components in the proof of
Proposition~\ref{prop:balanced-conditional-source} for $f(x,y)=xy$
on the rectangle $[\frac{5}{2},\frac{7}{2}]\times[1,2]$, where
$f_x-f_y=y-x$ does not vanish. The law $\mu_f$ on the left is spread
along the level curves $xy=c$, with $c$ in a Cantor-like set shown
at an early stage of its construction, so that $f(X,Y)$ has entropy
$sn+O(1)$, while the sum is transverse to these curves and has entropy
$n+O(1)$. The law $\mu_S$ on the right is spread along the lines
$x+y=c$, with $c$ in a similar set, and the roles of the two outputs
are reversed.
{For each of the two laws, both marginal densities are
bounded below by positive constants on their intervals, so every
full horizontal and vertical slice of the rectangle meets its support.}
The equal mixture of the two laws has both output entropies equal
to $\frac{1+s}{2}n+O(1)$.
}
\label{fig:level-sets}
\end{figure}

\begin{proof}
Admissibility implies that $f_x$, $f_y$ and $f_x-f_y$ are nonzero polynomials, since otherwise $f$ would depend only on $y$, only on $x$ or only on $x+y$, respectively. Their product is a nonzero polynomial, whose zero set has empty interior, so we can choose a compact rectangle inside the interior of $B$ that avoids all three zero sets.

On this rectangle, we apply Lemma~\ref{lem:cantor-level-source} twice. With base $f$ and transverse function $S=x+y$, whose determinant $f_xS_y-f_yS_x$ equals $f_x-f_y$, we obtain a fixed law $\mu_f$ for which $(H_n(S),H_n(f(X,Y)))=(n,sn)+O(1)$. With base $S$ and transverse function $f$, we obtain a fixed law $\mu_S$ for which these entropies are $(sn,n)+O(1)$. Both laws have marginal entropies $n+O(1)$ and joint entropy $(1+s)n+O(1)$, and Figure~\ref{fig:level-sets} shows the two components for $f(x,y)=xy$.

We take $\mu=\frac{\mu_f+\mu_S}{2}$, with a Frostman constant at least as large as the constants of both components. Lemma~\ref{lem:progressions-mixtures} then preserves both exact conditional bounds, and its entropy inequalities, applied separately to $S$, $f(X,Y)$, $X$, $Y$ and $(X,Y)$, place each entropy between the average of the two component entropies and that average plus one bit. This proves all the entropy identities, with errors uniform in $n$. The marginal densities remain bounded above and below on the common support intervals, and after exact conditioning, the mixture is supported on the union of two null sets, so both conditional laws remain singular.
\end{proof}

\subsubsection{Mixtures at small exponents}
\label{sec:conditional-small-exponents}

For $s\le\frac{1}{2}$, we combine the law of Lemma~\ref{lem:cantor-level-source} with a product of progressions, and we obtain an upper bound for the possible gain in Theorem~\ref{thm:conditional-main} at small exponents.

\begin{proposition}\label{prop:conditional-low-upper}
Fix an admissible polynomial $f$, an exponent $0<s\le\frac{1}{2}$ and a bounded rectangle $B$ with nonempty interior. There are constants $C,K<\infty$, depending only on $f$, $s$ and $B$, such that, for every integer $n\ge1$, there is a conditionally $(s,C;s,C)$-Frostman law $(X_n,Y_n)$ supported in $B$ with
$$\max\{H_n(X_n+Y_n),H_n(f(X_n,Y_n))\}\le(2s-s^2)n+K.$$
The law may depend on $n$, while $B$, $C$ and $K$ remain fixed.
\end{proposition}

\begin{proof}
We first construct the laws for all sufficiently large $n$, with a common support and a common Frostman constant, and then fill in the finitely many remaining scales. We fix such an $n$ and write $(X,Y)=(X_n,Y_n)$. The law will be a mixture of a fixed law for which $f(X,Y)$ has small entropy and a law built from progressions for which $X+Y$ has small entropy.

Admissibility makes $f_x$, $f_y$ and $f_x-f_y$ nonzero polynomials, so we can choose a compact rectangle in the interior of $B$ on which all three are bounded away from zero. Lemma~\ref{lem:cantor-level-source}, with base $f$ and transverse function $x+y$, provides a fixed law on this rectangle with output entropies $(n,sn)+O(1)$ and both exact conditional bounds, and neither the rectangle nor this law depends on $n$. For the second law, we take independent aligned smoothed progressions in smaller coordinate intervals of the same rectangle, which can be given the same step by choosing intervals of equal length. For all sufficiently large $n$, Lemma~\ref{lem:progressions-mixtures} provides the exact conditional bounds and the bounds $(sn,2sn)+O(1)$ for the output entropies, where we used $2s\le1$.

We give the law built from progressions the weight $1-s$ and the fixed law the weight $s$. The larger of the two Frostman constants serves as a common constant, the mixture lemma preserves both exact conditional bounds, and its upper bound for entropies shows that
$$\begin{aligned}
H_n(X+Y)&\le((1-s)s+s)n+O(1)=(2s-s^2)n+O(1),\\
H_n(f(X,Y))&\le((1-s)2s+s^2)n+O(1)=(2s-s^2)n+O(1).
\end{aligned}$$
At the finitely many remaining coarse scales, we use a fixed uniform product law on a smaller rectangle inside $B$. Its exact conditional laws are $s$-Frostman and its entropy at each of these scales is finite, so enlarging the common Frostman constant and the additive constant in the entropy bound covers every $n\ge1$. The weights and both computations of coefficients remain valid at $s=\frac{1}{2}$.
\end{proof}

\subsection{Prescribing the two output laws}
\label{sec:prescribed-output-laws}

We now prove the matching assertion in Proposition~\ref{thm:joint-models} by choosing independent laws for the two outputs and pulling their product back through a regular inverse branch.

\begin{proof}[Proof of the matching assertion in Proposition~\ref{thm:joint-models}]
Let $\mathcal T(x,y)=(x+y,f(x,y))$. We prescribe independent laws for the two outputs and pull their product back through a regular inverse branch of $\mathcal T$. This produces the required output entropies at once, and the remaining task is to verify the joint rectangle bound in the original coordinates.

The Jacobian determinant of $\mathcal T$ is $f_y-f_x$, which is not identically zero, and since the zero set of a nonzero polynomial has empty interior, there is a regular point of $\mathcal T$ in $(-R,R)^2$. By the inverse function theorem, $\mathcal T$ is invertible on a neighborhood of this point, and we choose a compact product rectangle with nonempty interior inside the image of this neighborhood, small enough that its inverse image lies in $(-R,R)^2$.

For $0<s\le1$, we apply Lemma~\ref{joint-models:scalar-source} with exponent $s$ to each of the two sides of this rectangle and call the resulting laws $\alpha$ and $\beta$, with interval constants $C_\alpha$ and $C_\beta$. Their Frostman bounds hold for intervals of every length, since for lengths at least one they follow from the total mass once the constants are at least one. Let $(X,Y)$ be the pullback of $\alpha\otimes\beta$ through the chosen inverse branch. Its outputs $S$ and $F$ are independent, with laws $\alpha$ and $\beta$, so that $H_n(S)=sn+O(1)$ and $H_n(F)=sn+O(1)$ for every $n\ge1$. The branch and both scalar laws are chosen independently of $n$, so this is a single fixed law in $[-R,R]^2$, also when $s=1$.

It remains to verify the joint rectangle bound. Let $I$ and $J$ be intervals of lengths $r_1,r_2>0$. If $(X,Y)\in I\times J$, then $S$ lies in $I+J$, an interval of length $r_1+r_2$. When the sum is fixed, say $S=u$, the first coordinate lies in the interval
$$I\cap(u-J)\cap[-R,R]\cap[u-R,u+R],$$
whose length is at most $\min\{r_1,r_2\}$, and an empty interval contributes no mass. Let $L_{f,R}\ge1$ be an upper bound for $|f_x-f_y|$ on $[-R,R]^2$. The derivative of $x\mapsto f(x,u-x)$ is $f_x-f_y$, so the possible values of $F$ lie in an interval of length at most $L_{f,R}\min\{r_1,r_2\}$, and restricting to the chosen inverse branch can only decrease this set. Since $(S,F)$ has the product law $\alpha\otimes\beta$, integrating first with respect to $\beta$ and then over $u\in I+J$, we obtain
$$\begin{aligned}
\mathbb P(X\in I,Y\in J)&\le C_\beta L_{f,R}^s\min\{r_1,r_2\}^s\alpha(I+J)\\
&\le C_\alpha C_\beta L_{f,R}^s(r_1+r_2)^s\min\{r_1,r_2\}^s\\
&\le2^sC_\alpha C_\beta L_{f,R}^sr_1^sr_2^s,
\end{aligned}$$
where the last step uses $(r_1+r_2)\min\{r_1,r_2\}\le2r_1r_2$. By continuity from above, the same bound holds when one of the intervals has length zero. All choices depend only on $s$, $f$ and $R$, so the constant in the last line is a finite number $C_0=C_0(s,f,R)$ independent of $n$.

Thus, the law is jointly $(s,s,C_0)$-Frostman and has the prescribed output laws at every scale, which proves \eqref{joint-models:fixed-upper}. Together with the lower bounds of Section~\ref{sec:joint-models}, these fixed laws establish the sharpness assertions and complete the proof of Proposition~\ref{thm:joint-models}.
\end{proof}

This construction says nothing about exact conditional Frostman bounds, and it does not provide matching upper bounds for $M_{\mathrm j,n}$.

\subsection{Three-image constructions}
\label{sec:three-image-constructions}

The following fixed law attains the coefficient in Corollary~\ref{cor:sharp-three-images} for every $\frac{3}{2}\le\alpha\le2$.

\begin{proposition}\label{prop:three-image-fixed-sharpness}
Let $\frac{3}{2}\le\alpha\le2$, and let $B_0$ and $B_1$ be compact product rectangles with nonempty interiors and $B_0\subset\operatorname{int}B_1$. Suppose that $g$ is real analytic near $B_1$, with
$$g_ug_v\partial_u\partial_v\ln\left|\frac{g_u}{g_v}\right|\ne0\text{ on }B_1.$$
There are a finite constant $C_0$ and a fixed probability law of $(U,V)$ supported in the interior of $B_0$ such that $\mathbb P((U,V)\in B(z,r))\le C_0r^\alpha$ for all $z\in\mathbb R^2$ and $r>0$, and such that, for every integer $n\ge1$,
\begin{equation}\label{eq:three-image-fixed-sharpness}
H_n(W)=\frac{\alpha+1}{3}n+O(1)\text{ for }W\in\{U,V,g(U,V)\}.
\end{equation}
The constant $C_0$ and the bounded errors depend only on $\alpha$, $g$, $B_0$ and $B_1$. In particular, the coefficient in Corollary~\ref{cor:sharp-three-images} is attained by a fixed law on each such chart, for every Frostman constant $C\ge C_0$.
\end{proposition}

\begin{proof}
We construct three fixed laws, in each of which a different output has entropy coefficient $\alpha-1$ while the other two have coefficient one, and an equal mixture then balances all three outputs. Let $t=\alpha-1\in[\frac{1}{2},1]$. On each compact nondegenerate interval $J$ used below, let $\eta_{t,J}$ be the fixed $t$-Frostman probability measure of Lemma~\ref{joint-models:scalar-source}, which satisfies $H_n(\eta_{t,J})=tn+O(1)$, and let $\lambda_J$ be the normalized Lebesgue measure on $J$. These choices are made before $n$.

We choose a compact nondegenerate rectangle $B_*=I\times J$ in the interior of $B_0$, on which $|g_u|$ and $|g_v|$ are bounded above and away from zero. The two product laws $\nu_1=\eta_{t,I}\otimes\lambda_J$ and $\nu_2=\lambda_I\otimes\eta_{t,J}$ have entropy vectors
$$(H_n(U),H_n(V),H_n(g(U,V)))=(t,1,1)n+O(1)\text{ and }(1,t,1)n+O(1),$$
respectively, with coordinatewise bounded errors. Indeed, under $\nu_1$, we condition on $U=u$. The map $v\mapsto g(u,v)$ is monotone, with derivative bounded away from zero uniformly in $u$, so it pushes $\lambda_J$ forward to a law with a uniformly bounded density, and integrating in $u$ preserves this bound. Under $\nu_2$, the same argument applies to $u\mapsto g(u,v)$. A probability law with bounded support and bounded density has entropy $n+O(1)$, since its atoms at level $n$ have mass $O(2^{-n})$ and its support meets $O(2^n)$ dyadic intervals, and these bounds hold for all $n\ge1$.

The third law must have small entropy for $g$ and full entropy for both coordinates, so we prescribe it in the coordinates $u$ and $g$. Let $T(u,v)=(u,g(u,v))$. Since $g_v\ne0$, the inverse function theorem provides an inverse branch $T^{-1}(u,z)=(u,\phi(u,z))$ on a compact nondegenerate product rectangle $I'\times K$ whose image under $T^{-1}$ lies in the interior of $B_*$. After shrinking this rectangle if necessary, the branch and its inverse are bi-Lipschitz, and
$$\partial_u\phi(u,z)=-\frac{g_u(u,\phi(u,z))}{g_v(u,\phi(u,z))}$$
is bounded above and away from zero in absolute value. We take independent random variables $U$ and $Z$ with laws $\lambda_{I'}$ and $\eta_{t,K}$, set $V=\phi(U,Z)$, and denote the law of $(U,V)$ by $\nu_3$. Then $g(U,V)=Z$. For each fixed $z$, the bound on the derivative shows that $\phi(U,z)$ has a uniformly bounded density, and integrating this bound against the law of $Z$ shows that $V$ has a bounded density as well. Thus, $\nu_3$ has entropy vector $(1,1,t)n+O(1)$.

Each product $\eta_{t,J_1}\otimes\lambda_{J_2}$, in either order of the coordinates, assigns mass at most $Cr^{t+1}=Cr^\alpha$ to every ball of radius $r$, as we see by enclosing the ball in a product of two intervals of length $2r$ and applying the two interval bounds. For $\nu_3$, we note that $T$ is Lipschitz, so it maps a ball of radius $r$ into a ball of radius at most a fixed multiple of $r$, and the product bound for $\lambda_{I'}\otimes\eta_{t,K}$ yields the same isotropic Frostman estimate. Enlarging the constants to a common $C_0\ge1$, we cover all three laws and all radii $r>0$.

Finally, we set $\nu=\frac{\nu_1+\nu_2+\nu_3}{3}$. Averaging the bounds for balls shows that $\nu$ is $\alpha$-Frostman with constant $C_0$. For each of the three observables, the mixture inequalities \eqref{eq:finite-mixture-upper} place its entropy between the average of the component entropies and that average plus $\log3$, and the average coefficient is $\frac{t+2}{3}=\frac{\alpha+1}{3}$ in every coordinate. Thus, this fixed law satisfies \eqref{eq:three-image-fixed-sharpness} for every $n\ge1$, including the case $\alpha=2$, in which all the scalar measures are uniform. Its support lies in the interior of $B_0$, and replacing $C_0$ by any larger constant preserves the Frostman bound.
\end{proof}

\subsection{Arithmetic constructions and counterexamples}
\label{sec:arithmetic-counterexamples}

In this subsection, we construct independent inputs with marginal entropies $sn+O(1)$ for which both outputs have entropy strictly below $\frac{3s}{2}n$. The same arithmetic source also yields a joint counterexample for some exponents above $\frac{1}{2}$.

\begin{theorem}\label{thm:arithmetic-counterexample}
For every $0<s\le\frac{1}{2}$, there are constants $C_s,K_s<\infty$, an unbounded sequence of integers $n_j$ and independent $(s,C_s)$-Frostman random variables $X_j$ and $Y_j$ supported in $[\frac{1}{8},\frac{15}{32}]$ such that
\begin{equation}\label{eq:arith-minimal-marginals}
H_{n_j}(X_j)=sn_j+O_s(1)\text{ and }H_{n_j}(Y_j)=sn_j+O_s(1),
\end{equation}
and
\begin{equation}\label{eq:arith-conditional-upper}
\max\{H_{n_j}(X_j+Y_j),H_{n_j}(X_j^2+Y_j)\}\le\frac{149259}{100000}sn_j+K_s.
\end{equation}
The Frostman bounds at all scales also hold for both exact conditional laws.
\end{theorem}

Since $\frac{149259}{100000}=1.49259<\frac{3}{2}$, independence together with a bounded excess of marginal entropy does not imply the coefficient $\frac{3s}{2}$. The coefficient $\min\{\frac{3s}{2},\frac{1+2s}{3}\}$ also fails for $0<s\le\frac{807}{2000}$, and at the endpoint $s=\frac{807}{2000}$, its gap from our upper bound exceeds $\frac{1}{20000}$.

\subsubsection{Arithmetic sources from number fields}

The only external arithmetic input is the following special case of Sawin's construction \cite[Lemma~11(4), equation~(9)]{Sawin}. For a number field $K$, we write $\mathcal O_K$ for its ring of integers and $\Delta_K$ for its discriminant.

\begin{lemma}\label{lem:arith-tower}
Let $\mathcal T$ be a set of $r$ odd primes, an odd number of which are congruent to $3$ modulo $4$, and let $A_0=\prod_{p\in\mathcal T}p$ and $K_0=\mathbb Q(\sqrt{A_0})$. Let $\mathcal S$ be a set of $t$ odd primes outside $\mathcal T$, all inert in $K_0$, and suppose that $r+t+1<\frac{(r-1)^2}{4}$. Then there are totally real, everywhere unramified extensions $K/K_0$ of unbounded degree in which the primes of $K_0$ above $\mathcal S$ split completely. Every $p\in\mathcal S$ has residue degree two in $K$ over $\mathbb Q$, and $|\Delta_K|^{\frac{1}{[K:\mathbb Q]}}=\sqrt{4A_0}$.
\end{lemma}

\begin{proof}
We apply Sawin's lemma to the quadratic field $K_0$, with $\mathcal S$ as the set of prescribed inert primes and with no prescribed split primes, and the hypothesis $r+t+1<\frac{(r-1)^2}{4}$ guarantees its numerical condition. It provides an infinite unramified tower over $K_0$ that is totally real and in which the primes above $\mathcal S$ split completely, and the finite subextensions of this tower have unbounded degree and are the required fields $K$. For each $p\in\mathcal S$, the residue degree in $K$ over $K_0$ is one, and multiplying by the residue degree two of $p$ in $K_0$, we obtain the absolute residue degree two.

The squarefree integer $A_0$ has an odd number of prime factors congruent to $3$ modulo $4$, so $A_0\equiv3\pmod4$ and the discriminant of $K_0$ is $\Delta_{K_0}=4A_0$. Since $K/K_0$ is unramified, its relative different is trivial, and therefore
$$|\Delta_K|=|\Delta_{K_0}|^{[K:K_0]}\text{ and }[K:\mathbb Q]=2[K:K_0].$$
Taking roots of order $[K:\mathbb Q]$, we obtain the assertion about the discriminant.
\end{proof}

We fix $\mathcal T$ to be the set of the first $96$ odd primes and $\mathcal S$ to be the set of the first $2159$ primes greater than $509$ that are inert in $K_0$. Then $\max\mathcal T=509$, exactly $51$ primes of $\mathcal T$ are congruent to $3$ modulo $4$, $\min\mathcal S=521$ and $\max\mathcal S=41849$. Moreover, $\frac{95^2}{4}-(96+2159+1)=\frac{1}{4}>0$, so Lemma~\ref{lem:arith-tower} applies. We set
\begin{equation}\label{eq:arith-fixed-parameters}
\begin{aligned}
R_0&=\sqrt{4A_0}, & m&=\prod_{p\in\mathcal S}p, & \rho&=\prod_{p\in\mathcal S}\sqrt{\frac{p^2+1}{2p^2}}.
\end{aligned}
\end{equation}
The construction proceeds from these fixed sets of primes to number fields of unbounded degree, then to weighted and uniform finite sources, and finally to smoothed Frostman laws and a mixture that balances the coefficients of the two outputs. The products of primes and the parameters $a$, $b$, $h$, ${\mathcal{U}}$ and ${\mathcal{V}}$ defined below remain fixed, while the field $K$, its degree $D$, the scale of observation $n_D=\lfloor\frac{hD}{s\ln2}\rfloor$ and the law of the source vary. The Frostman bounds hold at every scale, and the upper bounds for the entropies are evaluated along the scales $n_D$.

Lemma~\ref{lem:arith-tower}, applied to these sets of primes, provides fields $K$ of unbounded degree. We write $D=[K:\mathbb Q]$, let $\sigma_1,\dots,\sigma_D$ be the real embeddings of $K$, and set $\sigma=(\sigma_1,\dots,\sigma_D)$ and $\Lambda=\sigma(\mathcal O_K)\subset\mathbb R^D$. Let $\Lambda^*$ be the dual lattice of $\Lambda$ with respect to the Euclidean inner product. The covolume of $\Lambda$ is $|\Delta_K|^{\frac{1}{2}}=R_0^{\frac{D}{2}}$, and the algebraic norm yields
\begin{equation}\label{eq:arith-norm-products}
\prod_i|\sigma_i(\alpha)|\ge1\text{ for }0\ne\alpha\in\mathcal O_K\text{ and }\prod_i|\xi_i|\ge R_0^{-D}\text{ for }0\ne\xi\in\Lambda^*.
\end{equation}
The first inequality holds because the norm of a nonzero algebraic integer is a nonzero integer. For the second, we recall that $\Lambda^*$ is the image under $\sigma$ of the codifferent of $K$, the fractional ideal of all $\beta\in K$ for which $\operatorname{Tr}_{K/\mathbb Q}(\beta\alpha)$ is an integer for every $\alpha\in\mathcal O_K$. If $\xi=\sigma(\beta)$ is a nonzero vector of $\Lambda^*$, then the product of the principal ideal $(\beta)$ with the different of $K$ is an integral ideal, and taking norms, we obtain $|N_{K/\mathbb Q}(\beta)||\Delta_K|\ge1$, which is the second inequality because $|\Delta_K|=R_0^D$.

The arithmetic parameters are defined with natural logarithms. Let
$$\phi(t)=\left(\frac{\sin\pi t}{\pi t}\right)^6,$$
with $\phi(0)=1$, and let
$${c_\phi}=\int\phi=\frac{11}{20}\text{ and }m_4=\frac{1}{{c_\phi}}\int t^4\phi(t)dt=\frac{15}{22\pi^4}<\frac{1}{16}.$$
Both integrals can be computed from the expansion of $\sin^6$ as a combination of cosines, using $\int\frac{1-\cos(\tau t)}{t^2}dt=\pi|\tau|$ for real $\tau$, and ${c_\phi}$ is also the value at the origin of the sixfold convolution of the indicator function of $[-\frac{1}{2},\frac{1}{2}]$. We set
\begin{equation}\label{eq:arith-source-exponents}
\begin{aligned}
L&=\frac{m}{32}, & x_0&=\frac{3m}{16}, & h&=\ln\frac{m{c_\phi}}{32\sqrt{R_0}},\\
a&=\ln\frac{m}{96R_0}, & b&=\ln(2m).
\end{aligned}
\end{equation}
The finite certification described below shows that
\begin{equation}\label{eq:arith-transfer-margins}
\begin{aligned}
m&>192R_0, & h&>a>0, & h&>\frac{3}{5}b, & \frac{2}{5}a&>\frac{3}{5}(b-h).
\end{aligned}
\end{equation}

For the lattice arguments below, we normalize the Fourier transform as $\widehat\mu(\xi)=\int e^{-2\pi ix\cdot\xi}d\mu(x)$. We first construct a weighted law on the lattice with an exact normalization and exact moments. We give $\alpha\in\mathcal O_K$ the weight $\prod_i\phi(\frac{\sigma_i(\alpha)-x_0}{L})$ and let $Z_K$ be the sum of these weights. The Fourier transform of the continuous product weight $x\mapsto\prod_i\phi(\frac{x_i-x_0}{L})$ is supported in $[-\frac{3}{L},\frac{3}{L}]^D$, and every vector in this box has coordinate product at most $(\frac{3}{L})^D<R_0^{-D}$, because $L>3R_0$. Thus, \eqref{eq:arith-norm-products} excludes all nonzero dual vectors from this box, and only the zero term survives in the Poisson summation formula. This yields the normalization and the coordinate moments of the resulting probability measure $q_K$, namely
\begin{equation}\label{eq:arith-poisson-identities}
Z_K=\frac{(L{c_\phi})^D}{R_0^{\frac{D}{2}}}=e^{hD}\text{ and }\mathbb E_{q_K}|\sigma_i(\alpha)-x_0|^4=L^4m_4.
\end{equation}
For the moment identity, we note that multiplying the weight by the fourth power of a coordinate differentiates its Fourier transform without enlarging the support.

All these applications of Poisson summation are absolutely convergent. Indeed, an aligned box of volume $v$ contains at most $1+\lceil v\rceil$ points of $\Lambda$. To see this, we subdivide one side of the box into $1+\lceil v\rceil$ equal pieces, which produces boxes of volume less than one, and by the norm bound applied to differences, each of these contains at most one lattice point. Summing over dyadic boxes, we obtain absolute convergence for the product weights, including the weight multiplied by a fourth power, whose factor in that coordinate still decays like $(1+|x|)^{-2}$, and the same bounds justify the periodization underlying the Poisson summation formula.

To obtain a finite support and to keep the first coordinate in a fixed positive interval, we condition on the event
$$\sum_i|\sigma_i(\alpha)-x_0|^4\le\left(\frac{m}{2}\right)^4D\text{ and }|\sigma_1(\alpha)-x_0|\le\frac{m}{32}.$$
By \eqref{eq:arith-poisson-identities} and Markov's inequality, the two complementary events have probabilities at most $\frac{m_4}{16^4}$ and $m_4$, whose sum is less than $\frac{1}{8}$. We denote the conditional law by $\mu_K$ and its finite support by $A_K$. Conditioning increases masses by a factor of at most $\frac{8}{7}$, and Minkowski's inequality shows that $\|\sigma(\alpha)\|_4\le(\frac{m}{2}+x_0)D^{\frac{1}{4}}$ for every $\alpha\in A_K$. Hence,
\begin{equation}\label{eq:arith-truncated-source}
\begin{aligned}
\max_\alpha\mu_K(\alpha)&\le2e^{-hD}, & \|\sigma(\alpha)\|_4&\le\frac{11m}{16}D^{\frac{1}{4}}, & \frac{\sigma_1(\alpha)}{m}&\in\left[\frac{5}{32},\frac{7}{32}\right],
\end{aligned}
\end{equation}
and in particular $\|\sigma(\alpha)\|_2\le D^{\frac{1}{4}}\|\sigma(\alpha)\|_4\le m\sqrt D$.

\begin{lemma}\label{lem:arith-coarse-source}
The law of $\frac{\sigma_1(\alpha)}{m}$ under $\mu_K$ assigns mass at most $C(r+e^{-aD})$ to every interval of length $r$, and the distinct points of its support are separated by at least $2e^{-\frac{1}{2}}e^{-bD}$, where $C$ does not depend on $D$.
\end{lemma}

\begin{proof}
We prove the interval estimate first for the lattice law $q_K$ and then transfer it through the fixed rescaling and the conditioning, while the separation follows from the algebraic norm bound. Let $\nu$ be the law of $\sigma_1(\alpha)$ under $q_K$. Poisson summation with a character in the first coordinate shows that $\widehat\nu$ agrees with the Fourier transform of the probability density $\frac{1}{L{c_\phi}}\phi(\frac{x-x_0}{L})$ whenever
$$|t|\le T_D=\frac{1}{2}\left(\frac{L}{3}\right)^{D-1}R_0^{-D}.$$
Indeed, a nonzero dual vector contributing to this formula would have first coordinate of absolute value at most $|t|+\frac{3}{L}$ and all other coordinates of absolute value at most $\frac{3}{L}$, so its coordinate product would be at most $\frac{1}{2}R_0^{-D}+(\frac{3}{L})^D<R_0^{-D}$, since $L>6R_0$ and $D\ge2$, which contradicts \eqref{eq:arith-norm-products}.

For an interval of length $r_0\ge T_D^{-1}$ centered at $u$, we majorize its indicator function by $\frac{\pi^2}{4}\operatorname{sinc}^2(\frac{x-u}{r_0})$, where $\operatorname{sinc}(t)=\frac{\sin(\pi t)}{\pi t}$. The Fourier transform of this majorant is supported in $[-\frac{1}{r_0},\frac{1}{r_0}]$, which lies inside the range where the two Fourier transforms agree, so the expectation of the majorant under $\nu$ equals its expectation under the continuous comparison density. The Fourier transform of that density is supported in $[-\frac{3}{L},\frac{3}{L}]$ and bounded by one, and consequently
$$\nu\left(\left[u-\frac{r_0}{2},u+\frac{r_0}{2}\right]\right)\le Cr_0\int_{-\frac{1}{r_0}}^{\frac{1}{r_0}}|\widehat\nu(t)|dt\le\frac{Cr_0}{L}.$$
Enlarging shorter intervals to length $T_D^{-1}$, we obtain the bound $\frac{C(r_0+T_D^{-1})}{L}$ for every $r_0$. Now $T_D^{-1}=\frac{2L}{3}e^{-aD}$ and $\frac{m}{L}=32$, so after rescaling by $m$ and conditioning on an event of probability at least $\frac{7}{8}$, we obtain the bound $C(r+e^{-aD})$ for intervals of normalized length $r$.

For the separation, we take distinct $\alpha,\beta\in A_K$. By \eqref{eq:arith-truncated-source}, their difference satisfies $\|\sigma(\alpha-\beta)\|_2\le2m\sqrt D$, and the inequality between the arithmetic and geometric means of the squares of the last $D-1$ coordinates shows that
$$\prod_{i>1}|\sigma_i(\alpha-\beta)|\le\left(\frac{4m^2D}{D-1}\right)^{\frac{D-1}{2}}\le e^{\frac{1}{2}}(2m)^{D-1}.$$
The algebraic integer $\alpha-\beta$ has a nonzero integral norm, so its first coordinate, divided by $m$, has absolute value at least $2e^{-\frac{1}{2}}(2m)^{-D}=2e^{-\frac{1}{2}}e^{-bD}$.
\end{proof}

\subsubsection{Uniform thinning and output counts}

The law $\mu_K$ has atoms of mass $O(e^{-hD})$, but its support may be much larger than $e^{hD}$, so its entropy is not controlled from above. We therefore pass to the uniform law on a subset of size comparable to $e^{hD}$ that keeps the coarse interval bound.

\begin{lemma}\label{lem:arithmetic-uniform-thinning}
Let $h>a>0$ and $C$ be fixed, and suppose that $\mu_D$ are probability measures on finite sets $A_D\subset\mathbb R$ such that $|A_D|\le e^{CD}$ and $\max_x\mu_D(x)\le Ce^{-hD}$, and such that $\mu_D(J)\le C(|J|+e^{-aD})$ for every interval $J$. Then, for all sufficiently large $D$, there are subsets $B_D\subset A_D$ with $|B_D|\asymp e^{hD}$ whose uniform laws satisfy the same interval bound with a fixed larger constant.
\end{lemma}

\begin{proof}
The selection should have about $e^{hD}$ points and preserve all interval bounds at once. We retain each $x\in A_D$ independently with probability $\lambda_D\mu_D(x)$, where $\lambda_D=\frac{e^{hD}}{2C}$. These probabilities are at most $\frac{1}{2}$, and the expected number of retained points is $\lambda_D$. By Chernoff's inequality, this number lies between $\frac{\lambda_D}{2}$ and $2\lambda_D$, except with probability at most $2e^{-c\lambda_D}$.

For an interval $J$ of length $r$, let $N_J$ be the number of retained points in $J$, and let $\overline\mu=C\lambda_D(r+e^{-aD})$, which bounds the expectation of $N_J$. By independence, $\mathbb E2^{N_J}\le e^{\overline\mu}$, and the exponential Markov inequality yields
$$\mathbb P(N_J>2\overline\mu)\le e^{-(2\ln2-1)\overline\mu}.$$
To make this estimate simultaneous over all intervals, it suffices to consider the at most $(|A_D|+1)^2$ sets of consecutive points that an interval can cut out of $A_D$. For each nonempty one, we apply the estimate to the shortest closed interval containing it, where for a single point, we let the length tend to zero, and every other interval cutting out the same set is at least as long. Therefore, the probability that the bound on the cardinality or any of the interval bounds fails is at most
$$2e^{-c\lambda_D}+e^{2CD+O(1)}\exp(-c'e^{(h-a)D})=o(1),$$
since $h>a$ makes $e^{(h-a)D}$ grow much faster than $O(D)$, the logarithm of the number of sets. For sufficiently large $D$, the failure probability is less than one, so a successful realization $B_D$ exists. Dividing its interval counts by $|B_D|\ge\frac{\lambda_D}{2}$, we obtain the interval bound with a fixed larger constant, while the bounds on the cardinality show that $|B_D|\asymp e^{hD}$.
\end{proof}

We apply the thinning lemma to the normalized first coordinates of the points of $A_K$. The map $\alpha\mapsto\frac{\sigma_1(\alpha)}{m}$ is injective on $A_K$ by Lemma~\ref{lem:arith-coarse-source}. By the inequality between the arithmetic and geometric means and \eqref{eq:arith-norm-products}, distinct points of $\sigma(A_K)$ are at Euclidean distance at least $\sqrt D$, and packing disjoint balls of radius $\frac{\sqrt D}{2}$ into the ball of radius $(m+\frac{1}{2})\sqrt D$ shows that $|A_K|\le(2m+1)^D\le(3m)^D$. Together with \eqref{eq:arith-truncated-source}, Lemma~\ref{lem:arith-coarse-source} and $h-a=\frac{1}{2}\ln R_0+\ln\frac{33}{20}>0$, this verifies every hypothesis of Lemma~\ref{lem:arithmetic-uniform-thinning}. We obtain subsets $B_K\subset A_K$ with $|B_K|\asymp e^{hD}$ whose uniform laws satisfy the same coarse interval bound, and their points keep the separation of Lemma~\ref{lem:arith-coarse-source}. The atoms of these uniform laws have mass $O(e^{-hD})$.

We now convert these discrete laws into laws that are Frostman at all scales. We fix $0<s\le\frac{3}{5}$ and set
\begin{equation}\label{eq:arith-observation-scale}
n_D=\left\lfloor\frac{hD}{s\ln2}\right\rfloor\text{ and }\delta_D=2^{-n_D}.
\end{equation}
We choose a point of $B_K$ uniformly, add to its normalized first coordinate an independent uniform noise on an interval of length comparable to $\delta_D$, and denote the resulting law by $\mu$. Two thresholds are relevant, namely the noise scale $\delta_D$ and the scale $e^{-aD}$, below which the coarse interval bound alone is insufficient. The definition of $n_D$ shows that $e^{-hD}\le\delta_D^s$, and \eqref{eq:arith-transfer-margins} together with $s\le\frac{3}{5}$ shows that $h>sb$ and $s(b-h)-a(1-s)<0$, where for the last inequality, we use $b-h>0$ to compare with $\frac{3}{5}(b-h)-\frac{2}{5}a<0$. Moreover, $\frac{h}{s}>a$, so $\delta_D<e^{-aD}$ for all sufficiently large $D$, and since $h>sb$, the separation of the centers is much larger than $\delta_D$. The density of $\mu$ is therefore $O(\delta_D^{s-1})$, and for $r<\delta_D$, every interval $J$ of length $r$ satisfies
$$\mu(J)\le Cr\delta_D^{s-1}=Cr^s\left(\frac{r}{\delta_D}\right)^{1-s}\le Cr^s.$$
For $\delta_D\le r\le e^{-aD}$, we enlarge the interval by the length of the noise so that it contains every contributing center. The separation and the bound for atoms give the estimate $O(e^{-hD}+re^{(b-h)D})$, while the coarse interval bound gives $O(e^{-aD})$, and taking the smaller of the two, we find that
$$\mu(J)\le Ce^{-hD}+C\min\{re^{(b-h)D},e^{-aD}\}\le Cr^s(1+e^{(s(b-h)-a(1-s))D})\le C'r^s,$$
where the middle step uses $e^{-hD}\le\delta_D^s\le r^s$ and $\min\{u,v\}\le u^sv^{1-s}$. For $e^{-aD}\le r\le1$, the coarse bound is $O(r)\le O(r^s)$, and longer intervals are controlled by the total mass. This covers all lengths, including both thresholds, with a constant independent of $D$.

The positive square root is bi-Lipschitz on $[\frac{5}{32},\frac{7}{32}]$, so applying it to the centers preserves their separation and their coarse interval bound up to fixed factors, and adding the same independent noise yields the same estimates at all scales. In both constructions, the distinct noise intervals remain separated by much more than $\delta_D$, so the observed dyadic interval determines the center and each center contributes boundedly many dyadic intervals. Hence, both laws have entropy $\log|B_K|+O(1)=sn_D+O_s(1)$.

To control the entropies of the outputs, we bound the numbers of exact sums and of exact polynomial values, keeping track of the covolume of the lattice. Let $C_{\mathrm g}=\sqrt{2\pi e}$. For $T,q>0$ with $\sqrt{2\pi}T\ge qR_0$, every coset of $q\Lambda$ contains at most
\begin{equation}\label{eq:arith-gaussian-count}
\frac{2(C_{\mathrm g}T)^D}{q^DR_0^{\frac{D}{2}}}
\end{equation}
points in a ball of radius $T\sqrt D$. To see this, we majorize the indicator function of the ball by $e^{\frac{D}{2}}e^{-\frac{|x|^2}{2T^2}}$ and apply Poisson summation over the coset. The zero term is $\frac{(C_{\mathrm g}T)^D}{q^DR_0^{\frac{D}{2}}}$. The nonzero points of the dual lattice $(q\Lambda)^*=q^{-1}\Lambda^*$ have Euclidean norm at least $\frac{\sqrt D}{qR_0}$ and are separated by at least this distance, by the norm bound for the dual lattice and the inequality between the arithmetic and geometric means. We partition these points into annuli whose radii are consecutive integer multiples of $\frac{\sqrt D}{qR_0}$. Packing disjoint balls of half this radius bounds the number of points in the $j$th annulus by $(2j+3)^D$, and the condition $\sqrt{2\pi}T\ge qR_0$ bounds each Gaussian weight there by $e^{-\pi j^2D}$. Their total contribution, relative to the zero term, is therefore at most
$$\sum_{j\ge1}(2j+3)^De^{-\pi j^2D}<1\text{ for }D\ge2.$$
The phases coming from the translation of the coset have modulus one and do not affect the bound, which proves \eqref{eq:arith-gaussian-count}.

We write $A_K+A_K=\{\alpha+\beta:\alpha,\beta\in A_K\}$ and $A_K^2+mA_K=\{\alpha^2+m\beta:\alpha,\beta\in A_K\}$. For sums, $\sigma(A_K+A_K)$ lies in the ball of radius $2m\sqrt D$, and we use \eqref{eq:arith-gaussian-count} with $T=2m$ and $q=1$. For polynomial values, the bound for the fourth norm yield
$$\|\sigma(\alpha^2+m\beta)\|_2\le\|\sigma(\alpha)\|_4^2+m\|\sigma(\beta)\|_2\le2m^2\sqrt D.$$
Moreover, $\mathcal O_K/p\mathcal O_K\cong(\mathbb F_{p^2})^{\frac{D}{2}}$ for every $p\in\mathcal S$, because $p$ splits into $\frac{D}{2}$ primes of $K$, each of residue degree two. Each field $\mathbb F_{p^2}$ contains $\frac{p^2+1}{2}$ squares, including zero, so by the Chinese remainder theorem, there are exactly
$$\prod_{p\in\mathcal S}\left(\frac{p^2+1}{2}\right)^{\frac{D}{2}}=(m\rho)^D$$
classes of squares modulo $m\mathcal O_K$, and every value $\alpha^2+m\beta$ lies in one of these classes. We apply \eqref{eq:arith-gaussian-count} to each corresponding coset of $m\Lambda$ with $T=2m^2$ and $q=m$. In both applications, the condition on $T$ holds because $m>192R_0$, and we obtain
\begin{equation}\label{eq:arith-sum-count}
|A_K+A_K|\le2\left(\frac{2C_{\mathrm g}m}{\sqrt{R_0}}\right)^D
\end{equation}
and
\begin{equation}\label{eq:arith-polynomial-count}
|A_K^2+mA_K|\le2\left(\frac{2C_{\mathrm g}\rho m^2}{\sqrt{R_0}}\right)^D.
\end{equation}
Both inequalities remain valid when $A_K$ is replaced by $B_K$.

We define
\begin{equation}\label{eq:arith-output-coefficients}
{\mathcal{U}}=\frac{1}{h}\ln\frac{2C_{\mathrm g}m}{\sqrt{R_0}}\text{ and }{\mathcal{V}}=\frac{1}{h}\ln\frac{2C_{\mathrm g}\rho m^2}{\sqrt{R_0}}.
\end{equation}
The finite bounds we need are
\begin{equation}\label{eq:arith-coefficient-bounds}
1<{\mathcal{U}}<\overline {\mathcal{U}}=\frac{1000299}{10^6}\text{ and }\frac{9}{5}<{\mathcal{V}}<\overline {\mathcal{V}}=\frac{1970181}{10^6},
\end{equation}
and these bounds on ${\mathcal{U}}$ and ${\mathcal{V}}$ determine the weights of the mixtures and the gaps in the entropy. To certify them, together with the inequalities \eqref{eq:arith-transfer-margins}, it suffices to use the following rational intervals, whose decimal endpoints are exact.
$$\begin{array}{c|cc}
&\text{lower}&\text{upper}\\\hline
\ln m&20786.509452&20786.509453\\
\ln R_0&243.850314&243.850315\\
\ln\rho&-748.252319&-748.252318\\
\ln\frac{640}{11}&4.063572&4.063573\\
\ln C_{\mathrm g}&1.418925&1.418939
\end{array}$$
The primes are generated by a sieve, and inertness is tested by the criterion $A_0^{\frac{p-1}{2}}\equiv-1\pmod p$. The products $A_0$ and $m$ have $702$ and $29989$ binary digits, respectively. For the logarithms, we first multiply each positive rational argument by an integral power of two to bring it into $[1,2]$, and for the resulting argument $x$, we set $z=\frac{x-1}{x+1}$ and use
$$\ln x=2\sum_{j=0}^{17}\frac{z^{2j+1}}{2j+1}+E,\text{ where }0\le E\le\frac{2z^{37}}{37(1-z^2)}.$$
The same formula bounds $\ln2$, so the adjustment by powers of two introduces no uncontrolled error, and the rational bounds $\frac{333}{106}<\pi<\frac{355}{113}$ suffice for the last row of the table. Substituting these intervals into the definitions of $a$, $b$, $h$, ${\mathcal{U}}$ and ${\mathcal{V}}$ proves all the displayed inequalities between the parameters, with strict margins. 
\subsubsection{Proof of Theorem~\ref{thm:arithmetic-counterexample}}

\begin{proof}
We fix $0<s\le\frac{1}{2}$ and use the smoothed sources constructed above, which are available throughout $0<s\le\frac{3}{5}$. Let $Z_D$ and $Z'_D$ be independent uniform choices from $\frac{\sigma_1(B_K)}{m}$. We use $Z'_D$ as the center of $Y_D$, and, with an independent binary selector, either $Z_D$ or $\sqrt{Z_D}$ as the center of $X_D$, and we add independent uniform noise on $[0,\frac{\delta_D}{100}]$ to both coordinates. The selector and both possible centers of $X_D$ are independent of $Y_D$, so the inputs remain independent. Each branch is $s$-Frostman with a constant independent of $D$ and has entropy $sn_D+O_s(1)$, and Lemma~\ref{lem:progressions-mixtures}, applied to the binary mixture defining $X_D$, shows that $X_D$ has the same properties. By independence, both exact conditional laws satisfy these bounds as well. The centers lie in $[\frac{5}{32},\frac{7}{32}]$ or in its image under the square root, and both of these intervals lie strictly inside $[\frac{1}{8},\frac{15}{32}]$, so for sufficiently large $D$, the noisy supports are contained in this fixed interval.

It remains to balance the bounds for the two outputs. On the branch with center $Z_D$, the exact centers of the sum and of the polynomial belong to
$$\frac{\sigma_1(B_K+B_K)}{m}\text{ and }\frac{\sigma_1(B_K^2+mB_K)}{m^2},$$
respectively. By \eqref{eq:arith-sum-count}, \eqref{eq:arith-polynomial-count} and $\frac{hD}{\ln2}=sn_D+O_s(1)$, the entropies of $X_D+Y_D$ and $X_D^2+Y_D$ on this branch have coefficients at most $s{\mathcal{U}}$ and $s{\mathcal{V}}$, respectively. On the branch with center $\sqrt{Z_D}$, the centers of the polynomial belong to $\frac{\sigma_1(B_K+B_K)}{m}$, while the sum has at most $|B_K|^2$ centers, so the two coefficients are at most $2s$ and $s{\mathcal{U}}$. The noise moves each output by $O(\delta_D)$, uniformly on the fixed support, so every center contributes only boundedly many observed dyadic intervals, and all these entropy bounds have errors $O_s(1)$.

We choose the probability $\theta$ of the branch with center $Z_D$ so that the two averaged bounds agree, that is,
$$\theta {\mathcal{U}}+2(1-\theta)=\theta {\mathcal{V}}+(1-\theta){\mathcal{U}},\text{ which gives }\theta=\frac{2-{\mathcal{U}}}{{\mathcal{V}}+2-2{\mathcal{U}}}.$$
Since $1<{\mathcal{U}}<2$ and ${\mathcal{V}}>{\mathcal{U}}$, the denominator is positive and $0<\theta<1$. Both averaged coefficients then equal $sR_{\mathrm u}$, where
\begin{equation}\label{eq:arith-uniform-rate}
R_{\mathrm u}=\frac{2{\mathcal{V}}-{\mathcal{U}}^2}{{\mathcal{V}}+2-2{\mathcal{U}}}<\frac{2\overline {\mathcal{V}}-\overline {\mathcal{U}}^2}{\overline {\mathcal{V}}+2-2\overline {\mathcal{U}}}=\frac{2939763910599}{1969583000000}<\frac{149259}{100000}.
\end{equation}
The comparison with the rational bounds follows from
$$\partial_{\mathcal{U}}R_{\mathrm u}=\frac{2(2-{\mathcal{U}})({\mathcal{V}}-{\mathcal{U}})}{({\mathcal{V}}+2-2{\mathcal{U}})^2}>0\text{ and }\partial_{\mathcal{V}}R_{\mathrm u}=\frac{(2-{\mathcal{U}})^2}{({\mathcal{V}}+2-2{\mathcal{U}})^2}>0$$
on the range \eqref{eq:arith-coefficient-bounds}. The binary selector costs at most one bit. Taking the fields along any sequence of unbounded degrees, we obtain $n_D\to\infty$ with the same Frostman constants and bounded errors in the entropy, which proves the theorem.
\end{proof}

\subsubsection{A joint arithmetic counterexample}

The same arithmetic source, mixed with two laws built from progressions, shows that the coefficient $\frac{1+3s}{4}$, which is sharp for $s\ge\frac{3}{4}$, is not a valid lower bound for every $s>\frac{1}{2}$.

\begin{corollary}\label{cor:arithmetic-joint-counterexample}
Let $s=\frac{507}{1000}$ and $f(x,y)=x^2+y$. There are constants $C,K_s<\infty$, an unbounded sequence of integers $n_j$ and jointly $(s,s,C)$-Frostman laws $(X_j,Y_j)$ supported in $[0,1]^2$ such that
$$M_{\mathrm j,n_j}\le\left(\frac{1+3s}{4}-\frac{1}{20000}\right)n_j+K_s.$$
The laws may vary with $j$, but their support and the constants do not, and their Frostman bounds hold at every scale. More generally, there is $s_*>\frac{507}{1000}$ such that the coefficient $\frac{1+3s}{4}$ fails for each fixed $s\in(\frac{1}{2},s_*)$.
\end{corollary}

\begin{proof}
We take two independent copies of the smoothed arithmetic law with centers in $\frac{\sigma_1(B_K)}{m}$, which is available for $0<s\le\frac{3}{5}$, and we mix this law with two laws built from progressions so that all three entries of $M_{\mathrm j,n}$ have the same upper coefficient. Let $s_*=\frac{1}{{\mathcal{U}}+{\mathcal{V}}-1}$. By \eqref{eq:arith-coefficient-bounds}, $\frac{507}{1000}<s_*<\frac{5}{9}<\frac{3}{5}$. We fix $\frac{1}{2}<s<s_*$ and use the scales \eqref{eq:arith-observation-scale}. Listing the coefficients in the order $\frac{H_n(X)+H_n(Y)}{2}$, $H_n(S)$ and $H_n(F)$, the arithmetic law has the vector of upper coefficients $s(1,{\mathcal{U}},{\mathcal{V}})$, and by independence, it satisfies a joint Frostman bound with a constant independent of $D$.

Two laws built from progressions provide the remaining vectors $(s,1,s)$ and $(\frac{1+s}{2},s,s)$. For the first, we take independent aligned smoothed progressions as $X^2$ and $Y$, with $X$ on a fixed branch of the positive square root. The marginals and $F=X^2+Y$ have coefficient $s$, while the sum has coefficient at most one. For the second, we take independent aligned smoothed progressions $T$ and $S$, invert $T=X^2-X$ on a fixed branch near $X=\frac{1}{5}$ and set $Y=S-X$, which lies near $\frac{3}{10}$. We choose the supporting intervals so small that both laws remain in $[0,1]^2$. Then $X+Y=S$ and $F=T+S$, so both output coefficients equal $s$, and $X$ has coefficient $s$ while $Y$ has coefficient at most one. Lemma~\ref{lem:progressions-mixtures} provides all the marginal Frostman bounds used here, and for the second law, the independence of $X$ and $S$ also yields
$$\mathbb P(X\in I,Y\in J)=\int_I\mathbb P(S\in J+x)d\mathcal L(X)(x)\lesssim|I|^s|J|^s.$$
Thus, the three laws satisfy a common joint Frostman bound at every scale.

To balance the three entries, we set $E_s=1-s+s(3{\mathcal{V}}-{\mathcal{U}}-2)$ and mix the laws, in the order just described, with the weights $\vartheta_s=\frac{1-s}{E_s}$, $\frac{s({\mathcal{V}}-{\mathcal{U}})}{E_s}$ and $\frac{2s({\mathcal{V}}-1)}{E_s}$. Since ${\mathcal{V}}>{\mathcal{U}}>1$ and $s<1$, these weights are positive, and their numerators add up to $E_s$. The marginal and polynomial entries of the averaged vector both equal $s+\vartheta_ss({\mathcal{V}}-1)$, and the entry for the sum agrees with them because
$$\vartheta_ss({\mathcal{U}}-1)+\frac{s({\mathcal{V}}-{\mathcal{U}})}{E_s}(1-s)=\vartheta_ss({\mathcal{V}}-1).$$
The label of the mixture costs at most $\log3$ bits, so every term of the maximum has this upper coefficient with a bounded error, and its gap below $\frac{1+3s}{4}$ is
\begin{equation}\label{eq:arith-joint-gap}
\frac{1+3s}{4}-s-\vartheta_ss({\mathcal{V}}-1)=\frac{(1-s)(1-s({\mathcal{U}}+{\mathcal{V}}-1))}{4E_s}.
\end{equation}
The denominator is positive, and the numerator is positive exactly when $s<s_*$, which covers the interval $\frac{1}{2}<s<s_*$. At $s=\frac{507}{1000}$, we use ${\mathcal{U}}>1$, ${\mathcal{U}}<\overline {\mathcal{U}}$ and ${\mathcal{V}}<\overline {\mathcal{V}}$ to bound the numerator from below and the denominator from above, and we obtain the lower bound
$$\frac{(1-s)(1-s(\overline {\mathcal{U}}+\overline {\mathcal{V}}-1))}{4(1-s+s(3\overline {\mathcal{V}}-3))}=\frac{5956919}{98432265050}>\frac{1}{20000}.$$
Mixtures preserve the common rectangle bound and the support in the unit square, and the degrees of the number fields provide an unbounded sequence of scales, which completes the proof of both claims.
\end{proof}

\section{Some open questions}
\label{sec:conjectures-support}

We ask whether the gains in the main estimates can be strengthened to the formulas below. As in the main theorems, the polynomial, the support and the Frostman constants remain fixed, while the law may vary with $n$.

\begin{conjecture}\label{conj:joint-target}
Under the hypotheses of Theorem~\ref{thm:joint-main}, for every $\varepsilon>0$,
\begin{equation}\label{eq:joint-target}
M_{\mathrm j,n}\ge(s+\gamma(s,d)-\varepsilon)n-O_\varepsilon(1),
\end{equation}
where $d=\deg f$ and
\begin{equation}\label{eq:joint-target-gain}
\gamma(s,d)=\min\left\{\frac{s}{2d+1},\frac{1-s}{4}\right\}.
\end{equation}
\end{conjecture}

The two expressions in the minimum agree at $s=\frac{2d+1}{2d+5}$. Theorem~\ref{thm:joint-main} proves the conjecture for $s\ge\frac{3}{4}$, and for every $s>\frac{1}{2}$ when $d\ge8$. The examples of Proposition~\ref{prop:joint-high-mixture-upper} give matching upper bounds whenever the minimum in \eqref{eq:joint-target-gain} is attained by $\frac{1-s}{4}$. For $d=2$ and $s\le\frac{1}{2}$, part (ii) of Proposition~\ref{prop:joint-quadratic-mixture-upper} gives the proposed coefficient $\frac{6s}{5}$ as an upper bound, and the matching lower bound is open. For $d\ge3$ and $0<s\le\frac{1}{2}$, Proposition~\ref{prop:degree-upper-latest} gives the upper gain $\frac{s}{2d-1}$, which leaves a gap from $\frac{s}{2d+1}$. Here $d$ is the degree of the full polynomial, including any nonlinear outer composition. Finally, Corollary~\ref{cor:arithmetic-joint-counterexample} shows that the coefficient $\frac{1+3s}{4}$ cannot hold for every $s>\frac{1}{2}$.

\begin{question}\label{conj:conditional-target}
Under the hypotheses of Theorem~\ref{thm:conditional-main}, is it true that, for every $\varepsilon>0$,
\begin{equation}\label{eq:conditional-target}
M_n\ge\left(\min\left\{\frac{4s}{3},\frac{1+2s}{3}\right\}-\varepsilon\right)n-O_\varepsilon(1)?
\end{equation}
\end{question}

The two coefficients agree at $s=\frac{1}{2}$, and Theorem~\ref{thm:conditional-main} answers the question affirmatively for $s\ge\frac{2}{3}$. Proposition~\ref{upper-modular:prop} provides examples with the upper coefficient $\frac{1+2s}{3}$ for $s\ge\frac{1}{4}$, so if the proposed lower bound holds, then the coefficient $\frac{1+2s}{3}$ cannot be increased for $\frac{1}{2}\le s<1$ in an estimate valid for every admissible polynomial. We have no matching example for the coefficient $\frac{4s}{3}$ at smaller exponents.

Theorem~\ref{thm:arithmetic-counterexample} rules out the coefficient $\frac{3s}{2}$ even for independent inputs with $H_n(X)=H_n(Y)=sn+O(1)$, but its coefficient $1.49259s$ does not decide \eqref{eq:conditional-target}. The upper coefficient $2s-s^2$ of Proposition~\ref{prop:conditional-low-upper} also leaves the question open, although it rules out the lower coefficient $\frac{1+2s}{3}$ for every $0<s<\frac{1}{3}$.

We also ask for sharp versions of Theorems~\ref{thm:joint-main}--\ref{thm:conditional-classification} with unequal Frostman exponents.

\begin{question}\label{q:unequal-exponents}
Let $0<s_1,s_2<1$ with $s_1\ne s_2$, and let $f$ be admissible. What are the optimal leading coefficients in lower bounds for $M_{\mathrm j,n}$ under joint $(s_1,s_2,C)$-Frostman bounds, and for $M_n$ under conditional $(s_1,C_1;s_2,C_2)$-Frostman bounds, uniformly over laws with fixed support and fixed Frostman constants? How do these coefficients depend on $(s_1,s_2)$ and on $f$, and how does the conditional classification according to whether $f\in\mathfrak E_S$ extend to this setting?
\end{question}

\section*{Acknowledgements}
The authors would like to thank Chun-Yen Shen for many helpful discussions, particularly regarding possible extensions of Corollary~\ref{cor:polynomial-image-product} in terms of Hausdorff dimension via incidence bounds.

\end{document}